\documentclass{article}
\usepackage[utf8]{inputenc}
\usepackage{amsmath}
\usepackage{amssymb}
\usepackage{amsfonts}
\usepackage{mathtools}
\usepackage{thmtools}
\usepackage{pgf,tikz}
\usepackage{mathrsfs}
\usetikzlibrary{arrows}
\usepackage{graphicx}
\usepackage{textcomp}
\usepackage{natbib}
\usepackage{graphicx}
\usepackage{float}
\usepackage{lipsum}
\usepackage{mwe}
\DeclarePairedDelimiter{\ceil}{\lceil}{\rceil}
\usepackage{url}
\usepackage{hyperref}
\usepackage{xcolor}
\let\bibsection\relax
\usepackage[alphabetic]{amsrefs}
\usepackage{microtype}

\newtheorem{lemma}{Lemma}[section]

\newtheorem{Definition}[lemma]{Definition}

\newtheorem{Remark}{Remark}[section]

\title{\Large Construction of tropical morphisms from tropical modifications of nonhyperelliptic genus $3$ metric graphs with tree gonality $3$ to metric trees}

\author{Hamdi D\"ervodeli}
\date{ }

\begin{document}

\maketitle
\begin{abstract}

In  this  article,  we  look  into  the  tree gonality of  genus $3$ metric  graphs  $\Gamma$ which is defined as the minimum of degrees of all tropical morphisms from any tropical modification of $\Gamma$ to any metric tree.  It is denoted by $\textnormal{tgon} (\Gamma)$ and is  at  most $3$. We define hyperelliptic metric graphs in terms of tropical morphisms and tree gonality.  Let $\Gamma$ be a genus $3$ metric graph with $\textnormal{tgon} (\Gamma) = 3$ which is not hyperelliptic. In  this  paper, for such  metric  graphs  $\Gamma$,  we  construct  a tropical modification $\Gamma'$ of $\Gamma$, a metric tree $T$ and a tropical map $\varphi:\Gamma' \to T$ of degree $3$. 
\end{abstract}
\tableofcontents
\newpage 

\section{Introduction}
We look into the tree gonality of metric graphs. Its motivation comes from the striking interplay between graphs and algebraic curves discovered over the last two decades. For example, there exists a good theory of divisors on graphs (see \cite{BN1}) (also including such notions as linear systems, linear equivalences, canonical divisors, degrees, and ranks), and maps between metric graphs with suitable balancing conditions that behave similarly to morphisms between curves (see \cite{BN1}). \par
Recall that the gonality of an algebraic curve $C$ is the minimum of degrees of all non-constant morphisms from $C$ to the projective line $\mathbb{P}^1$. There are two notions of graph gonality in the literature, which are both inspired by the gonality of an algebraic curve. They are tree (or geometric) gonality and divisorial gonality e.g., studied for ordinary or metric graphs (see \cite{BM}). Yet another variant is stable gonality, which is the infimum of the divisorial gonality over all subdivisions of an ordinary graph (see \cite{CFK}).
\par 
We study a tropical version of gonality, where the roles of algebraic curves and the projective line are played by metric graphs and metric trees, respectively, and the morphisms are replaced by the tropical morphisms (see \cite{BN1}, \cite{LC}, \cite{GM}, \cite{BBG},  \cite{MCH}). 
The tree gonality of a metric graph $\Gamma$ is defined as minimum of degrees of all tropical morphisms from any tropical modification of $\Gamma$ to any metric tree. The tree gonality of any metric graph of genus $g$ is at most $\ceil*{\frac{g}{2}}+1$ (see Theorem 1, \cite{JA}). Its proof is entirely combinatorial and provides an explicit method to construct divisors with degree $\ceil*{\frac{g}{2}}+1$ and rank $1$ on genus-$g$ metric graphs. In this article, we are interested in constructing a degree-($\ceil*{\frac{g}{2}}+1$) tropical morphism from a tropical modification of $\Gamma$ to a metric tree, where $\Gamma$ is of genus $g$ with tree gonality $\ceil*{\frac{g}{2}}+1$. Interest for such method dates back to (\cite{BM}, Remark 3.13).  In this regard, our modest contribution is on the case where $g=3$ and $\Gamma$ is not hyperelliptic, i.e., given a nonhyperelliptic genus $3$ metric graph $\Gamma$ with tree gonality $3$, we construct a tropical modification $\Gamma'$, a metric tree $T$, and a degree $3$ tropical morphism $\phi: \Gamma' \to T$ (Problem 1). We emphasize that our constructions are more direct than in \cite{JA} in the sense that we avoid constructing divisors of certain degree and rank, but rather make explicit constructions of tropical morphisms from tropical modifications of metric graphs to metric trees.
\par 

\textbf{Problem 1}. Let $\Gamma$ be a genus $3$ metric graph with tree gonality $3$ which is not hyperelliptic. Construct a tropical modification $\Gamma'$ of $\Gamma$, a metric tree $T$ and a tropical morphism $\varphi:\Gamma' \to T$ of degree $3$.
\newpage 
\section{Preliminaries}
\subsection{Metric graphs.}
A \textit{graph} $G$ is defined by the following data: a set  $V$ called the \textit{vertex set}, a set $E$ called the \textit{edge set} and a map $\partial:E\to P(V)$  such that for any $e \in E$ we have $|\partial(e)| = 1$ or $|\partial(e)|=2$, where $P(V)$ is the power set of $V$. We write $G=(V,E,\partial)$. The elements of $V$ (resp. $E$) are called \textit{vertices} (resp. \textit{edges}) of $G$. An edge $e \in E$ with $|\partial(e)| = 1$ is called a \textit{loop}. Two or more edges $e_1,e_2,\ldots,e_n \in E$ are called \textit{multiple edges} if there exist $v_1,v_2 \in V$ such that $\partial(e_i)=\left\lbrace v_i, v_j \right\rbrace$ for all $i=1,2,\ldots,n$. The graph $G$ is said to be \textit{finite} if both $V$ and $E$ are finite sets. A length map on $G$ is any function $l:E \to (0, +\infty)$. In this article, unless stated otherwise, a graph is always assumed to be finite with multiple edges allowed. \par 
Let  $G = (V,E,\partial)$ be a graph. A \textit{path} in the graph $G$ is a sequence of edges $(e_1,e_2\ldots,e_{n-1})$ for which there exists a sequence of vertices $(v_1,v_2,\ldots,v_n)$ such that $\partial(e_i) = \left\lbrace v_i, v_{i+1} \right\rbrace$  for $i=1,2,\ldots,n-1$. If $w=(e_1,e_2,\ldots,e_{n-1})$ is a path in $G$ with vertex sequence $(v_1,v_2,\ldots,v_n)$ then $w$ is said to be a \textit{path from} $v_1$ \textit{to} $v_n$. A graph $G$ is said to be \textit{connected} if for any two vertices $v_1$ and $v_2$ there exists a path from $v_1$ to $v_2$.  Let  $e \in E$ with $\partial(e)=\left\lbrace v,w\right\rbrace$. \textit{Subdividing} the edge $e \in E$ with $\partial(e)=\left\lbrace v,w\right\rbrace$ into edges $e_1,e_2$ yields the graph $G'=(V',E',\partial')$ where $V'=V \cup \left\lbrace z \right\rbrace$, $E'= (E \setminus \left\lbrace e \right\rbrace) \cup \left\lbrace e_1,e_2 \right\rbrace$ and $\partial'$ is given by $\partial'|_{E \setminus \left\lbrace e \right\rbrace}=\partial$ and $\partial'(e_1)=\left\lbrace v,z \right\rbrace  , \partial'(e_2)= \left\lbrace z,w \right\rbrace$. \par

Let $G=(V,E, \partial)$ be a connected graph with no loops. An \textit{orientation} on $G$ is a map $\Vec{\partial}:E \to V \times V$ such that if we write $\Vec{\partial}(e)= (v_1, v_2)$ then $\partial(e)=  \left\lbrace v_1, v_2 \right\rbrace $. Note that giving an orientation $\Vec{\partial}$ on $G$ is equivalent to giving a map $(\partial_0,\partial_1):E \to V\times V$ where $\partial_0, \partial_1:E \to V$ are endpoint maps. \par  Fix an orientation $(\partial_0,\partial_1):E \to V \times V$ on $G$ and choose a length map $l$ on $G$. Let $(X,d)$ be the disjoint union of the real metric spaces $[0,l(e)]$ for $e \in E$ i.e., the set $$ X = \bigsqcup_{e \in E} [0,l(e)]:= \bigcup_{e \in E} [0,l(e)] \times \left\lbrace e \right\rbrace $$
together with the metric $d:X\times X \to [0, \infty]$ given by
\begin{equation*}
 d((x_1,e_1),(x_2,e_2))  =\begin{cases}
    |x_1-x_2|, & \text{if $e_1=e_2$}\\
    \infty, & \text{otherwise}.
  \end{cases}
\end{equation*}
Consider the relation $\sim_1$ on $X$ defined by $x \sim_1 y$ if there exists a vertex $v \in V$ such that $x , y \in \{ (0, e) \in X \mid \partial_0 (e) = v \} \cup \{ (l(e) , e) \in X \mid \partial_1 (e) = v \}$ and let $\sim$ be the equivalence relation on $X$ generated by $\sim_1$ i.e., $x \sim y$ if and only if $x=y$ or there exists a finite subset $\left\lbrace z_1,z_2,\ldots,z_n \right\rbrace \subset X$ such that $x=z_1, z_n=y$ and $z_i \sim_1 z_{i+1}$ for $i=1,2,\ldots,n-1$. Let $\Bar{X}:= X/_{\sim}$ be the quotient space of $X$ with respect to the equivalence relation $\sim$ and  $\Bar{d}: \Bar{X}\times \Bar{X} \to [0, \infty) $ be given by
$$\Bar{d}(\Bar{x},\Bar{y}):= \textnormal{inf}  \sum_{i=1}^{k} d(x_i, y_i)$$
where the infimum is taken over all $k \in \mathbb{N}$ and sequences $(x_1,y_1,x_2, \ldots, x_k, y_k)$ in $X$ such that $x_1 \in \Bar{x}, x_{i+1} \sim y_i$ for $i=1,2,\ldots,k-1$ and $y_k \in \Bar{y}$. Then, $\Gamma:=(\Bar{X},\Bar{d})$ is a metric space. In this case, we say that the metric space $(\Bar{X},\Bar{d})$ is obtained from $(G, l)$ by gluing intervals $[0, l(e)]$, one for each $e \in E$, along their endpoints in the manner prescribed by $G$. We often regard each edge $e \in E$ as a subset of $\Gamma$ and each vertex $v \in V$ as a point in $\Gamma$.

\begin{Definition} \textnormal{
A \textit{metric graph} is a metric space $\Gamma$ such that there exists a loopless connected graph $G$ with a length map $l$ such that $\Gamma$ is isometric to the metric space obtained from $(G, l)$ by gluing intervals $[0, l(e)]$, one for each $e \in E(G)$, along their endpoints in the manner prescribed by $G$.}
\end{Definition}
The pair $(G,l)$  is called a \textit{model} of $\Gamma$ whereas $\Gamma$ is called a \textit{realization} of the model $(G,l)$. The construction of a metric graph from a graph that may have loops will be given in the following way. Let $G=(V,E,\partial)$ be a connected graph with loops and $l$ a length function on $G$. Subdividing all the loops $e \in E$, say into $e_1,e_2,\ldots, e_n$, yields a graph $G'$ with no loops. The length map $l'$ on $G'$ is given by $l'=l$ on $E\setminus \left\lbrace e \in E \hspace{0.1cm} | \hspace{0.1cm} \partial(e)=1 \right\rbrace$ and $l'(e_1)+l'(e_2) + \ldots+l'(e_n)= l(e)$ for edges $e_1, e_2$ for which a loop $e \in E$ subdivided to $e_1,e_2$. Then  $\Gamma$ does not depend on the choice of the subdivision $(G' , l')$. Thus, we define $\Gamma$ to be the realization of $(G,l)$, and we also call $(G,l)$ a model (that may have loops) of $\Gamma$. \par

The first Betti number of $\Gamma$ is equal to $g(G):= |E(G)|-|V(G)|+1$. It is called the \textit{genus} of $\Gamma$ and it is denoted by $g(\Gamma)$. A metric graph $\Gamma$ of genus $g(\Gamma)=0$ is called a \textit{metric tree}. \par 
Let $\Gamma$ be a metric graph. A \textit{vertex set} of $\Gamma$ is a finite subset $S \subset \Gamma$ such that the subspace $\Gamma \setminus S$ is isometric to a disjoint union of finitely many real open intervals. Any vertex set $S \not= \emptyset$ of $\Gamma$ induces a model $(G_S,l_S)$ of $\Gamma$ in the following way. The graph $G_S = (V,E,\partial)$ is given by its vertex set $V:=S$,  its edge set $E$ defined to be the set of closures of finitely many connected components of $\Gamma \setminus V$ and the map $\partial: E \to P(V)$ given by $e \longmapsto \partial(\textnormal{int} (e))$, where $\textnormal{int} (e) = e \setminus S $ and $\partial(\textnormal{int} (e)) \subset V $ is its boundary in $\Gamma$. Each edge $e \in E$  is isometric to either a segment or a circle. The length map $l_S:E \to (0, \infty)$ assigns each edge $e \in E$ the length of the segment or circle isometric to it. \par 
We single out a particular model for $\Gamma$. A point $x \in \Gamma$ is called an \textit{essential} vertex if for any $\varepsilon > 0$, the open ball $ B(x,\varepsilon):= \left\lbrace y \in \Gamma \hspace{0.1cm} | \hspace{0.1cm} \Bar{d} (x,y) < \varepsilon \right\rbrace$ is not isometric to $(-\varepsilon, \varepsilon) \subset \mathbb{R}$. If $x \in \Gamma$ is an essential vertex, then for any model $(G,l)$ of $\Gamma$ and any edge $e \in E(G)$ we have $ x \notin \textnormal{int } (e)$, and so,  the set of essential vertices of $\Gamma$ is a subset of $E(G)$ for any model $(G, l)$ of $\Gamma$. Since $G$ is a finite graph, $\Gamma$ has only finitely many essential vertices.

\begin{lemma} \textnormal{Let $\Gamma$ be a metric graph, $\mathcal{E}$ the set of essential vertices of $\Gamma$, and  $S$ a finite nonempty subset of $\Gamma$. Then, the set $S$ is a vertex set of $\Gamma$ if and only if $\mathcal{E} \subseteq S$.
} 
\end{lemma}
\textit{Proof.} Suppose that $\emptyset \not=S $ is a vertex set in $\Gamma$. Then, $S$ induces a model $(G,l)$ of $\Gamma$ where $S=V(G)$ and, so
$$
\Gamma \setminus S  =  
\Gamma \setminus V(G)   \equiv   \bigsqcup_{e \in E(G)} (0,l(e)).
$$ 
If $\mathcal{E} \cap (\Gamma \setminus S) \not= \emptyset $ then there exists $x \in \mathcal{E}$ and an edge $e \in E(G) $ such that $x \in \textnormal{int } (e)$ which contradicts $x$ being an essential vertex. Therefore, $ \mathcal{E}\cap (\Gamma \setminus S) = \emptyset$ and $ \mathcal{E} \subseteq S$. Now, assume that $\mathcal{E} \subseteq S$. If  $\mathcal{E}= \emptyset$ then $\Gamma$ is isometric to a circle, and so, any non-empty finite subset of $\Gamma$ is a vertex set. Suppose that $ \mathcal{E} \not= \emptyset$. Let $(G,l)$ be a model of $\Gamma$, and $V$, $E$ be the set of vertices, edges of $G$ respectively. Then, the set $V= \bigcup_{e \in E} \partial (e)$, where $\partial(e)$ is the boundary set of $e \subset \Gamma$, is a vertex set of $\Gamma$. 
As $\mathcal{E}$ is the set of essential vertices, and $V$ is a vertex set, it follows, from what we have shown above, that $\mathcal{E} \subseteq V$. Now, if $\mathcal{E} = V$, then $\mathcal{E}$ is a vertex set. Assume that $\mathcal{E} \subsetneq V$. We know that the set $ V \setminus \mathcal{E}$ is always finite. If this is a one-element set i.e., $V \setminus \mathcal{E} = \left\lbrace x_1 \right\rbrace $, then there exist unique edges $e_1,f_1 \in E, e_1 \not=f_1$ such that $x_1$ is a common endpoint of $e_1$ and $f_1$.  Then, we obtain that 
\begin{align*}
     \Gamma \setminus \mathcal{E} = (\Gamma \setminus V) \cup \left\lbrace x_1 \right\rbrace 
   &\equiv  \bigsqcup_{e \in E} (0,l(e)) \cup \left\lbrace x_1 \right\rbrace \\ 
   &\equiv  \bigsqcup\limits_{\substack{ e \in E \\ e \not=  e_1, f_1}}   (0,l(e)) \sqcup (0,l(e_1)+l(f_1))
\end{align*}
which implies that $\mathcal{E}$ is a vertex set. If $V \setminus \mathcal{E} = \left\lbrace x_1,x_2 \right\rbrace $, then there exist unique edges $e_i,f_i \in E$ with $ e_i \not=f_i$ such that $x_i$ is a common endpoint of $e_i$ and $f_i$ for $i=1,2$. In the case when one of $e_1$ and $e_2$ is equal to one of $f_1$ and $f_2$, say, $f_1=e_2$, we have that
$$ \Gamma \setminus \mathcal{E} \equiv  \bigsqcup\limits_{\substack{ e \in E \\ e \not=  e_1, e_2, f_2}}   (0,l(e)) \hspace{0.1cm} \sqcup \hspace{0.1cm} (0,l(e_1)+l(e_2)+l(f_2)).$$
If both $e_1$ and $e_2$ are different to both $f_1$ and $f_2$, then
$$ \Gamma \setminus \mathcal{E} \equiv  \bigsqcup\limits_{\substack{ e \in E \\ e \not=  e_1, e_2, f_1, f_2}}   (0,l(e)) \hspace{0.1cm} \sqcup \hspace{0.1cm}  (0,l(e_1)+l(e_2)) \hspace{0.1cm} \sqcup \hspace{0.1cm} (0,l(f_1)+l(f_2))$$
and therefore, $\mathcal{E}$ is a vertex set. Similarly we get we get that $\mathcal{E}$ is a vertex set if $V \setminus \mathcal{E} = \left\lbrace x_1,x_2,\ldots,x_n \right\rbrace$, Thus, $\Gamma \setminus \mathcal{E}$  is isometric to a disjoint union of finitely many open real intervals. Since $\Gamma \setminus S \subset \Gamma \setminus \mathcal{E}$, we have that $\Gamma \setminus S$ is also is isometric to a disjoint union of finitely many open intervals, and therefore, $S$ is a vertex set. $\square$ \par 
A metric graph is said to be a \textit{metric loop} if it is isometric to a circle. If $\Gamma$ is not a metric loop, then $\mathcal{E} \not= \emptyset$ is a vertex set of $\Gamma$. The model $(G_{\mathcal{E}}, l_{\mathcal{E}})$ induced by the essential vertex set $\mathcal{E}$ is called the \textit{essential model} of $\Gamma$. From Lemma 2.1, the essential model $(G_{\mathcal{E}}, l_{\mathcal{E}})$ is minimal in the sense that any other model of $\Gamma$ can be obtained by a sequence of edge subdivisions of $G_{\mathcal{E}}$. Thus, all models are refinements of the essential model. In addition, this implies that the \textit{valence} of a point $x \in \Gamma$ defined as the valence of $x$ in $G_S$ for $S$ a vertex set of $\Gamma$ and $x \in S$, is well-defined notion. The valence of the point $x \in \Gamma$ is denoted by $\textnormal{val} (x)$. 
 
\subsection{Harmonic maps and tropical morphisms.}

\begin{Definition}
\textnormal{
Let $\Gamma_1$ and $\Gamma_2$ be metric graphs with loopless models $(G_1,l_1)$ and $(G_2,l_2)$ respectively, where $E(G_1) = \left\lbrace e_1 \right\rbrace $ and $E(G_2) = \left\lbrace e_2 \right\rbrace $. A map $\varphi: \Gamma_1 \to \Gamma_2$ is said to be \textit{linear} if there exist isometries $\rho_1:\Gamma_1 \to [0, l_1(e_1)]$ and $\rho_2:\Gamma_2 \to [0, l_2(e_2)]$ such that the map $\rho_2 \circ \varphi \circ \rho_1^{-1}: [0, l_1(e_1)] \to [0, l_2(e_2)] $ is an affine linear map. 
}
\end{Definition}

\begin{Definition}
\textnormal{Let $\Gamma_1$ and $\Gamma_2$ be two metric graphs. A continuous map $\varphi: \Gamma_1 \to \Gamma_2$ is said to be \textit{piecewise linear} if there exist loopless models $(G_1, l_1)$ and $(G_2,l_2)$ of $\Gamma_1$ and $\Gamma_2$ respectively, such that for any edge $e_1 \in E(G_1)$ there exists an edge $e_2 \in E(G_2)$ such that $\varphi(e_1) \subseteq e_2$ and $\varphi_{|e_1} :e_1 \to e_2$ is a linear map. 
}
\end{Definition}
Let $\varphi:\Gamma_1 \to \Gamma_2$ be a piecewise linear map of metric graphs, $v \in \Gamma_1$ and $w:=\varphi(v)$. Let $(G_1,l_1)$ (resp., $(G_2,l_2)$) be loopless models of $\Gamma_1$ (resp., $\Gamma_2$) such that for all $e_1 \in E(G_1)$ there exists $e_2 \in E(G_2)$ such that $\varphi(e_1) \subseteq e_2$, $\varphi_{|e_1} :e_1 \to e_2$ is a linear map, and assume that $v \in V(G_1)$ and $w \in V(G_2)$. Fix a direction $\Vec{w}$ at $w$ (i.e., a 'unit vector' starting at $w$ with direction of a path emanating from $w$), and let $e_2 \in E(G_2)$ such that $w$ is an endpoint of $e_2$ and $e_2$ is in the direction $\Vec{w}$.
Let $\left\lbrace e_{v1},e_{v2}, \ldots , e_{vr} \right\rbrace \subseteq E(G_1)$ be the set of edges emanating from $v$. Without loss of generality, assume that 
$$ \left\lbrace e_{v1},e_{v2}, \ldots , e_{vs} \right\rbrace = \left\lbrace e_{vj} \hspace{0.1cm} | \hspace{0.1cm} \varphi(e_{vj}) \subseteq e_2,  \hspace{0.1cm} j=1,2,\ldots, r \right\rbrace $$ 
for some $s$ such that $0\leqslant s \leqslant r$. Then, $\varphi|_{e_{vj}} : e_{vj} \to e$ is a linear map for $j=1,2,\ldots, s$ because of the choice of models $(G_1,l_1)$ and $(G_2,l_2)$. Denote by $m_{\varphi,\Vec{w}} (v)$ the sum of slopes of these linear maps $\varphi|_{e_{vj}}$, $j=1,2,\ldots,s$.  i.e.,
$$ m_{\varphi,\Vec{w}} (v) = \sum_{j=1}^s \textnormal{slope } (\rho \circ \varphi \circ \rho_{vj}^{-1}) $$
where $\rho : e_2 \to [0, l_2(e_2)]$ and  $\rho_{vj}: e_{vj} \to [0, l_1(e_{vj})]$ are the chosen isometries with unique parametrizations $\rho(w)=\rho_{vj}(v)=0$ for $i=1,2,\ldots,s$ i.e., that map initial endpoints of $e_2, e_{vj}, j=1,2,\ldots,s$ to $0$. This definition of the slope of the linear maps $\varphi|_{e_{vj}}, j=1,2,\ldots,s$, and their sum $m_{\varphi,\Vec{w}}(v)$ is independent of the choice of such models $(G_1,l_1)$ and $(G_2,l_2)$. \par 
\begin{Definition} \textnormal{A continuous map  $\varphi:\Gamma_1 \to \Gamma_2$ is said to be a \textit{harmonic map} of metric graphs if it is piecewise linear with integer slopes and satisfies the harmonicity condition: For any point $v \in \Gamma$ and any two directions $\Vec{w_1}, \Vec{w_2}$ emanating from $w:=\varphi(v)$ we have $m_{\varphi, \Vec{w_1}}(v) = m_{\varphi,\Vec{w_2}}(v).$
}
\end{Definition}
Let $\varphi: \Gamma_1 \to \Gamma_2$ be a harmonic map and $v \in \Gamma$. Then, $m_{\varphi}(v) := m_{\varphi, \Vec{w_1}}(v) = m_{\varphi,\Vec{w_2}}(v)$ for any two directions $\Vec{w_1}, \Vec{w_2}$ emanating from $\varphi(v)$ is said to be the \textit{local degree} of $\varphi$ at $v$. The \textit{degree} of a non-constant harmonic map $\varphi: \Gamma_1 \to \Gamma_2$ is defined to be the sum of all local degrees of $\varphi$ at the pre-images under $\varphi$ of any point $w \in \Gamma'$ i.e., 
 $$\textnormal{deg }\varphi := \sum_{v \in \Gamma, \varphi(v)=w} m_{\varphi}(w) $$
for any $w \in \Gamma'$. The degree of $\varphi$ is independent of the choice of $w$. (see Section 2.4, \cite{K}).

\begin{Definition} 
\textnormal{A non-constant harmonic map $\varphi: \Gamma \to \Gamma'$ of metric graphs is said to be a \textit{tropical morphism} between metric graphs if the slopes of $\varphi$ along the edges of linearity are nonzero and the following inequality
$$(k - 2)  \geqslant  m_{\varphi}(v) \cdot ( l-2)$$ 
holds for all points $v \in \Gamma$, where $k$ is the valence of $v$, and $l$ is the valence of $w:=\varphi(v)$. The above inequality is known as the \textit{Riemann-Hurwitz condition}.  
}
\end{Definition}
\subsection{Tree gonality.}
Let $\Gamma$ be a metric graph, $T$ a metric tree, and let $v \in \Gamma$, $w \in T$ be two points such that $\textnormal{val}(w)=1$. Denote by $\Gamma'$ the quotient space of $\Gamma \sqcup T$ with respect to the equivalence relation $\sim$ that identifies $v$ with $w$. The metric space $\Gamma'$ is a metric graph, and we say that  $\Gamma'$ is obtained by \textit{grafting} the metric tree $T$ onto the point $v \in \Gamma$. In this article, we allow the inverse operation of grafting a metric tree onto a point of a metric graph, and we call it \textit{deleting} a metric tree onto a point of the metric graph.  
\begin{Definition}
\textnormal{
A \textit{tropical modification} of a metric graph $\Gamma$ is another metric graph $\Gamma'$ that is obtained by grafting or deleting a finite number of metric trees onto points of $\Gamma$. 
}
\end{Definition}
Given a tropical modification $\Gamma'$ of $\Gamma$ and a tropical morphism $\varphi: \Gamma' \to T$ of metric graphs, then there exists a tropical modification $\Gamma''$ (resp. $T'$)  of  $\Gamma'$ (resp. $T$) respectively and a tropical morphism $\varphi':\Gamma'' \to T'$ that extends $\varphi$ and has the same degree as $\varphi$ (\cite{FJ}). The following definition is the key definition in this article.  
\begin{Definition} \textnormal{The \textit{tree gonality} of a metric graph $\Gamma$, denoted by $\textnormal{tgon}(\Gamma)$, is defined as the minimum of degrees of all tropical morphisms from any tropical modification of $\Gamma$ to any metric tree. }
\end{Definition}

In order to study tree gonality and tropical morphisms of metric graphs, we consider the equivalence relation on metric graphs under tropical modification called \textit{tropical equivalence}. Metric graphs under tropical equivalence are said to be \textit{tropically equivalent}. \par 
First, we recall the notions of contracting and deleting an edge of a graph. Let $G = (V,E, \partial)$ be a graph and $e \in E$ with $\partial(e)=\left\lbrace v,w\right\rbrace$. \textit{Contracting} $G$ at the edge $e \in E$ yields the graph $G_1 = (V_1,E_1,\partial_1)$ where $V_1:=V/\sim$ where $\sim$ identifies $v$ with $w$, $E_1:=E\setminus \left\lbrace e \right\rbrace$ and $\partial_1 : E_1 \to P(V_1)$ given as follows: for $e' \in E_1$ such that $\partial(e')  = \left\lbrace v',w' \right\rbrace $ we define  $\partial_1(e') = \left\lbrace p(v'),p(w') \right\rbrace$, where $p: V \to V_1 $ is the quotient map. \textit{Deleting} the edge $e \in E$ yields the graph $G':= (V, E \setminus \left\lbrace e \right\rbrace, \partial|_{E \setminus \left\lbrace e \right\rbrace})$.\par
Next, we work with the notion of dangling edges which is due to \cite{JA}. Note that we regard a singleton graph (a graph without an edge) as a tree. 
\begin{Definition}\textnormal{
Let $G$ be a connected graph. An edge $e \in E(G)$ is said to be \textit{dangling} if deleting $e$ gives a graph with two connected components and one of them is a tree. }
\end{Definition}
Let $\Gamma$ be a metric graph with model $(G,l)$. Assume that $g(\Gamma) \geq 2$. Denote by $\tilde{G}$ the graph obtained by successively contracting the dangling edges of $G$, and let $\tilde{l}$ be a length map on $\tilde{G}$ given as the restriction of $l$ on $E(\tilde{G})$. Let $\tilde{\Gamma}$ be metric graph which is the realization of $(\tilde{G}, \tilde{l})$. Then, the metric graph $\Gamma$ is a tropical modification of $\tilde{\Gamma}$, and note that by construction, $\tilde{\Gamma}$ satisfies the following property:  $\tilde{\Gamma}$ is the unique metric graph tropically equivalent to $\Gamma$ whose essential model $(\mathcal{E},l_{\mathcal{E}})$ has valency at least $3$ i.e., every vertex point has valence at least three. 
\subsection{Hyperelliptic metric graphs.}

We first recall the basic theory of divisors on metric graphs (\cite{MCH}, \cite{BN1}). Let $\Gamma$ be a metric graph. An element of the free abelian group $\textnormal{Div}(\Gamma)$ generated by points of $\Gamma$ is called a \textit{divisor} on $\Gamma$. If
$$ D = \sum_{v \in \Gamma} D(v) \cdot v$$ 
is a divisor in $\Gamma$, then define the \textit{degree} of $D$ to be
$$ \textnormal{deg} (D) := \sum_{x \in \Gamma} D(v) \in \mathbb{Z} $$
Denote by $\textnormal{Div}^{0}(\Gamma)$ the subgroup of divisors of degree $0$. A function $f : \Gamma \to \mathbb{R}$ is called \textit{rational function} on $\Gamma$ if it is continuous, piecewise-linear with integer slopes along its domains of linearity. We denote by $\textnormal{Rat}(\Gamma)$ the set of rational functions on $\Gamma$. For $f \in \textnormal{Rat}(\Gamma)$ and a point $v$ in $\Gamma$, the sum of the outgoing slopes of $f$ at $v$ is denoted by $\textnormal{ord}_{v}(f)$. This sum is $0$ except for all but finitely many
points of $\Gamma$, and therefore,
$$\textnormal{div} (f) := \sum_{v \in \Gamma} \textnormal{ord}_{v}(f)$$
is a divisor on $\Gamma$. The set of principal divisors on $\Gamma$ is defined to be $\textnormal{Prin} (\Gamma) := \left\lbrace \textnormal{div} (f) \hspace{0.1cm} | \hspace{0.1cm} f \in \textnormal{Rat}(\Gamma) \right\rbrace $. Note that $\textnormal{Prin} (\Gamma)$ is a subgroup of $\textnormal{Div}^{0} (\Gamma)$. Two divisors $D$ and $D'$ are said to be \textit{linearly equivalent}, and we write $D \sim D'$, if $D - D' \in \textnormal{Prin} (\Gamma)$. A divisor $D = \sum_{v \in \Gamma} D(v) \cdot v \in \textnormal{Div}(\Gamma)$ is said to be \textit{effective}, and we write $D \geqslant 0$,  if $D(v) \geqslant 0$ for all $v \in \Gamma$.  Denote by $\textnormal{Div}^{k}_{+}(\Gamma)$ the set of all effective divisors with degree $k$.  For  a divisor $D \in \textnormal{Div}(\Gamma)$ a \textit{complete linear system} $|D|$ is defined to be
 $|D|:= \left\lbrace D' \in \textnormal{Div}(\Gamma) \hspace{0.1cm} | \hspace{0.1cm} D' \geqslant 0, D' \sim D \right\rbrace.$
The \textit{rank} of a divisor $D$ is defined to be $-1$ if $|D|=\emptyset$, and 
$$\textnormal{max} \left\lbrace k \in \mathbb{Z} \hspace{0.1cm} | \hspace{0.1cm}  \forall D' \in \textnormal{Div}^{k}_{+} (\Gamma) \textnormal{ we have }  |D - D'| \not= \emptyset \right\rbrace $$
if $|D| \not= \emptyset$. The rank of the divisor $D$ is simply denoted by $r(D)$. In the literature, there exists a notion of a hyperelliptic metric graph. For example in \cite{MCH}, a metric graph $\Gamma$ is said to be hyperelliptic if there exists a divisor $D \in \textnormal{Div} (\Gamma)$ such that $\textnormal{deg} (D) = 2$ and $r(D)=1$. In this article, we give a definiton of hyperelliptic metric graphs in terms of tropical morphisms and their tree gonality and which is different to the one given in \cite{MCH}. 
\begin{Definition}
\textnormal{A metric graph $\Gamma$ is said to be \textit{hyperelliptic} if there exists a tropical morphism from $\Gamma$ to a metric tree with degree $\textnormal{tgon}(\Gamma)$.
}
\end{Definition}
One of our goals in this article is to investigate genus $3$ nonhyperelliptic metric graphs $\Gamma$ with tree gonality $3$. Note that if $\Gamma$ is hyperelliptic in the sense of Kawaguchi-Yamaki (\cite{KJ}) that does not imply that $\Gamma$ is hyperelliptic in our sense. For example, the metric graph $\Gamma$ in Figure 25 is hyperelliptic in the sense of Kawaguchi-Yamaki but is not hyperelliptic in our sense. This is because the harmonic map coming from the unique hyperelliptic involution $\iota$ (Theorem 3.5, \cite{KJ}) does not satisfy the Riemann-Hurwitz condition.
\newpage 
\section{Construction of tropical morphisms}
The main result in this article is the constructive solution given to the Problem 1 stated below. Before we do that, we give the following lemma, which will be useful to construct tropical morphisms. 
\begin{lemma} \textnormal{
Let $\Gamma=(G,l)$, $T=(H,m)$ be two metric graphs where $H$ does not have multiple edges and $\psi:V(G) \to V(H)$ a map on the set of vertices. Suppose that for any $v,w \in V(G)$ that are the endpoints of some non-loop edge $e \in E(G)$, we have $\psi(v) = \psi(w)$, or $\psi(v)$ and $\psi(w)$ are endpoints of some edge $e' \in H$. Then, there exists a unique continuous map $\varphi: \Gamma \to T$ such that $\varphi|_{V(G)} = \psi$ and $\varphi$ is linear over each edge $e$ in $G$.} \end{lemma} 
\textit{Proof. } If $e \in E(G)$ is an edge with endpoints $v, w$ such that $\psi(v)=\psi(w)$, then take $\varphi_e : e \to T$ to be the constant map on $e$ with image $\psi(v)=\psi(w)$. In the case when $e \in E(G)$ is an edge with endpoints $v,w$ such that $\psi(v)$ and $\psi(w)$ are endpoints of some edge $e' \in H$, then choose $\varphi_e:e \to e'$ to be the linear map with slope $m(e')/l(e)$. Now, we take  $\varphi:\Gamma \to T$ to be the unique continuous map such that $\varphi|_e=\varphi_e$ for all edges $e \in E(G)$. $\square$ 

\begin{flushleft}
    \textbf{Problem 1.} 
Let $\Gamma$ be a genus $3$ metric graph with tree gonality $3$ which is not hyperelliptic. Construct a tropical modification $\Gamma'$ of $\Gamma$, a metric tree $T$ and a tropical morphism $\varphi:\Gamma' \to T$ of degree $3$.
\end{flushleft} 
\textit{Solution of Problem 1.} 
Consider genus $3$ nonhyperelliptic metric graphs with tree gonality $3$ up to tropical equivalence. There is a complete list (up to tropical equivalence) of genus $3$ metric graphs (Figure 4, \cite{C}), and also a complete list of genus $3$ hyperelliptic metric graphs (the tropical hyperelliptic curves of genus $3$ with unmarked vertices in Figure 2, \cite{MCH}). Note that there is a hyperelliptic metric graph in the latter list, namely the one in Figure 25, which is not hyperelliptic in our sense. Based on this, now it is enough to make the constructions for the tropically equivalent metric graphs $\Gamma$ whose essential model $(G,l)$ has valency at least $3$. They are depicted in Figures 1,5,7,9,\ldots,25. We divide the constructions into four cases depending on the number bridges (edges of a connected graph whose deletion increases its number of connected components) that the essential model $(G,l)$ possesses. \par  
\textbf{Case 1.} If the metric graph $\Gamma$ has no bridges, then $\Gamma$ is one of the metric graphs given in Figure 1, 5, 7, 9, 11, or 13. \par
\textit{Solution of Case 1.} \par 
\textbf{Case 1.1.} Consider the metric graph $\Gamma$ whose essential model $(G,l)$ is given in Figure 1, where the graph $G=(V,E,\partial)$ is given by $V = \left\lbrace v_1,v_2,v_3,v_4 \right\rbrace$, $E=\left\lbrace e_1,e_2,\ldots, e_6  \right\rbrace$, and $\partial(e_1)= \left\lbrace v_1,v_2 \right\rbrace $,
$\partial(e_2)= \left\lbrace v_1,v_3 \right\rbrace $, $\partial(e_3)= \left\lbrace v_4,v_1 \right\rbrace $, $\partial(e_5)= \left\lbrace v_2,v_3 \right\rbrace $, and $\partial(e_6)= \left\lbrace v_2, v_4 \right\rbrace $. The length map $l$ on $E$ is defined by assigning $ e_1 \mapsto a,  e_2  \mapsto b,  e_3 \mapsto c,   e_4 \mapsto d,   e_5 \mapsto e$ and $  e_6 \mapsto f$, where $a,b,c,d,e$, and $f$ are real positive numbers. \par 
\begin{figure}[H]
    \centering
    \begin{tikzpicture}[line cap=round,line join=round,>=triangle 45,x=1.0cm,y=1.0cm]
\definecolor{wrwrwr}{rgb}{0.3803921568627451,0.3803921568627451,0.3803921568627451}
\clip(-1.5,-1.5) rectangle (3.,3.);
\draw (-1.0139485196985758,0.9984976653834488)-- (2.005135575749092,-0.9940978376120061);
\draw (2.005135575749092,-0.9940978376120061)-- (2.,1.);
\draw (2.,1.)-- (0.027635493230869612,2.0249862578356526);
\draw (0.027635493230869612,2.0249862578356526)-- (-1.0139485196985758,0.9984976653834488);
\draw [shift={(1.0163855344899808,0.5154442101118233)}] plot[domain=-0.990914194401082:2.150678459188711,variable=\t]({1.*1.8045342440464094*cos(\t r)+0.*1.8045342440464094*sin(\t r)},{0.*1.8045342440464094*cos(\t r)+1.*1.8045342440464094*sin(\t r)});
\draw (-1.0139485196985758,0.9984976653834488)-- (2.,1.);
\begin{scriptsize}
\draw [fill=black] (-1.0139485196985758,0.9984976653834488) circle (0.5pt);
\draw[color=black] (-1.257745383827365,0.9267194287657506) node {$v_1$};
\draw [fill=black] (2.005135575749092,-0.9940978376120061) circle (0.5pt);
\draw[color=black] (2.23949301797127,-1.1180368948665094) node {$v_2$};
\draw [fill=black] (2.,1.) circle (0.5pt);
\draw[color=black] (2.2535947857204586,1.1523477127527588) node {$v_3$};
\draw [fill=black] (0.027635493230869612,2.0249862578356526) circle (0.5pt);
\draw[color=black] (-0.12960396389232126,2.3933032746813026) node {$v_4$};
\draw[color=black] (0.4274158622006065,-0.4341011590308914) node {$a$};
\draw[color=black] (1.7247784951259066,0.11586778318744059) node {$e$};
\draw[color=black] (1.2735219271518892,1.6106551646013687) node {$d$};
\draw[color=black] (-0.686623789985249,1.5683498613538047) node {$c$};
\draw[color=black] (2.5849863278263774,1.7939781453408128) node {$f$};
\draw[color=black] (0.46972116544817066,1.2017038998749168) node {$b$};
\draw[color=wrwrwr] (-2.8582960233602077,4.289991036947089) node {$Curve$};
\end{scriptsize}
\end{tikzpicture}
    \begin{center}
        \textnormal{Figure 1. The essential model $(G,l)$ of $\Gamma$}
    \end{center}
\end{figure}

Choose any vertex, say $v_1 \in V(G)$, and without loss of generality, assume that $c\leqslant b\leqslant a$. It is enough to consider the following three subcases: (A) $c<b \leqslant a$, (B) $c=b<a$, and (C) $c=b=a$. We give the constructions for each subcase separately as follows.  
\par 
\textbf{Case 1.1.A.} Let $(G_1,l_1)$ be another model of $\Gamma$ given in Figure 1.1, where the graph $G_1=(V_1,E_1,\partial_1)$ is obtained by subdividing the edges $ e_i \in E$ into $ e_i',e_i'',e_i''' (i=1,2)$, and $e_j \in E $ into $e_j',e_j'' (j=4,5,6)$ with orientation given by:  
\begin{center}
\begin{tabular}{|c|c|c|}
\hline
$ \partial_1(e_1') = \left\lbrace v_1,v_6 \right\rbrace $ & $\partial_1(e_1'') = \left\lbrace v_6, v_4\right\rbrace$ & $\partial_1(e_1''') = \left\lbrace v_5,v_2 \right\rbrace $  \\
\hline 
$ \partial_1(e_2'') = \left\lbrace v_8, v_7 \right\rbrace$ & $\partial_1(e_2''') = \left\lbrace v_7,v_3 \right\rbrace$ & $ 
 \partial_1(e_4') =\left\lbrace v_4,v_9 \right\rbrace$  \\ 
 \hline 
$ \partial_1(e_4'') = \left\lbrace v_9,v_3 \right\rbrace $ & $\partial_1(e_5') = \left\lbrace v_2,v_{10} \right\rbrace $ & $\partial_1(e_5'') = \left\lbrace v_3,v_{10} \right\rbrace$ \\
\hline 
$ \partial_1(e_2') = \left\lbrace v_1,v_8 \right\rbrace$ & $\partial_1(e_6'') = \left\lbrace v_4,v_{11} \right\rbrace$ & $\partial_1(e_6') = \left\lbrace v_2,v_{11} \right\rbrace$  \\ 
\hline 
\end{tabular}
\end{center}
and length map $l_1$, which is equal to $l$ on $E \setminus \left\lbrace  e_1,e_2,e_4,e_5,e_6 \right\rbrace$, whereas on $\left\lbrace  e_1,e_2,e_4,e_5,e_6 \right\rbrace$  it is equal to 

\begin{center}
\begin{tabular}{ll}
$l_1(e_1') = l_1(e_1'') = (a-c)/2$ & $l_1(e_4') = l_1(e_4'') = d/2$  \\ 
$l_1(e_2') = l_1(e_2'') = (b-c)/2$ & $l_1(e_5') = l_1(e_5'') = e/2$  \\
 $l_1(e_6') = l_1(e_1''') = l_1(e_2''') = c$ & $l_1(e_6'') = f/2$. \\
\end{tabular} 
\end{center}

\begin{figure}[H]
    \centering
    \begin{tikzpicture}[line cap=round,line join=round,>=triangle 45,x=1.0cm,y=1.0cm]
 \definecolor{uuuuuu}{rgb}{0.26666666666666666,0.26666666666666666,0.26666666666666666}
\definecolor{wrwrwr}{rgb}{0.3803921568627451,0.3803921568627451,0.3803921568627451}
\definecolor{rvwvcq}{rgb}{0.08235294117647059,0.396078431372549,0.7529411764705882}
\definecolor{dtsfsf}{rgb}{0.8274509803921568,0.1843137254901961,0.1843137254901961}
\definecolor{sexdts}{rgb}{0.1803921568627451,0.49019607843137253,0.19607843137254902}
\clip(-3.,-0.8) rectangle (2.5,2.5);
\draw (-1.,0.)-- (1.,0.);
\draw (-1.,1.)-- (1.,1.);
\draw (-1.,2.)-- (1.,2.);
\draw (-1.,-1.5)-- (1.,-1.5);
\draw [color=sexdts] (-1.,-1.5)-- (-2.,-2.);
\draw [color=dtsfsf] (-1.,-1.5)-- (-1.5,-1.5);
\draw [color=rvwvcq] (1.,-1.5)-- (1.2894301310058747,-1.012638385870342);
\draw [color=sexdts] (1.,-1.5)-- (2.,-1.5);
\draw [color=dtsfsf] (1.,-1.5)-- (1.5,-2.);
\draw [shift={(1.,1.5)},color=dtsfsf]  plot[domain=-1.5707963267948966:1.5707963267948966,variable=\t]({1.*0.5*cos(\t r)+0.*0.5*sin(\t r)},{0.*0.5*cos(\t r)+1.*0.5*sin(\t r)});
\draw [shift={(1.,1.)},color=sexdts]  plot[domain=-1.5707963267948966:1.5707963267948966,variable=\t]({1.*1.*cos(\t r)+0.*1.*sin(\t r)},{0.*1.*cos(\t r)+1.*1.*sin(\t r)});
\draw [shift={(0.7237164077379609,0.5)},color=rvwvcq]  plot[domain=-1.0659842856832995:1.0659842856832993,variable=\t]({1.*0.5712553048797154*cos(\t r)+0.*0.5712553048797154*sin(\t r)},{0.*0.5712553048797154*cos(\t r)+1.*0.5712553048797154*sin(\t r)});
\draw (2.9059569101891936,-1.3536469272472165) node[anchor=north west] {$T=(T', t')$};
\draw (2.9059569101891936,1.3092357812953126) node[anchor=north west] {$\Gamma' = (G',l')$};
\draw [->,line width=0.4pt] (3.063701734239714,0.5091002477485305) -- (3.063701734239714,-0.9908997522514695);
\draw [shift={(-1.,1.5)},color=dtsfsf]  plot[domain=1.5707963267948966:4.71238898038469,variable=\t]({1.*0.5*cos(\t r)+0.*0.5*sin(\t r)},{0.*0.5*cos(\t r)+1.*0.5*sin(\t r)});
\draw [shift={(-1.,1.)},color=sexdts]  plot[domain=1.5707963267948966:4.71238898038469,variable=\t]({1.*1.*cos(\t r)+0.*1.*sin(\t r)},{0.*1.*cos(\t r)+1.*1.*sin(\t r)});
\begin{scriptsize}
\draw [fill=black] (-1.,0.) circle (0.5pt);
\draw[color=black] (-0.8921547425214919,-0.3258676362308018) node {$v_5$};
\draw [fill=black] (1.,0.) circle (0.5pt);
\draw[color=black] (0.8363831560061263,-0.3165241881306526) node {$v_2$};
\draw[color=black] (0.05620523964366082,0.17400683712718173) node {$c$};
\draw [fill=black] (1.,1.) circle (0.5pt);
\draw[color=black] (0.8083528117056784,1.243831644594268) node {$v_3$};
\draw[color=black] (0.0655486877438101,1.183099231943298) node {$c$};
\draw [fill=black] (-1.,2.) circle (0.5pt);
\draw[color=black] (-0.9855892235229847,2.2716109356106826) node {$v_1$};
\draw [fill=black] (1.,2.) circle (0.5pt);
\draw[color=black] (0.8457266041062756,2.2622674875105333) node {$v_4$};
\draw[color=black] (0.05620523964366082,2.1735047305591157) node {$c$};
\draw [fill=black] (-1.,-1.5) circle (0.5pt);
\draw[color=black] (-0.9762457754228354,-1.8208193322546777) node {$w_1$};
\draw [fill=black] (1.,-1.5) circle (0.5pt);
\draw[color=black] (0.8831003965068728,-1.79278898795423) node {$w_0$};
\draw [fill=sexdts] (-2.,-2.) circle (0.5pt);
\draw[color=sexdts] (-2.2562981651432876,-1.624606922151544) node {$w_6$};
\draw [fill=dtsfsf] (-1.5,-1.5) circle (0.5pt);
\draw[color=dtsfsf] (-1.639630590533435,-1.2041517576448288) node {$w_8$};
\draw [fill=rvwvcq] (1.2894301310058747,-1.012638385870342) circle (0.5pt);
\draw[color=rvwvcq] (1.5838590040180691,-0.9238483146403522) node {$w_{10}$};
\draw [fill=sexdts] (2.,-1.5) circle (0.5pt);
\draw[color=sexdts] (2.3033045077295644,-1.4003641677479626) node {$w_{11}$};
\draw [fill=dtsfsf] (1.5,-2.) circle (0.5pt);
\draw[color=dtsfsf] (1.4717376268162776,-2.203900704360796) node {$w_9$};
\draw [color=dtsfsf] (1.5,1.5)-- ++(-1.0pt,-1.0pt) -- ++(2.0pt,2.0pt) ++(-2.0pt,0) -- ++(2.0pt,-2.0pt);
\draw[color=dtsfsf] (1.6632783128693382,1.4447157787474763) node {$d$};
\draw [color=rvwvcq] (1.294962530147754,0.5032389790046032)-- ++(-1.0pt,-1.0pt) -- ++(2.0pt,2.0pt) ++(-2.0pt,0) -- ++(2.0pt,-2.0pt);
\draw[color=rvwvcq] (1.4857527989665018,0.3608757991301662) node {$e$};
\draw [color=sexdts] (2.,1.)-- ++(-1.0pt,-1.0pt) -- ++(2.0pt,2.0pt) ++(-2.0pt,0) -- ++(2.0pt,-2.0pt);
\draw[color=sexdts] (2.1491376140771012,0.9868868218401643) node {$f$};
\draw[color=dtsfsf] (1.2474948724126949,1.5521654318991924) node {$v_9$};
\draw[color=sexdts] (1.7426976217206072,1.0476192344911341) node {$v_{11}$};
\draw[color=rvwvcq] (0.9858783256085148,0.5243861408827777) node {$v_{10}$};
\draw [fill=rvwvcq] (1.294962530147754,2.498664332229594) circle (0.5pt);
\draw[color=black] (3.4011596594971056,-0.21841798307908572) node {$\varphi$};
\draw[color=wrwrwr] (-3.8586995143188907,3.8646688366861253) node {$Curve$};
\draw [fill=uuuuuu] (-1.,1.) circle (0.5pt);
\draw[color=uuuuuu] (-0.8921547425214919,0.7392854471862099) node {$v_7$};
\draw[color=dtsfsf] (-1.0463216361739551,1.4820895711480733) node {$b-c$};
\draw[color=sexdts] (-2.5225864359975425,0.9495130294395673) node {$a-c$};
\draw [color=dtsfsf] (-1.5,1.5)-- ++(-1.0pt,-1.0pt) -- ++(2.0pt,2.0pt) ++(-2.0pt,0) -- ++(2.0pt,-2.0pt);
\draw[color=dtsfsf] (-1.6583174867337334,1.3092357812953126) node {$v_8$};
\draw [color=sexdts] (-2.,1.)-- ++(-1.0pt,-1.0pt) -- ++(2.0pt,2.0pt) ++(-2.0pt,0) -- ++(2.0pt,-2.0pt);
\draw[color=sexdts] (-1.7984692082359728,0.9168109610890451) node {$v_6$};
\end{scriptsize}
\end{tikzpicture}
    \begin{center}
        \textnormal{Figure 1.1. The model $(G_1,l_1)$ of $\Gamma$}
    \end{center}
\end{figure}
\newpage 
Let $\Gamma'$ be the the tropical modification of $\Gamma$ with model $(G',l')$ in Figure 1.2, where the graph $G'$ is given by its vertex set $V(G') = V_1 \cup \left\lbrace v_6',v_8',v_9',v_{10}',v_{11}' \right\rbrace$, and edge set $E(G')=E_1 \cup \left\lbrace v_2v_9', v_3v_{11}' ,v_4v_{10}', v_5v_8', v_7v_6' \right\rbrace $. The length map $l'$ on $G'$ is defined by $l' = l_1$ on $E_1$, and

\begin{center}
\begin{equation*}
\begin{aligned}[c]
l'(v_2v_9')&= \frac{d}{2}\\
l'(v_3v_{11}')&= \frac{f}{2}\\
l'(v_4v_{10}')&= \frac{e}{2}
\end{aligned}
\qquad
\begin{aligned}[c]
l'(v_5v_8')&= \frac{b-c}{2}\\
l'(v_7v_6')&= \frac{a-c}{2}.\\
\end{aligned}
\end{equation*}
\end{center}
\begin{figure}[H]
    \centering
   \begin{tikzpicture}[line cap=round,line join=round,>=triangle 45,x=1.0cm,y=1.0cm]
 \definecolor{uuuuuu}{rgb}{0.26666666666666666,0.26666666666666666,0.26666666666666666}
\definecolor{wrwrwr}{rgb}{0.3803921568627451,0.3803921568627451,0.3803921568627451}
\definecolor{rvwvcq}{rgb}{0.08235294117647059,0.396078431372549,0.7529411764705882}
\definecolor{dtsfsf}{rgb}{0.8274509803921568,0.1843137254901961,0.1843137254901961}
\definecolor{sexdts}{rgb}{0.1803921568627451,0.49019607843137253,0.19607843137254902}
\clip(-3.,-0.8) rectangle (2.5,3.);
\draw (-1.,0.)-- (1.,0.);
\draw (-1.,1.)-- (1.,1.);
\draw (-1.,2.)-- (1.,2.);
\draw (-1.,-1.5)-- (1.,-1.5);
\draw [color=sexdts] (-1.,-1.5)-- (-2.,-2.);
\draw [color=dtsfsf] (-1.,-1.5)-- (-1.5,-1.5);
\draw [color=rvwvcq] (1.,-1.5)-- (1.2894301310058747,-1.012638385870342);
\draw [color=sexdts] (1.,-1.5)-- (2.,-1.5);
\draw [color=dtsfsf] (1.,-1.5)-- (1.5,-2.);
\draw [shift={(1.,1.5)},color=dtsfsf]  plot[domain=-1.5707963267948966:1.5707963267948966,variable=\t]({1.*0.5*cos(\t r)+0.*0.5*sin(\t r)},{0.*0.5*cos(\t r)+1.*0.5*sin(\t r)});
\draw [shift={(1.,1.)},color=sexdts]  plot[domain=-1.5707963267948966:1.5707963267948966,variable=\t]({1.*1.*cos(\t r)+0.*1.*sin(\t r)},{0.*1.*cos(\t r)+1.*1.*sin(\t r)});
\draw [shift={(0.7237164077379609,0.5)},color=rvwvcq]  plot[domain=-1.0659842856832995:1.0659842856832993,variable=\t]({1.*0.5712553048797154*cos(\t r)+0.*0.5712553048797154*sin(\t r)},{0.*0.5712553048797154*cos(\t r)+1.*0.5712553048797154*sin(\t r)});
\draw [color=sexdts] (1.,1.)-- (2.,0.5);
\draw [color=rvwvcq] (1.,2.)-- (1.294962530147754,2.498664332229594);
\draw [color=dtsfsf] (1.,0.)-- (1.5,-0.5);
\draw (2.9059569101891936,-1.3536469272472165) node[anchor=north west] {$T=(T', t')$};
\draw (2.9059569101891936,1.3092357812953126) node[anchor=north west] {$\Gamma' = (G',l')$};
\draw [->,line width=0.4pt] (3.063701734239714,0.5091002477485305) -- (3.063701734239714,-0.9908997522514695);
\draw [shift={(-1.,1.5)},color=dtsfsf]  plot[domain=1.5707963267948966:4.71238898038469,variable=\t]({1.*0.5*cos(\t r)+0.*0.5*sin(\t r)},{0.*0.5*cos(\t r)+1.*0.5*sin(\t r)});
\draw [shift={(-1.,1.)},color=sexdts]  plot[domain=1.5707963267948966:4.71238898038469,variable=\t]({1.*1.*cos(\t r)+0.*1.*sin(\t r)},{0.*1.*cos(\t r)+1.*1.*sin(\t r)});
\draw [color=dtsfsf] (-1.,0.)-- (-1.5,0.);
\draw [color=sexdts] (-1.,1.)-- (-2.,0.5);
\begin{scriptsize}
\draw [fill=black] (-1.,0.) circle (0.5pt);
\draw[color=black] (-0.8921547425214919,0.25342614597845015) node {$v_5$};
\draw [fill=black] (1.,0.) circle (0.5pt);
\draw[color=black] (0.8176962598058277,0.25342614597845015) node {$v_2$};
\draw [fill=black] (1.,1.) circle (0.5pt);
\draw[color=black] (0.8083528117056784,1.243831644594268) node {$v_3$};
\draw [fill=black] (-1.,2.) circle (0.5pt);
\draw[color=black] (-0.9855892235229847,2.2716109356106826) node {$v_1$};
\draw [fill=black] (1.,2.) circle (0.5pt);
\draw[color=black] (0.8457266041062756,2.2622674875105333) node {$v_4$};
\draw [fill=black] (-1.,-1.5) circle (0.5pt);
\draw[color=black] (-0.9762457754228354,-1.8208193322546777) node {$w_1$};
\draw [fill=black] (1.,-1.5) circle (0.5pt);
\draw[color=black] (0.8831003965068728,-1.79278898795423) node {$w_0$};
\draw [fill=sexdts] (-2.,-2.) circle (0.5pt);
\draw[color=sexdts] (-2.2562981651432876,-1.624606922151544) node {$w_6$};
\draw [fill=dtsfsf] (-1.5,-1.5) circle (0.5pt);
\draw[color=dtsfsf] (-1.639630590533435,-1.2041517576448288) node {$w_8$};
\draw [fill=rvwvcq] (1.2894301310058747,-1.012638385870342) circle (0.5pt);
\draw[color=rvwvcq] (1.5838590040180691,-0.9238483146403522) node {$w_{10}$};
\draw [fill=sexdts] (2.,-1.5) circle (0.5pt);
\draw[color=sexdts] (2.3033045077295644,-1.4003641677479626) node {$w_{11}$};
\draw [fill=dtsfsf] (1.5,-2.) circle (0.5pt);
\draw[color=dtsfsf] (1.4717376268162776,-2.203900704360796) node {$w_9$};
\draw [color=dtsfsf] (1.5,1.5)-- ++(-1.0pt,-1.0pt) -- ++(2.0pt,2.0pt) ++(-2.0pt,0) -- ++(2.0pt,-2.0pt);
\draw [color=rvwvcq] (1.294962530147754,0.5032389790046032)-- ++(-1.0pt,-1.0pt) -- ++(2.0pt,2.0pt) ++(-2.0pt,0) -- ++(2.0pt,-2.0pt);
\draw [color=sexdts] (2.,1.)-- ++(-1.0pt,-1.0pt) -- ++(2.0pt,2.0pt) ++(-2.0pt,0) -- ++(2.0pt,-2.0pt);
\draw[color=dtsfsf] (1.2755252167131428,1.5521654318991924) node {$v_9$};
\draw[color=sexdts] (1.7426976217206072,1.1036799230920296) node {$v_{11}$};
\draw[color=rvwvcq] (0.9858783256085148,0.5243861408827777) node {$v_{10}$};
\draw [fill=sexdts] (2.,0.5) circle (0.5pt);
\draw[color=sexdts] (2.2472438191286686,0.5056992446824792) node {$v_{11}'$};
\draw [fill=rvwvcq] (1.294962530147754,2.498664332229594) circle (0.5pt);
\draw[color=rvwvcq] (1.5558286597176214,2.5519143786151597) node {$v_{10}'$};
\draw [fill=dtsfsf] (1.5,-0.5) circle (0.5pt);
\draw[color=dtsfsf] (1.7520410698207565,-0.4566759096328909) node {$v_9'$};
\draw[color=black] (3.4011596594971056,-0.21841798307908572) node {$\varphi$};
\draw[color=wrwrwr] (-3.8586995143188907,3.8646688366861253) node {$Curve$};
\draw [fill=uuuuuu] (-1.,1.) circle (0.5pt);
\draw[color=uuuuuu] (-0.8921547425214919,1.2344881964941188) node {$v_7$};
\draw [color=dtsfsf] (-1.5,1.5)-- ++(-1.0pt,-1.0pt) -- ++(2.0pt,2.0pt) ++(-2.0pt,0) -- ++(2.0pt,-2.0pt);
\draw[color=dtsfsf] (-1.2565492184273142,1.5988826723999385) node {$v_8$};
\draw [color=sexdts] (-2.,1.)-- ++(-1.0pt,-1.0pt) -- ++(2.0pt,2.0pt) ++(-2.0pt,0) -- ++(2.0pt,-2.0pt);
\draw[color=sexdts] (-2.2189243727426904,1.3372661255957603) node {$v_6$};
\draw [fill=dtsfsf] (-1.5,0.) circle (0.5pt);
\draw[color=dtsfsf] (-1.6676609348338827,-0.1857159147285634) node {$v_8'$};
\draw [fill=sexdts] (-2.,0.5) circle (0.5pt);
\draw[color=sexdts] (-2.1908940284422425,0.7299419990860606) node {$v_6'$};
\end{scriptsize}
\end{tikzpicture}
    \begin{center}
        \textnormal{Figure 1.2. The model $(G',l')$ of $\Gamma'$}
    \end{center}
\end{figure}

Let $T$ be the metric tree with model $(T',t')$ in Figure 1.3, where the tree $T'$ is given with its vertex set $V(T')= \left\lbrace w_0,w_1,w_6,w_8,w_9,w_{10},w_{11}\right\rbrace$, and edge set $E(T')= \left\lbrace w_0w_1,w_0w_9,w_0w_{10}, w_0w_{11},w_1w_6,w_1w_8 \right\rbrace$, whereas the length map $t'$ on $T$ is defined by

\begin{center}
\begin{equation*}
\begin{aligned}[c]
t'(w_0w_9) &= \frac{d}{2} \\
t'(w_0w_{10}) &=  \frac{e}{2}\\ 
t'(w_0w_{11}) &=  \frac{f}{2}
\end{aligned}
\qquad
\begin{aligned}[c]
t'(w_0w_1) &=  c\\
t'(w_1w_8) &=  \frac{b-c}{2}\\
t'(w_1w_6) &=  \frac{a-c}{2}.\\
\end{aligned}
\end{equation*}
\end{center}
\begin{figure}[H]
    \centering
    \begin{tikzpicture}[line cap=round,line join=round,>=triangle 45,x=1.0cm,y=1.0cm]
\definecolor{uuuuuu}{rgb}{0.26666666666666666,0.26666666666666666,0.26666666666666666}
\definecolor{wrwrwr}{rgb}{0.3803921568627451,0.3803921568627451,0.3803921568627451}
\definecolor{rvwvcq}{rgb}{0.08235294117647059,0.396078431372549,0.7529411764705882}
\definecolor{dtsfsf}{rgb}{0.8274509803921568,0.1843137254901961,0.1843137254901961}
\definecolor{sexdts}{rgb}{0.1803921568627451,0.49019607843137253,0.19607843137254902}
\clip(-2.5,-2.5) rectangle (2.8,-0.5);
\draw (-1.,0.)-- (1.,0.);
\draw (-1.,1.)-- (1.,1.);
\draw (-1.,2.)-- (1.,2.);
\draw (-1.,-1.5)-- (1.,-1.5);
\draw [color=sexdts] (-1.,-1.5)-- (-2.,-2.);
\draw [color=dtsfsf] (-1.,-1.5)-- (-1.5,-1.5);
\draw [color=rvwvcq] (1.,-1.5)-- (1.2894301310058747,-1.012638385870342);
\draw [color=sexdts] (1.,-1.5)-- (2.,-1.5);
\draw [color=dtsfsf] (1.,-1.5)-- (1.5,-2.);
\draw [shift={(1.,1.5)},color=dtsfsf]  plot[domain=-1.5707963267948966:1.5707963267948966,variable=\t]({1.*0.5*cos(\t r)+0.*0.5*sin(\t r)},{0.*0.5*cos(\t r)+1.*0.5*sin(\t r)});
\draw [shift={(1.,1.)},color=sexdts]  plot[domain=-1.5707963267948966:1.5707963267948966,variable=\t]({1.*1.*cos(\t r)+0.*1.*sin(\t r)},{0.*1.*cos(\t r)+1.*1.*sin(\t r)});
\draw [shift={(0.7237164077379609,0.5)},color=rvwvcq]  plot[domain=-1.0659842856832995:1.0659842856832993,variable=\t]({1.*0.5712553048797154*cos(\t r)+0.*0.5712553048797154*sin(\t r)},{0.*0.5712553048797154*cos(\t r)+1.*0.5712553048797154*sin(\t r)});
\draw (2.9059569101891936,-1.3536469272472162) node[anchor=north west] {$T=(T', t')$};
\draw (2.9059569101891936,1.3092357812953133) node[anchor=north west] {$\Gamma' = (G',l')$};
\draw [->,line width=0.4pt] (3.063701734239714,0.5091002477485305) -- (3.063701734239714,-0.9908997522514695);
\draw [shift={(-1.,1.5)},color=dtsfsf]  plot[domain=1.5707963267948966:4.71238898038469,variable=\t]({1.*0.5*cos(\t r)+0.*0.5*sin(\t r)},{0.*0.5*cos(\t r)+1.*0.5*sin(\t r)});
\draw [shift={(-1.,1.)},color=sexdts]  plot[domain=1.5707963267948966:4.71238898038469,variable=\t]({1.*1.*cos(\t r)+0.*1.*sin(\t r)},{0.*1.*cos(\t r)+1.*1.*sin(\t r)});
\begin{scriptsize}
\draw [fill=black] (-1.,0.) circle (0.5pt);
\draw[color=black] (-0.8921547425214919,-0.3258676362308013) node {$v_5$};
\draw [fill=black] (1.,0.) circle (0.5pt);
\draw[color=black] (0.8363831560061263,-0.3165241881306521) node {$v_2$};
\draw[color=black] (0.05620523964366082,0.1740068371271823) node {$c$};
\draw [fill=black] (1.,1.) circle (0.5pt);
\draw[color=black] (0.8083528117056784,1.2438316445942688) node {$v_3$};
\draw[color=black] (0.0655486877438101,1.1830992319432987) node {$c$};
\draw [fill=black] (-1.,2.) circle (0.5pt);
\draw[color=black] (-0.9855892235229847,2.2716109356106835) node {$v_1$};
\draw [fill=black] (1.,2.) circle (0.5pt);
\draw[color=black] (0.8457266041062756,2.2622674875105346) node {$v_4$};
\draw[color=black] (0.05620523964366082,2.1735047305591166) node {$c$};
\draw [fill=black] (-1.,-1.5) circle (0.5pt);
\draw[color=black] (-0.9855892235229847,-1.7367282993533344) node {$w_1$};
\draw [fill=black] (1.,-1.5) circle (0.5pt);
\draw[color=black] (0.9111307408073206,-1.7367282993533344) node {$w_0$};
\draw [fill=sexdts] (-2.,-2.) circle (0.5pt);
\draw[color=sexdts] (-2.2562981651432876,-2.0170317423578115) node {$w_6$};
\draw [fill=dtsfsf] (-1.5,-1.5) circle (0.5pt);
\draw[color=dtsfsf] (-1.7984692082359728,-1.381677271547664) node {$w_8$};
\draw [fill=rvwvcq] (1.2894301310058747,-1.012638385870342) circle (0.5pt);
\draw[color=rvwvcq] (1.5838590040180691,-0.9238483146403518) node {$w_{10}$};
\draw [fill=sexdts] (2.,-1.5) circle (0.5pt);
\draw[color=sexdts] (2.3033045077295644,-1.5311724411500516) node {$w_{11}$};
\draw [fill=dtsfsf] (1.5,-2.) circle (0.5pt);
\draw[color=dtsfsf] (1.5277983154171735,-2.203900704360796) node {$w_9$};
\draw [color=dtsfsf] (1.5,1.5)-- ++(-1.0pt,-1.0pt) -- ++(2.0pt,2.0pt) ++(-2.0pt,0) -- ++(2.0pt,-2.0pt);
\draw[color=dtsfsf] (1.6632783128693382,1.4447157787474771) node {$d$};
\draw [color=rvwvcq] (1.294962530147754,0.5032389790046032)-- ++(-1.0pt,-1.0pt) -- ++(2.0pt,2.0pt) ++(-2.0pt,0) -- ++(2.0pt,-2.0pt);
\draw[color=rvwvcq] (1.4857527989665018,0.3608757991301668) node {$e$};
\draw [color=sexdts] (2.,1.)-- ++(-1.0pt,-1.0pt) -- ++(2.0pt,2.0pt) ++(-2.0pt,0) -- ++(2.0pt,-2.0pt);
\draw[color=sexdts] (2.1491376140771012,0.9868868218401651) node {$f$};
\draw[color=dtsfsf] (1.2755252167131428,1.5521654318991933) node {$v_9'$};
\draw[color=sexdts] (1.7426976217206072,1.047619234491135) node {$v_{11}$};
\draw[color=rvwvcq] (0.9858783256085148,0.5243861408827782) node {$v_{10}$};
\draw [fill=rvwvcq] (1.294962530147754,2.498664332229594) circle (0.5pt);
\draw[color=black] (3.4011596594971056,-0.21841798307908522) node {$\varphi$};
\draw[color=wrwrwr] (-3.8586995143188907,3.8646688366861266) node {$Curve$};
\draw [fill=uuuuuu] (-1.,1.) circle (0.5pt);
\draw[color=uuuuuu] (-0.8921547425214919,1.2344881964941194) node {$v_7$};
\draw[color=dtsfsf] (-1.0463216361739551,1.482089571148074) node {$b-c$};
\draw[color=sexdts] (-2.5225864359975425,0.9495130294395682) node {$a-c$};
\draw [color=dtsfsf] (-1.5,1.5)-- ++(-1.0pt,-1.0pt) -- ++(2.0pt,2.0pt) ++(-2.0pt,0) -- ++(2.0pt,-2.0pt);
\draw[color=dtsfsf] (-1.6583174867337334,1.5241350875987456) node {$v_8$};
\draw [color=sexdts] (-2.,1.)-- ++(-1.0pt,-1.0pt) -- ++(2.0pt,2.0pt) ++(-2.0pt,0) -- ++(2.0pt,-2.0pt);
\draw[color=sexdts] (-1.7984692082359728,0.9168109610890458) node {$v_6$};
\end{scriptsize}
\end{tikzpicture}
    \begin{center}
        \textnormal{Figure 1.3. The model $(T',t')$ of $T$}
    \end{center}
\end{figure}

Let $\psi:V(G') \to V(T)$ the map on the set of vertices given by $v_2, v_3,v_4 \mapsto w_0$, $v_1,v_7,v_5 \mapsto w_1$, and $v_i, v_i' \mapsto w_i $ for $i=6,8,9,10,11$. Then, the map $\psi$ satisfies the condition in Lemma 3.1, and so, there exist a unique continuous map $\varphi: \Gamma' \to T$ such that $\varphi|_{V(G')} = \psi$, and $\varphi$ is linear on each edge $e' \in E(G')$ with slope $t'(e)/l'(e')$, where $e = \varphi(e') \in E(T)$ with endpoints $\psi(v)$ and $\psi(w)$. The map $\varphi$ given in Figure 2. By construction, the models $(G',l')$ and $(T',t')$ satisfy the condition in Definition 2.4, and therefore, $\varphi$ is a piecewise linear function. From our choice of length maps $t', l'$, the slope of $\varphi|{e'}$ is equal to $1$ for all edges $e'$. Thus, $\varphi$ has non-zero integer slopes along its edges of linearity. It is remaining to show that the map $\varphi$ satisfies (i) the harmonicity condition and (ii) the Riemann-Hurtswitz condition on every point $v \in \Gamma'$.  
\begin{itemize}
\item[\textbf{(i)}]
Assume that $v \in \Gamma^{'}$ is a vertex point, say $v= v_1 \in V(G')$. Then, for all the directions $\Vec{w}$ at $\varphi(v)=w_1 $, we have that $m_{\varphi,\Vec{w}}(v_1)=1$. We check that the harmonicity condition holds on every other vertex point in a similar fashion, and this checking process terminates because the vertex set is finite. Whenever $v$ is not a vertex point, say $v \in \textnormal{int}(e)$ for some edge $e \in E(G')$, we have that $ \varphi(v)\in \textnormal{int}(e') $ where $e'=\varphi(e)$. Consider the new vertex sets on $\Gamma'$ and $T$ by adding $v$ and $w$ respectively. There are only two directions $\Vec{w}_1$ and $\Vec{w}_2$ at $\varphi(v)$ because $\textnormal{val}(\varphi(v))=2$. The slopes of $\varphi$ at $v$ with directions $\Vec{w}_1$ and $\Vec{w}_2$ at $\varphi(v)$ are equal to the slope of the same linear map $\varphi|_{e}$ i.e.,  $m_{\varphi, \Vec{w}_1}(v) = m_{\varphi, \Vec{w}_2}(v)$, and so, we get that $\varphi$ is a harmonic map. Its degree is $3$ because for a fixed $w \in T$, say $w_1$, the degree of $\varphi$ is given by
\begin{align*}
    \deg(\varphi) &= \sum_{v \in \Gamma', \varphi(v) =w_1} m_{\varphi}(v) \\
    & = m_{\varphi}(v_1)+ m_{\varphi}(v_7)+ m_{\varphi}(v_5) \\
    & = 3.
\end{align*}
 
\item [\textbf{(ii)}] Assume that $v \in \Gamma'$ is a vertex point, say $v=v_8 \in V(G')$. Then, $m_{\varphi}(v_8) = 2$, $\textnormal{val} (v_8) = 2$, and $\textnormal{val} (\varphi(v_8)) = 1 $. Therefore, $$(\textnormal{val} (v_8)-2) - m_{\varphi}(v_8) \cdot \big(\textnormal{val } (\varphi(v_8)) - 2\big) = 2 > 0.$$
Similarly, we check that the Riemann-Hurwitz condition holds on every other vertex point. Now, assume that $v$ is not a vertex point.  Consider the new vertex sets on $\Gamma'$ and $T$ (just like in part (i)) by adding $v$ and $w$ respectively.  Then, we have that $\textnormal{val}(v) = \textnormal{val}(\varphi (v))= 2$, and so the Riemann-Hurwitz condition holds. 
\end{itemize}

From (i) and (ii), we obtain that the map $\varphi:\Gamma' \to T$ is a tropical morphism of metric graphs of degree $3$, and so, the solution for the Case 1.1.A is finished.

\begin{figure}[H]
    \centering
    \begin{tikzpicture}[line cap=round,line join=round,>=triangle 45,x=1.0cm,y=1.0cm]
\definecolor{eqeqeq}{rgb}{0.8784313725490196,0.8784313725490196,0.8784313725490196}
\definecolor{uuuuuu}{rgb}{0.26666666666666666,0.26666666666666666,0.26666666666666666}
\definecolor{wrwrwr}{rgb}{0.3803921568627451,0.3803921568627451,0.3803921568627451}
\definecolor{cqcqcq}{rgb}{0.7529411764705882,0.7529411764705882,0.7529411764705882}
\definecolor{rvwvcq}{rgb}{0.08235294117647059,0.396078431372549,0.7529411764705882}
\definecolor{dtsfsf}{rgb}{0.8274509803921568,0.1843137254901961,0.1843137254901961}
\definecolor{sexdts}{rgb}{0.1803921568627451,0.49019607843137253,0.19607843137254902}
\clip(-3.,-3.) rectangle (4.,3.);
\draw (-1.,0.)-- (1.,0.);
\draw (-1.,1.)-- (1.,1.);
\draw (-1.,2.)-- (1.,2.);
\draw (-1.,-1.5)-- (1.,-1.5);
\draw [color=sexdts] (-1.,-1.5)-- (-2.,-2.);
\draw [color=dtsfsf] (-1.,-1.5)-- (-1.5,-1.5);
\draw [color=rvwvcq] (1.,-1.5)-- (1.2894301310058747,-1.012638385870342);
\draw [color=sexdts] (1.,-1.5)-- (2.,-1.5);
\draw [color=dtsfsf] (1.,-1.5)-- (1.5,-2.);
\draw [shift={(1.,1.5)},color=dtsfsf]  plot[domain=-1.5707963267948966:1.5707963267948966,variable=\t]({1.*0.5*cos(\t r)+0.*0.5*sin(\t r)},{0.*0.5*cos(\t r)+1.*0.5*sin(\t r)});
\draw [shift={(1.,1.)},color=sexdts]  plot[domain=-1.5707963267948966:1.5707963267948966,variable=\t]({1.*1.*cos(\t r)+0.*1.*sin(\t r)},{0.*1.*cos(\t r)+1.*1.*sin(\t r)});
\draw [shift={(0.7237164077379609,0.5)},color=rvwvcq]  plot[domain=-1.0659842856832995:1.0659842856832993,variable=\t]({1.*0.5712553048797154*cos(\t r)+0.*0.5712553048797154*sin(\t r)},{0.*0.5712553048797154*cos(\t r)+1.*0.5712553048797154*sin(\t r)});
\draw [color=sexdts] (1.,1.)-- (2.,0.5);
\draw [color=rvwvcq] (1.,2.)-- (1.294962530147754,2.498664332229594);
\draw [color=dtsfsf] (1.,0.)-- (1.5,-0.5);
\draw [line width=0.4pt,dash pattern=on 2pt off 2pt,color=cqcqcq] (-1.,-1.5)-- (-1.,0.);
\draw [line width=0.4pt,dash pattern=on 2pt off 2pt,color=cqcqcq] (-1.,0.)-- (-1.,1.);
\draw [line width=0.4pt,dash pattern=on 2pt off 2pt,color=cqcqcq] (-1.,1.)-- (-1.,2.);
\draw [line width=0.4pt,dash pattern=on 2pt off 2pt,color=cqcqcq] (1.,-1.5)-- (1.,0.);
\draw [line width=0.4pt,dash pattern=on 2pt off 2pt,color=cqcqcq] (1.,0.)-- (1.,1.);
\draw [line width=0.4pt,dash pattern=on 2pt off 2pt,color=cqcqcq] (1.,1.)-- (1.,2.);
\draw [line width=0.4pt,dash pattern=on 2pt off 2pt,color=cqcqcq] (1.2894301310058747,-1.012638385870342)-- (1.294962530147754,0.5032389790046032);
\draw [line width=0.4pt,dash pattern=on 2pt off 2pt,color=cqcqcq] (1.294962530147754,0.5032389790046032)-- (1.294962530147754,2.498664332229594);
\draw [line width=0.4pt,dash pattern=on 2pt off 2pt,color=cqcqcq] (1.5,1.5)-- (1.5,-0.5);
\draw [line width=0.4pt,dash pattern=on 2pt off 2pt,color=cqcqcq] (1.5,-0.5)-- (1.5,-2.);
\draw [line width=0.4pt,dash pattern=on 2pt off 2pt,color=cqcqcq] (2.,0.5)-- (2.,-1.5);
\draw [line width=0.4pt,dash pattern=on 2pt off 2pt,color=cqcqcq] (2.,0.5)-- (2.,1.);
\draw (2.905956910189192,-1.3536469272472158) node[anchor=north west] {$T$};
\draw (2.905956910189192,1.3092357812953128) node[anchor=north west] {$\Gamma'$};
\draw [->,line width=0.4pt] (3.063701734239714,0.5091002477485305) -- (3.063701734239714,-0.9908997522514695);
\draw [shift={(-1.,1.5)},color=dtsfsf]  plot[domain=1.5707963267948966:4.71238898038469,variable=\t]({1.*0.5*cos(\t r)+0.*0.5*sin(\t r)},{0.*0.5*cos(\t r)+1.*0.5*sin(\t r)});
\draw [shift={(-1.,1.)},color=sexdts]  plot[domain=1.5707963267948966:4.71238898038469,variable=\t]({1.*1.*cos(\t r)+0.*1.*sin(\t r)},{0.*1.*cos(\t r)+1.*1.*sin(\t r)});
\draw [color=dtsfsf] (-1.,0.)-- (-1.5,0.);
\draw [color=sexdts] (-1.,1.)-- (-2.,0.5);
\draw [line width=0.4pt,dash pattern=on 2pt off 2pt,color=cqcqcq] (-1.5,-1.5)-- (-1.5,1.5);
\draw [line width=0.4pt,dash pattern=on 2pt off 2pt,color=eqeqeq] (-2.,-2.)-- (-2.,0.5);
\begin{scriptsize}
\draw [fill=black] (-1.,0.) circle (0.5pt);
\draw [fill=black] (1.,0.) circle (0.5pt);
\draw[color=black] (0.8270397079059764,0.30948683457934595) node {$v_2$};
\draw[color=black] (0.05620523964366028,0.17400683712718223) node {$c$};
\draw [fill=black] (1.,1.) circle (0.5pt);
\draw[color=black] (0.8083528117056777,1.2438316445942683) node {$v_3$};
\draw[color=black] (0.06554868774380956,1.1830992319432985) node {$c$};
\draw [fill=black] (-1.,2.) circle (0.5pt);
\draw[color=black] (-0.985589223522985,2.271610935610683) node {$v_1$};
\draw [fill=black] (1.,2.) circle (0.5pt);
\draw[color=black] (0.8457266041062749,2.2622674875105337) node {$v_4$};
\draw[color=black] (0.05620523964366028,2.1735047305591157) node {$c$};
\draw [fill=black] (-1.,-1.5) circle (0.5pt);
\draw [fill=black] (1.,-1.5) circle (0.5pt);
\draw[color=black] (0.04686179154351099,-1.3116014107965441) node {$c$};
\draw [fill=sexdts] (-2.,-2.) circle (0.5pt);
\draw[color=sexdts] (-2.3263740258944074,-2.0310469145080345) node {$\frac{a-c}{2}$};
\draw [fill=dtsfsf] (-1.5,-1.5) circle (0.5pt);
\draw[color=dtsfsf] (-1.8124843803861967,-1.4797834765992302) node {$\frac{b-c}{2}$};
\draw [fill=rvwvcq] (1.2894301310058747,-1.012638385870342) circle (0.5pt);
\draw[color=rvwvcq] (1.392318317965008,-0.7696814209878893) node {$\frac{e}{2}$};
\draw [fill=sexdts] (2.,-1.5) circle (0.5pt);
\draw[color=sexdts] (2.1584810621772497,-1.4237227879983347) node {$\frac{f}{2}$};
\draw [fill=dtsfsf] (1.5,-2.) circle (0.5pt);
\draw[color=dtsfsf] (1.644591416669039,-1.9282689854063928) node {$\frac{d}{2}$};
\draw [color=dtsfsf] (1.5,1.5)-- ++(-1.0pt,-1.0pt) -- ++(2.0pt,2.0pt) ++(-2.0pt,0) -- ++(2.0pt,-2.0pt);
\draw[color=dtsfsf] (1.6632783128693374,1.4447157787474767) node {$d$};
\draw [color=rvwvcq] (1.294962530147754,0.5032389790046032)-- ++(-1.0pt,-1.0pt) -- ++(2.0pt,2.0pt) ++(-2.0pt,0) -- ++(2.0pt,-2.0pt);
\draw[color=rvwvcq] (1.420348662265456,0.38890614343061436) node {$e$};
\draw [color=sexdts] (2.,1.)-- ++(-1.0pt,-1.0pt) -- ++(2.0pt,2.0pt) ++(-2.0pt,0) -- ++(2.0pt,-2.0pt);
\draw[color=sexdts] (2.1491376140771,0.9868868218401646) node {$f$};
\draw [fill=sexdts] (2.,0.5) circle (0.5pt);
\draw[color=sexdts] (2.1584810621772497,0.5664316573334496) node {$\frac{f}{2}$};
\draw [fill=rvwvcq] (1.294962530147754,2.498664332229594) circle (0.5pt);
\draw[color=rvwvcq] (1.3082272850636647,2.7527985127683676) node {$\frac{e}{2}$};
\draw [fill=dtsfsf] (1.5,-0.5) circle (0.5pt);
\draw[color=dtsfsf] (1.6632783128693376,-0.423973841282368) node {$\frac{d}{2}$};
\draw[color=black] (3.4011596594971047,-0.2184179830790851) node {$\varphi$};
\draw[color=wrwrwr] (-4.951882942036356,4.1823460720912) node {$Curve$};
\draw [fill=uuuuuu] (-1.,1.) circle (0.5pt);
\draw[color=dtsfsf] (-1.0556650842741047,1.463402674947775) node {$b-c$};
\draw[color=sexdts] (-2.531929884097692,0.9308261332392692) node {$a-c$};
\draw [color=dtsfsf] (-1.5,1.5)-- ++(-1.0pt,-1.0pt) -- ++(2.0pt,2.0pt) ++(-2.0pt,0) -- ++(2.0pt,-2.0pt);
\draw [color=sexdts] (-2.,1.)-- ++(-1.0pt,-1.0pt) -- ++(2.0pt,2.0pt) ++(-2.0pt,0) -- ++(2.0pt,-2.0pt);
\draw [fill=dtsfsf] (-1.5,0.) circle (0.5pt);
\draw[color=dtsfsf] (-1.8405147246866447,-0.012862124875802225) node {$\frac{b-c}{2}$};
\draw [fill=sexdts] (-2.,0.5) circle (0.5pt);
\draw[color=sexdts] (-2.3263740258944074,0.43562338393136046) node {$\frac{a-c}{2}$};
\end{scriptsize}
\end{tikzpicture}
    \begin{center}
        \textnormal{Figure 2. The tropical morphism $\varphi:\Gamma' \to T$}
    \end{center}
\end{figure}

\begin{Remark} \textnormal{Let $\varphi:\Gamma' \to T$ be non-constant piecewise linear map  with nonzero integer slopes (as in Case 1.1.A), where the models $(G',l')$, $(T',t')$ of $\Gamma'$, $T$ respectively, are taken so that the condition in the Definition 2.4 is satisfied. In order to show that $\varphi$ satisfies the harmonicity and the Riemann-Hurwitz condition on $\Gamma'$, it is enough to check those conditions on vertex points. This is due to the parts (i) and (ii) above.}
\end{Remark} 
\par 

\textbf{Case 1.1.B.}  Let $\Gamma_1'$ be the tropical modification of $\Gamma$ with model $(G_1',l_1')$, where the graph $G_1'$ is obtained by contracting the edges $v_1v_8, v_8v_7$, and $v_5v_8'$ of $G'$ in Figure 1.2. Let $T_1$ be the metric tree with model $(T_1', t_1')$, where the tree $T_1'$ is obtained by contracting the edge $w_1w_8$ of $T'$ in Figure 1.3. Next, let $\psi_1:V(G_1') \to V(T_1)$ the map on the set of vertices given by $v_2, v_3,v_4 \mapsto w_0$, $v_1,v_5 \mapsto w_1$, and $v_i, v_i' \mapsto w_i$ for $i=6,9,10,11$. This map $\psi$ satisfies the condition in Lemma 3.1, and so, there exist a unique continuous map $\varphi_1: \Gamma_1' \to T_1$, given in Figure 3, such that $\varphi_1|_{V(G_1')} = \psi_1$ and $\varphi_1$ is linear on each edge $e' \in E(G_1')$ with slope $t_1'(e)/l_1'(e')$, where $e = \varphi_1(e') \in E(T_1)$ with endpoints $\psi_1(v)$ and $\psi_1(w)$. Following the reasoning in (i) and (ii), we get that $\varphi_1$ is a tropical map of degree $3$, and thus, the solution of Case 1.1.B is done. 
\newpage 
\begin{figure}[H]
\centering
    \begin{tikzpicture}[line cap=round,line join=round,>=triangle 45,x=1.0cm,y=1.0cm]
\definecolor{eqeqeq}{rgb}{0.8784313725490196,0.8784313725490196,0.8784313725490196}
\definecolor{wrwrwr}{rgb}{0.3803921568627451,0.3803921568627451,0.3803921568627451}
\definecolor{cqcqcq}{rgb}{0.7529411764705882,0.7529411764705882,0.7529411764705882}
\definecolor{dtsfsf}{rgb}{0.8274509803921568,0.1843137254901961,0.1843137254901961}
\definecolor{rvwvcq}{rgb}{0.08235294117647059,0.396078431372549,0.7529411764705882}
\definecolor{sexdts}{rgb}{0.1803921568627451,0.49019607843137253,0.19607843137254902}
\clip(-3.,-3.) rectangle (4.,3.);
\draw (-1.,0.)-- (1.,0.);
\draw (-1.,2.)-- (1.,2.);
\draw (-1.,-1.5)-- (1.,-1.5);
\draw [color=sexdts] (-1.,-1.5)-- (-2.,-2.);
\draw [color=rvwvcq] (1.,-1.5)-- (1.2894301310058747,-1.012638385870342);
\draw [color=sexdts] (1.,-1.5)-- (2.,-1.5);
\draw [color=dtsfsf] (1.,-1.5)-- (1.5,-2.);
\draw [shift={(1.,1.5)},color=dtsfsf]  plot[domain=-1.5707963267948966:1.5707963267948966,variable=\t]({1.*0.5*cos(\t r)+0.*0.5*sin(\t r)},{0.*0.5*cos(\t r)+1.*0.5*sin(\t r)});
\draw [shift={(1.,1.)},color=sexdts]  plot[domain=-1.5707963267948966:1.5707963267948966,variable=\t]({1.*1.*cos(\t r)+0.*1.*sin(\t r)},{0.*1.*cos(\t r)+1.*1.*sin(\t r)});
\draw [shift={(0.7237164077379609,0.5)},color=rvwvcq]  plot[domain=-1.0659842856832995:1.0659842856832993,variable=\t]({1.*0.5712553048797154*cos(\t r)+0.*0.5712553048797154*sin(\t r)},{0.*0.5712553048797154*cos(\t r)+1.*0.5712553048797154*sin(\t r)});
\draw [color=sexdts] (1.,1.)-- (2.,0.5);
\draw [color=rvwvcq] (1.,2.)-- (1.294962530147754,2.498664332229594);
\draw [color=dtsfsf] (1.,0.)-- (1.5,-0.5);
\draw [line width=0.4pt,dash pattern=on 2pt off 2pt,color=cqcqcq] (-1.,-1.5)-- (-1.,0.);
\draw [line width=0.4pt,dash pattern=on 2pt off 2pt,color=cqcqcq] (-1.,0.)-- (-1.,1.);
\draw [line width=0.4pt,dash pattern=on 2pt off 2pt,color=cqcqcq] (-1.,1.)-- (-1.,2.);
\draw [line width=0.4pt,dash pattern=on 2pt off 2pt,color=cqcqcq] (1.,-1.5)-- (1.,0.);
\draw [line width=0.4pt,dash pattern=on 2pt off 2pt,color=cqcqcq] (1.,0.)-- (1.,1.);
\draw [line width=0.4pt,dash pattern=on 2pt off 2pt,color=cqcqcq] (1.,1.)-- (1.,2.);
\draw [line width=0.4pt,dash pattern=on 2pt off 2pt,color=cqcqcq] (1.2894301310058747,-1.012638385870342)-- (1.294962530147754,0.5032389790046032);
\draw [line width=0.4pt,dash pattern=on 2pt off 2pt,color=cqcqcq] (1.294962530147754,0.5032389790046032)-- (1.294962530147754,2.498664332229594);
\draw [line width=0.4pt,dash pattern=on 2pt off 2pt,color=cqcqcq] (1.5,1.5)-- (1.5,-0.5);
\draw [line width=0.4pt,dash pattern=on 2pt off 2pt,color=cqcqcq] (1.5,-0.5)-- (1.5,-2.);
\draw [line width=0.4pt,dash pattern=on 2pt off 2pt,color=cqcqcq] (2.,0.5)-- (2.,-1.5);
\draw [line width=0.4pt,dash pattern=on 2pt off 2pt,color=cqcqcq] (2.,0.5)-- (2.,1.);
\draw (2.9059569101891944,-1.3536469272472165) node[anchor=north west] {$T_1$};
\draw (2.9059569101891944,1.3092357812953115) node[anchor=north west] {$\Gamma_1'$};
\draw [->,line width=0.4pt] (3.063701734239714,0.5091002477485305) -- (3.063701734239714,-0.9908997522514695);
\draw [shift={(-1.,1.)},color=sexdts]  plot[domain=1.5707963267948966:4.71238898038469,variable=\t]({1.*1.*cos(\t r)+0.*1.*sin(\t r)},{0.*1.*cos(\t r)+1.*1.*sin(\t r)});
\draw [line width=0.4pt,dash pattern=on 2pt off 2pt,color=eqeqeq] (-2.,-2.)-- (-2.,0.5);
\draw (1.,1.)-- (-1.,2.);
\draw [color=sexdts] (-1.,2.)-- (-2.,0.5);
\begin{scriptsize}
\draw [fill=black] (-1.,0.) circle (0.5pt);
\draw [fill=black] (1.,0.) circle (0.5pt);
\draw[color=black] (0.8270397079059781,0.30948683457934484) node {$v_2$};
\draw[color=black] (0.05620523964366188,0.17400683712718112) node {$c$};
\draw [fill=black] (1.,1.) circle (0.5pt);
\draw[color=black] (0.8083528117056795,0.9448413053894918) node {$v_3$};
\draw [fill=black] (-1.,2.) circle (0.5pt);
\draw[color=black] (-0.9855892235229836,2.2716109356106813) node {$v_1$};
\draw [fill=black] (1.,2.) circle (0.5pt);
\draw[color=black] (0.8457266041062766,2.262267487510532) node {$v_4$};
\draw[color=black] (0.05620523964366188,2.1735047305591144) node {$c$};
\draw [fill=black] (-1.,-1.5) circle (0.5pt);
\draw [fill=black] (1.,-1.5) circle (0.5pt);
\draw[color=black] (0.046861791543512596,-1.311601410796545) node {$c$};
\draw [fill=sexdts] (-2.,-2.) circle (0.5pt);
\draw[color=sexdts] (-2.3450609220947047,-1.984329674007289) node {$\frac{a-c}{2}$};
\draw [fill=rvwvcq] (1.2894301310058747,-1.012638385870342) circle (0.5pt);
\draw[color=rvwvcq] (1.39231831796501,-0.7696814209878903) node {$\frac{e}{2}$};
\draw [fill=sexdts] (2.,-1.5) circle (0.5pt);
\draw[color=sexdts] (2.1024203735763556,-1.433066236098485) node {$\frac{f}{2}$};
\draw [fill=dtsfsf] (1.5,-2.) circle (0.5pt);
\draw[color=dtsfsf] (1.8127734824717279,-1.9189255373062444) node {$\frac{d}{2}$};
\draw [color=dtsfsf] (1.5,1.5)-- ++(-1.0pt,-1.0pt) -- ++(2.0pt,2.0pt) ++(-2.0pt,0) -- ++(2.0pt,-2.0pt);
\draw[color=dtsfsf] (1.6632783128693394,1.4447157787474754) node {$d$};
\draw [color=rvwvcq] (1.294962530147754,0.5032389790046032)-- ++(-1.0pt,-1.0pt) -- ++(2.0pt,2.0pt) ++(-2.0pt,0) -- ++(2.0pt,-2.0pt);
\draw[color=rvwvcq] (1.429692110365607,0.42627993583121015) node {$e$};
\draw [color=sexdts] (2.,1.)-- ++(-1.0pt,-1.0pt) -- ++(2.0pt,2.0pt) ++(-2.0pt,0) -- ++(2.0pt,-2.0pt);
\draw[color=sexdts] (2.1491376140771026,0.9868868218401634) node {$f$};
\draw [fill=sexdts] (2.,0.5) circle (0.5pt);
\draw[color=sexdts] (2.1024203735763556,0.5757751054335977) node {$\frac{f}{2}$};
\draw [fill=rvwvcq] (1.294962530147754,2.498664332229594) circle (0.5pt);
\draw[color=rvwvcq] (1.3082272850636665,2.752798512768366) node {$\frac{e}{2}$};
\draw [fill=dtsfsf] (1.5,-0.5) circle (0.5pt);
\draw[color=dtsfsf] (1.8127734824717279,-0.42397384128236903) node {$\frac{d}{2}$};
\draw[color=black] (3.4245182797474802,-0.19505936282871308) node {$\varphi_1$};
\draw[color=wrwrwr] (-3.550365727013963,3.7712343556846313) node {$Curve$};
\draw[color=sexdts] (-2.531929884097691,0.930826133239268) node {$a-c$};
\draw [color=sexdts] (-2.,1.)-- ++(-1.0pt,-1.0pt) -- ++(2.0pt,2.0pt) ++(-2.0pt,0) -- ++(2.0pt,-2.0pt);
\draw [fill=sexdts] (-2.,0.5) circle (0.5pt);
\draw[color=sexdts] (-2.3450609220947047,0.5290578649328516) node {$\frac{a-c}{2}$};
\draw[color=black] (0.046861791543512596,1.089664750941805) node {$c$};
\end{scriptsize}
\end{tikzpicture}
    \begin{center}
        \textnormal{Figure 3. The tropical map $\varphi_1:\Gamma_1' \to T_1$}
    \end{center}
    
\end{figure}
\begin{figure}[H]
    \centering
    \begin{tikzpicture}[line cap=round,line join=round,>=triangle 45,x=1.0cm,y=1.0cm]
\definecolor{eqeqeq}{rgb}{0.8784313725490196,0.8784313725490196,0.8784313725490196}
\definecolor{wrwrwr}{rgb}{0.3803921568627451,0.3803921568627451,0.3803921568627451}
\definecolor{cqcqcq}{rgb}{0.7529411764705882,0.7529411764705882,0.7529411764705882}
\definecolor{dtsfsf}{rgb}{0.8274509803921568,0.1843137254901961,0.1843137254901961}
\definecolor{sexdts}{rgb}{0.1803921568627451,0.49019607843137253,0.19607843137254902}
\definecolor{rvwvcq}{rgb}{0.08235294117647059,0.396078431372549,0.7529411764705882}
\clip(-2.,-3.) rectangle (4.,3.);
\draw (-1.,2.)-- (1.,2.);
\draw (-1.,-1.5)-- (1.,-1.5);
\draw [color=rvwvcq] (1.,-1.5)-- (1.2894301310058747,-1.012638385870342);
\draw [color=sexdts] (1.,-1.5)-- (2.,-1.5);
\draw [color=dtsfsf] (1.,-1.5)-- (1.5,-2.);
\draw [shift={(1.,1.5)},color=dtsfsf]  plot[domain=-1.5707963267948966:1.5707963267948966,variable=\t]({1.*0.5*cos(\t r)+0.*0.5*sin(\t r)},{0.*0.5*cos(\t r)+1.*0.5*sin(\t r)});
\draw [shift={(1.,1.)},color=sexdts]  plot[domain=-1.5707963267948966:1.5707963267948966,variable=\t]({1.*1.*cos(\t r)+0.*1.*sin(\t r)},{0.*1.*cos(\t r)+1.*1.*sin(\t r)});
\draw [shift={(0.7237164077379609,0.5)},color=rvwvcq]  plot[domain=-1.0659842856832995:1.0659842856832993,variable=\t]({1.*0.5712553048797154*cos(\t r)+0.*0.5712553048797154*sin(\t r)},{0.*0.5712553048797154*cos(\t r)+1.*0.5712553048797154*sin(\t r)});
\draw [color=sexdts] (1.,1.)-- (2.,0.5);
\draw [color=rvwvcq] (1.,2.)-- (1.294962530147754,2.498664332229594);
\draw [color=dtsfsf] (1.,0.)-- (1.5,-0.5);
\draw [line width=0.4pt,dash pattern=on 2pt off 2pt,color=cqcqcq] (-1.,-1.5)-- (-1.,0.);
\draw [line width=0.4pt,dash pattern=on 2pt off 2pt,color=cqcqcq] (-1.,0.)-- (-1.,1.);
\draw [line width=0.4pt,dash pattern=on 2pt off 2pt,color=cqcqcq] (-1.,1.)-- (-1.,2.);
\draw [line width=0.4pt,dash pattern=on 2pt off 2pt,color=cqcqcq] (1.,-1.5)-- (1.,0.);
\draw [line width=0.4pt,dash pattern=on 2pt off 2pt,color=cqcqcq] (1.,0.)-- (1.,1.);
\draw [line width=0.4pt,dash pattern=on 2pt off 2pt,color=cqcqcq] (1.,1.)-- (1.,2.);
\draw [line width=0.4pt,dash pattern=on 2pt off 2pt,color=cqcqcq] (1.2894301310058747,-1.012638385870342)-- (1.294962530147754,0.5032389790046032);
\draw [line width=0.4pt,dash pattern=on 2pt off 2pt,color=cqcqcq] (1.294962530147754,0.5032389790046032)-- (1.294962530147754,2.498664332229594);
\draw [line width=0.4pt,dash pattern=on 2pt off 2pt,color=cqcqcq] (1.5,1.5)-- (1.5,-0.5);
\draw [line width=0.4pt,dash pattern=on 2pt off 2pt,color=cqcqcq] (1.5,-0.5)-- (1.5,-2.);
\draw [line width=0.4pt,dash pattern=on 2pt off 2pt,color=cqcqcq] (2.,0.5)-- (2.,-1.5);
\draw [line width=0.4pt,dash pattern=on 2pt off 2pt,color=cqcqcq] (2.,0.5)-- (2.,1.);
\draw (2.9059569101891944,-1.3536469272472165) node[anchor=north west] {$T_2'$};
\draw (2.9059569101891944,1.3092357812953115) node[anchor=north west] {$\Gamma_2'$};
\draw [->,line width=0.4pt] (3.063701734239714,0.5091002477485305) -- (3.063701734239714,-0.9908997522514695);
\draw [line width=0.4pt,dash pattern=on 2pt off 2pt,color=eqeqeq] (-2.,-2.)-- (-2.,0.5);
\draw (1.,1.)-- (-1.,2.);
\draw (1.,0.)-- (-1.,2.);
\begin{scriptsize}
\draw [fill=black] (1.,0.) circle (0.5pt);
\draw[color=black] (0.8083528117056795,-0.12965522612766855) node {$v_2$};
\draw [fill=black] (1.,1.) circle (0.5pt);
\draw[color=black] (0.8083528117056795,0.9448413053894918) node {$v_3$};
\draw [fill=black] (-1.,2.) circle (0.5pt);
\draw[color=black] (-0.9855892235229836,2.2716109356106813) node {$v_1$};
\draw [fill=black] (1.,2.) circle (0.5pt);
\draw[color=black] (0.8457266041062766,2.262267487510532) node {$v_4$};
\draw[color=black] (0.05620523964366188,2.1735047305591144) node {$c$};
\draw [fill=black] (-1.,-1.5) circle (0.5pt);
\draw [fill=black] (1.,-1.5) circle (0.5pt);
\draw[color=black] (0.046861791543512596,-1.311601410796545) node {$c$};
\draw [fill=rvwvcq] (1.2894301310058747,-1.012638385870342) circle (0.5pt);
\draw[color=rvwvcq] (1.39231831796501,-0.7696814209878903) node {$\frac{e}{2}$};
\draw [fill=sexdts] (2.,-1.5) circle (0.5pt);
\draw[color=sexdts] (2.1024203735763556,-1.433066236098485) node {$\frac{f}{2}$};
\draw [fill=dtsfsf] (1.5,-2.) circle (0.5pt);
\draw[color=dtsfsf] (1.8127734824717279,-1.9189255373062444) node {$\frac{d}{2}$};
\draw [color=dtsfsf] (1.5,1.5)-- ++(-1.0pt,-1.0pt) -- ++(2.0pt,2.0pt) ++(-2.0pt,0) -- ++(2.0pt,-2.0pt);
\draw[color=dtsfsf] (1.6632783128693394,1.4447157787474754) node {$d$};
\draw [color=rvwvcq] (1.294962530147754,0.5032389790046032)-- ++(-1.0pt,-1.0pt) -- ++(2.0pt,2.0pt) ++(-2.0pt,0) -- ++(2.0pt,-2.0pt);
\draw[color=rvwvcq] (1.429692110365607,0.42627993583121015) node {$e$};
\draw [color=sexdts] (2.,1.)-- ++(-1.0pt,-1.0pt) -- ++(2.0pt,2.0pt) ++(-2.0pt,0) -- ++(2.0pt,-2.0pt);
\draw[color=sexdts] (2.1491376140771026,0.9868868218401634) node {$f$};
\draw [fill=sexdts] (2.,0.5) circle (0.5pt);
\draw[color=sexdts] (2.1024203735763556,0.5757751054335977) node {$\frac{f}{2}$};
\draw [fill=rvwvcq] (1.294962530147754,2.498664332229594) circle (0.5pt);
\draw[color=rvwvcq] (1.3082272850636665,2.752798512768366) node {$\frac{e}{2}$};
\draw [fill=dtsfsf] (1.5,-0.5) circle (0.5pt);
\draw[color=dtsfsf] (1.8127734824717279,-0.42397384128236903) node {$\frac{d}{2}$};
\draw[color=black] (3.4245182797474802,-0.19505936282871308) node {$\varphi_2$};
\draw[color=wrwrwr] (-3.550365727013963,3.7712343556846313) node {$Curve$};
\draw[color=black] (0.03751834344336331,1.6689585331510566) node {$c$};
\draw[color=black] (0.03751834344336331,0.46365372823180695) node {$c$};
\end{scriptsize}
\end{tikzpicture}
    \begin{center}
        \textnormal{Figure 4. The tropical map $\varphi_2: \Gamma_2' \to T_2$}
    \end{center}
\end{figure}

\textbf{Case 1.1.C.} Let $\Gamma_2'$ be the tropical modification of $\Gamma$ with model $(G_2',l_2')$, where $G_2'$ is obtained by contracting the edges $v_1v_6, v_6v_5, v_7v_6', v_1v_8, v_8v_7$ and $v_5v_8'$ of the graph $G'$ as in Figure 1.2. Let $T_2$ be the metric tree with model $(T_2', t_2')$, where the tree $T_2'$ is obtained by contracting the edges $w_1w_6, w_1w_8$ of the tree $T'$ as in Figure 1.3. Next, let $\psi_2:V(G_2') \to V(T_2)$ the map on the set of vertices given by $v_2, v_3,v_4 \mapsto w_0$, $v_1 \mapsto w_1$ and $v_i, v_i' \mapsto w_i$ for $i=9,10,11$. The function $\psi_2$ satisfies the condition in Lemma 3.1 and so, there exist a unique continuous map $\varphi_2: \Gamma_2' \to T_2$, given in Figure 4, such that $\varphi_2|_{V(G_2')} = \psi_2$ and $\varphi_2$ is linear on each edge $e' \in E(G_2')$ with slope $t_2'(e)/l_2'(e')$, where $e = \varphi_2(e') \in E(T_2)$ with endpoints $\psi_2(v)$ and $\psi_2(w)$. Following the reasoning in (i) and (ii), we conclude that $\varphi_2$ is a tropical map of degree $3$, and therefore, the solution of Case 1.1.C is finished. \par
\textbf{Case 1.2.}  Consider the metric graph $\Gamma$ with essential model $(G,l)$ in Figure 5. The graph $G$ is given by its vertex set $V(G) = \left\lbrace v_1, v_2, v_3, v_4 \right\rbrace$, and edge set $E(G)=\left\lbrace v_1v_2,v_3v_4,e_1,e_2, e_3,e_4 \right\rbrace$, where $e_1,e_2$ (resp., $e_3,e_4$) are two edges with endpoints $v_1,v_4$ (resp., $v_2,v_3$).  The length map $l:E(G) \to (0, \infty)$ is defined by assigning $v_1v_2 \mapsto a, v_3v_4 \mapsto b, e_1 \mapsto c, e_2 \mapsto d, e_3 \mapsto e$ and $e_4 \mapsto f$ where $a,b,c,d,e$, and $f$ are real positive numbers such that $a<b$. Note that if $a=b$,  then $\Gamma$ is a hyperelliptic metric graph.

\begin{figure}[H]
    \centering
    \begin{tikzpicture}[line cap=round,line join=round,>=triangle 45,x=1.0cm,y=1.0cm]
\definecolor{wrwrwr}{rgb}{0.3803921568627451,0.3803921568627451,0.3803921568627451}
\clip(-2.,-2.) rectangle (3.2,3.);
\draw (0.,-1.)-- (1.,-1.);
\draw (1.,-1.)-- (1.,2.);
\draw (1.,2.)-- (0.,1.);
\draw (0.,1.)-- (0.,-1.);
\draw [shift={(0.,0.)}] plot[domain=1.5707963267948966:4.71238898038469,variable=\t]({1.*1.*cos(\t r)+0.*1.*sin(\t r)},{0.*1.*cos(\t r)+1.*1.*sin(\t r)});
\draw [shift={(1.,0.5)}] plot[domain=-1.5707963267948966:1.5707963267948966,variable=\t]({1.*1.5*cos(\t r)+0.*1.5*sin(\t r)},{0.*1.5*cos(\t r)+1.*1.5*sin(\t r)});
\begin{scriptsize}
\draw [fill=black] (0.,-1.) circle (0.5pt);
\draw[color=black] (-0.0412632622492525,-1.2356951217874457) node {$v_1$};
\draw [fill=black] (1.,-1.) circle (0.5pt);
\draw[color=black] (1.027050956857946,-1.2613346630460183) node {$v_2$};
\draw [fill=black] (0.,1.) circle (0.5pt);
\draw[color=black] (-0.07544931726068285,1.2769799215526718) node {$v_4$};
\draw [fill=black] (1.,2.) circle (0.5pt);
\draw[color=black] (0.9928649018465157,2.285468544389862) node {$v_3$};
\draw[color=black] (0.5014403610572042,-1.282700947428162) node {$a$};
\draw[color=black] (1.2877196263201023,0.4351483168962039) node {$e$};
\draw[color=black] (0.47580081979863154,1.8025905173534107) node {$b$};
\draw[color=black] (-0.3190249592171241,-0.1032820495338212) node {$d$};
\draw[color=black] (-1.3275135820543194,-0.12892159079239382) node {$c$};
\draw[color=black] (2.8004525605758954,0.4522413444019191) node {$f$};
\draw[color=wrwrwr] (-4.156409634250181,3.7341026254992156) node {$Curve$};
\end{scriptsize}
\end{tikzpicture}
    \begin{center}
        \textnormal{Figure 5. The essential model $(G,l)$ of $\Gamma$}
    \end{center}
\end{figure}

Let $(G_1,l_1)$ be another model of $\Gamma$  as in Figure 5.1. The graph $G_1$ is obtained from $G$ by subdividing the following edges: $v_3v_4 \in E(G)$ into $v_3v_6$, $v_6v_5, v_5v_4$; $e_1 \in E(G) $ into $v_1v_7,v_7v_4 $; $e_2 \in E(G) $ into $v_1v_8,v_8v_4 $; $e_3\in E(G)$ into $v_2v_9,v_9v_3 $, and $e_4 \in E(G)$ into $v_2v_{10},v_{10}v_3 $, such that

\begin{equation*}
\begin{aligned}[c]
l_1(v_3v_6) = l_1(v_6v_5) &= \frac{b-a}{2}\\
l_1(v_1v_7) = l_1(v_7v_4) &= \frac{c}{2}\\
l_1(v_1v_8) = l_1(v_8v_4) &= \frac{d}{2} 
\end{aligned}
\qquad
\begin{aligned}[c]
l_1(v_4v_5) &= a \\
l_1(v_2v_9) = l_1(v_9v_3) &= \frac{e}{2} \\
l_1(v_2v_{10}) = l_1(v_{10}v_3) &= \frac{f}{2}.
\end{aligned}
\end{equation*}
\begin{figure}[H]
    \centering
    \begin{tikzpicture}[line cap=round,line join=round,>=triangle 45,x=1.0cm,y=1.0cm]
\definecolor{wrwrwr}{rgb}{0.3803921568627451,0.3803921568627451,0.3803921568627451}
\definecolor{sexdts}{rgb}{0.1803921568627451,0.49019607843137253,0.19607843137254902}
\definecolor{ffqqqq}{rgb}{1.,0.,0.}
\definecolor{rvwvcq}{rgb}{0.08235294117647059,0.396078431372549,0.7529411764705882}
\clip(-2.,-0.9) rectangle (2.,4.);
\draw (-0.5,0.)-- (0.5,0.);
\draw (-0.5,3.)-- (0.5,3.);
\draw [shift={(-0.25,1.)},color=rvwvcq]  plot[domain=-0.9272952180016123:0.9272952180016122,variable=\t]({1.*1.25*cos(\t r)+0.*1.25*sin(\t r)},{0.*1.25*cos(\t r)+1.*1.25*sin(\t r)});
\draw [shift={(0.5,1.)},color=ffqqqq]  plot[domain=-1.5707963267948966:1.5707963267948966,variable=\t]({1.*1.*cos(\t r)+0.*1.*sin(\t r)},{0.*1.*cos(\t r)+1.*1.*sin(\t r)});
\draw (-0.5,-2.)-- (0.5,-2.);
\draw [shift={(0.2333333333333335,2.5)},color=sexdts]  plot[domain=-1.080839000541168:1.0808390005411685,variable=\t]({1.*0.5666666666666665*cos(\t r)+0.*0.5666666666666665*sin(\t r)},{0.*0.5666666666666665*cos(\t r)+1.*0.5666666666666665*sin(\t r)});
\draw [shift={(1.4657697822512776,1.5)},color=ffqqqq]  plot[domain=2.4897854997887947:3.7933998073907915,variable=\t]({1.*2.4727011216101786*cos(\t r)+0.*2.4727011216101786*sin(\t r)},{0.*2.4727011216101786*cos(\t r)+1.*2.4727011216101786*sin(\t r)});
\draw [shift={(0.14435667947216052,1.5)},color=rvwvcq]  plot[domain=1.9765323738923735:4.3066529332872125,variable=\t]({1.*1.632542658058401*cos(\t r)+0.*1.632542658058401*sin(\t r)},{0.*1.632542658058401*cos(\t r)+1.*1.632542658058401*sin(\t r)});
\draw [color=ffqqqq] (-0.5,-2.)-- (-1.,-3.);
\draw [color=rvwvcq] (-0.5,-2.)-- (-1.5,-1.);
\draw [color=ffqqqq] (0.5,-2.)-- (1.5,-2.);
\draw [color=rvwvcq] (0.5,-2.)-- (1.,-3.);
\draw [color=sexdts] (0.5,-2.)-- (0.8,-1.0123047916035581);
\draw (2.302153011953767,1.7453380218400287) node[anchor=north west] {$\Gamma'$};
\draw (2.3131450563766167,-1.8913527146761007) node[anchor=north west] {$T$};
\draw [->] (2.398333400653702,0.5662062130719998) -- (2.398333400653702,-1.2211093707448013);
\begin{scriptsize}
\draw [fill=black] (-0.5,0.) circle (0.5pt);
\draw[color=black] (-0.49807030476718667,-0.17850861351833328) node {$v_1$};
\draw [fill=black] (0.5,0.) circle (0.5pt);
\draw[color=black] (0.4967097155007086,-0.21574435484784996) node {$v_2$};
\draw[color=black] (0.010311749789610594,-0.24677413928911385) node {$a$};
\draw [fill=black] (0.5,2.) circle (0.5pt);
\draw[color=black] (0.4967097155007086,1.794985676946051) node {$v_3$};
\draw [fill=black] (-0.5,3.) circle (0.5pt);
\draw[color=black] (-0.5035663269786116,3.22235576124419) node {$v_4$};
\draw [fill=black] (0.5,3.) circle (0.5pt);
\draw[color=black] (0.41976540454076083,3.197531933691179) node {$v_5$};
\draw[color=black] (-0.022664383478938413,2.7444970808487263) node {$a$};
\draw [color=rvwvcq] (1.,1.)-- ++(-1.0pt,-1.0pt) -- ++(2.0pt,2.0pt) ++(-2.0pt,0) -- ++(2.0pt,-2.0pt);
\draw[color=rvwvcq] (1.255160780677336,1.0999185054617395) node {$v_9$};
\draw[color=rvwvcq] (1.0545559699603295,1.4536580480921482) node {$e$};
\draw [color=ffqqqq] (1.5,1.)-- ++(-1.0pt,-1.0pt) -- ++(2.0pt,2.0pt) ++(-2.0pt,0) -- ++(2.0pt,-2.0pt);
\draw[color=ffqqqq] (1.8047630018198193,1.0750946779087285) node {$v_{10}$};
\draw[color=ffqqqq] (1.268900836205898,1.875663116493337) node {$f$};
\draw [fill=black] (-0.5,-2.) circle (0.5pt);
\draw [fill=black] (0.5,-2.) circle (0.5pt);
\draw[color=black] (0.0267998164238851,-2.257504171083015) node {$a$};
\draw [color=sexdts] (0.8,2.5)-- ++(-1.0pt,-1.0pt) -- ++(2.0pt,2.0pt) ++(-2.0pt,0) -- ++(2.0pt,-2.0pt);
\draw[color=sexdts] (0.5956381153063556,2.5024647622068676) node {$v_6$};
\draw[color=sexdts] (1.2139406140916496,2.396963495106571) node {$b-a$};
\draw [color=ffqqqq] (-1.006924325949443,1.5058893192882263)-- ++(-1.0pt,-1.0pt) -- ++(2.0pt,2.0pt) ++(-2.0pt,0) -- ++(2.0pt,-2.0pt);
\draw[color=ffqqqq] (-1.256521369943814,1.4970997463099176) node {$v_8$};
\draw[color=ffqqqq] (-0.7206592043298925,1.0688887210204756) node {$d$};
\draw [color=rvwvcq] (-1.4879256850319909,1.4708484492843676)-- ++(-1.0pt,-1.0pt) -- ++(2.0pt,2.0pt) ++(-2.0pt,0) -- ++(2.0pt,-2.0pt);
\draw[color=rvwvcq] (-1.7456673467606243,1.484687832533412) node {$v_7$};
\draw[color=rvwvcq] (-1.237285292203827,1.0192410659144533) node {$c$};
\draw [fill=black] (-1.,-3.) circle (0.5pt);
\draw[color=ffqqqq] (-0.5282984269300233,-2.6919211532607092) node {$\frac{d}{2}$};
\draw [fill=black] (-1.5,-1.) circle (0.5pt);
\draw[color=rvwvcq] (-1.1218688257639056,-1.5872608271517146) node {$\frac{c}{2}$};
\draw [fill=black] (1.5,-2.) circle (0.5pt);
\draw[color=ffqqqq] (1.3788212804343947,-2.245092257306509) node {$\frac{f}{2}$};
\draw [fill=black] (1.,-3.) circle (0.5pt);
\draw[color=rvwvcq] (0.6863224817948654,-2.8036283772492596) node {$\frac{e}{2}$};
\draw [fill=black] (0.8,-1.0123047916035581) circle (0.5pt);
\draw[color=sexdts] (1.0820360810174536,-1.5003774307161755) node {$\frac{b-a}{2}$};
\draw[color=black] (2.6648904779078064,-0.22195031173610275) node {$\varphi$};
\draw[color=wrwrwr] (-2.2705374679516965,4.978641560619727) node {$Curve$};
\end{scriptsize}
\end{tikzpicture}
    \begin{center}
        \textnormal{Figure 5.1. The model $(G_1,l_1)$ of $\Gamma$}
    \end{center}
\end{figure}
Let $\Gamma' $ be the tropical modification of $\Gamma$ with model $(G',l')$ in Figure 5.2, where the graph $G'$ is given with its vertex set $V(G') = V(G_1) \cup \left\lbrace v_6',v_7',\ldots,v_{11}' \right\rbrace$, and edge set $E(G')=  \left\lbrace v_2v_6', v_3v_{11}',v_{11}'v_7', v_{11}'v_8',v_5v_9', v_5v_{10}' \right\rbrace \cup E(G_1)$. The length map $l'$ on $G'$ is given by $l' = l_1$ on $E(G_1)$, and 
\begin{align*}
 l'(v_1v_7) = l'(v_7v_4) = l'(v_{11}v_{7}') &= \frac{c}{2} & l'(v_5v_{10}') = l'(v_2v_{10})= l'(v_{10},v_{3}) &= \frac{f}{2} \\ 
l'(v_{11}v_{8}')=l'(v_{1}v_{8})=l'(v_8v_4)&= \frac{d}{2}  & l'(v_2v_6')=l'(v_3v_6)=l'(v_6v_5)&= \frac{b-a}{2} \\
l'(v_5v_9')=l'(v_2v_{9})= l'(v_{9},v_{3})&= \frac{e}{2} & l'(v_1v_2)=l'(v_4v_5)=l'(v_3v_{11}) &= a.
\end{align*}

\begin{figure}[ht]
    \centering
    \begin{tikzpicture}[line cap=round,line join=round,>=triangle 45,x=1.0cm,y=1.0cm]
\definecolor{wrwrwr}{rgb}{0.3803921568627451,0.3803921568627451,0.3803921568627451}
\definecolor{sexdts}{rgb}{0.1803921568627451,0.49019607843137253,0.19607843137254902}
\definecolor{ffqqqq}{rgb}{1.,0.,0.}
\definecolor{rvwvcq}{rgb}{0.08235294117647059,0.396078431372549,0.7529411764705882}
\clip(-2.,-0.9) rectangle (2.,4.);
\draw (-0.5,0.)-- (0.5,0.);
\draw (-0.5,2.)-- (0.5,2.);
\draw (-0.5,3.)-- (0.5,3.);
\draw [shift={(-0.25,1.)},color=rvwvcq]  plot[domain=-0.9272952180016123:0.9272952180016122,variable=\t]({1.*1.25*cos(\t r)+0.*1.25*sin(\t r)},{0.*1.25*cos(\t r)+1.*1.25*sin(\t r)});
\draw [shift={(0.5,1.)},color=ffqqqq]  plot[domain=-1.5707963267948966:1.5707963267948966,variable=\t]({1.*1.*cos(\t r)+0.*1.*sin(\t r)},{0.*1.*cos(\t r)+1.*1.*sin(\t r)});
\draw (-0.5,-2.)-- (0.5,-2.);
\draw [shift={(0.2333333333333335,2.5)},color=sexdts]  plot[domain=-1.080839000541168:1.0808390005411685,variable=\t]({1.*0.5666666666666665*cos(\t r)+0.*0.5666666666666665*sin(\t r)},{0.*0.5666666666666665*cos(\t r)+1.*0.5666666666666665*sin(\t r)});
\draw [shift={(1.4657697822512776,1.5)},color=ffqqqq]  plot[domain=2.4897854997887947:3.7933998073907915,variable=\t]({1.*2.4727011216101786*cos(\t r)+0.*2.4727011216101786*sin(\t r)},{0.*2.4727011216101786*cos(\t r)+1.*2.4727011216101786*sin(\t r)});
\draw [shift={(0.14435667947216052,1.5)},color=rvwvcq]  plot[domain=1.9765323738923735:4.3066529332872125,variable=\t]({1.*1.632542658058401*cos(\t r)+0.*1.632542658058401*sin(\t r)},{0.*1.632542658058401*cos(\t r)+1.*1.632542658058401*sin(\t r)});
\draw [color=ffqqqq] (-0.5,-2.)-- (-1.,-3.);
\draw [color=rvwvcq] (-0.5,-2.)-- (-1.5,-1.);
\draw [color=ffqqqq] (-0.5,2.)-- (-1.,3.);
\draw [color=rvwvcq] (-0.5,2.)-- (-1.5,3.);
\draw [color=ffqqqq] (0.5,-2.)-- (1.5,-2.);
\draw [color=rvwvcq] (0.5,-2.)-- (1.,-3.);
\draw [color=sexdts] (0.5,-2.)-- (0.8,-1.0123047916035581);
\draw [color=sexdts] (0.5,0.)-- (0.7962429260233621,0.8144444886021116);
\draw [color=ffqqqq] (0.5,3.)-- (1.5,3.);
\draw [color=rvwvcq] (0.5,3.)-- (0.9949411571371803,3.473667079314102);
\draw (2.302153011953767,1.7453380218400287) node[anchor=north west] {$\Gamma'$};
\draw (2.3131450563766167,-1.8913527146761007) node[anchor=north west] {$T$};
\draw [->] (2.398333400653702,0.5662062130719998) -- (2.398333400653702,-1.2211093707448013);
\begin{scriptsize}
\draw [fill=black] (-0.5,0.) circle (0.5pt);
\draw[color=black] (-0.49807030476718667,-0.17850861351833328) node {$v_1$};
\draw [fill=black] (0.5,0.) circle (0.5pt);
\draw[color=black] (0.4967097155007086,-0.21574435484784996) node {$v_2$};
\draw[color=black] (0.010311749789610594,-0.24677413928911385) node {$a$};
\draw [fill=black] (-0.5,2.) circle (0.5pt);
\draw [fill=black] (0.5,2.) circle (0.5pt);
\draw[color=black] (0.4967097155007086,1.794985676946051) node {$v_3$};
\draw[color=black] (-0.017168361267513588,1.7639558925047871) node {$a$};
\draw [fill=black] (-0.5,3.) circle (0.5pt);
\draw[color=black] (-0.5035663269786116,3.22235576124419) node {$v_4$};
\draw [fill=black] (0.5,3.) circle (0.5pt);
\draw[color=black] (0.41976540454076083,3.197531933691179) node {$v_5$};
\draw[color=black] (-0.022664383478938413,2.7444970808487263) node {$a$};
\draw [color=rvwvcq] (1.,1.)-- ++(-1.0pt,-1.0pt) -- ++(2.0pt,2.0pt) ++(-2.0pt,0) -- ++(2.0pt,-2.0pt);
\draw[color=rvwvcq] (1.255160780677336,1.0999185054617395) node {$v_9$};
\draw[color=rvwvcq] (1.0545559699603295,1.4536580480921482) node {$e$};
\draw [color=ffqqqq] (1.5,1.)-- ++(-1.0pt,-1.0pt) -- ++(2.0pt,2.0pt) ++(-2.0pt,0) -- ++(2.0pt,-2.0pt);
\draw[color=ffqqqq] (1.8047630018198193,1.0750946779087285) node {$v_{10}$};
\draw[color=ffqqqq] (1.268900836205898,1.875663116493337) node {$f$};
\draw [fill=black] (-0.5,-2.) circle (0.5pt);
\draw [fill=black] (0.5,-2.) circle (0.5pt);
\draw[color=black] (0.0267998164238851,-2.257504171083015) node {$a$};
\draw [color=sexdts] (0.8,2.5)-- ++(-1.0pt,-1.0pt) -- ++(2.0pt,2.0pt) ++(-2.0pt,0) -- ++(2.0pt,-2.0pt);
\draw[color=sexdts] (0.5956381153063556,2.5024647622068676) node {$v_6$};
\draw[color=sexdts] (1.2139406140916496,2.396963495106571) node {$b-a$};
\draw [color=ffqqqq] (-1.006924325949443,1.5058893192882263)-- ++(-1.0pt,-1.0pt) -- ++(2.0pt,2.0pt) ++(-2.0pt,0) -- ++(2.0pt,-2.0pt);
\draw[color=ffqqqq] (-1.256521369943814,1.571571228968951) node {$v_8$};
\draw[color=ffqqqq] (-0.7206592043298925,1.0688887210204756) node {$d$};
\draw [color=rvwvcq] (-1.4879256850319909,1.4708484492843676)-- ++(-1.0pt,-1.0pt) -- ++(2.0pt,2.0pt) ++(-2.0pt,0) -- ++(2.0pt,-2.0pt);
\draw[color=rvwvcq] (-1.7511633689720492,1.596395056521962) node {$v_7$};
\draw[color=rvwvcq] (-1.237285292203827,1.0192410659144533) node {$c$};
\draw [fill=black] (-1.,-3.) circle (0.5pt);
\draw[color=ffqqqq] (-0.5282984269300233,-2.6919211532607092) node {$\frac{d}{2}$};
\draw [fill=black] (-1.5,-1.) circle (0.5pt);
\draw[color=rvwvcq] (-1.1218688257639056,-1.5872608271517146) node {$\frac{c}{2}$};
\draw [fill=black] (-1.,3.) circle (0.5pt);
\draw[color=black] (-0.9817202593725721,3.22235576124419) node {$v_{8}'$};
\draw [fill=black] (-1.5,3.) circle (0.5pt);
\draw[color=black] (-1.492850325035082,3.234767675020696) node {$v_{7}'$};
\draw [fill=black] (1.5,-2.) circle (0.5pt);
\draw[color=ffqqqq] (1.3788212804343947,-2.245092257306509) node {$\frac{f}{2}$};
\draw [fill=black] (1.,-3.) circle (0.5pt);
\draw[color=rvwvcq] (0.6863224817948654,-2.8036283772492596) node {$\frac{e}{2}$};
\draw [fill=black] (0.8,-1.0123047916035581) circle (0.5pt);
\draw[color=sexdts] (1.0820360810174536,-1.5003774307161755) node {$\frac{b-a}{2}$};
\draw [fill=black] (0.7962429260233621,0.8144444886021116) circle (0.5pt);
\draw[color=black] (0.716550603957702,1.1371542467912563) node {$v_6'$};
\draw[color=sexdts] (0.5104497710292706,0.7089432215018143) node {$\frac{b-a}{2}$};
\draw [fill=black] (1.5,3.) circle (0.5pt);
\draw[color=black] (1.815755046242669,3.0734127959261235) node {$v_{10}'$};
\draw[color=ffqqqq] (1.5217178579314403,3.439564252333038) node {$\frac{f}{2}$};
\draw [fill=black] (0.9949411571371803,3.473667079314102) circle (0.5pt);
\draw[color=black] (1.0243278477974929,3.731244226080918) node {$v_9'$};
\draw[color=rvwvcq] (0.7577707705433883,3.5885072176511046) node {$\frac{e}{2}$};
\draw[color=black] (2.6648904779078064,-0.22195031173610275) node {$\varphi$};
\draw[color=wrwrwr] (-2.2705374679516965,4.978641560619727) node {$Curve$};
\end{scriptsize}
\end{tikzpicture}
    \begin{center}
        \textnormal{Figure 5.2. The model $(G',l')$ of $\Gamma'$}
    \end{center}
\end{figure}
Choose $T$ to be the metric tree with model $(T',t')$ in Figure 5.3, where the tree $T'$ is given by its vertex set $V(T')= \left\lbrace w_1,w_2,w_6,w_7,\ldots, w_{10} \right\rbrace$, and edge set $E(T')= \left\lbrace w_1w_2,w_1w_7,w_1w_8,w_2w_6,w_2w_9,w_2w_{10} \right\rbrace$. The length map $t'$ on $T'$ is given by
\begin{equation*}
\begin{aligned}[c]
t'(w_2w_1) &= a\\
t'(w_2w_6) &=  \frac{b-a}{2}\\ 
t'(w_1w_7) &=  \frac{c}{2}
\end{aligned}
\qquad
\begin{aligned}[c]
t'(w_1w_8) &=  \frac{d}{2}\\
t'(w_2w_9) &=  \frac{e}{2}\\
t'(w_2w_{10}) &=  \frac{f}{2}.\\
\end{aligned}
\end{equation*}
\begin{figure}[H]
    \centering
   \begin{tikzpicture}[line cap=round,line join=round,>=triangle 45,x=1.0cm,y=1.0cm]
\definecolor{wrwrwr}{rgb}{0.3803921568627451,0.3803921568627451,0.3803921568627451}
\definecolor{sexdts}{rgb}{0.1803921568627451,0.49019607843137253,0.19607843137254902}
\definecolor{ffqqqq}{rgb}{1.,0.,0.}
\definecolor{rvwvcq}{rgb}{0.08235294117647059,0.396078431372549,0.7529411764705882}
\clip(-2.,-4.) rectangle (2.2,-0.5);
\draw (-0.5,0.)-- (0.5,0.);
\draw (-0.5,3.)-- (0.5,3.);
\draw [shift={(-0.25,1.)},color=rvwvcq]  plot[domain=-0.9272952180016123:0.9272952180016122,variable=\t]({1.*1.25*cos(\t r)+0.*1.25*sin(\t r)},{0.*1.25*cos(\t r)+1.*1.25*sin(\t r)});
\draw [shift={(0.5,1.)},color=ffqqqq]  plot[domain=-1.5707963267948966:1.5707963267948966,variable=\t]({1.*1.*cos(\t r)+0.*1.*sin(\t r)},{0.*1.*cos(\t r)+1.*1.*sin(\t r)});
\draw (-0.5,-2.)-- (0.5,-2.);
\draw [shift={(0.2333333333333335,2.5)},color=sexdts]  plot[domain=-1.080839000541168:1.0808390005411685,variable=\t]({1.*0.5666666666666665*cos(\t r)+0.*0.5666666666666665*sin(\t r)},{0.*0.5666666666666665*cos(\t r)+1.*0.5666666666666665*sin(\t r)});
\draw [shift={(1.4657697822512776,1.5)},color=ffqqqq]  plot[domain=2.4897854997887947:3.7933998073907915,variable=\t]({1.*2.4727011216101786*cos(\t r)+0.*2.4727011216101786*sin(\t r)},{0.*2.4727011216101786*cos(\t r)+1.*2.4727011216101786*sin(\t r)});
\draw [shift={(0.14435667947216052,1.5)},color=rvwvcq]  plot[domain=1.9765323738923735:4.3066529332872125,variable=\t]({1.*1.632542658058401*cos(\t r)+0.*1.632542658058401*sin(\t r)},{0.*1.632542658058401*cos(\t r)+1.*1.632542658058401*sin(\t r)});
\draw [color=ffqqqq] (-0.5,-2.)-- (-1.,-3.);
\draw [color=rvwvcq] (-0.5,-2.)-- (-1.5,-1.);
\draw [color=ffqqqq] (0.5,-2.)-- (1.5,-2.);
\draw [color=rvwvcq] (0.5,-2.)-- (1.,-3.);
\draw [color=sexdts] (0.5,-2.)-- (0.8,-1.0123047916035581);
\draw (2.302153011953767,1.7453380218400287) node[anchor=north west] {$\Gamma'$};
\draw (2.3131450563766167,-1.8913527146761007) node[anchor=north west] {$T$};
\draw [->] (2.398333400653702,0.5662062130719998) -- (2.398333400653702,-1.2211093707448013);
\begin{scriptsize}
\draw [fill=black] (-0.5,0.) circle (0.5pt);
\draw[color=black] (-0.49807030476718667,-0.17850861351833328) node {$v_1$};
\draw [fill=black] (0.5,0.) circle (0.5pt);
\draw[color=black] (0.4967097155007086,-0.21574435484784996) node {$v_2$};
\draw[color=black] (0.010311749789610594,-0.24677413928911385) node {$a$};
\draw [fill=black] (0.5,2.) circle (0.5pt);
\draw[color=black] (0.4967097155007086,1.794985676946051) node {$v_3$};
\draw [fill=black] (-0.5,3.) circle (0.5pt);
\draw[color=black] (-0.5035663269786116,3.22235576124419) node {$v_4$};
\draw [fill=black] (0.5,3.) circle (0.5pt);
\draw[color=black] (0.41976540454076083,3.197531933691179) node {$v_5$};
\draw[color=black] (-0.022664383478938413,2.7444970808487263) node {$a$};
\draw [color=rvwvcq] (1.,1.)-- ++(-1.0pt,-1.0pt) -- ++(2.0pt,2.0pt) ++(-2.0pt,0) -- ++(2.0pt,-2.0pt);
\draw[color=rvwvcq] (1.255160780677336,1.0999185054617395) node {$v_9$};
\draw[color=rvwvcq] (1.0545559699603295,1.4536580480921482) node {$e$};
\draw [color=ffqqqq] (1.5,1.)-- ++(-1.0pt,-1.0pt) -- ++(2.0pt,2.0pt) ++(-2.0pt,0) -- ++(2.0pt,-2.0pt);
\draw[color=ffqqqq] (1.8047630018198193,1.0750946779087285) node {$v_{10}$};
\draw[color=ffqqqq] (1.268900836205898,1.875663116493337) node {$f$};
\draw [fill=black] (-0.5,-2.) circle (0.5pt);
\draw[color=black] (-0.3991419049615397,-1.8292931457935728) node {$w_1$};
\draw [fill=black] (0.5,-2.) circle (0.5pt);
\draw[color=black] (0.32083700473511384,-1.8292931457935728) node {$w_2$};
\draw[color=black] (0.0267998164238851,-2.2078565159769927) node {$a$};
\draw [color=sexdts] (0.8,2.5)-- ++(-1.0pt,-1.0pt) -- ++(2.0pt,2.0pt) ++(-2.0pt,0) -- ++(2.0pt,-2.0pt);
\draw[color=sexdts] (0.5956381153063556,2.5024647622068676) node {$v_6$};
\draw[color=sexdts] (1.2139406140916496,2.396963495106571) node {$b-a$};
\draw [color=ffqqqq] (-1.006924325949443,1.5058893192882263)-- ++(-1.0pt,-1.0pt) -- ++(2.0pt,2.0pt) ++(-2.0pt,0) -- ++(2.0pt,-2.0pt);
\draw[color=ffqqqq] (-1.256521369943814,1.4970997463099176) node {$v_8$};
\draw[color=ffqqqq] (-0.7206592043298925,1.0688887210204756) node {$d$};
\draw [color=rvwvcq] (-1.4879256850319909,1.4708484492843676)-- ++(-1.0pt,-1.0pt) -- ++(2.0pt,2.0pt) ++(-2.0pt,0) -- ++(2.0pt,-2.0pt);
\draw[color=rvwvcq] (-1.7456673467606243,1.484687832533412) node {$v_7$};
\draw[color=rvwvcq] (-1.237285292203827,1.0192410659144533) node {$c$};
\draw [fill=black] (-1.,-3.) circle (0.5pt);
\draw[color=black] (-1.1411049035038925,-3.1697798336561736) node {$w_8$};
\draw[color=ffqqqq] (-0.5282984269300233,-2.6919211532607092) node {$\frac{d}{2}$};
\draw [fill=black] (-1.5,-1.) circle (0.5pt);
\draw[color=black] (-1.723683257914925,-0.9356353538851725) node {$w_7$};
\draw[color=rvwvcq] (-1.1218688257639056,-1.5872608271517146) node {$\frac{c}{2}$};
\draw [fill=black] (1.5,-2.) circle (0.5pt);
\draw[color=black] (1.8762112905683421,-1.9782361111116396) node {$w_{10}$};
\draw[color=ffqqqq] (1.1150122142860026,-2.2326803435300038) node {$\frac{f}{2}$};
\draw [fill=black] (1.,-3.) circle (0.5pt);
\draw[color=black] (1.1287522698145647,-3.219427488762196) node {$w_9$};
\draw[color=rvwvcq] (0.6808264595834406,-2.7167449808137207) node {$\frac{e}{2}$};
\draw [fill=black] (0.8,-1.0123047916035581) circle (0.5pt);
\draw[color=black] (0.9253994479918457,-0.8115162161201168) node {$w_6$};
\draw[color=sexdts] (1.0820360810174536,-1.5003774307161755) node {$\frac{b-a}{2}$};
\draw[color=black] (2.6648904779078064,-0.22195031173610275) node {$\varphi$};
\draw[color=wrwrwr] (-2.2705374679516965,4.978641560619727) node {$Curve$};
\end{scriptsize}
\end{tikzpicture}
    \begin{center}
        \textnormal{Figure 5.3. The model $(T',t')$ of $T$}
    \end{center}
\end{figure}
Let $\psi:V(G') \to V(T)$ the map on the set of vertices given by $v_1, v_4,v_{11}' \mapsto w_1$, $v_2,v_3,v_5 \mapsto w_2$, and $v_i, v_i' \mapsto w_i$ for $i=6,8,9,10,11$. The function $\psi$ satisfies the condition in Lemma 3.1, and so, there exist a unique continuous map $\varphi: \Gamma' \to T$, shown in Figure 6, such that $\varphi|_{V(G')} = \psi$, and $\varphi$ is linear on each edge $e' \in E(G')$ with slope $t'(e)/l'(e')$, where $e = \varphi(e') \in E(T)$ with endpoints $\psi(v)$ and $\psi(w)$. The tropical morphism $\varphi:\Gamma' \to T$ is of degree $3$ essentially because of the reasoning in (i) and (ii). \par 
\begin{Remark} \textnormal{The constructions of tropical morphisms of the remaining metric graphs are done similarly as for the metric graph in the Case 1.1.A. In order to avoid tedious writing, we give the construction of a model, a tropical modification, a metric tree, and a tropical morphism, using only figures from now on. The vertices labeled with a small $\times$ are the 'midpoints' of the edges i.e., when subdividing an edge $e$ into $e_1$ and $e_2$ then both lengths of $e_1$ and $e_2$ are equal to the half of the length of edge $e$.}
\end{Remark}
\begin{figure}[H]
    \centering
    
\begin{tikzpicture}[line cap=round,line join=round,>=triangle 45,x=1.0cm,y=1.0cm]
\definecolor{wrwrwr}{rgb}{0.3803921568627451,0.3803921568627451,0.3803921568627451}
\definecolor{cqcqcq}{rgb}{0.7529411764705882,0.7529411764705882,0.7529411764705882}
\definecolor{sexdts}{rgb}{0.1803921568627451,0.49019607843137253,0.19607843137254902}
\definecolor{ffqqqq}{rgb}{1.,0.,0.}
\definecolor{rvwvcq}{rgb}{0.08235294117647059,0.396078431372549,0.7529411764705882}
\definecolor{wrwrwr}{rgb}{0.3803921568627451,0.3803921568627451,0.3803921568627451}
\definecolor{cqcqcq}{rgb}{0.7529411764705882,0.7529411764705882,0.7529411764705882}
\definecolor{sexdts}{rgb}{0.1803921568627451,0.49019607843137253,0.19607843137254902}
\definecolor{ffqqqq}{rgb}{1.,0.,0.}
\definecolor{rvwvcq}{rgb}{0.08235294117647059,0.396078431372549,0.7529411764705882}
\clip(-2.,-3.5) rectangle (3.,4.5);
\draw (-0.5,0.)-- (0.5,0.);
\draw (-0.5,2.)-- (0.5,2.);
\draw (-0.5,3.)-- (0.5,3.);
\draw [shift={(-0.25,1.)},color=rvwvcq]  plot[domain=-0.9272952180016123:0.9272952180016122,variable=\t]({1.*1.25*cos(\t r)+0.*1.25*sin(\t r)},{0.*1.25*cos(\t r)+1.*1.25*sin(\t r)});
\draw [shift={(0.5,1.)},color=ffqqqq]  plot[domain=-1.5707963267948966:1.5707963267948966,variable=\t]({1.*1.*cos(\t r)+0.*1.*sin(\t r)},{0.*1.*cos(\t r)+1.*1.*sin(\t r)});
\draw (-0.5,-2.)-- (0.5,-2.);
\draw [shift={(0.2333333333333335,2.5)},color=sexdts]  plot[domain=-1.080839000541168:1.0808390005411685,variable=\t]({1.*0.5666666666666665*cos(\t r)+0.*0.5666666666666665*sin(\t r)},{0.*0.5666666666666665*cos(\t r)+1.*0.5666666666666665*sin(\t r)});
\draw [shift={(1.4657697822512776,1.5)},color=ffqqqq]  plot[domain=2.4897854997887947:3.7933998073907915,variable=\t]({1.*2.4727011216101786*cos(\t r)+0.*2.4727011216101786*sin(\t r)},{0.*2.4727011216101786*cos(\t r)+1.*2.4727011216101786*sin(\t r)});
\draw [shift={(0.14435667947216052,1.5)},color=rvwvcq]  plot[domain=1.9765323738923735:4.3066529332872125,variable=\t]({1.*1.632542658058401*cos(\t r)+0.*1.632542658058401*sin(\t r)},{0.*1.632542658058401*cos(\t r)+1.*1.632542658058401*sin(\t r)});
\draw [color=ffqqqq] (-0.5,-2.)-- (-1.,-3.);
\draw [color=rvwvcq] (-0.5,-2.)-- (-1.5,-1.);
\draw [color=ffqqqq] (-0.5,2.)-- (-1.,3.);
\draw [color=rvwvcq] (-0.5,2.)-- (-1.5,3.);
\draw [color=ffqqqq] (0.5,-2.)-- (1.5,-2.);
\draw [color=rvwvcq] (0.5,-2.)-- (1.,-3.);
\draw [color=sexdts] (0.5,-2.)-- (0.8,-1.0123047916035581);
\draw [line width=0.4pt,dash pattern=on 2pt off 2pt,color=cqcqcq] (0.5,-2.)-- (0.5,3.);
\draw [line width=0.4pt,dash pattern=on 2pt off 2pt,color=cqcqcq] (-0.5,-2.)-- (-0.5,3.);
\draw [line width=0.4pt,dash pattern=on 2pt off 2pt,color=cqcqcq] (-1.,-3.)-- (-1.,3.);
\draw [line width=0.4pt,dash pattern=on 2pt off 2pt,color=cqcqcq] (-1.5,-1.)-- (-1.5,3.);
\draw [line width=0.4pt,dash pattern=on 2pt off 2pt,color=cqcqcq] (0.8,-1.0123047916035581)-- (0.8,2.5);
\draw [color=sexdts] (0.5,0.)-- (0.7962429260233621,0.8144444886021116);
\draw [color=ffqqqq] (0.5,3.)-- (1.5,3.);
\draw [color=rvwvcq] (0.5,3.)-- (0.9949411571371803,3.473667079314102);
\draw [line width=0.4pt,dash pattern=on 2pt off 2pt,color=cqcqcq] (1.,-3.)-- (0.9949411571371803,3.473667079314102);
\draw [line width=0.4pt,dash pattern=on 2pt off 2pt,color=cqcqcq] (1.5,-2.)-- (1.5,3.);
\draw (2.302153011953767,1.7453380218400287) node[anchor=north west] {$\Gamma'$};
\draw (2.3131450563766167,-1.8913527146761007) node[anchor=north west] {$T$};
\draw [->] (2.398333400653702,0.5662062130719998) -- (2.398333400653702,-1.2211093707448013);
\begin{scriptsize}
\draw [fill=black] (-0.5,0.) circle (0.5pt);
\draw[color=black] (-0.49807030476718667,-0.24056818240086106) node {$v_1$};
\draw [fill=black] (0.5,0.) circle (0.5pt);
\draw[color=black] (0.4967097155007086,-0.21574435484784996) node {$v_2$};
\draw[color=black] (0.010311749789610594,-0.24677413928911385) node {$a$};
\draw [fill=black] (-0.5,2.) circle (0.5pt);
\draw [fill=black] (0.5,2.) circle (0.5pt);
\draw[color=black] (0.4967097155007086,1.794985676946051) node {$v_3$};
\draw[color=black] (-0.017168361267513588,1.7639558925047871) node {$a$};
\draw [fill=black] (-0.5,3.) circle (0.5pt);
\draw[color=black] (-0.5035663269786116,3.309239157679729) node {$v_4$};
\draw [fill=black] (0.5,3.) circle (0.5pt);
\draw[color=black] (-0.022664383478938413,2.7444970808487263) node {$a$};
\draw [color=rvwvcq] (1.,1.)-- ++(-1.0pt,-1.0pt) -- ++(2.0pt,2.0pt) ++(-2.0pt,0) -- ++(2.0pt,-2.0pt);
\draw[color=rvwvcq] (1.1095161920745777,1.3543627378801035) node {$e$};
\draw [color=ffqqqq] (1.5,1.)-- ++(-1.0pt,-1.0pt) -- ++(2.0pt,2.0pt) ++(-2.0pt,0) -- ++(2.0pt,-2.0pt);
\draw[color=ffqqqq] (1.268900836205898,1.875663116493337) node {$f$};
\draw [fill=black] (-0.5,-2.) circle (0.5pt);
\draw [fill=black] (0.5,-2.) circle (0.5pt);
\draw[color=black] (0.0267998164238851,-2.257504171083015) node {$a$};
\draw [color=sexdts] (0.8,2.5)-- ++(-1.0pt,-1.0pt) -- ++(2.0pt,2.0pt) ++(-2.0pt,0) -- ++(2.0pt,-2.0pt);
\draw[color=sexdts] (1.2139406140916496,2.396963495106571) node {$b-a$};
\draw [color=ffqqqq] (-1.006924325949443,1.5058893192882263)-- ++(-1.0pt,-1.0pt) -- ++(2.0pt,2.0pt) ++(-2.0pt,0) -- ++(2.0pt,-2.0pt);
\draw[color=ffqqqq] (-0.7206592043298925,1.0688887210204756) node {$d$};
\draw [color=rvwvcq] (-1.4879256850319909,1.4708484492843676)-- ++(-1.0pt,-1.0pt) -- ++(2.0pt,2.0pt) ++(-2.0pt,0) -- ++(2.0pt,-2.0pt);
\draw[color=rvwvcq] (-1.237285292203827,1.0192410659144533) node {$c$};
\draw [fill=black] (-1.,-3.) circle (0.5pt);
\draw[color=ffqqqq] (-0.5282984269300233,-2.6919211532607092) node {$\frac{d}{2}$};
\draw [fill=black] (-1.5,-1.) circle (0.5pt);
\draw[color=rvwvcq] (-1.1053807591296312,-0.6563672939137973) node {$\frac{c}{2}$};
\draw [fill=black] (-1.,3.) circle (0.5pt);
\draw [fill=black] (-1.5,3.) circle (0.5pt);
\draw [fill=black] (1.5,-2.) circle (0.5pt);
\draw[color=ffqqqq] (1.3788212804343947,-2.245092257306509) node {$\frac{f}{2}$};
\draw [fill=black] (1.,-3.) circle (0.5pt);
\draw[color=rvwvcq] (0.6863224817948654,-2.8036283772492596) node {$\frac{e}{2}$};
\draw [fill=black] (0.8,-1.0123047916035581) circle (0.5pt);
\draw[color=sexdts] (1.0820360810174536,-1.5003774307161755) node {$\frac{b-a}{2}$};
\draw [fill=black] (0.7962429260233621,0.8144444886021116) circle (0.5pt);
\draw[color=sexdts] (0.5049537488178458,0.6841193939488033) node {$\frac{b-a}{2}$};
\draw [fill=black] (1.5,3.) circle (0.5pt);
\draw[color=ffqqqq] (1.5437019467771396,3.3650927696740047) node {$\frac{f}{2}$};
\draw [fill=black] (0.9949411571371803,3.473667079314102) circle (0.5pt);
\draw[color=rvwvcq] (0.900667348040434,3.662978700310138) node {$\frac{e}{2}$};
\draw[color=black] (2.6648904779078064,-0.22195031173610275) node {$\varphi$};
\draw[color=wrwrwr] (-2.2705374679516965,4.978641560619727) node {$Curve$};
\end{scriptsize}
\end{tikzpicture}
    \begin{center}
        \textnormal{Figure 6. The tropical morphism $\varphi:\Gamma' \to T$ }
    \end{center}
\end{figure}
\begin{figure}[H]
    \centering
    \begin{tikzpicture}[line cap=round,line join=round,>=triangle 45,x=1.0cm,y=1.0cm]
\definecolor{wrwrwr}{rgb}{0.3803921568627451,0.3803921568627451,0.3803921568627451}
\clip(-1.5,-1.5) rectangle (3.,2.5);
\draw (0.,2.)-- (0.,0.);
\draw (0.,0.)-- (2.5,0.);
\draw [shift={(0.,1.)}] plot[domain=1.5707963267948966:4.71238898038469,variable=\t]({1.*1.*cos(\t r)+0.*1.*sin(\t r)},{0.*1.*cos(\t r)+1.*1.*sin(\t r)});
\draw [shift={(1.25,0.)}] plot[domain=3.141592653589793:6.283185307179586,variable=\t]({1.*1.25*cos(\t r)+0.*1.25*sin(\t r)},{0.*1.25*cos(\t r)+1.*1.25*sin(\t r)});
\draw (0.,2.)-- (2.5,0.);
\begin{scriptsize}
\draw [fill=black] (0.,0.) circle (0.5pt);
\draw[color=black] (-0.26153711690148784,-0.18535025360631566) node {$v_1$};
\draw [fill=black] (0.,2.) circle (0.5pt);
\draw[color=black] (0.005022400770458867,2.2228771129471316) node {$v_4$};
\draw [fill=black] (2.5,0.) circle (0.5pt);
\draw[color=black] (2.753342945043289,-0.09343317854702377) node {$v_3$};
\draw[color=black] (-0.29370809317224,0.9774007458937266) node {$d$};
\draw[color=black] (1.250498767823865,0.18691390038381644) node {$e$};
\draw[color=black] (-1.2588373812948057,0.9222505008581515) node {$c$};
\draw[color=black] (1.241307060317936,-0.9528578303514028) node {$f$};
\draw[color=black] (1.2688821828357235,1.3083022161071773) node {$b$};
\draw[color=wrwrwr] (-2.996070099915424,4.984985218478853) node {$Curve$};
\end{scriptsize}
\end{tikzpicture}
    \begin{center}
        \textnormal{Figure 7. The essential model $(G,l)$ of $\Gamma_1$}
    \end{center}
\end{figure}
\textbf{Case 1.3.}  Consider the metric graph $\Gamma_1$ with essential model $(G,l)$ in Figure 7, where $b, c,d,e$, and $f$ are real positive numbers. The model $(G_1,l_1)$ which is obtained by subdividing $(G,l)$ is shown in Figure 7.1. The tropical modification $\Gamma_1'$, the metric tree $T_1$ with models $(G_1',l_1')$, $(T_1',t_1')$ is given in Figure 7.2, 7.3, respectively. The construction of the tropical morphism $\varphi_1: \Gamma_1' \to T_1$ of degree $3$ is depicted in Figure 8. 
\begin{figure}[H]
    \centering

    \begin{center}
        \textnormal{Figure 8. The tropical morphism $\varphi_1: \Gamma_1' \to T_1$}
    \end{center}
\end{figure}

\textbf{Case 1.4.}  Consider the metric graph $\Gamma_2$ with essential model in Figure 9, where $a,b,c,d$, and $e$ are real positive numbers such that $b>a$. Note that if $a=b$,  then $\Gamma_2$ is a hyperelliptic metric graph. The model $(G_1,l_1)$ that is obtained by subdividing $(G,l)$ is shown in Figure 9.1.
\begin{figure}[H]
    \centering
   %
    \begin{center}
        \textnormal{Figure 9. The essential model $(G,l)$ of $\Gamma_2$}
    \end{center}
\end{figure}
\begin{figure}[H]
    \centering
    %
    \begin{center}
        \textnormal{Figure 9.1. The model $(G_1,l_1)$ of $\Gamma_2$}
    \end{center}
\end{figure}
The tropical modification $\Gamma_2'$, the metric tree $T_2$ with models $(G_2',l_2')$, $(T_2',t_2')$ is given in Figure 9.2, 9.3 respectively. 
\begin{figure}[H]
    \centering
    %
    \begin{center}
        \textnormal{Figure 9.3. The model $(T_2',t_2')$ of $T_2$}
    \end{center}
\end{figure}
The construction of the tropical morphism $\varphi_2: \Gamma_2' \to T_2$ of degree 3 is depicted in Figure 10.
\begin{figure}[H]
    \centering
    %
    \begin{center}
        \textnormal{Figure 10. The tropical morphism $\varphi_2:\Gamma_2' \to T_2$}
    \end{center}
\end{figure}
\textbf{Case 1.5.}  Consider the metric graph $\Gamma_3$ with essential model $(G,l)$ in Figure 11, where $a,b,c$, and $e$ are real positive numbers such that $b>a$. Note that if $a=b$, then $\Gamma_2$ is a hyperelliptic metric graph. The model $(G_1,l_1)$ which is obtained by subdividing $(G,l)$ is shown in Figure 11.1. The tropical modification $\Gamma_3'$, the metric tree $T_3$ with models $(G_3',l_3')$, $(T_3',t_3')$ is given in Figure 11.2, 11.3, respectively. The construction of the tropical morphism $\varphi_3: \Gamma_3' \to T_3$ of degree $3$ is depicted in Figure 12.
\begin{figure}[H]
    \centering
    \begin{tikzpicture}[line cap=round,line join=round,>=triangle 45,x=1.0cm,y=1.0cm]
\definecolor{uuuuuu}{rgb}{0.26666666666666666,0.26666666666666666,0.26666666666666666}
\clip(-2.,-1.) rectangle (3.,1.5);
\draw(-0.6227249882937179,0.) circle (0.6227249882937179cm);
\draw(2.,0.) circle (0.5cm);
\draw [shift={(0.75,-0.08626957576535084)}] plot[domain=0.11452278703110289:3.02706986655869,variable=\t]({1.*0.7549453223265465*cos(\t r)+0.*0.7549453223265465*sin(\t r)},{0.*0.7549453223265465*cos(\t r)+1.*0.7549453223265465*sin(\t r)});
\draw (0.,0.)-- (1.5,0.);
\begin{scriptsize}
\draw [fill=uuuuuu] (0.,0.) circle (0.5pt);
\draw[color=uuuuuu] (-0.17461940465933834,0.07273129131192979) node {$v_1$};
\draw[color=black] (-0.9738180228938504,-0.04275983849074533) node {$c$};
\draw [fill=black] (1.5,0.) circle (0.5pt);
\draw[color=black] (1.737913704872962,0.06349200092771579) node {$v_2$};
\draw[color=black] (2.29689077311791,-0.09819558079602939) node {$e$};
\draw[color=black] (0.7262114078015278,0.8996477806990835) node {$b$};
\draw[color=black] (0.7077328270330998,-0.23678493655923952) node {$a$};
\end{scriptsize}
\end{tikzpicture}
 
    \begin{center}
        \textnormal{Figure 11. The essential model $(G,l)$ of $\Gamma_3$}
    \end{center}
\end{figure}
\begin{figure}[H]
    \centering
    \begin{tikzpicture}[line cap=round,line join=round,>=triangle 45,x=1.0cm,y=1.0cm]
\definecolor{qqqqff}{rgb}{0.,0.,1.}
\definecolor{wrwrwr}{rgb}{0.3803921568627451,0.3803921568627451,0.3803921568627451}
\definecolor{rvwvcq}{rgb}{0.08235294117647059,0.396078431372549,0.7529411764705882}
\definecolor{ttzzqq}{rgb}{0.2,0.6,0.}
\clip(-2.5,-1.) rectangle (3.,3.5);
\draw (-0.5,-2.)-- (1.,-2.);
\draw [color=ttzzqq] (1.,-2.)-- (1.5,-3.);
\draw [color=rvwvcq] (1.,-2.)-- (2.,-1.5);
\draw [color=rvwvcq] (-0.5,-2.)-- (-1.5439241262940946,-1.99957780037612);
\draw (-0.5,2.5)-- (1.,2.5);
\draw [shift={(0.25,1.5)},color=ttzzqq]  plot[domain=-0.9272952180016123:0.9272952180016122,variable=\t]({1.*1.25*cos(\t r)+0.*1.25*sin(\t r)},{0.*1.25*cos(\t r)+1.*1.25*sin(\t r)});
\draw (1.,0.5)-- (-0.5061173037923797,0.49515013832992044);
\draw [color=rvwvcq] (-0.5061173037923797,0.49515013832992044)-- (-1.5,0.5);
\draw [color=rvwvcq] (1.,2.5)-- (2.,3.);
\draw [color=ttzzqq] (1.,0.5)-- (1.5,-0.5);
\draw [color=rvwvcq] (1.5026573768370155,0.5039366400034262) circle (0.5026727918071416cm);
\draw (3.4754684863824683,0.9442011672970076) node[anchor=north west] {$\Gamma_3'$};
\draw (3.505405221646941,-1.9496832416019991) node[anchor=north west] {$T_3'$};
\draw [->,line width=0.4pt] (3.615173250950007,0.1359093151562506) -- (3.635131074459655,-1.4008430950866702);
\draw (3.954456250614029,-0.32312062556566085) node[anchor=north west] {$\varphi_3$};
\draw (1.,0.5)-- (-0.5,2.5);
\draw [color=qqqqff] (-0.9970147837599868,2.5054451564022004) circle (0.49704461067818734cm);
\begin{scriptsize}
\draw [fill=black] (-0.5,-2.) circle (0.5pt);
\draw [fill=black] (1.,-2.) circle (0.5pt);
\draw[color=black] (0.2722378130839064,-1.8449046681763455) node {$a$};
\draw [fill=black] (1.5,-3.) circle (0.5pt);
\draw[color=ttzzqq] (1.6593065470044674,-2.114335285556598) node {$\frac{b-a}{2}$};
\draw [fill=black] (2.,-1.5) circle (0.5pt);
\draw[color=rvwvcq] (1.68924328226894,-1.4257903744737306) node {$\frac{e}{2}$};
\draw [fill=black] (-1.5439241262940946,-1.99957780037612) circle (0.5pt);
\draw[color=rvwvcq] (-1.074915273817358,-1.695220991853983) node {$\frac{c}{2}$};
\draw [fill=black] (-0.5,2.5) circle (0.5pt);
\draw[color=black] (-0.40133873036672585,2.770342018429829) node {$v_1$};
\draw [fill=black] (1.,2.5) circle (0.5pt);
\draw[color=black] (0.856004150741121,2.7603631066750047) node {$v_5$};
\draw[color=black] (0.3121534601032031,2.2763858865660334) node {$a$};
\draw [fill=black] (1.,0.5) circle (0.5pt);
\draw[color=black] (0.806109591967,0.34546646200755804) node {$v_2$};
\draw [color=rvwvcq] (-1.4940295675199735,2.510890312804401)-- ++(-1.0pt,-1.0pt) -- ++(2.0pt,2.0pt) ++(-2.0pt,0) -- ++(2.0pt,-2.0pt);
\draw[color=rvwvcq] (-1.7484918172679902,2.640616165617115) node {$v_7$};
\draw [color=ttzzqq] (1.5,1.5)-- ++(-1.0pt,-1.0pt) -- ++(2.0pt,2.0pt) ++(-2.0pt,0) -- ++(2.0pt,-2.0pt);
\draw[color=ttzzqq] (1.734148385165649,1.5728726078509299) node {$v_6$};
\draw[color=ttzzqq] (1.0505929299601924,1.458115122670452) node {$b-a$};
\draw [color=rvwvcq] (2.005314753674031,0.5078732800068523)-- ++(-1.0pt,-1.0pt) -- ++(2.0pt,2.0pt) ++(-2.0pt,0) -- ++(2.0pt,-2.0pt);
\draw[color=rvwvcq] (2.26303070817133,0.6049181676329862) node {$v_9$};
\draw [fill=black] (-0.5061173037923797,0.49515013832992044) circle (0.5pt);
\draw[color=black] (-0.28159178930883566,0.33548755025273386) node {$v_{11}$};
\draw[color=black] (0.2722378130839064,0.7396334763231124) node {$a$};
\draw [fill=black] (-1.5,0.5) circle (0.5pt);
\draw[color=black] (-1.7584707290228143,0.5550236088588654) node {$v_7'$};
\draw[color=rvwvcq] (-1.1148309208366547,0.1408987710336627) node {$\frac{c}{2}$};
\draw [fill=black] (2.,3.) circle (0.5pt);
\draw[color=black] (2.26303070817133,3.1096250180938507) node {$v_9'$};
\draw[color=rvwvcq] (2.128315399481204,2.685521268513824) node {$\frac{e}{2}$};
\draw [fill=black] (1.5,-0.5) circle (0.5pt);
\draw[color=black] (1.7241694734108248,-0.46282539013319907) node {$v_6'$};
\draw[color=ttzzqq] (1.1603609592632584,-0.46781484601061113) node {$\frac{b-a}{2}$};
\draw[color=rvwvcq] (1.7990113115720059,0.5600130647362774) node {$e$};
\draw[color=wrwrwr] (-6.383696327383822,4.491704296136997) node {$Curve$};
\draw[color=black] (0.3121534601032031,1.7275457400507044) node {$a$};
\draw[color=qqqqff] (-0.9252315974949954,3.1545301209905596) node {$c$};
\end{scriptsize}
\end{tikzpicture}
    \begin{center}
        \textnormal{Figure 11.1. The model $(G_3',l_3')$ of $\Gamma_3'$ }
    \end{center}
\end{figure}
 \begin{figure}[H]
    \centering
    \begin{tikzpicture}[line cap=round,line join=round,>=triangle 45,x=1.0cm,y=1.0cm]
\definecolor{qqqqff}{rgb}{0.,0.,1.}
\definecolor{wrwrwr}{rgb}{0.3803921568627451,0.3803921568627451,0.3803921568627451}
\definecolor{rvwvcq}{rgb}{0.08235294117647059,0.396078431372549,0.7529411764705882}
\definecolor{ttzzqq}{rgb}{0.2,0.6,0.}
\clip(-2.5,-3.5) rectangle (3.,-1.);
\draw (-0.5,-2.)-- (1.,-2.);
\draw [color=ttzzqq] (1.,-2.)-- (1.5,-3.);
\draw [color=rvwvcq] (1.,-2.)-- (2.,-1.5);
\draw [color=rvwvcq] (-0.5,-2.)-- (-1.5439241262940946,-1.99957780037612);
\draw (-0.5,2.5)-- (1.,2.5);
\draw [shift={(0.25,1.5)},color=ttzzqq]  plot[domain=-0.9272952180016123:0.9272952180016122,variable=\t]({1.*1.25*cos(\t r)+0.*1.25*sin(\t r)},{0.*1.25*cos(\t r)+1.*1.25*sin(\t r)});
\draw (1.,0.5)-- (-0.5360540390568522,0.5051290500847446);
\draw [color=rvwvcq] (-0.5360540390568522,0.5051290500847446)-- (-1.5,0.5);
\draw [color=rvwvcq] (1.,2.5)-- (2.,3.);
\draw [color=ttzzqq] (1.,0.5)-- (1.5,-0.5);
\draw [color=rvwvcq] (1.5026573768370155,0.5039366400034262) circle (0.5026727918071416cm);
\draw (3.4754684863824683,0.9442011672970076) node[anchor=north west] {$\Gamma_3'$};
\draw (3.505405221646941,-1.9496832416019991) node[anchor=north west] {$T_3'$};
\draw [->,line width=0.4pt] (3.615173250950007,0.1359093151562506) -- (3.635131074459655,-1.4008430950866702);
\draw (3.954456250614029,-0.32312062556566085) node[anchor=north west] {$\varphi_3$};
\draw (1.,0.5)-- (-0.5,2.5);
\draw [color=qqqqff] (-0.9970147837599868,2.5054451564022004) circle (0.49704461067818734cm);
\begin{scriptsize}
\draw [fill=black] (-0.5,-2.) circle (0.5pt);
\draw[color=black] (-0.4911489361601435,-2.2789873295111964) node {$w_1$};
\draw [fill=black] (1.,-2.) circle (0.5pt);
\draw[color=black] (0.856004150741121,-2.2889662412660208) node {$w_2$};
\draw[color=black] (0.2722378130839064,-1.8449046681763455) node {$a$};
\draw [fill=black] (1.5,-3.) circle (0.5pt);
\draw[color=black] (1.7241694734108246,-3.296836328503261) node {$w_6$};
\draw[color=ttzzqq] (1.6593065470044674,-2.114335285556598) node {$\frac{b-a}{2}$};
\draw [fill=black] (2.,-1.5) circle (0.5pt);
\draw[color=black] (2.312925266945451,-1.4607165656156151) node {$w_9$};
\draw[color=rvwvcq] (1.68924328226894,-1.4257903744737306) node {$\frac{e}{2}$};
\draw [fill=black] (-1.5439241262940946,-1.99957780037612) circle (0.5pt);
\draw[color=black] (-1.6686605232293967,-2.2490505942467243) node {$w_7$};
\draw[color=rvwvcq] (-1.074915273817358,-1.695220991853983) node {$\frac{c}{2}$};
\draw [fill=black] (-0.5,2.5) circle (0.5pt);
\draw[color=black] (-0.40133873036672585,2.770342018429829) node {$v_1$};
\draw [fill=black] (1.,2.5) circle (0.5pt);
\draw[color=black] (0.856004150741121,2.7603631066750047) node {$v_5$};
\draw[color=black] (0.3121534601032031,2.2763858865660334) node {$a$};
\draw [fill=black] (1.,0.5) circle (0.5pt);
\draw[color=black] (0.806109591967,0.34546646200755804) node {$v_2$};
\draw [color=rvwvcq] (-1.4940295675199735,2.510890312804401)-- ++(-1.0pt,-1.0pt) -- ++(2.0pt,2.0pt) ++(-2.0pt,0) -- ++(2.0pt,-2.0pt);
\draw[color=rvwvcq] (-1.7484918172679902,2.640616165617115) node {$v_7$};
\draw [color=ttzzqq] (1.5,1.5)-- ++(-1.0pt,-1.0pt) -- ++(2.0pt,2.0pt) ++(-2.0pt,0) -- ++(2.0pt,-2.0pt);
\draw[color=ttzzqq] (1.734148385165649,1.5728726078509299) node {$v_6$};
\draw[color=ttzzqq] (1.0505929299601924,1.458115122670452) node {$b-a$};
\draw [color=rvwvcq] (2.005314753674031,0.5078732800068523)-- ++(-1.0pt,-1.0pt) -- ++(2.0pt,2.0pt) ++(-2.0pt,0) -- ++(2.0pt,-2.0pt);
\draw[color=rvwvcq] (2.26303070817133,0.6049181676329862) node {$v_9$};
\draw [fill=black] (-0.5360540390568522,0.5051290500847446) circle (0.5pt);
\draw[color=black] (-0.3115285245733082,0.34546646200755804) node {$v_{11}$};
\draw[color=black] (0.22234325430978547,0.7496123880779366) node {$a$};
\draw [fill=black] (-1.5,0.5) circle (0.5pt);
\draw[color=black] (-1.7584707290228143,0.5550236088588654) node {$v_7'$};
\draw[color=rvwvcq] (-1.1347887443463032,0.1408987710336627) node {$\frac{c}{2}$};
\draw [fill=black] (2.,3.) circle (0.5pt);
\draw[color=black] (2.26303070817133,3.1096250180938507) node {$v_9'$};
\draw[color=rvwvcq] (2.128315399481204,2.685521268513824) node {$\frac{e}{2}$};
\draw [fill=black] (1.5,-0.5) circle (0.5pt);
\draw[color=black] (1.7241694734108248,-0.46282539013319907) node {$v_6'$};
\draw[color=ttzzqq] (1.1603609592632584,-0.46781484601061113) node {$\frac{b-a}{2}$};
\draw[color=rvwvcq] (1.7990113115720059,0.5600130647362774) node {$e$};
\draw[color=wrwrwr] (-6.383696327383822,4.491704296136997) node {$Curve$};
\draw[color=black] (0.2622589013290822,1.6077987989928144) node {$a$};
\draw[color=qqqqff] (-0.9252315974949954,3.1545301209905596) node {$c$};
\end{scriptsize}
\end{tikzpicture}
    \begin{center}
        \textnormal{Figure 11.2. The model $(T_3',t_3')$ of $T_3$}
    \end{center}
\end{figure}
\begin{figure}[H]
    \centering
    \begin{tikzpicture}[line cap=round,line join=round,>=triangle 45,x=1.0cm,y=1.0cm]
\definecolor{eqeqeq}{rgb}{0.8784313725490196,0.8784313725490196,0.8784313725490196}
\definecolor{qqqqff}{rgb}{0.,0.,1.}
\definecolor{wrwrwr}{rgb}{0.3803921568627451,0.3803921568627451,0.3803921568627451}
\definecolor{rvwvcq}{rgb}{0.08235294117647059,0.396078431372549,0.7529411764705882}
\definecolor{ttzzqq}{rgb}{0.2,0.6,0.}
\clip(-2.5,-3.5) rectangle (4.5,3.5);
\draw (-0.5,-2.)-- (1.,-2.);
\draw [color=ttzzqq] (1.,-2.)-- (1.5,-3.);
\draw [color=rvwvcq] (1.,-2.)-- (2.,-1.5);
\draw [color=rvwvcq] (-0.5,-2.)-- (-1.5439241262940946,-1.99957780037612);
\draw (-0.5,2.5)-- (1.,2.5);
\draw [shift={(0.25,1.5)},color=ttzzqq]  plot[domain=-0.9272952180016123:0.9272952180016122,variable=\t]({1.*1.25*cos(\t r)+0.*1.25*sin(\t r)},{0.*1.25*cos(\t r)+1.*1.25*sin(\t r)});
\draw (1.,0.5)-- (-0.5061173037923797,0.49515013832992044);
\draw [color=rvwvcq] (-0.5061173037923797,0.49515013832992044)-- (-1.5,0.5);
\draw [color=rvwvcq] (1.,2.5)-- (2.,3.);
\draw [color=ttzzqq] (1.,0.5)-- (1.5,-0.5);
\draw [color=rvwvcq] (1.5026573768370155,0.5039366400034262) circle (0.5026727918071416cm);
\draw (3.4754684863824683,0.9442011672970076) node[anchor=north west] {$\Gamma_3'$};
\draw (3.505405221646941,-1.9496832416019991) node[anchor=north west] {$T_3$};
\draw [->,line width=0.4pt] (3.615173250950007,0.1359093151562506) -- (3.635131074459655,-1.4008430950866702);
\draw (3.954456250614029,-0.32312062556566085) node[anchor=north west] {$\varphi_3$};
\draw (1.,0.5)-- (-0.5,2.5);
\draw [color=qqqqff] (-0.9970147837599868,2.5054451564022004) circle (0.49704461067818734cm);
\draw [line width=0.4pt,dash pattern=on 2pt off 2pt,color=eqeqeq] (-1.4940295675199735,2.510890312804401)-- (-1.5439241262940946,-1.99957780037612);
\draw [line width=0.4pt,dash pattern=on 2pt off 2pt,color=eqeqeq] (-0.5,2.5)-- (-0.5,-2.);
\draw [line width=0.4pt,dash pattern=on 2pt off 2pt,color=eqeqeq] (1.,2.5)-- (1.,-2.);
\draw [line width=0.4pt,dash pattern=on 2pt off 2pt,color=eqeqeq] (2.,3.)-- (2.,-1.5);
\draw [line width=0.4pt,dash pattern=on 2pt off 2pt,color=eqeqeq] (1.5,1.5)-- (1.5,-3.);
\begin{scriptsize}
\draw [fill=black] (-0.5,-2.) circle (0.5pt);
\draw [fill=black] (1.,-2.) circle (0.5pt);
\draw[color=black] (0.2722378130839064,-1.8449046681763455) node {$a$};
\draw [fill=black] (1.5,-3.) circle (0.5pt);
\draw[color=ttzzqq] (1.6593065470044674,-2.114335285556598) node {$\frac{b-a}{2}$};
\draw [fill=black] (2.,-1.5) circle (0.5pt);
\draw[color=rvwvcq] (1.68924328226894,-1.4257903744737306) node {$\frac{e}{2}$};
\draw [fill=black] (-1.5439241262940946,-1.99957780037612) circle (0.5pt);
\draw[color=rvwvcq] (-1.074915273817358,-1.695220991853983) node {$\frac{c}{2}$};
\draw [fill=black] (-0.5,2.5) circle (0.5pt);
\draw[color=black] (-0.40133873036672585,2.770342018429829) node {$v_1$};
\draw [fill=black] (1.,2.5) circle (0.5pt);
\draw[color=black] (0.3121534601032031,2.2763858865660334) node {$a$};
\draw [fill=black] (1.,0.5) circle (0.5pt);
\draw[color=black] (0.806109591967,0.34546646200755804) node {$v_2$};
\draw [color=rvwvcq] (-1.4940295675199735,2.510890312804401)-- ++(-1.0pt,-1.0pt) -- ++(2.0pt,2.0pt) ++(-2.0pt,0) -- ++(2.0pt,-2.0pt);
\draw [color=ttzzqq] (1.5,1.5)-- ++(-1.0pt,-1.0pt) -- ++(2.0pt,2.0pt) ++(-2.0pt,0) -- ++(2.0pt,-2.0pt);
\draw[color=ttzzqq] (1.0505929299601924,1.458115122670452) node {$b-a$};
\draw [color=rvwvcq] (2.005314753674031,0.5078732800068523)-- ++(-1.0pt,-1.0pt) -- ++(2.0pt,2.0pt) ++(-2.0pt,0) -- ++(2.0pt,-2.0pt);
\draw [fill=black] (-0.5061173037923797,0.49515013832992044) circle (0.5pt);
\draw[color=black] (0.2722378130839064,0.7396334763231124) node {$a$};
\draw [fill=black] (-1.5,0.5) circle (0.5pt);
\draw[color=rvwvcq] (-1.1148309208366547,0.1408987710336627) node {$\frac{c}{2}$};
\draw [fill=black] (2.,3.) circle (0.5pt);
\draw[color=rvwvcq] (2.128315399481204,2.685521268513824) node {$\frac{e}{2}$};
\draw [fill=black] (1.5,-0.5) circle (0.5pt);
\draw[color=ttzzqq] (1.1603609592632584,-0.46781484601061113) node {$\frac{b-a}{2}$};
\draw[color=rvwvcq] (1.7990113115720059,0.5600130647362774) node {$e$};
\draw[color=wrwrwr] (-6.383696327383822,4.491704296136997) node {$Curve$};
\draw[color=black] (0.3121534601032031,1.7275457400507044) node {$a$};
\draw[color=qqqqff] (-0.9252315974949954,3.1545301209905596) node {$c$};
\end{scriptsize}
\end{tikzpicture}
    \begin{center}
        \textnormal{Figure 12. The tropical morphism $\varphi_3:\Gamma_3' \to T_3$ }
    \end{center}
\end{figure}

\textbf{Case 1.6.} Consider the metric graph $\Gamma_4$ with essential model $(G,l)$ in Figure 13, where $b,c,d$, and $e$ are real positive numbers. The model $(G_1,l_1)$ that is obtained by subdividing $(G,l)$ is shown in Figure 13.1. The tropical modification $\Gamma_4'$, the metric tree $T_4$ with models $(G_4',l_4')$, $(T_4',t_4')$ is given in Figure 13.2, 13.3, respectively. The construction of the tropical morphism $\varphi_4: \Gamma_4' \to T_4$ of degree $3$ is depicted in Figure 14.
\begin{figure}[H]
    \centering
    \begin{tikzpicture}[line cap=round,line join=round,>=triangle 45,x=1.0cm,y=1.0cm]

\definecolor{wrwrwr}{rgb}{0.3803921568627451,0.3803921568627451,0.3803921568627451}
\clip(-1.5,-1.5) rectangle (3.,2.);
\draw (1.4949220900329574,0.)-- (-1.,0.);
\draw(1.9929146880375257,0.) circle (0.4979925980045683cm);
\draw [shift={(0.2474610450164787,0.)}] plot[domain=3.141592653589793:6.283185307179586,variable=\t]({1.*1.2474610450164787*cos(\t r)+0.*1.2474610450164787*sin(\t r)},{0.*1.2474610450164787*cos(\t r)+1.*1.2474610450164787*sin(\t r)});
\draw [shift={(0.2474610450164787,0.)}] plot[domain=0.:3.141592653589793,variable=\t]({1.*1.2474610450164787*cos(\t r)+0.*1.2474610450164787*sin(\t r)},{0.*1.2474610450164787*cos(\t r)+1.*1.2474610450164787*sin(\t r)});
\begin{scriptsize}
\draw [fill=black] (1.4949220900329574,0.) circle (0.5pt);
\draw[color=black] (1.7345223022804426,0.07227638556236948) node {$v_3$};
\draw [fill=black] (-1.,0.) circle (0.5pt);
\draw[color=black] (-1.2064528519729523,0.08167247231717266) node {$v_1$};
\draw[color=black] (0.2640347251537452,0.28368833754544087) node {$d$};
\draw[color=black] (2.3029855509460346,-0.09215513264668602) node {$e$};
\draw[color=black] (0.2170542913797293,-1.012971634617397) node {$c$};
\draw[color=black] (0.3298073324373675,1.646120916991901) node {$b$};
\draw[color=wrwrwr] (-4.997773857536034,3.6380913090101736) node {$Curve$};
\end{scriptsize}
\end{tikzpicture}
    \begin{center}
        \textnormal{Figure 13. The essential model $(G,l)$ of $\Gamma_4$}
    \end{center}
\end{figure}
\begin{figure}[H]
    \centering
    \begin{tikzpicture}[line cap=round,line join=round,>=triangle 45,x=1.0cm,y=1.0cm]
\definecolor{qqttzz}{rgb}{0.,0.2,0.6}
\definecolor{wrwrwr}{rgb}{0.3803921568627451,0.3803921568627451,0.3803921568627451}
\definecolor{uuuuuu}{rgb}{0.26666666666666666,0.26666666666666666,0.26666666666666666}
\definecolor{sexdts}{rgb}{0.1803921568627451,0.49019607843137253,0.19607843137254902}
\definecolor{rvwvcq}{rgb}{0.08235294117647059,0.396078431372549,0.7529411764705882}
\definecolor{qqqqff}{rgb}{0.,0.,1.}
\definecolor{dtsfsf}{rgb}{0.8274509803921568,0.1843137254901961,0.1843137254901961}
\clip(-1.5,-0.5) rectangle (1.5,3.);
\draw [shift={(0.7427138764381784,1.5036430617809107)},color=dtsfsf]  plot[domain=2.211365906928214:4.070969994016963,variable=\t]({1.*1.242719216312013*cos(\t r)+0.*1.242719216312013*sin(\t r)},{0.*1.242719216312013*cos(\t r)+1.*1.242719216312013*sin(\t r)});
\draw [shift={(-0.003960176597837122,1.5039601765978372)},color=rvwvcq]  plot[domain=1.566820425794195:4.7155154751509825,variable=\t]({1.*0.9960476960475828*cos(\t r)+0.*0.9960476960475828*sin(\t r)},{0.*0.9960476960475828*cos(\t r)+1.*0.9960476960475828*sin(\t r)});
\draw [shift={(-0.7464501113912964,1.5063728582396283)},color=sexdts]  plot[domain=-0.9287282318863292:0.9287282318863297,variable=\t]({1.*1.2464664028803656*cos(\t r)+0.*1.2464664028803656*sin(\t r)},{0.*1.2464664028803656*cos(\t r)+1.*1.2464664028803656*sin(\t r)});
\draw [color=sexdts] (0.,-2.)-- (0.5,-2.5);
\draw [color=dtsfsf] (0.,-2.)-- (-0.5,-1.5);
\draw [color=rvwvcq] (0.,-2.)-- (-1.,-2.5);
\draw (1.926513264302514,1.0516805380025631) node[anchor=north west] {$\Gamma_4'$};
\draw (1.8765830064975002,-1.7527356087123767) node[anchor=north west] {$T_4$};
\draw [color=qqttzz] (0.,-2.)-- (1.,-2.);
\draw [color=qqttzz] (0.4884961611223245,0.49686081579665) circle (0.48946709838537134cm);
\begin{scriptsize}
\draw [fill=black] (-8.460437625978602E-4,0.5113549601832207) circle (0.5pt);
\draw [fill=black] (-8.460437625978602E-4,0.5079173487142253) circle (0.5pt);
\draw[color=black] (-0.17471841832515247,0.25279641312234286) node {$v_1$};
\draw [color=dtsfsf] (-0.5,1.5)-- ++(-1.0pt,-1.0pt) -- ++(2.0pt,2.0pt) ++(-2.0pt,0) -- ++(2.0pt,-2.0pt);
\draw[color=dtsfsf] (-0.7489163830828118,1.484409438979349) node {$v_3'$};
\draw[color=dtsfsf] (-0.5866430452165168,1.9046557755048816) node {$d$};
\draw [color=qqqqff] (-1.,1.5)-- ++(-1.0pt,-1.0pt) -- ++(2.0pt,2.0pt) ++(-2.0pt,0) -- ++(2.0pt,-2.0pt);
\draw[color=qqqqff] (-1.2648623804012884,1.492731148613518) node {$v_3''$};
\draw[color=rvwvcq] (-1.1025890425349933,1.8880123562365436) node {$c$};
\draw [color=sexdts] (0.5,1.5)-- ++(-1.0pt,-1.0pt) -- ++(2.0pt,2.0pt) ++(-2.0pt,0) -- ++(2.0pt,-2.0pt);
\draw[color=sexdts] (0.24136706338329633,1.4511226004426732) node {$v_2'$};
\draw[color=sexdts] (0.6865785288113366,1.4552834552597578) node {$b$};
\draw [color=rvwvcq] (0.9778383660072469,0.48580428287907457)-- ++(-1.0pt,-1.0pt) -- ++(2.0pt,2.0pt) ++(-2.0pt,0) -- ++(2.0pt,-2.0pt);
\draw[color=rvwvcq] (1.2566156387519114,0.4941259925132427) node {$v_4'$};
\draw [fill=black] (0.,-2.) circle (0.5pt);
\draw [fill=black] (0.5,-2.5) circle (0.5pt);
\draw[color=black] (0.6657742547259141,-2.626515120300118) node {$w_2'$};
\draw[color=sexdts] (0.27049304710288774,-2.522493749873006) node {$\frac{b}{2}$};
\draw [fill=black] (1.,-2.) circle (0.5pt);
\draw[color=black] (1.0152860593610111,-2.135534251884149) node {$w_4'$};
\draw [fill=black] (-0.5,-1.5) circle (0.5pt);
\draw[color=black] (-0.7322729638144738,-1.5613362871264904) node {$w_3'$};
\draw[color=dtsfsf] (-0.170557563508068,-1.5238885937727302) node {$\frac{d}{2}$};
\draw [fill=black] (-1.,-2.5) circle (0.5pt);
\draw[color=black] (-1.2066104129621056,-2.6514802492026246) node {$w_3''$};
\draw[color=rvwvcq] (-0.5866430452165168,-2.0731214296278817) node {$\frac{c}{2}$};
\draw [fill=uuuuuu] (0.,2.5) circle (0.5pt);
\draw[color=uuuuuu] (-0.024927644910110895,2.724344174470524) node {$v_3$};
\draw[color=wrwrwr] (-5.513095148644551,3.67717992758287) node {$Curve$};
\draw[color=qqttzz] (0.7947607540555333,-1.6820010768219404) node {$\frac{e}{2}$};
\draw[color=qqttzz] (0.7947607540555333,0.39010462208613067) node {$e$};
\end{scriptsize}
\end{tikzpicture}
    \begin{center}
        \textnormal{Figure 13.1.  The model $(G_1,l_1)$ of $\Gamma_4$}
    \end{center}
\end{figure}
\begin{figure}[H]
    \centering
    \begin{tikzpicture}[line cap=round,line join=round,>=triangle 45,x=1.0cm,y=1.0cm]
\definecolor{qqttzz}{rgb}{0.,0.2,0.6}
\definecolor{wrwrwr}{rgb}{0.3803921568627451,0.3803921568627451,0.3803921568627451}
\definecolor{uuuuuu}{rgb}{0.26666666666666666,0.26666666666666666,0.26666666666666666}
\definecolor{sexdts}{rgb}{0.1803921568627451,0.49019607843137253,0.19607843137254902}
\definecolor{rvwvcq}{rgb}{0.08235294117647059,0.396078431372549,0.7529411764705882}
\definecolor{qqqqff}{rgb}{0.,0.,1.}
\definecolor{dtsfsf}{rgb}{0.8274509803921568,0.1843137254901961,0.1843137254901961}
\clip(-1.5,-1.) rectangle (1.5,3.);
\draw [shift={(0.7427138764381784,1.5036430617809107)},color=dtsfsf]  plot[domain=2.211365906928214:4.070969994016963,variable=\t]({1.*1.242719216312013*cos(\t r)+0.*1.242719216312013*sin(\t r)},{0.*1.242719216312013*cos(\t r)+1.*1.242719216312013*sin(\t r)});
\draw [shift={(-0.003960176597837122,1.5039601765978372)},color=rvwvcq]  plot[domain=1.566820425794195:4.7155154751509825,variable=\t]({1.*0.9960476960475828*cos(\t r)+0.*0.9960476960475828*sin(\t r)},{0.*0.9960476960475828*cos(\t r)+1.*0.9960476960475828*sin(\t r)});
\draw [shift={(-0.7464501113912964,1.5063728582396283)},color=sexdts]  plot[domain=-0.9287282318863292:0.9287282318863297,variable=\t]({1.*1.2464664028803656*cos(\t r)+0.*1.2464664028803656*sin(\t r)},{0.*1.2464664028803656*cos(\t r)+1.*1.2464664028803656*sin(\t r)});
\draw [color=sexdts] (0.,-2.)-- (0.5,-2.5);
\draw [color=dtsfsf] (-8.460437625978602E-4,0.5113549601832207)-- (-0.5130501526429195,-0.13835360745691855);
\draw [color=rvwvcq] (-8.460437625978602E-4,0.5113549601832207)-- (-1.0149414271162547,0.23634604266358505);
\draw [color=dtsfsf] (0.,-2.)-- (-0.5,-1.5);
\draw [color=rvwvcq] (0.,-2.)-- (-1.,-2.5);
\draw (1.926513264302514,1.0516805380025631) node[anchor=north west] {$\Gamma_4'$};
\draw (1.8765830064975002,-1.7527356087123767) node[anchor=north west] {$T_4$};
\draw [color=sexdts] (-8.460437625978602E-4,0.5079173487142253)-- (0.5201443361279533,-0.5294442924895384);
\draw [color=qqttzz] (0.,2.5)-- (1.,2.5);
\draw [color=qqttzz] (0.,-2.)-- (1.,-2.);
\draw [color=qqttzz] (0.4884961611223245,0.49686081579665) circle (0.48946709838537134cm);
\begin{scriptsize}
\draw [fill=black] (-8.460437625978602E-4,0.5113549601832207) circle (0.5pt);
\draw [fill=black] (-8.460437625978602E-4,0.5079173487142253) circle (0.5pt);
\draw[color=black] (-0.09982303161763167,0.06971880117062572) node {$v_1$};
\draw [color=dtsfsf] (-0.5,1.5)-- ++(-1.0pt,-1.0pt) -- ++(2.0pt,2.0pt) ++(-2.0pt,0) -- ++(2.0pt,-2.0pt);
\draw[color=dtsfsf] (-0.7489163830828118,1.484409438979349) node {$v_3'$};
\draw[color=dtsfsf] (-0.5866430452165168,1.9046557755048816) node {$d$};
\draw [color=qqqqff] (-1.,1.5)-- ++(-1.0pt,-1.0pt) -- ++(2.0pt,2.0pt) ++(-2.0pt,0) -- ++(2.0pt,-2.0pt);
\draw[color=qqqqff] (-1.2648623804012884,1.492731148613518) node {$v_3''$};
\draw[color=rvwvcq] (-1.1025890425349933,1.8880123562365436) node {$c$};
\draw [color=sexdts] (0.5,1.5)-- ++(-1.0pt,-1.0pt) -- ++(2.0pt,2.0pt) ++(-2.0pt,0) -- ++(2.0pt,-2.0pt);
\draw[color=sexdts] (0.24136706338329633,1.4511226004426732) node {$v_2'$};
\draw[color=sexdts] (0.6865785288113366,1.4552834552597578) node {$b$};
\draw [fill=black] (1.,2.5) circle (0.5pt);
\draw[color=black] (1.0069643497268421,2.716022464836355) node {$x_4'$};
\draw [color=rvwvcq] (0.9778383660072469,0.48580428287907457)-- ++(-1.0pt,-1.0pt) -- ++(2.0pt,2.0pt) ++(-2.0pt,0) -- ++(2.0pt,-2.0pt);
\draw[color=rvwvcq] (1.2566156387519114,0.4941259925132427) node {$v_4'$};
\draw [fill=black] (0.5201443361279533,-0.5294442924895384) circle (0.5pt);
\draw[color=black] (0.715704512530928,-0.7208436140754255) node {$x_2'$};
\draw [fill=black] (0.,-2.) circle (0.5pt);
\draw [fill=black] (0.5,-2.5) circle (0.5pt);
\draw[color=black] (0.6657742547259141,-2.626515120300118) node {$w_2'$};
\draw[color=sexdts] (0.27049304710288774,-2.522493749873006) node {$\frac{b}{2}$};
\draw [fill=black] (1.,-2.) circle (0.5pt);
\draw[color=black] (1.0152860593610111,-2.135534251884149) node {$w_4'$};
\draw [fill=black] (-0.5130501526429195,-0.13835360745691855) circle (0.5pt);
\draw[color=black] (-0.6989861252777979,-0.27147129383030166) node {$x_3'$};
\draw[color=dtsfsf] (-0.5866430452165168,0.115488204158555) node {$\frac{d}{2}$};
\draw [fill=black] (-1.0149414271162547,0.23634604266358505) circle (0.5pt);
\draw[color=black] (-1.2149321225962744,0.2777615420248497) node {$x_3''$};
\draw[color=rvwvcq] (-0.7364338186315583,0.5815039436720169) node {$\frac{c}{2}$};
\draw [fill=black] (-0.5,-1.5) circle (0.5pt);
\draw[color=black] (-0.7322729638144738,-1.5613362871264904) node {$w_3'$};
\draw[color=dtsfsf] (-0.170557563508068,-1.5238885937727302) node {$\frac{d}{2}$};
\draw [fill=black] (-1.,-2.5) circle (0.5pt);
\draw[color=black] (-1.2066104129621056,-2.6514802492026246) node {$w_3''$};
\draw[color=rvwvcq] (-0.5866430452165168,-2.0731214296278817) node {$\frac{c}{2}$};
\draw [fill=uuuuuu] (0.,2.5) circle (0.5pt);
\draw[color=uuuuuu] (-0.024927644910110895,2.724344174470524) node {$v_3$};
\draw[color=wrwrwr] (-5.513095148644551,3.67717992758287) node {$Curve$};
\draw[color=sexdts] (0.2538496278345498,-0.3754926642574137) node {$\frac{b}{2}$};
\draw[color=qqttzz] (0.7614739155188572,2.2208807416033016) node {$\frac{e}{2}$};
\draw[color=qqttzz] (0.7947607540555333,-1.6820010768219404) node {$\frac{e}{2}$};
\draw[color=qqttzz] (0.7947607540555333,0.39010462208613067) node {$e$};
\end{scriptsize}
\end{tikzpicture}
    \begin{center}
        \textnormal{Figure 13.2. The model $(G_4',l_4')$ of $\Gamma_4'$}
    \end{center}
\end{figure}
\begin{figure}[H]
    \centering
    \begin{tikzpicture}[line cap=round,line join=round,>=triangle 45,x=1.0cm,y=1.0cm]
\definecolor{qqttzz}{rgb}{0.,0.2,0.6}
\definecolor{wrwrwr}{rgb}{0.3803921568627451,0.3803921568627451,0.3803921568627451}
\definecolor{uuuuuu}{rgb}{0.26666666666666666,0.26666666666666666,0.26666666666666666}
\definecolor{cqcqcq}{rgb}{0.7529411764705882,0.7529411764705882,0.7529411764705882}
\definecolor{sexdts}{rgb}{0.1803921568627451,0.49019607843137253,0.19607843137254902}
\definecolor{rvwvcq}{rgb}{0.08235294117647059,0.396078431372549,0.7529411764705882}
\definecolor{qqqqff}{rgb}{0.,0.,1.}
\definecolor{dtsfsf}{rgb}{0.8274509803921568,0.1843137254901961,0.1843137254901961}
\clip(-1.5,-3.) rectangle (1.5,-1.);
\draw [shift={(0.7427138764381784,1.5036430617809107)},color=dtsfsf]  plot[domain=2.211365906928214:4.070969994016963,variable=\t]({1.*1.242719216312013*cos(\t r)+0.*1.242719216312013*sin(\t r)},{0.*1.242719216312013*cos(\t r)+1.*1.242719216312013*sin(\t r)});
\draw [shift={(-0.003960176597837122,1.5039601765978372)},color=rvwvcq]  plot[domain=1.566820425794195:4.7155154751509825,variable=\t]({1.*0.9960476960475828*cos(\t r)+0.*0.9960476960475828*sin(\t r)},{0.*0.9960476960475828*cos(\t r)+1.*0.9960476960475828*sin(\t r)});
\draw [shift={(-0.7464501113912964,1.5063728582396283)},color=sexdts]  plot[domain=-0.9287282318863292:0.9287282318863297,variable=\t]({1.*1.2464664028803656*cos(\t r)+0.*1.2464664028803656*sin(\t r)},{0.*1.2464664028803656*cos(\t r)+1.*1.2464664028803656*sin(\t r)});
\draw [color=sexdts] (0.,-2.)-- (0.5,-2.5);
\draw [color=dtsfsf] (-8.460437625978602E-4,0.5113549601832207)-- (-0.5130501526429195,-0.13835360745691855);
\draw [color=rvwvcq] (-8.460437625978602E-4,0.5113549601832207)-- (-1.0149414271162547,0.23634604266358505);
\draw [color=dtsfsf] (0.,-2.)-- (-0.5,-1.5);
\draw [color=rvwvcq] (0.,-2.)-- (-1.,-2.5);
\draw [line width=0.4pt,dash pattern=on 1pt off 1pt,color=cqcqcq] (-0.5,1.5)-- (-0.5130501526429195,-0.13835360745691855);
\draw [line width=0.4pt,dash pattern=on 1pt off 1pt,color=cqcqcq] (-1.,1.5)-- (-1.0149414271162547,0.23634604266358505);
\draw (1.9265132643025131,1.0516805380025631) node[anchor=north west] {$\Gamma_4'$};
\draw (1.876583006497499,-1.7527356087123767) node[anchor=north west] {$T_4$};
\draw [->] (2.2109114979186684,-0.07267187425176508) -- (2.220844272249744,-1.0560165330282543);
\draw [color=sexdts] (-8.460437625978602E-4,0.5079173487142253)-- (0.5201443361279533,-0.5294442924895384);
\draw [color=qqttzz] (0.,2.5)-- (1.,2.5);
\draw [color=qqttzz] (0.,-2.)-- (1.,-2.);
\draw [color=qqttzz] (0.4884961611223245,0.49686081579665) circle (0.48946709838537134cm);
\begin{scriptsize}
\draw [fill=black] (-8.460437625978602E-4,0.5113549601832207) circle (0.5pt);
\draw [fill=black] (-8.460437625978602E-4,0.5079173487142253) circle (0.5pt);
\draw[color=black] (-0.09982303161763265,0.06971880117062572) node {$v_1$};
\draw [color=dtsfsf] (-0.5,1.5)-- ++(-1.0pt,-1.0pt) -- ++(2.0pt,2.0pt) ++(-2.0pt,0) -- ++(2.0pt,-2.0pt);
\draw[color=dtsfsf] (-0.5866430452165177,1.9046557755048816) node {$d$};
\draw [color=qqqqff] (-1.,1.5)-- ++(-1.0pt,-1.0pt) -- ++(2.0pt,2.0pt) ++(-2.0pt,0) -- ++(2.0pt,-2.0pt);
\draw[color=rvwvcq] (-1.1025890425349942,1.8880123562365436) node {$c$};
\draw [color=sexdts] (0.5,1.5)-- ++(-1.0pt,-1.0pt) -- ++(2.0pt,2.0pt) ++(-2.0pt,0) -- ++(2.0pt,-2.0pt);
\draw[color=sexdts] (0.6865785288113356,1.4552834552597578) node {$b$};
\draw [fill=black] (1.,2.5) circle (0.5pt);
\draw [color=rvwvcq] (0.9778383660072469,0.48580428287907457)-- ++(-1.0pt,-1.0pt) -- ++(2.0pt,2.0pt) ++(-2.0pt,0) -- ++(2.0pt,-2.0pt);
\draw [fill=black] (0.5201443361279533,-0.5294442924895384) circle (0.5pt);
\draw [fill=black] (0.,-2.) circle (0.5pt);
\draw [fill=black] (0.5,-2.5) circle (0.5pt);
\draw[color=black] (0.6657742547259132,-2.626515120300118) node {$w_2'$};
\draw[color=sexdts] (0.2704930471028868,-2.522493749873006) node {$\frac{b}{2}$};
\draw [fill=black] (1.,-2.) circle (0.5pt);
\draw[color=black] (1.0152860593610102,-2.135534251884149) node {$w_4'$};
\draw [fill=black] (-0.5130501526429195,-0.13835360745691855) circle (0.5pt);
\draw[color=dtsfsf] (-0.5866430452165178,0.115488204158555) node {$\frac{d}{2}$};
\draw [fill=black] (-1.0149414271162547,0.23634604266358505) circle (0.5pt);
\draw[color=rvwvcq] (-0.7364338186315594,0.5815039436720169) node {$\frac{c}{2}$};
\draw [fill=black] (-0.5,-1.5) circle (0.5pt);
\draw[color=black] (-0.7322729638144748,-1.5613362871264904) node {$w_3'$};
\draw[color=dtsfsf] (-0.17055756350806894,-1.5238885937727302) node {$\frac{d}{2}$};
\draw [fill=black] (-1.,-2.5) circle (0.5pt);
\draw[color=black] (-1.2066104129621065,-2.6514802492026246) node {$w_3''$};
\draw[color=rvwvcq] (-0.5866430452165178,-2.0731214296278817) node {$\frac{c}{2}$};
\draw [fill=uuuuuu] (0.,2.5) circle (0.5pt);
\draw[color=uuuuuu] (-0.02492764491011184,2.724344174470524) node {$v_3$};
\draw[color=wrwrwr] (-5.513095148644552,3.67717992758287) node {$Curve$};
\draw[color=black] (2.6213760187556225,-0.5544094213920463) node {$\varphi_4$};
\draw[color=sexdts] (0.25384962783454884,-0.3754926642574137) node {$\frac{b}{2}$};
\draw[color=qqttzz] (0.7614739155188563,2.2208807416033016) node {$\frac{e}{2}$};
\draw[color=qqttzz] (0.7947607540555323,-1.6820010768219404) node {$\frac{e}{2}$};
\draw[color=qqttzz] (0.7947607540555324,0.39010462208613067) node {$e$};
\end{scriptsize}
\end{tikzpicture}
    \begin{center}
        \textnormal{Figure 13.3.  The model $(T_4',t_4')$ of $T_4$}
    \end{center}
\end{figure}
\begin{figure}[H]
    \centering
    \begin{tikzpicture}[line cap=round,line join=round,>=triangle 45,x=1.0cm,y=1.0cm]

\definecolor{qqttzz}{rgb}{0.,0.2,0.6}
\definecolor{wrwrwr}{rgb}{0.3803921568627451,0.3803921568627451,0.3803921568627451}
\definecolor{uuuuuu}{rgb}{0.26666666666666666,0.26666666666666666,0.26666666666666666}
\definecolor{cqcqcq}{rgb}{0.7529411764705882,0.7529411764705882,0.7529411764705882}
\definecolor{sexdts}{rgb}{0.1803921568627451,0.49019607843137253,0.19607843137254902}
\definecolor{rvwvcq}{rgb}{0.08235294117647059,0.396078431372549,0.7529411764705882}
\definecolor{qqqqff}{rgb}{0.,0.,1.}
\definecolor{dtsfsf}{rgb}{0.8274509803921568,0.1843137254901961,0.1843137254901961}
\clip(-2.5,-3.) rectangle (3.,3.);
\draw [shift={(0.7427138764381784,1.5036430617809107)},color=dtsfsf]  plot[domain=2.211365906928214:4.070969994016963,variable=\t]({1.*1.242719216312013*cos(\t r)+0.*1.242719216312013*sin(\t r)},{0.*1.242719216312013*cos(\t r)+1.*1.242719216312013*sin(\t r)});
\draw [shift={(-0.003960176597837122,1.5039601765978372)},color=rvwvcq]  plot[domain=1.566820425794195:4.7155154751509825,variable=\t]({1.*0.9960476960475828*cos(\t r)+0.*0.9960476960475828*sin(\t r)},{0.*0.9960476960475828*cos(\t r)+1.*0.9960476960475828*sin(\t r)});
\draw [shift={(-0.7464501113912964,1.5063728582396283)},color=sexdts]  plot[domain=-0.9287282318863292:0.9287282318863297,variable=\t]({1.*1.2464664028803656*cos(\t r)+0.*1.2464664028803656*sin(\t r)},{0.*1.2464664028803656*cos(\t r)+1.*1.2464664028803656*sin(\t r)});
\draw [color=sexdts] (0.,-2.)-- (0.5,-2.5);
\draw [color=dtsfsf] (-8.460437625978602E-4,0.5113549601832207)-- (-0.5130501526429195,-0.13835360745691855);
\draw [color=rvwvcq] (-8.460437625978602E-4,0.5113549601832207)-- (-1.0149414271162547,0.23634604266358505);
\draw [color=dtsfsf] (0.,-2.)-- (-0.5,-1.5);
\draw [color=rvwvcq] (0.,-2.)-- (-1.,-2.5);
\draw [line width=0.4pt,dash pattern=on 2pt off 2pt,color=cqcqcq] (0.,2.5)-- (0.,-2.);
\draw [line width=0.4pt,dash pattern=on 2pt off 2pt,color=cqcqcq] (-0.5,1.5)-- (-0.5130501526429195,-0.13835360745691855);
\draw [line width=0.4pt,dash pattern=on 2pt off 2pt,color=cqcqcq] (-0.5130501526429195,-0.13835360745691855)-- (-0.5,-1.5);
\draw [line width=0.4pt,dash pattern=on 2pt off 2pt,color=cqcqcq] (-1.,1.5)-- (-1.0149414271162547,0.23634604266358505);
\draw [line width=0.4pt,dash pattern=on 2pt off 2pt,color=cqcqcq] (-1.0149414271162547,0.23634604266358505)-- (-1.,-2.5);
\draw [line width=0.4pt,dash pattern=on 2pt off 2pt,color=cqcqcq] (0.,2.5046166785127313)-- (0.,-2.);
\draw [line width=0.4pt,dash pattern=on 2pt off 2pt,color=cqcqcq] (0.5,1.5)-- (0.5,-2.5);
\draw [line width=0.4pt,dash pattern=on 2pt off 2pt,color=cqcqcq] (1.,2.5)-- (1.,-2.);
\draw (1.9265132643025131,1.0516805380025631) node[anchor=north west] {$\Gamma_4'$};
\draw (1.876583006497499,-1.7527356087123767) node[anchor=north west] {$T_4$};
\draw [->] (2.2109114979186684,-0.07267187425176508) -- (2.220844272249744,-1.0560165330282543);
\draw [color=sexdts] (-8.460437625978602E-4,0.5079173487142253)-- (0.5201443361279533,-0.5294442924895384);
\draw [color=qqttzz] (0.,2.5)-- (1.,2.5);
\draw [color=qqttzz] (0.,-2.)-- (1.,-2.);
\draw [color=qqttzz] (0.4884961611223245,0.49686081579665) circle (0.48946709838537134cm);
\begin{scriptsize}
\draw [fill=black] (-8.460437625978602E-4,0.5113549601832207) circle (0.5pt);
\draw [fill=black] (-8.460437625978602E-4,0.5079173487142253) circle (0.5pt);
\draw[color=black] (-0.09982303161763265,0.06971880117062572) node {$v_1$};
\draw [color=dtsfsf] (-0.5,1.5)-- ++(-1.0pt,-1.0pt) -- ++(2.0pt,2.0pt) ++(-2.0pt,0) -- ++(2.0pt,-2.0pt);
\draw[color=dtsfsf] (-0.5866430452165177,1.9046557755048816) node {$d$};
\draw [color=qqqqff] (-1.,1.5)-- ++(-1.0pt,-1.0pt) -- ++(2.0pt,2.0pt) ++(-2.0pt,0) -- ++(2.0pt,-2.0pt);
\draw[color=rvwvcq] (-1.1025890425349942,1.8880123562365436) node {$c$};
\draw [color=sexdts] (0.5,1.5)-- ++(-1.0pt,-1.0pt) -- ++(2.0pt,2.0pt) ++(-2.0pt,0) -- ++(2.0pt,-2.0pt);
\draw[color=sexdts] (0.6865785288113356,1.4552834552597578) node {$b$};
\draw [fill=black] (1.,2.5) circle (0.5pt);
\draw [color=rvwvcq] (0.9778383660072469,0.48580428287907457)-- ++(-1.0pt,-1.0pt) -- ++(2.0pt,2.0pt) ++(-2.0pt,0) -- ++(2.0pt,-2.0pt);
\draw [fill=black] (0.5201443361279533,-0.5294442924895384) circle (0.5pt);
\draw [fill=black] (0.,-2.) circle (0.5pt);
\draw [fill=black] (0.5,-2.5) circle (0.5pt);
\draw[color=sexdts] (0.2704930471028868,-2.522493749873006) node {$\frac{b}{2}$};
\draw [fill=black] (1.,-2.) circle (0.5pt);
\draw [fill=black] (-0.5130501526429195,-0.13835360745691855) circle (0.5pt);
\draw[color=dtsfsf] (-0.5866430452165178,0.115488204158555) node {$\frac{d}{2}$};
\draw [fill=black] (-1.0149414271162547,0.23634604266358505) circle (0.5pt);
\draw[color=rvwvcq] (-0.7364338186315594,0.5815039436720169) node {$\frac{c}{2}$};
\draw [fill=black] (-0.5,-1.5) circle (0.5pt);
\draw[color=dtsfsf] (-0.17055756350806894,-1.5238885937727302) node {$\frac{d}{2}$};
\draw [fill=black] (-1.,-2.5) circle (0.5pt);
\draw[color=rvwvcq] (-0.5866430452165178,-2.0731214296278817) node {$\frac{c}{2}$};
\draw [fill=uuuuuu] (0.,2.5) circle (0.5pt);
\draw[color=uuuuuu] (-0.02492764491011184,2.724344174470524) node {$v_3$};
\draw[color=wrwrwr] (-5.513095148644552,3.67717992758287) node {$Curve$};
\draw[color=black] (2.6213760187556225,-0.5544094213920463) node {$\varphi_4$};
\draw[color=sexdts] (0.25384962783454884,-0.3754926642574137) node {$\frac{b}{2}$};
\draw[color=qqttzz] (0.7614739155188563,2.2208807416033016) node {$\frac{e}{2}$};
\draw[color=qqttzz] (0.7947607540555323,-1.6820010768219404) node {$\frac{e}{2}$};
\draw[color=qqttzz] (0.7947607540555324,0.39010462208613067) node {$e$};
\end{scriptsize}
\end{tikzpicture}
    \begin{center}
        \textnormal{Figure 14. The tropical morphism $\varphi_4:\Gamma_4' \to T_4$}
    \end{center}
\end{figure}
\textbf{Case 2.} If the metric graph $\Gamma$ has  $1$ bridge, then $\Gamma$ is one of the metric graphs given in Figure 15, 17, or 19. \par 
\textit{Solution of Case 2.} \par 
\textbf{Case 2.1.} Consider the metric graph $\Gamma $ with essential model $(G,l)$ in Figure 15, where $a,b,c,d,e$, and $f$ are real positive numbers such that $b>a$. Note that if $a=b$, then $\Gamma$ is a  hyperelliptic metric graph. The model $(G_1,l_1)$ which is obtained by subdividing $(G,l)$ is shown in Figure 15.1. The tropical modification $\Gamma'$, the metric tree $T$ with model $(G',l')$, $(T',t')$ is given in Figure 15.2, 15.3, respectively. The construction of the tropical morphism $\varphi: \Gamma' \to T$ of degree $3$ is depicted in Figure 16. 
\begin{figure}[H]
    \centering
    \begin{tikzpicture}[line cap=round,line join=round,>=triangle 45,x=1.0cm,y=1.0cm]
\definecolor{wrwrwr}{rgb}{0.3803921568627451,0.3803921568627451,0.3803921568627451}
\definecolor{wrwrwr}{rgb}{0.3803921568627451,0.3803921568627451,0.3803921568627451}
\clip(-1.5,-1.5) rectangle (3.,2.5);
\draw(2.,0.5) circle (0.4937131202354966cm);
\draw [shift={(2.,0.5)}] plot[domain=2.498091544796509:3.7850937623830774,variable=\t]({1.*2.5*cos(\t r)+0.*2.5*sin(\t r)},{0.*2.5*cos(\t r)+1.*2.5*sin(\t r)});
\draw [shift={(0.625,0.5)}] plot[domain=1.965587446494658:4.317597860684928,variable=\t]({1.*1.625*cos(\t r)+0.*1.625*sin(\t r)},{0.*1.625*cos(\t r)+1.*1.625*sin(\t r)});
\draw (0.4124309592374233,0.5341627203576358)-- (1.5070018965129957,0.5265615332932221);
\draw (0.,2.)-- (0.4124309592374233,0.5341627203576358);
\draw (0.4124309592374233,0.5341627203576358)-- (0.,-1.);
\begin{scriptsize}
\draw [fill=black] (0.,2.) circle (0.5pt);
\draw[color=black] (-0.03223848403077798,2.206423874528652) node {$v_3$};
\draw [fill=black] (0.,-1.) circle (0.5pt);
\draw[color=black] (0.0057674512912904985,-1.2521162397795862) node {$v_1$};
\draw [fill=black] (0.4124309592374233,0.5341627203576358) circle (0.5pt);
\draw[color=black] (0.6062612293799726,0.2757223601675696) node {$v_2$};
\draw [fill=black] (1.5070018965129957,0.5265615332932221) circle (0.5pt);
\draw[color=black] (1.7616416631708545,0.5265615332932221) node {$v_4$};
\draw[color=black] (1.902263623862508,1.1004511566564572) node {$e$};
\draw[color=black] (-0.2716758765598094,0.49995737856777406) node {$d$};
\draw[color=black] (-1.2142230725471077,0.47715381737453294) node {$c$};
\draw[color=black] (0.9293116796175547,0.7127906163713579) node {$f$};
\draw[color=wrwrwr] (-3.935448041607212,4.262544975452561) node {$Curve$};
\draw[color=black] (0.3668238368509411,1.419701013361833) node {$b$};
\draw[color=black] (0.3364190885932863,-0.48059575274159455) node {$a$};
\end{scriptsize}
\end{tikzpicture}
    \begin{center}
        \textnormal{Figure 15. The essential model $(G,l)$ of $\Gamma$}
    \end{center}
\end{figure}
\begin{figure}[H]
    \centering
    \begin{tikzpicture}[line cap=round,line join=round,>=triangle 45,x=1.0cm,y=1.0cm]
\definecolor{wrwrwr}{rgb}{0.3803921568627451,0.3803921568627451,0.3803921568627451}
\definecolor{wvvxds}{rgb}{0.396078431372549,0.3411764705882353,0.8235294117647058}
\definecolor{sexdts}{rgb}{0.1803921568627451,0.49019607843137253,0.19607843137254902}
\definecolor{rvwvcq}{rgb}{0.08235294117647059,0.396078431372549,0.7529411764705882}
\definecolor{qqqqff}{rgb}{0.,0.,1.}
\definecolor{dtsfsf}{rgb}{0.8274509803921568,0.1843137254901961,0.1843137254901961}
\clip(-3.,-1.5) rectangle (3.,3.);
\draw (-1.5,0.5)-- (0.,0.5);
\draw (-1.5,2.5)-- (0.,2.5);
\draw [shift={(-0.75,1.5)},color=dtsfsf]  plot[domain=2.214297435588181:4.068887871591405,variable=\t]({1.*1.25*cos(\t r)+0.*1.25*sin(\t r)},{0.*1.25*cos(\t r)+1.*1.25*sin(\t r)});
\draw [shift={(-1.5,1.5)},color=rvwvcq]  plot[domain=1.5707963267948966:4.71238898038469,variable=\t]({1.*1.*cos(\t r)+0.*1.*sin(\t r)},{0.*1.*cos(\t r)+1.*1.*sin(\t r)});
\draw [shift={(-0.75,1.5)},color=sexdts]  plot[domain=-0.9272952180016123:0.9272952180016122,variable=\t]({1.*1.25*cos(\t r)+0.*1.25*sin(\t r)},{0.*1.25*cos(\t r)+1.*1.25*sin(\t r)});
\draw (0.,0.5)-- (1.,0.5);
\draw [color=wvvxds] (1.5,0.5) circle (0.5cm);
\draw [color=sexdts] (0.,-2.)-- (0.5,-2.5);
\draw (0.,-2.)-- (1.,-2.);
\draw [color=wvvxds] (1.,-2.)-- (2.,-2.);
\draw (0.,-2.)-- (-1.5,-2.);
\draw [color=rvwvcq] (-1.5,-2.)-- (-2.5,-2.5);
\draw (3.6294648150576125,1.4536427996091146) node[anchor=north west] {$\Gamma'$};
\draw (3.5626499389332755,-1.864829381232944) node[anchor=north west] {$T$};
\draw [->] (3.746390848275196,0.34562943721386136) -- (3.768662473649975,-1.3024708405197774);
\begin{scriptsize}
\draw [fill=black] (-1.5,0.5) circle (0.5pt);
\draw[color=black] (-1.3649471752365672,0.35119734355755816) node {$v_1$};
\draw [fill=black] (0.,0.5) circle (0.5pt);
\draw[color=black] (0.13838753756101208,0.32892571818277927) node {$v_2$};
\draw[color=black] (-0.7914528218360091,0.067234120029127) node {$a$};
\draw [fill=black] (-1.5,2.5) circle (0.5pt);
\draw[color=black] (-1.4985769274852407,2.812211947470628) node {$v_3$};
\draw [fill=black] (0.,2.5) circle (0.5pt);
\draw[color=black] (0.13838753756101205,2.77880450940846) node {$v_2''$};
\draw[color=black] (-0.7691811964612302,2.0828162164466186) node {$a$};
\draw [color=dtsfsf] (-2.,1.5)-- ++(-1.0pt,-1.0pt) -- ++(2.0pt,2.0pt) ++(-2.0pt,0) -- ++(2.0pt,-2.0pt);
\draw[color=dtsfsf] (-1.6656141177960828,1.5872725518577881) node {$v_3'$};
\draw[color=dtsfsf] (-1.9718489666992935,2.1273594671961766) node {$d$};
\draw [color=qqqqff] (-2.5,1.5)-- ++(-1.0pt,-1.0pt) -- ++(2.0pt,2.0pt) ++(-2.0pt,0) -- ++(2.0pt,-2.0pt);
\draw[color=qqqqff] (-2.7903311992224196,1.6206799899199564) node {$v_3''$};
\draw[color=rvwvcq] (-2.4618247249444303,2.105087841821398) node {$c$};
\draw [color=sexdts] (0.5,1.5)-- ++(-1.0pt,-1.0pt) -- ++(2.0pt,2.0pt) ++(-2.0pt,0) -- ++(2.0pt,-2.0pt);
\draw[color=sexdts] (0.7842646734296016,1.5984083645451774) node {$v_2'$};
\draw[color=sexdts] (-0.0787608098430827,1.526025582077146) node {$b-a$};
\draw [fill=black] (1.,0.5) circle (0.5pt);
\draw[color=black] (1.2965120570495174,0.6073210353675157) node {$v_4$};
\draw[color=black] (0.7230177036489596,0.8467410081463891) node {$f$};
\draw [color=rvwvcq] (2.,0.5)-- ++(-1.0pt,-1.0pt) -- ++(2.0pt,2.0pt) ++(-2.0pt,0) -- ++(2.0pt,-2.0pt);
\draw[color=rvwvcq] (2.3989575131010756,0.6073210353675157) node {$v_4'$};
\draw[color=wvvxds] (1.7363766582014017,0.40130850065081075) node {$e$};
\draw [fill=black] (0.,-2.) circle (0.5pt);
\draw[color=black] (-0.1066003415615564,-2.176632136479849) node {$w_2$};
\draw [fill=black] (0.5,-2.5) circle (0.5pt);
\draw[color=black] (0.7063139846178753,-2.599793018600648) node {$w_2'$};
\draw[color=sexdts] (0.17736288196687522,-2.761262302567795) node {$\frac{b-a}{2}$};
\draw [fill=black] (1.,-2.) circle (0.5pt);
\draw[color=black] (1.0403883652395596,-2.221175387229407) node {$w_4$};
\draw[color=black] (0.7452893290237385,-1.9149405383261966) node {$f$};
\draw [fill=black] (2.,-2.) circle (0.5pt);
\draw[color=black] (2.3321426369767386,-1.9205084446698912) node {$w_4'$};
\draw[color=wvvxds] (1.7141050328266227,-2.349237233134385) node {$\frac{e}{2}$};
\draw [fill=black] (-1.5,-2.) circle (0.5pt);
\draw[color=black] (-1.5097127401726302,-2.198903761854628) node {$w_1$};
\draw[color=black] (-0.7246379457116723,-2.371508858509164) node {$a$};
\draw [fill=black] (-2.5,-2.5) circle (0.5pt);
\draw[color=black] (-2.7569237611602513,-2.5775213932258696) node {$w_3''$};
\draw[color=rvwvcq] (-1.8716266525127885,-1.9594837890757542) node {$\frac{c}{2}$};
\draw[color=black] (4.342156827050538,-0.467284888965567) node {$\varphi$};
\draw[color=wrwrwr] (-5.3794076490404725,4.744275448732699) node {$Curve$};
\end{scriptsize}
\end{tikzpicture}
    \begin{center}
        \textnormal{Figure 15.1. The model $(G_1,l_1)$ of $\Gamma$ }
    \end{center}
\end{figure}
\begin{figure}[H]
    \centering
    \begin{tikzpicture}[line cap=round,line join=round,>=triangle 45,x=1.0cm,y=1.0cm]
\definecolor{wrwrwr}{rgb}{0.3803921568627451,0.3803921568627451,0.3803921568627451}
\definecolor{wvvxds}{rgb}{0.396078431372549,0.3411764705882353,0.8235294117647058}
\definecolor{sexdts}{rgb}{0.1803921568627451,0.49019607843137253,0.19607843137254902}
\definecolor{rvwvcq}{rgb}{0.08235294117647059,0.396078431372549,0.7529411764705882}
\definecolor{qqqqff}{rgb}{0.,0.,1.}
\definecolor{dtsfsf}{rgb}{0.8274509803921568,0.1843137254901961,0.1843137254901961}
\clip(-3.,-1.5) rectangle (3.,3.);
\draw (-1.5,0.5)-- (0.,0.5);
\draw (-1.5,2.5)-- (0.,2.5);
\draw (-1.5,-0.5)-- (0.,-0.5);
\draw [shift={(-0.75,1.5)},color=dtsfsf]  plot[domain=2.214297435588181:4.068887871591405,variable=\t]({1.*1.25*cos(\t r)+0.*1.25*sin(\t r)},{0.*1.25*cos(\t r)+1.*1.25*sin(\t r)});
\draw [shift={(-1.5,1.5)},color=rvwvcq]  plot[domain=1.5707963267948966:4.71238898038469,variable=\t]({1.*1.*cos(\t r)+0.*1.*sin(\t r)},{0.*1.*cos(\t r)+1.*1.*sin(\t r)});
\draw [shift={(-0.75,1.5)},color=sexdts]  plot[domain=-0.9272952180016123:0.9272952180016122,variable=\t]({1.*1.25*cos(\t r)+0.*1.25*sin(\t r)},{0.*1.25*cos(\t r)+1.*1.25*sin(\t r)});
\draw (0.,2.5)-- (1.,2.5);
\draw (0.,0.5)-- (1.,0.5);
\draw (1.,0.5)-- (0.,-0.5);
\draw [color=wvvxds] (1.,2.5)-- (2.,2.5);
\draw [color=wvvxds] (1.5,0.5) circle (0.5cm);
\draw [color=sexdts] (0.,-0.5)-- (0.5,-1.);
\draw [color=sexdts] (0.,-2.)-- (0.5,-2.5);
\draw (0.,-2.)-- (1.,-2.);
\draw [color=wvvxds] (1.,-2.)-- (2.,-2.);
\draw (0.,-2.)-- (-1.5,-2.);
\draw [color=dtsfsf] (-1.5,-0.5)-- (-2.,0.);
\draw [color=rvwvcq] (-1.5,-0.5)-- (-2.5,-1.);
\draw [color=rvwvcq] (-1.5,-2.)-- (-2.5,-2.5);
\draw (3.6294648150576125,1.4536427996091146) node[anchor=north west] {$\Gamma'$};
\draw (3.5626499389332755,-1.864829381232944) node[anchor=north west] {$T$};
\draw [->] (3.746390848275196,0.34562943721386136) -- (3.768662473649975,-1.3024708405197774);
\begin{scriptsize}
\draw [fill=black] (-1.5,-0.5) circle (0.5pt);
\draw[color=black] (-1.3538113625491777,-0.6621616109948825) node {$x_1$};
\draw [fill=black] (-1.5,0.5) circle (0.5pt);
\draw[color=black] (-1.3649471752365672,0.35119734355755816) node {$v_1$};
\draw [fill=black] (0.,0.5) circle (0.5pt);
\draw[color=black] (0.13838753756101208,0.32892571818277927) node {$v_2$};
\draw[color=black] (-0.7914528218360091,0.067234120029127) node {$a$};
\draw [fill=black] (-1.5,2.5) circle (0.5pt);
\draw[color=black] (-1.4985769274852407,2.812211947470628) node {$v_3$};
\draw [fill=black] (0.,2.5) circle (0.5pt);
\draw[color=black] (0.13838753756101205,2.8233477601580175) node {$v_2''$};
\draw[color=black] (-0.7691811964612302,2.0828162164466186) node {$a$};
\draw [fill=black] (0.,-0.5) circle (0.5pt);
\draw[color=black] (-0.20682265574806166,-0.12764260200018848) node {$x_2$};
\draw[color=black] (-0.7691811964612302,-0.8459025203368085) node {$a$};
\draw [color=dtsfsf] (-2.,1.5)-- ++(-1.0pt,-1.0pt) -- ++(2.0pt,2.0pt) ++(-2.0pt,0) -- ++(2.0pt,-2.0pt);
\draw[color=dtsfsf] (-1.6656141177960828,1.5872725518577881) node {$v_3'$};
\draw[color=dtsfsf] (-1.9718489666992935,2.1273594671961766) node {$d$};
\draw [color=qqqqff] (-2.5,1.5)-- ++(-1.0pt,-1.0pt) -- ++(2.0pt,2.0pt) ++(-2.0pt,0) -- ++(2.0pt,-2.0pt);
\draw[color=qqqqff] (-2.7903311992224196,1.6206799899199564) node {$v_3''$};
\draw[color=rvwvcq] (-2.4618247249444303,2.105087841821398) node {$c$};
\draw [color=sexdts] (0.5,1.5)-- ++(-1.0pt,-1.0pt) -- ++(2.0pt,2.0pt) ++(-2.0pt,0) -- ++(2.0pt,-2.0pt);
\draw[color=sexdts] (0.7842646734296016,1.5984083645451774) node {$v_2'$};
\draw[color=sexdts] (-0.0787608098430827,1.526025582077146) node {$b-a$};
\draw [fill=black] (1.,2.5) circle (0.5pt);
\draw[color=black] (1.1294748667386754,2.812211947470628) node {$x_4$};
\draw[color=black] (0.7230177036489596,2.105087841821398) node {$f$};
\draw [fill=black] (1.,0.5) circle (0.5pt);
\draw[color=black] (1.2965120570495174,0.6073210353675157) node {$v_4$};
\draw[color=black] (0.7230177036489596,0.8467410081463891) node {$f$};
\draw[color=black] (0.7118818909615701,-0.3336551367168935) node {$f$};
\draw [fill=black] (2.,2.5) circle (0.5pt);
\draw[color=black] (2.220784510102844,2.812211947470628) node {$x_4'$};
\draw[color=wvvxds] (1.7809199089509598,2.0716804037592293) node {$\frac{e}{2}$};
\draw [color=rvwvcq] (2.,0.5)-- ++(-1.0pt,-1.0pt) -- ++(2.0pt,2.0pt) ++(-2.0pt,0) -- ++(2.0pt,-2.0pt);
\draw[color=rvwvcq] (2.3989575131010756,0.6073210353675157) node {$v_4'$};
\draw[color=wvvxds] (1.7363766582014017,0.40130850065081075) node {$e$};
\draw [fill=black] (0.5,-1.) circle (0.5pt);
\draw[color=black] (0.7842646734296016,-0.8960136774300611) node {$x_2'$};
\draw[color=sexdts] (0.12168381852992788,-1.0574829613972083) node {$\frac{b-a}{2}$};
\draw [fill=black] (0.,-2.) circle (0.5pt);
\draw[color=black] (-0.1066003415615564,-2.176632136479849) node {$w_2$};
\draw [fill=black] (0.5,-2.5) circle (0.5pt);
\draw[color=black] (0.7063139846178753,-2.599793018600648) node {$w_2'$};
\draw[color=sexdts] (0.17736288196687522,-2.761262302567795) node {$\frac{b-a}{2}$};
\draw [fill=black] (1.,-2.) circle (0.5pt);
\draw[color=black] (1.0403883652395596,-2.221175387229407) node {$w_4$};
\draw[color=black] (0.7452893290237385,-1.9149405383261966) node {$f$};
\draw [fill=black] (2.,-2.) circle (0.5pt);
\draw[color=black] (2.3321426369767386,-1.9205084446698912) node {$w_4'$};
\draw[color=wvvxds] (1.7141050328266227,-2.349237233134385) node {$\frac{e}{2}$};
\draw [fill=black] (-1.5,-2.) circle (0.5pt);
\draw[color=black] (-1.5097127401726302,-2.198903761854628) node {$w_1$};
\draw[color=black] (-0.7246379457116723,-2.371508858509164) node {$a$};
\draw [fill=black] (-2.,0.) circle (0.5pt);
\draw[color=black] (-1.72129318123303,0.2843824674332215) node {$x_3'$};
\draw[color=dtsfsf] (-2.1388861570101354,0.11177737077868484) node {$\frac{d}{2}$};
\draw [fill=black] (-2.5,-1.) circle (0.5pt);
\draw[color=black] (-2.7791953865350303,-0.9071494901174506) node {$x_3''$};
\draw[color=rvwvcq] (-1.8604908398253988,-0.5006923270277354) node {$\frac{c}{2}$};
\draw [fill=black] (-2.5,-2.5) circle (0.5pt);
\draw[color=black] (-2.7569237611602513,-2.5775213932258696) node {$w_3''$};
\draw[color=rvwvcq] (-1.8716266525127885,-1.9594837890757542) node {$\frac{c}{2}$};
\draw[color=black] (4.342156827050538,-0.467284888965567) node {$\varphi$};
\draw[color=wrwrwr] (-5.3794076490404725,4.744275448732699) node {$Curve$};
\end{scriptsize}
\end{tikzpicture}
    \begin{center}
        \textnormal{Figure 15.2. The model $(G',l')$ of $\Gamma'$}
    \end{center}
\end{figure}
\begin{figure}[H]
    \centering
    \begin{tikzpicture}[line cap=round,line join=round,>=triangle 45,x=1.0cm,y=1.0cm]
\definecolor{wrwrwr}{rgb}{0.3803921568627451,0.3803921568627451,0.3803921568627451}
\definecolor{cqcqcq}{rgb}{0.7529411764705882,0.7529411764705882,0.7529411764705882}
\definecolor{wvvxds}{rgb}{0.396078431372549,0.3411764705882353,0.8235294117647058}
\definecolor{sexdts}{rgb}{0.1803921568627451,0.49019607843137253,0.19607843137254902}
\definecolor{rvwvcq}{rgb}{0.08235294117647059,0.396078431372549,0.7529411764705882}
\definecolor{qqqqff}{rgb}{0.,0.,1.}
\definecolor{dtsfsf}{rgb}{0.8274509803921568,0.1843137254901961,0.1843137254901961}
\clip(-3.,-3.5) rectangle (3.,-1.1);
\draw (-1.5,0.5)-- (0.,0.5);
\draw (-1.5,2.5)-- (0.,2.5);
\draw (-1.5,-0.5)-- (0.,-0.5);
\draw [shift={(-0.75,1.5)},color=dtsfsf]  plot[domain=2.214297435588181:4.068887871591405,variable=\t]({1.*1.25*cos(\t r)+0.*1.25*sin(\t r)},{0.*1.25*cos(\t r)+1.*1.25*sin(\t r)});
\draw [shift={(-1.5,1.5)},color=rvwvcq]  plot[domain=1.5707963267948966:4.71238898038469,variable=\t]({1.*1.*cos(\t r)+0.*1.*sin(\t r)},{0.*1.*cos(\t r)+1.*1.*sin(\t r)});
\draw [shift={(-0.75,1.5)},color=sexdts]  plot[domain=-0.9272952180016123:0.9272952180016122,variable=\t]({1.*1.25*cos(\t r)+0.*1.25*sin(\t r)},{0.*1.25*cos(\t r)+1.*1.25*sin(\t r)});
\draw (0.,2.5)-- (1.,2.5);
\draw (0.,0.5)-- (1.,0.5);
\draw (1.,0.5)-- (0.,-0.5);
\draw [color=wvvxds] (1.,2.5)-- (2.,2.5);
\draw [color=wvvxds] (1.5,0.5) circle (0.5cm);
\draw [color=sexdts] (0.,-0.5)-- (0.5,-1.);
\draw [color=sexdts] (0.,-2.)-- (0.5,-2.5);
\draw (0.,-2.)-- (1.,-2.);
\draw [color=wvvxds] (1.,-2.)-- (2.,-2.);
\draw (0.,-2.)-- (-1.5,-2.);
\draw [color=dtsfsf] (-1.5,-0.5)-- (-2.,0.);
\draw [color=rvwvcq] (-1.5,-0.5)-- (-2.5,-1.);
\draw [color=dtsfsf] (-1.5,-2.)-- (-2.,-1.5);
\draw [color=rvwvcq] (-1.5,-2.)-- (-2.5,-2.5);
\draw [line width=0.4pt,dash pattern=on 1pt off 1pt,color=cqcqcq] (-2.,1.5)-- (-2.,0.);
\draw [line width=0.4pt,dash pattern=on 1pt off 1pt,color=cqcqcq] (-2.5,1.5)-- (-2.5,-1.);
\draw (3.6294648150576134,1.453642799609114) node[anchor=north west] {$\Gamma'$};
\draw (3.5626499389332764,-1.8648293812329446) node[anchor=north west] {$T$};
\draw [->] (3.746390848275196,0.34562943721386136) -- (3.768662473649975,-1.3024708405197774);
\begin{scriptsize}
\draw [fill=black] (-1.5,-0.5) circle (0.5pt);
\draw[color=black] (-1.353811362549177,-0.6621616109948831) node {$x_1$};
\draw [fill=black] (-1.5,0.5) circle (0.5pt);
\draw[color=black] (-1.3649471752365665,0.35119734355755755) node {$v_1$};
\draw [fill=black] (0.,0.5) circle (0.5pt);
\draw[color=black] (0.13838753756101269,0.32892571818277866) node {$v_2$};
\draw[color=black] (-0.7914528218360085,0.06723412002912638) node {$a$};
\draw [fill=black] (-1.5,2.5) circle (0.5pt);
\draw[color=black] (-1.49857692748524,2.8122119474706277) node {$v_3$};
\draw [fill=black] (0.,2.5) circle (0.5pt);
\draw[color=black] (0.13838753756101269,2.823347760158017) node {$v_2''$};
\draw[color=black] (-0.7691811964612295,2.082816216446618) node {$a$};
\draw [fill=black] (0.,-0.5) circle (0.5pt);
\draw[color=black] (-0.20682265574806105,-0.12764260200018915) node {$x_2$};
\draw[color=black] (-0.7691811964612295,-0.8459025203368091) node {$a$};
\draw [color=dtsfsf] (-2.,1.5)-- ++(-1.0pt,-1.0pt) -- ++(2.0pt,2.0pt) ++(-2.0pt,0) -- ++(2.0pt,-2.0pt);
\draw[color=dtsfsf] (-1.6656141177960824,1.5872725518577875) node {$v_3'$};
\draw[color=dtsfsf] (-1.9718489666992927,2.127359467196176) node {$d$};
\draw [color=qqqqff] (-2.5,1.5)-- ++(-1.0pt,-1.0pt) -- ++(2.0pt,2.0pt) ++(-2.0pt,0) -- ++(2.0pt,-2.0pt);
\draw[color=qqqqff] (-2.790331199222419,1.6206799899199558) node {$v_3''$};
\draw[color=rvwvcq] (-2.46182472494443,2.105087841821397) node {$c$};
\draw [color=sexdts] (0.5,1.5)-- ++(-1.0pt,-1.0pt) -- ++(2.0pt,2.0pt) ++(-2.0pt,0) -- ++(2.0pt,-2.0pt);
\draw[color=sexdts] (0.7842646734296023,1.598408364545177) node {$v_2'$};
\draw[color=sexdts] (-0.0787608098430827,1.5260255820771453) node {$b-a$};
\draw [fill=black] (1.,2.5) circle (0.5pt);
\draw[color=black] (1.129474866738676,2.8122119474706277) node {$x_4$};
\draw[color=black] (0.7230177036489601,2.105087841821397) node {$f$};
\draw [fill=black] (1.,0.5) circle (0.5pt);
\draw[color=black] (1.296512057049518,0.6073210353675151) node {$v_4$};
\draw[color=black] (0.7230177036489601,0.8467410081463884) node {$f$};
\draw[color=black] (0.7118818909615707,-0.3336551367168941) node {$f$};
\draw [fill=black] (2.,2.5) circle (0.5pt);
\draw[color=black] (2.2207845101028445,2.8122119474706277) node {$x_4'$};
\draw[color=wvvxds] (1.7809199089509595,2.071680403759229) node {$\frac{e}{2}$};
\draw [color=rvwvcq] (2.,0.5)-- ++(-1.0pt,-1.0pt) -- ++(2.0pt,2.0pt) ++(-2.0pt,0) -- ++(2.0pt,-2.0pt);
\draw[color=rvwvcq] (2.398957513101076,0.6073210353675151) node {$v_4'$};
\draw[color=wvvxds] (1.7363766582014024,0.4013085006508101) node {$e$};
\draw [fill=black] (0.5,-1.) circle (0.5pt);
\draw [fill=black] (0.,-2.) circle (0.5pt);
\draw[color=black] (-0.10660034156155576,-2.1766321364798498) node {$w_2$};
\draw [fill=black] (0.5,-2.5) circle (0.5pt);
\draw[color=black] (0.706313984617876,-2.599793018600649) node {$w_2'$};
\draw[color=sexdts] (0.17736288196687522,-2.761262302567796) node {$\frac{b-a}{2}$};
\draw [fill=black] (1.,-2.) circle (0.5pt);
\draw[color=black] (1.0403883652395602,-2.2211753872294073) node {$w_4$};
\draw[color=black] (0.7230177036489601,-1.714495909953187) node {$f$};
\draw [fill=black] (2.,-2.) circle (0.5pt);
\draw[color=black] (2.332142636976739,-1.9205084446698917) node {$w_4'$};
\draw[color=wvvxds] (1.7141050328266227,-2.349237233134386) node {$\frac{e}{2}$};
\draw [fill=black] (-1.5,-2.) circle (0.5pt);
\draw[color=black] (-1.376082987923956,-2.232311199916797) node {$w_1$};
\draw[color=black] (-0.7691811964612295,-2.2935581696974388) node {$a$};
\draw [fill=black] (-2.,0.) circle (0.5pt);
\draw[color=black] (-1.7212931812330297,0.2843824674332208) node {$x_3'$};
\draw[color=dtsfsf] (-2.138886157010136,0.1117773707786842) node {$\frac{d}{2}$};
\draw [fill=black] (-2.5,-1.) circle (0.5pt);
\draw[color=rvwvcq] (-1.860490839825399,-0.5006923270277359) node {$\frac{c}{2}$};
\draw [fill=black] (-2.,-1.5) circle (0.5pt);
\draw[color=black] (-1.7435648066078087,-1.3080387468634715) node {$w_3'$};
\draw[color=dtsfsf] (-2.138886157010136,-1.4138289673936715) node {$\frac{d}{2}$};
\draw [fill=black] (-2.5,-2.5) circle (0.5pt);
\draw[color=black] (-2.756923761160251,-2.57752139322587) node {$w_3''$};
\draw[color=rvwvcq] (-1.8716266525127883,-1.959483789075755) node {$\frac{c}{2}$};
\draw[color=black] (4.342156827050539,-0.4672848889655676) node {$\varphi$};
\draw[color=wrwrwr] (-5.401679274415252,4.755411261420088) node {$Curve$};
\end{scriptsize}
\end{tikzpicture}

    \begin{center}
        \textnormal{Figure 15.3. The model $(T',t')$ of  $T$}
    \end{center}
\end{figure}
\begin{figure}[H]
    \centering
    \begin{tikzpicture}[line cap=round,line join=round,>=triangle 45,x=1.0cm,y=1.0cm]
\definecolor{wrwrwr}{rgb}{0.3803921568627451,0.3803921568627451,0.3803921568627451}
\definecolor{cqcqcq}{rgb}{0.7529411764705882,0.7529411764705882,0.7529411764705882}
\definecolor{wvvxds}{rgb}{0.396078431372549,0.3411764705882353,0.8235294117647058}
\definecolor{sexdts}{rgb}{0.1803921568627451,0.49019607843137253,0.19607843137254902}
\definecolor{rvwvcq}{rgb}{0.08235294117647059,0.396078431372549,0.7529411764705882}
\definecolor{qqqqff}{rgb}{0.,0.,1.}
\definecolor{dtsfsf}{rgb}{0.8274509803921568,0.1843137254901961,0.1843137254901961}
\clip(-3.,-3.5) rectangle (5.,3.);
\draw (-1.5,0.5)-- (0.,0.5);
\draw (-1.5,2.5)-- (0.,2.5);
\draw (-1.5,-0.5)-- (0.,-0.5);
\draw [shift={(-0.75,1.5)},color=dtsfsf]  plot[domain=2.214297435588181:4.068887871591405,variable=\t]({1.*1.25*cos(\t r)+0.*1.25*sin(\t r)},{0.*1.25*cos(\t r)+1.*1.25*sin(\t r)});
\draw [shift={(-1.5,1.5)},color=rvwvcq]  plot[domain=1.5707963267948966:4.71238898038469,variable=\t]({1.*1.*cos(\t r)+0.*1.*sin(\t r)},{0.*1.*cos(\t r)+1.*1.*sin(\t r)});
\draw [shift={(-0.75,1.5)},color=sexdts]  plot[domain=-0.9272952180016123:0.9272952180016122,variable=\t]({1.*1.25*cos(\t r)+0.*1.25*sin(\t r)},{0.*1.25*cos(\t r)+1.*1.25*sin(\t r)});
\draw (0.,2.5)-- (1.,2.5);
\draw (0.,0.5)-- (1.,0.5);
\draw (1.,0.5)-- (0.,-0.5);
\draw [color=wvvxds] (1.,2.5)-- (2.,2.5);
\draw [color=wvvxds] (1.5,0.5) circle (0.5cm);
\draw [color=sexdts] (0.,-0.5)-- (0.5,-1.);
\draw [color=sexdts] (0.,-2.)-- (0.5,-2.5);
\draw (0.,-2.)-- (1.,-2.);
\draw [color=wvvxds] (1.,-2.)-- (2.,-2.);
\draw (0.,-2.)-- (-1.5,-2.);
\draw [color=dtsfsf] (-1.5,-0.5)-- (-2.,0.);
\draw [color=rvwvcq] (-1.5,-0.5)-- (-2.5,-1.);
\draw [color=dtsfsf] (-1.5,-2.)-- (-2.,-1.5);
\draw [color=rvwvcq] (-1.5,-2.)-- (-2.5,-2.5);
\draw [line width=0.4pt,dash pattern=on 2pt off 2pt,color=cqcqcq] (-1.5,2.5)-- (-1.5,-2.);
\draw [line width=0.4pt,dash pattern=on 2pt off 2pt,color=cqcqcq] (-2.,1.5)-- (-2.,0.);
\draw [line width=0.4pt,dash pattern=on 2pt off 2pt,color=cqcqcq] (-2.,0.)-- (-2.,-1.5);
\draw [line width=0.4pt,dash pattern=on 2pt off 2pt,color=cqcqcq] (-2.5,1.5)-- (-2.5,-1.);
\draw [line width=0.4pt,dash pattern=on 2pt off 2pt,color=cqcqcq] (-2.5,-1.)-- (-2.5,-2.5);
\draw [line width=0.4pt,dash pattern=on 2pt off 2pt,color=cqcqcq] (0.,2.5)-- (0.,-2.);
\draw [line width=0.4pt,dash pattern=on 2pt off 2pt,color=cqcqcq] (0.5,1.5)-- (0.5,-2.5);
\draw [line width=0.4pt,dash pattern=on 2pt off 2pt,color=cqcqcq] (1.,2.5)-- (1.,-2.);
\draw [line width=0.4pt,dash pattern=on 2pt off 2pt,color=cqcqcq] (2.,2.5)-- (2.,-2.);
\draw (3.640600627744997,1.4536427996091132) node[anchor=north west] {$\Gamma'$};
\draw (3.5737857516206604,-1.8648293812329453) node[anchor=north west] {$T$};
\draw [->] (3.746390848275196,0.34562943721386136) -- (3.768662473649975,-1.3024708405197774);
\begin{scriptsize}
\draw [fill=black] (-1.5,-0.5) circle (0.5pt);
\draw [fill=black] (-1.5,0.5) circle (0.5pt);
\draw[color=black] (-1.3649471752365645,0.35119734355755694) node {$v_1$};
\draw [fill=black] (0.,0.5) circle (0.5pt);
\draw[color=black] (0.1383875375610125,0.32892571818277805) node {$v_2$};
\draw[color=black] (-0.7914528218360073,0.06723412002912574) node {$a$};
\draw [fill=black] (-1.5,2.5) circle (0.5pt);
\draw[color=black] (-1.4985769274852379,2.812211947470627) node {$v_3$};
\draw [fill=black] (0.,2.5) circle (0.5pt);
\draw[color=black] (-0.7691811964612283,2.0828162164466177) node {$a$};
\draw [fill=black] (0.,-0.5) circle (0.5pt);
\draw[color=black] (-0.7691811964612283,-0.8459025203368098) node {$a$};
\draw [color=dtsfsf] (-2.,1.5)-- ++(-1.0pt,-1.0pt) -- ++(2.0pt,2.0pt) ++(-2.0pt,0) -- ++(2.0pt,-2.0pt);
\draw[color=dtsfsf] (-1.7714043383262796,1.459210705952808) node {$d$};
\draw [color=qqqqff] (-2.5,1.5)-- ++(-1.0pt,-1.0pt) -- ++(2.0pt,2.0pt) ++(-2.0pt,0) -- ++(2.0pt,-2.0pt);
\draw[color=rvwvcq] (-2.227972658509247,1.4035316425158608) node {$c$};
\draw [color=sexdts] (0.5,1.5)-- ++(-1.0pt,-1.0pt) -- ++(2.0pt,2.0pt) ++(-2.0pt,0) -- ++(2.0pt,-2.0pt);
\draw[color=sexdts] (-0.07876080984308192,1.5260255820771447) node {$b-a$};
\draw [fill=black] (1.,2.5) circle (0.5pt);
\draw[color=black] (0.7230177036489591,2.1050878418213967) node {$f$};
\draw [fill=black] (1.,0.5) circle (0.5pt);
\draw[color=black] (1.296512057049516,0.6073210353675145) node {$v_4$};
\draw[color=black] (0.7230177036489591,0.8467410081463879) node {$f$};
\draw[color=black] (0.7118818909615696,-0.3336551367168948) node {$f$};
\draw [fill=black] (2.,2.5) circle (0.5pt);
\draw[color=wvvxds] (1.7809199089509569,2.071680403759228) node {$\frac{e}{2}$};
\draw [color=rvwvcq] (2.,0.5)-- ++(-1.0pt,-1.0pt) -- ++(2.0pt,2.0pt) ++(-2.0pt,0) -- ++(2.0pt,-2.0pt);
\draw[color=wvvxds] (1.7363766582013997,0.4013085006508095) node {$e$};
\draw [fill=black] (0.5,-1.) circle (0.5pt);
\draw[color=sexdts] (0.12168381852992766,-1.0574829613972094) node {$\frac{b-a}{2}$};
\draw [fill=black] (0.,-2.) circle (0.5pt);
\draw [fill=black] (0.5,-2.5) circle (0.5pt);
\draw[color=sexdts] (0.17736288196687494,-2.7612623025677965) node {$\frac{b-a}{2}$};
\draw [fill=black] (1.,-2.) circle (0.5pt);
\draw[color=black] (0.7230177036489591,-1.7144959099531876) node {$f$};
\draw [fill=black] (2.,-2.) circle (0.5pt);
\draw[color=wvvxds] (1.7141050328266203,-2.3492372331343865) node {$\frac{e}{2}$};
\draw [fill=black] (-1.5,-2.) circle (0.5pt);
\draw[color=black] (-0.7246379457116705,-2.3715088585091655) node {$a$};
\draw [fill=black] (-2.,0.) circle (0.5pt);
\draw[color=dtsfsf] (-2.1388861570101323,0.11177737077868358) node {$\frac{d}{2}$};
\draw [fill=black] (-2.5,-1.) circle (0.5pt);
\draw[color=rvwvcq] (-1.8604908398253959,-0.5006923270277366) node {$\frac{c}{2}$};
\draw [fill=black] (-2.,-1.5) circle (0.5pt);
\draw[color=dtsfsf] (-2.1388861570101323,-1.4138289673936721) node {$\frac{d}{2}$};
\draw [fill=black] (-2.5,-2.5) circle (0.5pt);
\draw[color=rvwvcq] (-1.8716266525127851,-1.9594837890757555) node {$\frac{c}{2}$};
\draw[color=black] (4.342156827050533,-0.46728488896556825) node {$\varphi$};
\draw[color=wrwrwr] (-5.2012346460422325,4.5438308203596876) node {$Curve$};
\end{scriptsize}
\end{tikzpicture}
    \begin{center}
        \textnormal{Figure 16. The tropical morphism $\varphi:\Gamma' \to T$}
    \end{center}
\end{figure}
\textbf{Case 2.2.} Consider the metric graph $\Gamma_1$ with essential model $(G,l)$ in Figure 17, where $a,b,c,d$, and $e$ are real positive numbers such that $b>a$. Note that if $a=b$, then $\Gamma_2$ is a hyperelliptic metric graph. The model $(G_1,l_1)$ that obtained by subdividing $(G,l)$ is shown in Figure 17.1. The tropical modification $\Gamma_1'$, the metric tree $T_1$ with model $(G_1',l_1')$, $(T_1',t_1')$ is given in Figure 17.2, 17.3, respectively. The construction of the tropical morphism $\varphi: \Gamma_1' \to T_1$ of degree $3$ is depicted in Figure 18.
\begin{figure}[H]
    \centering
    \begin{tikzpicture}[line cap=round,line join=round,>=triangle 45,x=1.0cm,y=1.0cm]
\definecolor{wrwrwr}{rgb}{0.3803921568627451,0.3803921568627451,0.3803921568627451}
\clip(-1.5,-1.) rectangle (3.,1.);
\draw (1.,0.)-- (1.5,0.);
\draw(2.,0.) circle (0.5cm);
\draw (0.,0.)-- (0.,0.);
\draw(-0.5,0.) circle (0.5cm);
\draw (0.,0.)-- (1.,0.);
\draw [shift={(0.5,0.)}] plot[domain=0.:3.141592653589793,variable=\t]({1.*0.5*cos(\t r)+0.*0.5*sin(\t r)},{0.*0.5*cos(\t r)+1.*0.5*sin(\t r)});
\begin{scriptsize}
\draw [fill=black] (0.,0.) circle (0.5pt);
\draw[color=black] (-0.20031206938676055,-0.040462119414003364) node {$v_1$};
\draw [fill=black] (1.,0.) circle (0.5pt);
\draw[color=black] (0.9928824007267678,-0.18961142817819357) node {$v_2$};
\draw [fill=black] (1.5,0.) circle (0.5pt);
\draw[color=black] (1.7645679547675823,0.004931148470750171) node {$v_4$};
\draw[color=black] (1.2490301266478785,0.24162461672696503) node {$d$};
\draw[color=black] (1.9364138974741503,0.676103037909606) node {$e$};
\draw[color=wrwrwr] (-3.108723590288486,2.0184468167873177) node {$Curve$};
\draw[color=black] (-0.4370055376429768,0.7214963057943595) node {$c$};
\draw[color=black] (0.4967988302719585,0.12489907073759879) node {$a$};
\draw[color=black] (0.5292225930467827,0.7214963057943595) node {$b$};
\end{scriptsize}
\end{tikzpicture}
    \begin{center}
        \textnormal{Figure 17. The essential model $(G,l)$ of $\Gamma_1$}
    \end{center}
\end{figure}
\begin{figure}[H]
    \centering
    \begin{tikzpicture}[line cap=round,line join=round,>=triangle 45,x=1.0cm,y=1.0cm]
\definecolor{wrwrwr}{rgb}{0.3803921568627451,0.3803921568627451,0.3803921568627451}
\definecolor{wvvxds}{rgb}{0.396078431372549,0.3411764705882353,0.8235294117647058}
\definecolor{rvwvcq}{rgb}{0.08235294117647059,0.396078431372549,0.7529411764705882}
\definecolor{sexdts}{rgb}{0.1803921568627451,0.49019607843137253,0.19607843137254902}
\definecolor{qqqqff}{rgb}{0.,0.,1.}
\clip(-3.,-0.5) rectangle (3.,3.);
\draw (-1.5,0.5)-- (0.,0.5);
\draw (-1.504144833828935,0.5182345338683995)-- (0.,2.5);
\draw [shift={(-0.75,1.5)},color=sexdts]  plot[domain=-0.9272952180016123:0.9272952180016122,variable=\t]({1.*1.25*cos(\t r)+0.*1.25*sin(\t r)},{0.*1.25*cos(\t r)+1.*1.25*sin(\t r)});
\draw (0.,0.5)-- (1.,0.5);
\draw [color=wvvxds] (1.5,0.5) circle (0.5cm);
\draw [color=sexdts] (0.,-2.)-- (0.5,-2.5);
\draw (0.,-2.)-- (1.,-2.);
\draw [color=wvvxds] (1.,-2.)-- (2.,-2.);
\draw (0.,-2.)-- (-1.5,-2.);
\draw [color=rvwvcq] (-1.5,-2.)-- (-2.5,-2.5);
\draw (3.6294648150576134,1.453642799609114) node[anchor=north west] {$\Gamma_1'$};
\draw (3.5626499389332764,-1.8648293812329446) node[anchor=north west] {$T_1$};
\draw [->] (3.746390848275196,0.34562943721386136) -- (3.768662473649975,-1.3024708405197774);
\draw [color=rvwvcq] (-1.9941205920740719,0.50709872118101) circle (0.4901022852345269cm);
\begin{scriptsize}
\draw [fill=black] (-1.5,0.5) circle (0.5pt);
\draw[color=black] (-1.2647248610500612,0.26211084205844193) node {$v_1$};
\draw [fill=black] (0.,0.5) circle (0.5pt);
\draw[color=black] (0.038165223374507416,0.2732466547458313) node {$v_2$};
\draw[color=black] (-0.6578230695873348,0.17859224690302095) node {$a$};
\draw [fill=black] (0.,2.5) circle (0.5pt);
\draw[color=black] (0.13838753756101269,2.823347760158017) node {$v_2''$};
\draw[color=black] (-0.6244156315251663,1.1696795760806828) node {$a$};
\draw [color=qqqqff] (-2.484096350319209,0.4959629084936205)-- ++(-1.0pt,-1.0pt) -- ++(2.0pt,2.0pt) ++(-2.0pt,0) -- ++(2.0pt,-2.0pt);
\draw[color=qqqqff] (-2.77919538653503,0.6184568480549045) node {$v_3''$};
\draw [color=sexdts] (0.5,1.5)-- ++(-1.0pt,-1.0pt) -- ++(2.0pt,2.0pt) ++(-2.0pt,0) -- ++(2.0pt,-2.0pt);
\draw[color=sexdts] (0.7842646734296023,1.598408364545177) node {$v_2'$};
\draw[color=sexdts] (-0.0787608098430827,1.5260255820771453) node {$b-a$};
\draw [fill=black] (1.,0.5) circle (0.5pt);
\draw[color=black] (1.2853762443621286,0.4291480323692838) node {$v_4$};
\draw[color=black] (0.7230177036489601,0.8467410081463884) node {$f$};
\draw [color=rvwvcq] (2.,0.5)-- ++(-1.0pt,-1.0pt) -- ++(2.0pt,2.0pt) ++(-2.0pt,0) -- ++(2.0pt,-2.0pt);
\draw[color=rvwvcq] (2.398957513101076,0.6073210353675151) node {$v_4'$};
\draw[color=wvvxds] (1.7363766582014024,0.4013085006508101) node {$e$};
\draw [fill=black] (0.5,-1.) circle (0.5pt);
\draw[color=black] (0.7842646734296023,-0.8960136774300618) node {$x_2'$};
\draw [fill=black] (0.,-2.) circle (0.5pt);
\draw[color=black] (-0.10660034156155576,-2.1766321364798498) node {$w_2$};
\draw [fill=black] (0.5,-2.5) circle (0.5pt);
\draw[color=black] (0.706313984617876,-2.599793018600649) node {$w_2'$};
\draw[color=sexdts] (0.17736288196687522,-2.761262302567796) node {$\frac{b-a}{2}$};
\draw [fill=black] (1.,-2.) circle (0.5pt);
\draw[color=black] (1.0403883652395602,-2.2211753872294073) node {$w_4$};
\draw[color=black] (0.7230177036489601,-1.714495909953187) node {$f$};
\draw [fill=black] (2.,-2.) circle (0.5pt);
\draw[color=black] (2.332142636976739,-1.9205084446698917) node {$w_4'$};
\draw[color=wvvxds] (1.7141050328266227,-2.349237233134386) node {$\frac{e}{2}$};
\draw [fill=black] (-1.5,-2.) circle (0.5pt);
\draw[color=black] (-1.5097127401726296,-2.1989037618546283) node {$w_1$};
\draw[color=black] (-0.7246379457116716,-2.371508858509165) node {$a$};
\draw [fill=black] (-2.5,-1.) circle (0.5pt);
\draw[color=black] (-2.77919538653503,-0.9071494901174513) node {$x_3''$};
\draw [fill=black] (-2.5,-2.5) circle (0.5pt);
\draw[color=black] (-2.756923761160251,-2.57752139322587) node {$w_3''$};
\draw[color=rvwvcq] (-1.8716266525127883,-1.959483789075755) node {$\frac{c}{2}$};
\draw[color=wrwrwr] (-5.401679274415252,4.755411261420088) node {$Curve$};
\draw[color=rvwvcq] (-2.2613800965714193,0.39017268796342064) node {$c$};
\end{scriptsize}
\end{tikzpicture}
    \begin{center}
        \textnormal{Figure 17.1. The model $(G_1,l_1)$ of $\Gamma_1$}
    \end{center}
\end{figure}
\begin{figure}[H]
    \centering
    \begin{tikzpicture}[line cap=round,line join=round,>=triangle 45,x=1.0cm,y=1.0cm]
\definecolor{wrwrwr}{rgb}{0.3803921568627451,0.3803921568627451,0.3803921568627451}
\definecolor{rvwvcq}{rgb}{0.08235294117647059,0.396078431372549,0.7529411764705882}
\definecolor{wvvxds}{rgb}{0.396078431372549,0.3411764705882353,0.8235294117647058}
\definecolor{sexdts}{rgb}{0.1803921568627451,0.49019607843137253,0.19607843137254902}
\definecolor{qqqqff}{rgb}{0.,0.,1.}
\clip(-3.,-1.5) rectangle (3.,3.);
\draw (-1.5,0.5)-- (0.,0.5);
\draw (-1.504144833828935,0.5182345338683995)-- (0.,2.5);
\draw (-1.5,-0.5)-- (0.,-0.5);
\draw [shift={(-0.75,1.5)},color=sexdts]  plot[domain=-0.9272952180016123:0.9272952180016122,variable=\t]({1.*1.25*cos(\t r)+0.*1.25*sin(\t r)},{0.*1.25*cos(\t r)+1.*1.25*sin(\t r)});
\draw (0.,2.5)-- (1.,2.5);
\draw (0.,0.5)-- (1.,0.5);
\draw (1.,0.5)-- (0.,-0.5);
\draw [color=wvvxds] (1.,2.5)-- (2.,2.5);
\draw [color=wvvxds] (1.5,0.5) circle (0.5cm);
\draw [color=sexdts] (0.,-0.5)-- (0.5,-1.);
\draw [color=sexdts] (0.,-2.)-- (0.5,-2.5);
\draw (0.,-2.)-- (1.,-2.);
\draw [color=wvvxds] (1.,-2.)-- (2.,-2.);
\draw (0.,-2.)-- (-1.5,-2.);
\draw [color=rvwvcq] (-1.5,-0.5)-- (-2.5,-1.);
\draw [color=rvwvcq] (-1.5,-2.)-- (-2.5,-2.5);
\draw (3.6294648150576134,1.453642799609114) node[anchor=north west] {$\Gamma_1'$};
\draw (3.5626499389332764,-1.8648293812329446) node[anchor=north west] {$T_1$};
\draw [->] (3.746390848275196,0.34562943721386136) -- (3.768662473649975,-1.3024708405197774);
\draw [color=rvwvcq] (-1.9941205920740719,0.50709872118101) circle (0.4901022852345269cm);
\begin{scriptsize}
\draw [fill=black] (-1.5,-0.5) circle (0.5pt);
\draw[color=black] (-1.353811362549177,-0.6621616109948831) node {$x_1$};
\draw [fill=black] (-1.5,0.5) circle (0.5pt);
\draw[color=black] (-1.3649471752365665,0.35119734355755755) node {$v_1$};
\draw [fill=black] (0.,0.5) circle (0.5pt);
\draw[color=black] (0.13838753756101269,0.32892571818277866) node {$v_2$};
\draw[color=black] (-0.7914528218360085,0.06723412002912638) node {$a$};
\draw [fill=black] (0.,2.5) circle (0.5pt);
\draw[color=black] (0.13838753756101269,2.823347760158017) node {$v_2''$};
\draw[color=black] (-0.6244156315251663,1.1696795760806828) node {$a$};
\draw [fill=black] (0.,-0.5) circle (0.5pt);
\draw[color=black] (-0.20682265574806105,-0.12764260200018915) node {$x_2$};
\draw[color=black] (-0.7691811964612295,-0.8459025203368091) node {$a$};
\draw [color=qqqqff] (-2.484096350319209,0.4959629084936205)-- ++(-1.0pt,-1.0pt) -- ++(2.0pt,2.0pt) ++(-2.0pt,0) -- ++(2.0pt,-2.0pt);
\draw[color=qqqqff] (-2.77919538653503,0.6184568480549045) node {$v_3''$};
\draw [color=sexdts] (0.5,1.5)-- ++(-1.0pt,-1.0pt) -- ++(2.0pt,2.0pt) ++(-2.0pt,0) -- ++(2.0pt,-2.0pt);
\draw[color=sexdts] (0.7842646734296023,1.598408364545177) node {$v_2'$};
\draw[color=sexdts] (-0.0787608098430827,1.5260255820771453) node {$b-a$};
\draw [fill=black] (1.,2.5) circle (0.5pt);
\draw[color=black] (1.129474866738676,2.8122119474706277) node {$x_4$};
\draw[color=black] (0.7230177036489601,2.105087841821397) node {$f$};
\draw [fill=black] (1.,0.5) circle (0.5pt);
\draw[color=black] (1.296512057049518,0.6073210353675151) node {$v_4$};
\draw[color=black] (0.7230177036489601,0.8467410081463884) node {$f$};
\draw[color=black] (0.7118818909615707,-0.3336551367168941) node {$f$};
\draw [fill=black] (2.,2.5) circle (0.5pt);
\draw[color=black] (2.2207845101028445,2.8122119474706277) node {$x_4'$};
\draw[color=wvvxds] (1.7809199089509595,2.071680403759229) node {$\frac{e}{2}$};
\draw [color=rvwvcq] (2.,0.5)-- ++(-1.0pt,-1.0pt) -- ++(2.0pt,2.0pt) ++(-2.0pt,0) -- ++(2.0pt,-2.0pt);
\draw[color=rvwvcq] (2.398957513101076,0.6073210353675151) node {$v_4'$};
\draw[color=wvvxds] (1.7363766582014024,0.4013085006508101) node {$e$};
\draw [fill=black] (0.5,-1.) circle (0.5pt);
\draw[color=black] (0.7842646734296023,-0.8960136774300618) node {$x_2'$};
\draw[color=sexdts] (0.12168381852992782,-1.0574829613972088) node {$\frac{b-a}{2}$};
\draw [fill=black] (0.,-2.) circle (0.5pt);
\draw[color=black] (-0.10660034156155576,-2.1766321364798498) node {$w_2$};
\draw [fill=black] (0.5,-2.5) circle (0.5pt);
\draw[color=black] (0.706313984617876,-2.599793018600649) node {$w_2'$};
\draw[color=sexdts] (0.17736288196687522,-2.761262302567796) node {$\frac{b-a}{2}$};
\draw [fill=black] (1.,-2.) circle (0.5pt);
\draw[color=black] (1.0403883652395602,-2.2211753872294073) node {$w_4$};
\draw[color=black] (0.7230177036489601,-1.714495909953187) node {$f$};
\draw [fill=black] (2.,-2.) circle (0.5pt);
\draw[color=black] (2.332142636976739,-1.9205084446698917) node {$w_4'$};
\draw[color=wvvxds] (1.7141050328266227,-2.349237233134386) node {$\frac{e}{2}$};
\draw [fill=black] (-1.5,-2.) circle (0.5pt);
\draw[color=black] (-1.5097127401726296,-2.1989037618546283) node {$w_1$};
\draw[color=black] (-0.7246379457116716,-2.371508858509165) node {$a$};
\draw [fill=black] (-2.5,-1.) circle (0.5pt);
\draw[color=black] (-2.77919538653503,-0.9071494901174513) node {$x_3''$};
\draw[color=rvwvcq] (-1.860490839825399,-0.5006923270277359) node {$\frac{c}{2}$};
\draw [fill=black] (-2.5,-2.5) circle (0.5pt);
\draw[color=black] (-2.756923761160251,-2.57752139322587) node {$w_3''$};
\draw[color=rvwvcq] (-1.8716266525127883,-1.959483789075755) node {$\frac{c}{2}$};
\draw[color=wrwrwr] (-5.401679274415252,4.755411261420088) node {$Curve$};
\draw[color=rvwvcq] (-2.2613800965714193,0.39017268796342064) node {$c$};
\end{scriptsize}
\end{tikzpicture}
    \begin{center}
        \textnormal{Figure 17.2. The model $(G_1',l_1')$ of $\Gamma_1'$}
    \end{center}
\end{figure}
\begin{figure}[H]
    \centering
    \begin{tikzpicture}[line cap=round,line join=round,>=triangle 45,x=1.0cm,y=1.0cm]
\definecolor{wrwrwr}{rgb}{0.3803921568627451,0.3803921568627451,0.3803921568627451}
\definecolor{cqcqcq}{rgb}{0.7529411764705882,0.7529411764705882,0.7529411764705882}
\definecolor{rvwvcq}{rgb}{0.08235294117647059,0.396078431372549,0.7529411764705882}
\definecolor{wvvxds}{rgb}{0.396078431372549,0.3411764705882353,0.8235294117647058}
\definecolor{sexdts}{rgb}{0.1803921568627451,0.49019607843137253,0.19607843137254902}
\definecolor{qqqqff}{rgb}{0.,0.,1.}
\clip(-3.,-3.5) rectangle (3.,-1.5);
\draw (-1.5,0.5)-- (0.,0.5);
\draw (-1.504144833828935,0.5182345338683995)-- (0.,2.5);
\draw (-1.5,-0.5)-- (0.,-0.5);
\draw [shift={(-0.75,1.5)},color=sexdts]  plot[domain=-0.9272952180016123:0.9272952180016122,variable=\t]({1.*1.25*cos(\t r)+0.*1.25*sin(\t r)},{0.*1.25*cos(\t r)+1.*1.25*sin(\t r)});
\draw (0.,2.5)-- (1.,2.5);
\draw (0.,0.5)-- (1.,0.5);
\draw (1.,0.5)-- (0.,-0.5);
\draw [color=wvvxds] (1.,2.5)-- (2.,2.5);
\draw [color=wvvxds] (1.5,0.5) circle (0.5cm);
\draw [color=sexdts] (0.,-0.5)-- (0.5,-1.);
\draw [color=sexdts] (0.,-2.)-- (0.5,-2.5);
\draw (0.,-2.)-- (1.,-2.);
\draw [color=wvvxds] (1.,-2.)-- (2.,-2.);
\draw (0.,-2.)-- (-1.5,-2.);
\draw [color=rvwvcq] (-1.5,-0.5)-- (-2.5,-1.);
\draw [color=rvwvcq] (-1.5,-2.)-- (-2.5,-2.5);
\draw [line width=0.4pt,dash pattern=on 1pt off 1pt,color=cqcqcq] (-2.484096350319209,0.4959629084936205)-- (-2.5,-1.);
\draw (3.6294648150576134,1.453642799609114) node[anchor=north west] {$\Gamma_1'$};
\draw (3.5626499389332764,-1.8648293812329446) node[anchor=north west] {$T_1$};
\draw [->] (3.746390848275196,0.34562943721386136) -- (3.768662473649975,-1.3024708405197774);
\draw [color=rvwvcq] (-1.9941205920740719,0.50709872118101) circle (0.4901022852345269cm);
\begin{scriptsize}
\draw [fill=black] (-1.5,-0.5) circle (0.5pt);
\draw[color=black] (-1.353811362549177,-0.6621616109948831) node {$x_1$};
\draw [fill=black] (-1.5,0.5) circle (0.5pt);
\draw[color=black] (-1.3649471752365665,0.35119734355755755) node {$v_1$};
\draw [fill=black] (0.,0.5) circle (0.5pt);
\draw[color=black] (0.13838753756101269,0.32892571818277866) node {$v_2$};
\draw[color=black] (-0.7914528218360085,0.06723412002912638) node {$a$};
\draw [fill=black] (0.,2.5) circle (0.5pt);
\draw[color=black] (0.13838753756101269,2.823347760158017) node {$v_2''$};
\draw[color=black] (-0.6244156315251663,1.1696795760806828) node {$a$};
\draw [fill=black] (0.,-0.5) circle (0.5pt);
\draw[color=black] (-0.20682265574806105,-0.12764260200018915) node {$x_2$};
\draw[color=black] (-0.7691811964612295,-0.8459025203368091) node {$a$};
\draw [color=qqqqff] (-2.484096350319209,0.4959629084936205)-- ++(-1.0pt,-1.0pt) -- ++(2.0pt,2.0pt) ++(-2.0pt,0) -- ++(2.0pt,-2.0pt);
\draw[color=qqqqff] (-2.77919538653503,0.6184568480549045) node {$v_3''$};
\draw [color=sexdts] (0.5,1.5)-- ++(-1.0pt,-1.0pt) -- ++(2.0pt,2.0pt) ++(-2.0pt,0) -- ++(2.0pt,-2.0pt);
\draw[color=sexdts] (0.7842646734296023,1.598408364545177) node {$v_2'$};
\draw[color=sexdts] (-0.0787608098430827,1.5260255820771453) node {$b-a$};
\draw [fill=black] (1.,2.5) circle (0.5pt);
\draw[color=black] (1.129474866738676,2.8122119474706277) node {$x_4$};
\draw[color=black] (0.7230177036489601,2.105087841821397) node {$f$};
\draw [fill=black] (1.,0.5) circle (0.5pt);
\draw[color=black] (1.296512057049518,0.6073210353675151) node {$v_4$};
\draw[color=black] (0.7230177036489601,0.8467410081463884) node {$f$};
\draw[color=black] (0.7118818909615707,-0.3336551367168941) node {$f$};
\draw [fill=black] (2.,2.5) circle (0.5pt);
\draw[color=black] (2.2207845101028445,2.8122119474706277) node {$x_4'$};
\draw[color=wvvxds] (1.7809199089509595,2.071680403759229) node {$\frac{e}{2}$};
\draw [color=rvwvcq] (2.,0.5)-- ++(-1.0pt,-1.0pt) -- ++(2.0pt,2.0pt) ++(-2.0pt,0) -- ++(2.0pt,-2.0pt);
\draw[color=rvwvcq] (2.398957513101076,0.6073210353675151) node {$v_4'$};
\draw[color=wvvxds] (1.7363766582014024,0.4013085006508101) node {$e$};
\draw [fill=black] (0.5,-1.) circle (0.5pt);
\draw[color=black] (0.7842646734296023,-0.8960136774300618) node {$x_2'$};
\draw[color=sexdts] (0.12168381852992782,-1.0574829613972088) node {$\frac{b-a}{2}$};
\draw [fill=black] (0.,-2.) circle (0.5pt);
\draw[color=black] (-0.10660034156155576,-2.1766321364798498) node {$w_2$};
\draw [fill=black] (0.5,-2.5) circle (0.5pt);
\draw[color=black] (0.706313984617876,-2.599793018600649) node {$w_2'$};
\draw[color=sexdts] (0.17736288196687522,-2.761262302567796) node {$\frac{b-a}{2}$};
\draw [fill=black] (1.,-2.) circle (0.5pt);
\draw[color=black] (1.0403883652395602,-2.2211753872294073) node {$w_4$};
\draw[color=black] (0.7230177036489601,-1.714495909953187) node {$f$};
\draw [fill=black] (2.,-2.) circle (0.5pt);
\draw[color=black] (2.332142636976739,-1.9205084446698917) node {$w_4'$};
\draw[color=wvvxds] (1.7141050328266227,-2.349237233134386) node {$\frac{e}{2}$};
\draw [fill=black] (-1.5,-2.) circle (0.5pt);
\draw[color=black] (-1.5097127401726296,-2.1989037618546283) node {$w_1$};
\draw[color=black] (-0.7246379457116716,-2.371508858509165) node {$a$};
\draw [fill=black] (-2.5,-1.) circle (0.5pt);
\draw[color=black] (-2.77919538653503,-0.9071494901174513) node {$x_3''$};
\draw[color=rvwvcq] (-1.860490839825399,-0.5006923270277359) node {$\frac{c}{2}$};
\draw [fill=black] (-2.5,-2.5) circle (0.5pt);
\draw[color=black] (-2.756923761160251,-2.57752139322587) node {$w_3''$};
\draw[color=rvwvcq] (-1.8716266525127883,-1.959483789075755) node {$\frac{c}{2}$};
\draw[color=wrwrwr] (-5.401679274415252,4.755411261420088) node {$Curve$};
\draw[color=rvwvcq] (-2.2613800965714193,0.39017268796342064) node {$c$};
\end{scriptsize}
\end{tikzpicture}
    \begin{center}
        \textnormal{Figure 17.3. The model $(T_1',t_1')$ of $T_1$}
    \end{center}
\end{figure}
\begin{figure}[H]
    \centering
    
\begin{tikzpicture}[line cap=round,line join=round,>=triangle 45,x=1.0cm,y=1.0cm]
\definecolor{wrwrwr}{rgb}{0.3803921568627451,0.3803921568627451,0.3803921568627451}
\definecolor{cqcqcq}{rgb}{0.7529411764705882,0.7529411764705882,0.7529411764705882}
\definecolor{rvwvcq}{rgb}{0.08235294117647059,0.396078431372549,0.7529411764705882}
\definecolor{wvvxds}{rgb}{0.396078431372549,0.3411764705882353,0.8235294117647058}
\definecolor{sexdts}{rgb}{0.1803921568627451,0.49019607843137253,0.19607843137254902}
\definecolor{qqqqff}{rgb}{0.,0.,1.}
\clip(-3.5,-3.5) rectangle (4.5,3.5);
\draw (-1.5,0.5)-- (0.,0.5);
\draw (-1.504144833828935,0.5182345338683995)-- (0.,2.5);
\draw (-1.5,-0.5)-- (0.,-0.5);
\draw [shift={(-0.75,1.5)},color=sexdts]  plot[domain=-0.9272952180016123:0.9272952180016122,variable=\t]({1.*1.25*cos(\t r)+0.*1.25*sin(\t r)},{0.*1.25*cos(\t r)+1.*1.25*sin(\t r)});
\draw (0.,2.5)-- (1.,2.5);
\draw (0.,0.5)-- (1.,0.5);
\draw (1.,0.5)-- (0.,-0.5);
\draw [color=wvvxds] (1.,2.5)-- (2.,2.5);
\draw [color=wvvxds] (1.5,0.5) circle (0.5cm);
\draw [color=sexdts] (0.,-0.5)-- (0.5,-1.);
\draw [color=sexdts] (0.,-2.)-- (0.5,-2.5);
\draw (0.,-2.)-- (1.,-2.);
\draw [color=wvvxds] (1.,-2.)-- (2.,-2.);
\draw (0.,-2.)-- (-1.5,-2.);
\draw [color=rvwvcq] (-1.5,-0.5)-- (-2.5,-1.);
\draw [color=rvwvcq] (-1.5,-2.)-- (-2.5,-2.5);
\draw [line width=0.4pt,dash pattern=on 2pt off 2pt,color=cqcqcq] (-1.504144833828935,0.5182345338683995)-- (-1.5,-2.);
\draw [line width=0.4pt,dash pattern=on 2pt off 2pt,color=cqcqcq] (-2.484096350319209,0.4959629084936205)-- (-2.5,-1.);
\draw [line width=0.4pt,dash pattern=on 2pt off 2pt,color=cqcqcq] (-2.5,-1.)-- (-2.5,-2.5);
\draw [line width=0.4pt,dash pattern=on 2pt off 2pt,color=cqcqcq] (0.,2.5)-- (0.,-2.);
\draw [line width=0.4pt,dash pattern=on 2pt off 2pt,color=cqcqcq] (0.5,1.5)-- (0.5,-2.5);
\draw [line width=0.4pt,dash pattern=on 2pt off 2pt,color=cqcqcq] (1.,2.5)-- (1.,-2.);
\draw [line width=0.4pt,dash pattern=on 2pt off 2pt,color=cqcqcq] (2.,2.5)-- (2.,-2.);
\draw (3.6294648150576134,1.453642799609114) node[anchor=north west] {$\Gamma_1'$};
\draw (3.5626499389332764,-1.8648293812329446) node[anchor=north west] {$T_1$};
\draw [->] (3.746390848275196,0.34562943721386136) -- (3.768662473649975,-1.3024708405197774);
\draw [color=rvwvcq] (-1.9941205920740719,0.50709872118101) circle (0.4901022852345269cm);
\draw (3.952403382991908,-0.17218585274974696) node[anchor=north west] {$\varphi_1$};
\begin{scriptsize}
\draw [fill=black] (-1.5,-0.5) circle (0.5pt);
\draw[color=black] (-1.353811362549177,-0.6621616109948831) node {$x_1$};
\draw [fill=black] (-1.5,0.5) circle (0.5pt);
\draw[color=black] (-1.3649471752365665,0.35119734355755755) node {$v_1$};
\draw [fill=black] (0.,0.5) circle (0.5pt);
\draw[color=black] (0.13838753756101269,0.32892571818277866) node {$v_2$};
\draw[color=black] (-0.7914528218360085,0.06723412002912638) node {$a$};
\draw [fill=black] (0.,2.5) circle (0.5pt);
\draw[color=black] (0.13838753756101269,2.823347760158017) node {$v_2''$};
\draw[color=black] (-0.6244156315251663,1.1696795760806828) node {$a$};
\draw [fill=black] (0.,-0.5) circle (0.5pt);
\draw[color=black] (-0.20682265574806105,-0.12764260200018915) node {$x_2$};
\draw[color=black] (-0.7691811964612295,-0.8459025203368091) node {$a$};
\draw [color=qqqqff] (-2.484096350319209,0.4959629084936205)-- ++(-1.0pt,-1.0pt) -- ++(2.0pt,2.0pt) ++(-2.0pt,0) -- ++(2.0pt,-2.0pt);
\draw[color=qqqqff] (-2.77919538653503,0.6184568480549045) node {$v_3''$};
\draw [color=sexdts] (0.5,1.5)-- ++(-1.0pt,-1.0pt) -- ++(2.0pt,2.0pt) ++(-2.0pt,0) -- ++(2.0pt,-2.0pt);
\draw[color=sexdts] (0.7842646734296023,1.598408364545177) node {$v_2'$};
\draw[color=sexdts] (-0.0787608098430827,1.5260255820771453) node {$b-a$};
\draw [fill=black] (1.,2.5) circle (0.5pt);
\draw[color=black] (1.129474866738676,2.8122119474706277) node {$x_4$};
\draw[color=black] (0.7230177036489601,2.105087841821397) node {$f$};
\draw [fill=black] (1.,0.5) circle (0.5pt);
\draw[color=black] (1.296512057049518,0.6073210353675151) node {$v_4$};
\draw[color=black] (0.7230177036489601,0.8467410081463884) node {$f$};
\draw[color=black] (0.7118818909615707,-0.3336551367168941) node {$f$};
\draw [fill=black] (2.,2.5) circle (0.5pt);
\draw[color=black] (2.2207845101028445,2.8122119474706277) node {$x_4'$};
\draw[color=wvvxds] (1.7809199089509595,2.071680403759229) node {$\frac{e}{2}$};
\draw [color=rvwvcq] (2.,0.5)-- ++(-1.0pt,-1.0pt) -- ++(2.0pt,2.0pt) ++(-2.0pt,0) -- ++(2.0pt,-2.0pt);
\draw[color=rvwvcq] (2.398957513101076,0.6073210353675151) node {$v_4'$};
\draw[color=wvvxds] (1.7363766582014024,0.4013085006508101) node {$e$};
\draw [fill=black] (0.5,-1.) circle (0.5pt);
\draw[color=black] (0.7842646734296023,-0.8960136774300618) node {$x_2'$};
\draw[color=sexdts] (0.12168381852992782,-1.0574829613972088) node {$\frac{b-a}{2}$};
\draw [fill=black] (0.,-2.) circle (0.5pt);
\draw[color=black] (-0.10660034156155576,-2.1766321364798498) node {$w_2$};
\draw [fill=black] (0.5,-2.5) circle (0.5pt);
\draw[color=black] (0.706313984617876,-2.599793018600649) node {$w_2'$};
\draw[color=sexdts] (0.17736288196687522,-2.761262302567796) node {$\frac{b-a}{2}$};
\draw [fill=black] (1.,-2.) circle (0.5pt);
\draw[color=black] (1.0403883652395602,-2.2211753872294073) node {$w_4$};
\draw[color=black] (0.7230177036489601,-1.714495909953187) node {$f$};
\draw [fill=black] (2.,-2.) circle (0.5pt);
\draw[color=black] (2.332142636976739,-1.9205084446698917) node {$w_4'$};
\draw[color=wvvxds] (1.7141050328266227,-2.349237233134386) node {$\frac{e}{2}$};
\draw [fill=black] (-1.5,-2.) circle (0.5pt);
\draw[color=black] (-1.5097127401726296,-2.1989037618546283) node {$w_1$};
\draw[color=black] (-0.7246379457116716,-2.371508858509165) node {$a$};
\draw [fill=black] (-2.5,-1.) circle (0.5pt);
\draw[color=black] (-2.77919538653503,-0.9071494901174513) node {$x_3''$};
\draw[color=rvwvcq] (-1.860490839825399,-0.5006923270277359) node {$\frac{c}{2}$};
\draw [fill=black] (-2.5,-2.5) circle (0.5pt);
\draw[color=black] (-2.756923761160251,-2.57752139322587) node {$w_3''$};
\draw[color=rvwvcq] (-1.8716266525127883,-1.959483789075755) node {$\frac{c}{2}$};
\draw[color=wrwrwr] (-5.401679274415252,4.755411261420088) node {$Curve$};
\draw[color=rvwvcq] (-2.2613800965714193,0.39017268796342064) node {$c$};
\end{scriptsize}
\end{tikzpicture}
    \begin{center}
        \textnormal{Figure 18. The tropical morphism $\varphi_1:\Gamma_1 \to T_1$}
    \end{center}
\end{figure}
\textbf{Case 2.3.} 
 Consider the metric graph $\Gamma_2 $ with essential model in Figure 19, where $a,b,c,d,e$, and $f$ are real positive numbers. The model $(G_1,l_1)$ which obtained by subdividing $(G,l)$ is shown in Figure 19.1. The tropical modification $\Gamma_2'$, the metric tree $T_2$ with model $(G_2',l_2')$, $(T_2',t_2')$ 7is given in Figure 19.2, 19.3, respectively. The construction of the tropical morphism $\varphi_2: \Gamma_2' \to T_2$ of degree $3$ is depicted in Figure 20. 
 \begin{figure}[H]
    \centering
    \begin{tikzpicture}[line cap=round,line join=round,>=triangle 45,x=1.0cm,y=1.0cm]
\definecolor{wrwrwr}{rgb}{0.3803921568627451,0.3803921568627451,0.3803921568627451}
\clip(-1.5,-1.5) rectangle (3.,2.5);
\draw (0.5,-0.5)-- (0.5,2.);
\draw [shift={(1.25,0.75)}] plot[domain=2.1112158270654806:4.171969480114106,variable=\t]({1.*1.457737973711325*cos(\t r)+0.*1.457737973711325*sin(\t r)},{0.*1.457737973711325*cos(\t r)+1.*1.457737973711325*sin(\t r)});
\draw [shift={(0.17605729651367463,0.75)}] plot[domain=1.3172207104540103:4.965964596725576,variable=\t]({1.*1.2912934891580727*cos(\t r)+0.*1.2912934891580727*sin(\t r)},{0.*1.2912934891580727*cos(\t r)+1.*1.2912934891580727*sin(\t r)});
\draw(2.,-0.5) circle (0.5cm);
\draw (1.5,-0.5)-- (0.5,-0.5);
\begin{scriptsize}
\draw [fill=black] (0.5,-0.5) circle (0.5pt);
\draw[color=black] (0.5155239218507216,-0.8017657742098742) node {$v_1$};
\draw [fill=black] (0.5,2.) circle (0.5pt);
\draw[color=black] (0.7804053224005073,2.0551693317199593) node {$v_3$};
\draw[color=black] (0.7756752973906896,0.7449524040004827) node {$b$};
\draw[color=black] (0.04725144587877921,0.7449524040004827) node {$d$};
\draw[color=black] (-0.7663128558098481,0.7260323039612122) node {$c$};
\draw [fill=black] (1.5,-0.5) circle (0.5pt);
\draw[color=black] (1.8210108245603793,-0.4328238234441009) node {$v_4$};
\draw[color=black] (2.024401899982536,0.2057295528812757) node {$e$};
\draw[color=black] (1.0405566979404752,-0.21997269800230873) node {$f$};
\draw[color=wrwrwr] (-3.3394464611506227,3.261325709223448) node {$Curve$};
\end{scriptsize}
\end{tikzpicture}
    \begin{center}
        \textnormal{Figure 19. The essential model $(G,l)$ of $\Gamma_2$}
    \end{center}
\end{figure}
\begin{figure}[H]
    \centering
    \begin{tikzpicture}[line cap=round,line join=round,>=triangle 45,x=1.0cm,y=1.0cm]
\definecolor{wrwrwr}{rgb}{0.3803921568627451,0.3803921568627451,0.3803921568627451}
\definecolor{uuuuuu}{rgb}{0.26666666666666666,0.26666666666666666,0.26666666666666666}
\definecolor{wvvxds}{rgb}{0.396078431372549,0.3411764705882353,0.8235294117647058}
\definecolor{sexdts}{rgb}{0.1803921568627451,0.49019607843137253,0.19607843137254902}
\definecolor{rvwvcq}{rgb}{0.08235294117647059,0.396078431372549,0.7529411764705882}
\definecolor{qqqqff}{rgb}{0.,0.,1.}
\definecolor{dtsfsf}{rgb}{0.8274509803921568,0.1843137254901961,0.1843137254901961}
\clip(-2.,-0.2) rectangle (2.5,3.);
\draw [shift={(0.75,1.5)},color=dtsfsf]  plot[domain=2.214297435588181:4.068887871591405,variable=\t]({1.*1.25*cos(\t r)+0.*1.25*sin(\t r)},{0.*1.25*cos(\t r)+1.*1.25*sin(\t r)});
\draw [shift={(0.,1.5)},color=rvwvcq]  plot[domain=1.5707963267948966:4.71238898038469,variable=\t]({1.*1.*cos(\t r)+0.*1.*sin(\t r)},{0.*1.*cos(\t r)+1.*1.*sin(\t r)});
\draw [shift={(-0.7464501113912964,1.5063728582396283)},color=sexdts]  plot[domain=-0.9287282318863292:0.9287282318863297,variable=\t]({1.*1.2464664028803656*cos(\t r)+0.*1.2464664028803656*sin(\t r)},{0.*1.2464664028803656*cos(\t r)+1.*1.2464664028803656*sin(\t r)});
\draw (0.,0.5081290379665255)-- (1.,0.5);
\draw [color=wvvxds] (1.5,0.5) circle (0.5cm);
\draw [color=sexdts] (0.,-0.5)-- (0.5,-1.);
\draw [color=sexdts] (0.,-2.)-- (0.5,-2.5);
\draw (0.,-2.)-- (1.,-2.);
\draw [color=wvvxds] (1.,-2.)-- (2.,-2.);
\draw [color=dtsfsf] (0.,-2.)-- (-0.5,-1.5);
\draw [color=rvwvcq] (0.,-2.)-- (-1.,-2.5);
\draw (3.3786274559977745,1.0891963363344506) node[anchor=north west] {$\Gamma_2'$};
\draw (3.3885602303288502,-1.761509896684261) node[anchor=north west] {$T_2$};
\draw [->] (3.567350168288212,-0.12260213205677878) -- (3.577282942619288,-1.105946790833268);
\begin{scriptsize}
\draw [fill=black] (0.,0.5) circle (0.5pt);
\draw[color=black] (0.15544218556372277,0.3343054871727011) node {$v_1$};
\draw [color=dtsfsf] (-0.5,1.5)-- ++(-1.0pt,-1.0pt) -- ++(2.0pt,2.0pt) ++(-2.0pt,0) -- ++(2.0pt,-2.0pt);
\draw[color=dtsfsf] (-0.7385075042330866,1.4765745352464008) node {$v_3'$};
\draw[color=dtsfsf] (-0.6044150507635653,1.9781796389657211) node {$d$};
\draw [color=qqqqff] (-1.,1.5)-- ++(-1.0pt,-1.0pt) -- ++(2.0pt,2.0pt) ++(-2.0pt,0) -- ++(2.0pt,-2.0pt);
\draw[color=qqqqff] (-1.314608415435475,1.4865073095774766) node {$v_3''$};
\draw[color=rvwvcq] (-1.1010537673173482,1.9583140903035698) node {$c$};
\draw [color=sexdts] (0.5,1.5)-- ++(-1.0pt,-1.0pt) -- ++(2.0pt,2.0pt) ++(-2.0pt,0) -- ++(2.0pt,-2.0pt);
\draw[color=sexdts] (0.18524050855694976,1.4368434379220982) node {$v_2'$};
\draw[color=sexdts] (0.7265767096005732,1.4517425994187116) node {$b$};
\draw [fill=black] (1.,0.5) circle (0.5pt);
\draw[color=black] (1.2977112336374237,0.40383490749023065) node {$v_4$};
\draw[color=black] (0.7067111609384219,0.8160450422298701) node {$f$};
\draw [color=rvwvcq] (2.,0.5)-- ++(-1.0pt,-1.0pt) -- ++(2.0pt,2.0pt) ++(-2.0pt,0) -- ++(2.0pt,-2.0pt);
\draw[color=rvwvcq] (2.300921441076065,0.5826248454495924) node {$v_4'$};
\draw[color=wvvxds] (1.7099213683770635,0.4088012946557685) node {$e$};
\draw [fill=black] (0.5,-1.) circle (0.5pt);
\draw[color=black] (0.7315430967661111,-1.0066190475225116) node {$x_2'$};
\draw[color=sexdts] (0.2100724443846389,-0.9917198860258982) node {$\frac{b}{2}$};
\draw [fill=black] (0.,-2.) circle (0.5pt);
\draw[color=black] (-0.033280526726714776,-2.238283064575892) node {$w_1$};
\draw [fill=black] (0.5,-2.5) circle (0.5pt);
\draw[color=black] (0.6918119994418085,-2.645526812149994) node {$w_2'$};
\draw[color=sexdts] (0.2696690903710929,-2.5710310046669265) node {$\frac{b}{2}$};
\draw [fill=black] (1.,-2.) circle (0.5pt);
\draw[color=black] (0.9897952293740782,-2.158820869927287) node {$w_4$};
\draw[color=black] (0.6967783866073463,-1.7664762838497987) node {$f$};
\draw [fill=black] (2.,-2.) circle (0.5pt);
\draw[color=black] (2.0426693084680982,-2.168753644258363) node {$w_4'$};
\draw[color=wvvxds] (1.7595852400324419,-1.6671485405390423) node {$\frac{e}{2}$};
\draw [fill=black] (-1.,-1.) circle (0.5pt);
\draw[color=black] (-1.2450789951179453,-0.9470224015360578) node {$x_3''$};
\draw [fill=black] (-0.5,-1.5) circle (0.5pt);
\draw[color=black] (-0.7683058272263136,-1.5727871843938235) node {$w_3'$};
\draw[color=dtsfsf] (-0.1673729801962362,-1.7565435095187232) node {$\frac{d}{2}$};
\draw [fill=black] (-1.,-2.5) circle (0.5pt);
\draw[color=black] (-1.2450789951179453,-2.675325135143221) node {$w_3''$};
\draw[color=rvwvcq] (-0.6044150507635653,-2.0346611907888414) node {$\frac{c}{2}$};
\draw [fill=uuuuuu] (0.,2.5) circle (0.5pt);
\draw[color=uuuuuu] (-0.023347752395639104,2.7678351982862353) node {$v_3$};
\draw[color=wrwrwr] (-4.64705420351136,3.845541213207943) node {$Curve$};
\draw[color=black] (4.049089723345381,-0.5993752999484099) node {$\varphi_2$};
\end{scriptsize}
\end{tikzpicture}
    \begin{center}
        \textnormal{Figure 19.1.  The model $(G_1,l_1)$ of $\Gamma_2$}
    \end{center}
\end{figure}
\begin{figure}[H]
    \centering
    \begin{tikzpicture}[line cap=round,line join=round,>=triangle 45,x=1.0cm,y=1.0cm]
\definecolor{wrwrwr}{rgb}{0.3803921568627451,0.3803921568627451,0.3803921568627451}
\definecolor{uuuuuu}{rgb}{0.26666666666666666,0.26666666666666666,0.26666666666666666}
\definecolor{wvvxds}{rgb}{0.396078431372549,0.3411764705882353,0.8235294117647058}
\definecolor{sexdts}{rgb}{0.1803921568627451,0.49019607843137253,0.19607843137254902}
\definecolor{rvwvcq}{rgb}{0.08235294117647059,0.396078431372549,0.7529411764705882}
\definecolor{qqqqff}{rgb}{0.,0.,1.}
\definecolor{dtsfsf}{rgb}{0.8274509803921568,0.1843137254901961,0.1843137254901961}
\clip(-2.,-1.5) rectangle (2.5,3.);
\draw [shift={(0.75,1.5)},color=dtsfsf]  plot[domain=2.214297435588181:4.068887871591405,variable=\t]({1.*1.25*cos(\t r)+0.*1.25*sin(\t r)},{0.*1.25*cos(\t r)+1.*1.25*sin(\t r)});
\draw [shift={(0.,1.5)},color=rvwvcq]  plot[domain=1.5707963267948966:4.71238898038469,variable=\t]({1.*1.*cos(\t r)+0.*1.*sin(\t r)},{0.*1.*cos(\t r)+1.*1.*sin(\t r)});
\draw [shift={(-0.7464501113912964,1.5063728582396283)},color=sexdts]  plot[domain=-0.9287282318863292:0.9287282318863297,variable=\t]({1.*1.2464664028803656*cos(\t r)+0.*1.2464664028803656*sin(\t r)},{0.*1.2464664028803656*cos(\t r)+1.*1.2464664028803656*sin(\t r)});
\draw (0.,2.5046166785127313)-- (1.,2.5);
\draw (0.,0.5081290379665255)-- (1.,0.5);
\draw (1.,0.5)-- (0.,-0.5);
\draw [color=wvvxds] (1.,2.5)-- (2.,2.5);
\draw [color=wvvxds] (1.5,0.5) circle (0.5cm);
\draw [color=sexdts] (0.,-0.5)-- (0.5,-1.);
\draw [color=sexdts] (0.,-2.)-- (0.5,-2.5);
\draw (0.,-2.)-- (1.,-2.);
\draw [color=wvvxds] (1.,-2.)-- (2.,-2.);
\draw [color=dtsfsf] (0.,-0.5)-- (-0.5,0.);
\draw [color=rvwvcq] (0.,-0.5)-- (-1.,-1.);
\draw [color=dtsfsf] (0.,-2.)-- (-0.5,-1.5);
\draw [color=rvwvcq] (0.,-2.)-- (-1.,-2.5);
\draw (3.3786274559977745,1.0891963363344506) node[anchor=north west] {$\Gamma_2'$};
\draw (3.3885602303288502,-1.761509896684261) node[anchor=north west] {$T_2$};
\draw [->] (3.567350168288212,-0.12260213205677878) -- (3.577282942619288,-1.105946790833268);
\begin{scriptsize}
\draw [fill=black] (0.,-0.5) circle (0.5pt);
\draw[color=black] (0.2448371545434037,-0.49011478230657785) node {$x_1$};
\draw [fill=black] (0.,0.5) circle (0.5pt);
\draw[color=black] (0.15544218556372277,0.3343054871727011) node {$v_1$};
\draw [color=dtsfsf] (-0.5,1.5)-- ++(-1.0pt,-1.0pt) -- ++(2.0pt,2.0pt) ++(-2.0pt,0) -- ++(2.0pt,-2.0pt);
\draw[color=dtsfsf] (-0.7881713758884649,1.4765745352464008) node {$v_3'$};
\draw[color=dtsfsf] (-0.6044150507635653,1.9781796389657211) node {$d$};
\draw [color=qqqqff] (-1.,1.5)-- ++(-1.0pt,-1.0pt) -- ++(2.0pt,2.0pt) ++(-2.0pt,0) -- ++(2.0pt,-2.0pt);
\draw[color=qqqqff] (-1.314608415435475,1.4865073095774766) node {$v_3''$};
\draw[color=rvwvcq] (-1.0315243469998185,2.0179107362900237) node {$c$};
\draw [color=sexdts] (0.5,1.5)-- ++(-1.0pt,-1.0pt) -- ++(2.0pt,2.0pt) ++(-2.0pt,0) -- ++(2.0pt,-2.0pt);
\draw[color=sexdts] (0.18524050855694976,1.4368434379220982) node {$v_2'$};
\draw[color=sexdts] (0.7265767096005732,1.4517425994187116) node {$b$};
\draw [fill=black] (1.,2.5) circle (0.5pt);
\draw[color=black] (0.9798624550430025,2.7579024239551595) node {$x_4$};
\draw[color=black] (0.7067111609384219,2.156969576925083) node {$f$};
\draw [fill=black] (1.,0.5) circle (0.5pt);
\draw[color=black] (1.1884507159955915,0.4634315534766845) node {$v_4$};
\draw[color=black] (0.7067111609384219,0.8160450422298701) node {$f$};
\draw[color=black] (0.6967783866073463,-0.25669458552630003) node {$f$};
\draw [fill=black] (2.,2.5) circle (0.5pt);
\draw[color=black] (2.2611903437517626,2.5592469373336466) node {$x_4'$};
\draw[color=wvvxds] (1.7496524657013663,2.1271712539318557) node {$\frac{e}{2}$};
\draw [color=rvwvcq] (2.,0.5)-- ++(-1.0pt,-1.0pt) -- ++(2.0pt,2.0pt) ++(-2.0pt,0) -- ++(2.0pt,-2.0pt);
\draw[color=rvwvcq] (2.300921441076065,0.5826248454495924) node {$v_4'$};
\draw[color=wvvxds] (1.7099213683770635,0.4088012946557685) node {$e$};
\draw [fill=black] (0.5,-1.) circle (0.5pt);
\draw[color=black] (0.7315430967661111,-1.0066190475225116) node {$x_2'$};
\draw[color=sexdts] (0.2100724443846389,-0.9917198860258982) node {$\frac{b}{2}$};
\draw [fill=black] (0.,-2.) circle (0.5pt);
\draw[color=black] (-0.033280526726714776,-2.238283064575892) node {$w_1$};
\draw [fill=black] (0.5,-2.5) circle (0.5pt);
\draw[color=black] (0.6918119994418085,-2.645526812149994) node {$w_2'$};
\draw[color=sexdts] (0.2696690903710929,-2.5710310046669265) node {$\frac{b}{2}$};
\draw [fill=black] (1.,-2.) circle (0.5pt);
\draw[color=black] (0.9897952293740782,-2.158820869927287) node {$w_4$};
\draw[color=black] (0.6967783866073463,-1.7664762838497987) node {$f$};
\draw [fill=black] (2.,-2.) circle (0.5pt);
\draw[color=black] (2.0426693084680982,-2.168753644258363) node {$w_4'$};
\draw[color=wvvxds] (1.7595852400324419,-1.6671485405390423) node {$\frac{e}{2}$};
\draw [fill=black] (-0.5,0.) circle (0.5pt);
\draw[color=black] (-0.7385075042330866,0.06612058023365855) node {$x_3'$};
\draw[color=dtsfsf] (-0.14750743153408485,-0.058039098904787076) node {$\frac{d}{2}$};
\draw [fill=black] (-1.,-1.) circle (0.5pt);
\draw[color=black] (-1.2450789951179453,-0.9470224015360578) node {$x_3''$};
\draw[color=rvwvcq] (-0.5448184047771112,-0.5546778154585695) node {$\frac{c}{2}$};
\draw [fill=black] (-0.5,-1.5) circle (0.5pt);
\draw[color=black] (-0.7683058272263136,-1.5727871843938235) node {$w_3'$};
\draw[color=dtsfsf] (-0.1673729801962362,-1.7565435095187232) node {$\frac{d}{2}$};
\draw [fill=black] (-1.,-2.5) circle (0.5pt);
\draw[color=black] (-1.2450789951179453,-2.675325135143221) node {$w_3''$};
\draw[color=rvwvcq] (-0.6044150507635653,-2.0346611907888414) node {$\frac{c}{2}$};
\draw [fill=uuuuuu] (0.,2.5) circle (0.5pt);
\draw[color=uuuuuu] (-0.023347752395639104,2.7678351982862353) node {$v_3$};
\draw[color=wrwrwr] (-4.64705420351136,3.845541213207943) node {$Curve$};
\draw[color=black] (4.049089723345381,-0.5993752999484099) node {$\varphi_2$};
\end{scriptsize}
\end{tikzpicture}
    \begin{center}
        \textnormal{Figure 19.2. The model $(G_2',l_2')$ of $\Gamma_2'$}
    \end{center}
\end{figure}
\begin{figure}[H]
    \centering
    \begin{tikzpicture}[line cap=round,line join=round,>=triangle 45,x=1.0cm,y=1.0cm]
\definecolor{wrwrwr}{rgb}{0.3803921568627451,0.3803921568627451,0.3803921568627451}
\definecolor{uuuuuu}{rgb}{0.26666666666666666,0.26666666666666666,0.26666666666666666}
\definecolor{cqcqcq}{rgb}{0.7529411764705882,0.7529411764705882,0.7529411764705882}
\definecolor{wvvxds}{rgb}{0.396078431372549,0.3411764705882353,0.8235294117647058}
\definecolor{sexdts}{rgb}{0.1803921568627451,0.49019607843137253,0.19607843137254902}
\definecolor{rvwvcq}{rgb}{0.08235294117647059,0.396078431372549,0.7529411764705882}
\definecolor{qqqqff}{rgb}{0.,0.,1.}
\definecolor{dtsfsf}{rgb}{0.8274509803921568,0.1843137254901961,0.1843137254901961}
\clip(-2.,-3.) rectangle (2.5,-1.2);
\draw [shift={(0.75,1.5)},color=dtsfsf]  plot[domain=2.214297435588181:4.068887871591405,variable=\t]({1.*1.25*cos(\t r)+0.*1.25*sin(\t r)},{0.*1.25*cos(\t r)+1.*1.25*sin(\t r)});
\draw [shift={(0.,1.5)},color=rvwvcq]  plot[domain=1.5707963267948966:4.71238898038469,variable=\t]({1.*1.*cos(\t r)+0.*1.*sin(\t r)},{0.*1.*cos(\t r)+1.*1.*sin(\t r)});
\draw [shift={(-0.7464501113912964,1.5063728582396283)},color=sexdts]  plot[domain=-0.9287282318863292:0.9287282318863297,variable=\t]({1.*1.2464664028803656*cos(\t r)+0.*1.2464664028803656*sin(\t r)},{0.*1.2464664028803656*cos(\t r)+1.*1.2464664028803656*sin(\t r)});
\draw (0.,2.5046166785127313)-- (1.,2.5);
\draw (0.,0.5081290379665255)-- (1.,0.5);
\draw (1.,0.5)-- (0.,-0.5);
\draw [color=wvvxds] (1.,2.5)-- (2.,2.5);
\draw [color=wvvxds] (1.5,0.5) circle (0.5cm);
\draw [color=sexdts] (0.,-0.5)-- (0.5,-1.);
\draw [color=sexdts] (0.,-2.)-- (0.5,-2.5);
\draw (0.,-2.)-- (1.,-2.);
\draw [color=wvvxds] (1.,-2.)-- (2.,-2.);
\draw [color=dtsfsf] (0.,-0.5)-- (-0.5,0.);
\draw [color=rvwvcq] (0.,-0.5)-- (-1.,-1.);
\draw [color=dtsfsf] (0.,-2.)-- (-0.5,-1.5);
\draw [color=rvwvcq] (0.,-2.)-- (-1.,-2.5);
\draw [line width=0.4pt,dash pattern=on 1pt off 1pt,color=cqcqcq] (-0.5,1.5)-- (-0.5,0.);
\draw [line width=0.4pt,dash pattern=on 1pt off 1pt,color=cqcqcq] (-1.,1.5)-- (-1.,-1.);
\draw (3.3786274559977745,1.0891963363344506) node[anchor=north west] {$\Gamma_2'$};
\draw (3.3885602303288502,-1.761509896684261) node[anchor=north west] {$T_2$};
\draw [->] (3.567350168288212,-0.12260213205677878) -- (3.577282942619288,-1.105946790833268);
\begin{scriptsize}
\draw [fill=black] (0.,-0.5) circle (0.5pt);
\draw[color=black] (0.2448371545434037,-0.49011478230657785) node {$x_1$};
\draw [fill=black] (0.,0.5) circle (0.5pt);
\draw[color=black] (0.15544218556372277,0.3343054871727011) node {$v_1$};
\draw [color=dtsfsf] (-0.5,1.5)-- ++(-1.0pt,-1.0pt) -- ++(2.0pt,2.0pt) ++(-2.0pt,0) -- ++(2.0pt,-2.0pt);
\draw[color=dtsfsf] (-0.7881713758884649,1.4765745352464008) node {$v_3'$};
\draw[color=dtsfsf] (-0.6044150507635653,1.9781796389657211) node {$d$};
\draw [color=qqqqff] (-1.,1.5)-- ++(-1.0pt,-1.0pt)-- ++(2.0pt,2.0pt) ++(-2.0pt,0) -- ++(2.0pt,-2.0pt);
\draw[color=qqqqff] (-1.314608415435475,1.4865073095774766) node {$v_3''$};
\draw[color=rvwvcq] (-1.0315243469998185,2.0179107362900237) node {$c$};
\draw [color=sexdts] (0.5,1.5)-- ++(-1.0pt,-1.0pt) -- ++(2.0pt,2.0pt) ++(-2.0pt,0) -- ++(2.0pt,-2.0pt);
\draw[color=sexdts] (0.18524050855694976,1.4368434379220982) node {$v_2'$};
\draw[color=sexdts] (0.7265767096005732,1.4517425994187116) node {$b$};
\draw [fill=black] (1.,2.5) circle (0.5pt);
\draw[color=black] (0.9798624550430025,2.7579024239551595) node {$x_4$};
\draw[color=black] (0.7067111609384219,2.156969576925083) node {$f$};
\draw [fill=black] (1.,0.5) circle (0.5pt);
\draw[color=black] (1.1884507159955915,0.4634315534766845) node {$v_4$};
\draw[color=black] (0.7067111609384219,0.8160450422298701) node {$f$};
\draw[color=black] (0.6967783866073463,-0.25669458552630003) node {$f$};
\draw [fill=black] (2.,2.5) circle (0.5pt);
\draw[color=black] (2.2611903437517626,2.5592469373336466) node {$x_4'$};
\draw[color=wvvxds] (1.7496524657013663,2.1271712539318557) node {$\frac{e}{2}$};
\draw [color=rvwvcq] (2.,0.5)-- ++(-1.0pt,-1.0pt) -- ++(2.0pt,2.0pt) ++(-2.0pt,0) -- ++(2.0pt,-2.0pt);
\draw[color=rvwvcq] (2.300921441076065,0.5826248454495924) node {$v_4'$};
\draw[color=wvvxds] (1.7099213683770635,0.4088012946557685) node {$e$};
\draw [fill=black] (0.5,-1.) circle (0.5pt);
\draw[color=black] (0.7315430967661111,-1.0066190475225116) node {$x_2'$};
\draw[color=sexdts] (1.0742238111882214,-0.514946718134267) node {$\frac{b}{2}$};
\draw [fill=black] (0.,-2.) circle (0.5pt);
\draw[color=black] (-0.033280526726714776,-2.238283064575892) node {$w_1$};
\draw [fill=black] (0.5,-2.5) circle (0.5pt);
\draw[color=black] (0.6918119994418085,-2.645526812149994) node {$w_2'$};
\draw[color=sexdts] (0.2696690903710929,-2.5710310046669265) node {$\frac{b}{2}$};
\draw [fill=black] (1.,-2.) circle (0.5pt);
\draw[color=black] (0.9897952293740782,-2.158820869927287) node {$w_4$};
\draw[color=black] (0.6967783866073463,-1.7664762838497987) node {$f$};
\draw [fill=black] (2.,-2.) circle (0.5pt);
\draw[color=black] (2.0426693084680982,-2.168753644258363) node {$w_4'$};
\draw[color=wvvxds] (1.7595852400324419,-1.6671485405390423) node {$\frac{e}{2}$};
\draw [fill=black] (-0.5,0.) circle (0.5pt);
\draw[color=black] (-0.7385075042330866,0.06612058023365855) node {$x_3'$};
\draw[color=dtsfsf] (-0.14750743153408485,-0.058039098904787076) node {$\frac{d}{2}$};
\draw [fill=black] (-1.,-1.) circle (0.5pt);
\draw[color=black] (-1.2450789951179453,-0.9470224015360578) node {$x_3''$};
\draw[color=rvwvcq] (-0.5448184047771112,-0.5546778154585695) node {$\frac{c}{2}$};
\draw [fill=black] (-0.5,-1.5) circle (0.5pt);
\draw[color=black] (-0.7683058272263136,-1.5727871843938235) node {$w_3'$};
\draw[color=dtsfsf] (-0.1673729801962362,-1.55788802289721) node {$\frac{d}{2}$};
\draw [fill=black] (-1.,-2.5) circle (0.5pt);
\draw[color=black] (-1.2450789951179453,-2.675325135143221) node {$w_3''$};
\draw[color=rvwvcq] (-0.6044150507635653,-2.0346611907888414) node {$\frac{c}{2}$};
\draw [fill=uuuuuu] (0.,2.5) circle (0.5pt);
\draw[color=uuuuuu] (-0.023347752395639104,2.7678351982862353) node {$v_3$};
\draw[color=wrwrwr] (-4.64705420351136,3.845541213207943) node {$Curve$};
\draw[color=black] (4.049089723345381,-0.5993752999484099) node {$\varphi_2$};
\end{scriptsize}
\end{tikzpicture}
    \begin{center}
        \textnormal{Figure 19.3. The model $(T_2',l_2')$ of $T_2$}
    \end{center}
\end{figure}
\begin{figure}[H]
    \centering
    \begin{tikzpicture}[line cap=round,line join=round,>=triangle 45,x=1.0cm,y=1.0cm] \definecolor{wrwrwr}{rgb}{0.3803921568627451,0.3803921568627451,0.3803921568627451}
\definecolor{uuuuuu}{rgb}{0.26666666666666666,0.26666666666666666,0.26666666666666666}
\definecolor{cqcqcq}{rgb}{0.7529411764705882,0.7529411764705882,0.7529411764705882}
\definecolor{wvvxds}{rgb}{0.396078431372549,0.3411764705882353,0.8235294117647058}
\definecolor{sexdts}{rgb}{0.1803921568627451,0.49019607843137253,0.19607843137254902}
\definecolor{rvwvcq}{rgb}{0.08235294117647059,0.396078431372549,0.7529411764705882}
\definecolor{qqqqff}{rgb}{0.,0.,1.}
\definecolor{dtsfsf}{rgb}{0.8274509803921568,0.1843137254901961,0.1843137254901961}
\clip(-1.5,-3.) rectangle (4.5,3.);
\draw [shift={(0.75,1.5)},color=dtsfsf]  plot[domain=2.214297435588181:4.068887871591405,variable=\t]({1.*1.25*cos(\t r)+0.*1.25*sin(\t r)},{0.*1.25*cos(\t r)+1.*1.25*sin(\t r)});
\draw [shift={(0.,1.5)},color=rvwvcq]  plot[domain=1.5707963267948966:4.71238898038469,variable=\t]({1.*1.*cos(\t r)+0.*1.*sin(\t r)},{0.*1.*cos(\t r)+1.*1.*sin(\t r)});
\draw [shift={(-0.7464501113912964,1.5063728582396283)},color=sexdts]  plot[domain=-0.9287282318863292:0.9287282318863297,variable=\t]({1.*1.2464664028803656*cos(\t r)+0.*1.2464664028803656*sin(\t r)},{0.*1.2464664028803656*cos(\t r)+1.*1.2464664028803656*sin(\t r)});
\draw (0.,2.5046166785127313)-- (1.,2.5);
\draw (0.,0.5081290379665255)-- (1.,0.5);
\draw (1.,0.5)-- (0.,-0.5);
\draw [color=wvvxds] (1.,2.5)-- (2.,2.5);
\draw [color=wvvxds] (1.5,0.5) circle (0.5cm);
\draw [color=sexdts] (0.,-0.5)-- (0.5,-1.);
\draw [color=sexdts] (0.,-2.)-- (0.5,-2.5);
\draw (0.,-2.)-- (1.,-2.);
\draw [color=wvvxds] (1.,-2.)-- (2.,-2.);
\draw [color=dtsfsf] (0.,-0.5)-- (-0.5,0.);
\draw [color=rvwvcq] (0.,-0.5)-- (-1.,-1.);
\draw [color=dtsfsf] (0.,-2.)-- (-0.5,-1.5);
\draw [color=rvwvcq] (0.,-2.)-- (-1.,-2.5);
\draw [line width=0.4pt,dash pattern=on 2pt off 2pt,color=cqcqcq] (0.,2.5)-- (0.,-2.);
\draw [line width=0.4pt,dash pattern=on 2pt off 2pt,color=cqcqcq] (-0.5,1.5)-- (-0.5,0.);
\draw [line width=0.4pt,dash pattern=on 2pt off 2pt,color=cqcqcq] (-0.5,0.)-- (-0.5,-1.5);
\draw [line width=0.4pt,dash pattern=on 2pt off 2pt,color=cqcqcq] (-1.,1.5)-- (-1.,-1.);
\draw [line width=0.4pt,dash pattern=on 2pt off 2pt,color=cqcqcq] (-1.,-1.)-- (-1.,-2.5);
\draw [line width=0.4pt,dash pattern=on 2pt off 2pt,color=cqcqcq] (0.,2.5046166785127313)-- (0.,-2.);
\draw [line width=0.4pt,dash pattern=on 2pt off 2pt,color=cqcqcq] (0.5,1.5)-- (0.5,-2.5);
\draw [line width=0.4pt,dash pattern=on 2pt off 2pt,color=cqcqcq] (1.,2.5)-- (1.,-2.);
\draw [line width=0.4pt,dash pattern=on 2pt off 2pt,color=cqcqcq] (2.,2.5)-- (2.,-2.);
\draw (3.3786274559977745,1.0891963363344506) node[anchor=north west] {$\Gamma_2'$};
\draw (3.3885602303288502,-1.761509896684261) node[anchor=north west] {$T_2$};
\draw [->] (3.567350168288212,-0.12260213205677878) -- (3.577282942619288,-1.105946790833268);
\begin{scriptsize}
\draw [fill=black] (0.,-0.5) circle (0.5pt);
\draw [fill=black] (0.,0.5) circle (0.5pt);
\draw[color=black] (0.15544218556372277,0.3343054871727011) node {$v_1$};
\draw [color=dtsfsf] (-0.5,1.5)-- ++(-1.0pt,-1.0pt) -- ++(2.0pt,2.0pt) ++(-2.0pt,0) -- ++(2.0pt,-2.0pt);
\draw[color=dtsfsf] (-0.6044150507635653,1.9781796389657211) node {$d$};
\draw [color=qqqqff] (-1.,1.5)-- ++(-1.0pt,-1.0pt) -- ++(2.0pt,2.0pt) ++(-2.0pt,0) -- ++(2.0pt,-2.0pt);
\draw[color=rvwvcq] (-1.0315243469998185,2.0179107362900237) node {$c$};
\draw [color=sexdts] (0.5,1.5)-- ++(-1.0pt,-1.0pt) -- ++(2.0pt,2.0pt) ++(-2.0pt,0) -- ++(2.0pt,-2.0pt);
\draw[color=sexdts] (0.7265767096005732,1.4517425994187116) node {$b$};
\draw [fill=black] (1.,2.5) circle (0.5pt);
\draw[color=black] (0.7067111609384219,2.156969576925083) node {$f$};
\draw [fill=black] (1.,0.5) circle (0.5pt);
\draw[color=black] (1.1884507159955915,0.4634315534766845) node {$v_4$};
\draw[color=black] (0.7067111609384219,0.8160450422298701) node {$f$};
\draw[color=black] (0.6967783866073463,-0.25669458552630003) node {$f$};
\draw [fill=black] (2.,2.5) circle (0.5pt);
\draw[color=wvvxds] (1.7496524657013663,2.1271712539318557) node {$\frac{e}{2}$};
\draw [color=rvwvcq] (2.,0.5)-- ++(-1.0pt,-1.0pt) -- ++(2.0pt,2.0pt) ++(-2.0pt,0) -- ++(2.0pt,-2.0pt);
\draw[color=wvvxds] (1.7099213683770635,0.4088012946557685) node {$e$};
\draw [fill=black] (0.5,-1.) circle (0.5pt);
\draw[color=sexdts] (0.2796018647021685,-1.0115854346880495) node {$\frac{b}{2}$};
\draw [fill=black] (0.,-2.) circle (0.5pt);
\draw [fill=black] (0.5,-2.5) circle (0.5pt);
\draw[color=sexdts] (0.2696690903710929,-2.5710310046669265) node {$\frac{b}{2}$};
\draw [fill=black] (1.,-2.) circle (0.5pt);
\draw[color=black] (0.6967783866073463,-1.7664762838497987) node {$f$};
\draw [fill=black] (2.,-2.) circle (0.5pt);
\draw[color=wvvxds] (1.7595852400324419,-1.6671485405390423) node {$\frac{e}{2}$};
\draw [fill=black] (-0.5,0.) circle (0.5pt);
\draw[color=dtsfsf] (-0.14750743153408485,-0.058039098904787076) node {$\frac{d}{2}$};
\draw [fill=black] (-1.,-1.) circle (0.5pt);
\draw[color=rvwvcq] (-0.5448184047771112,-0.5546778154585695) node {$\frac{c}{2}$};
\draw [fill=black] (-0.5,-1.5) circle (0.5pt);
\draw[color=dtsfsf] (-0.1673729801962362,-1.55788802289721) node {$\frac{d}{2}$};
\draw [fill=black] (-1.,-2.5) circle (0.5pt);
\draw[color=rvwvcq] (-0.6044150507635653,-2.0346611907888414) node {$\frac{c}{2}$};
\draw [fill=uuuuuu] (0.,2.5) circle (0.5pt);
\draw[color=uuuuuu] (-0.023347752395639104,2.7678351982862353) node {$v_3$};
\draw[color=wrwrwr] (-4.64705420351136,3.845541213207943) node {$Curve$};
\draw[color=black] (4.049089723345381,-0.5993752999484099) node {$\varphi_2$};
\end{scriptsize}
\end{tikzpicture}
    \begin{center}
        \textnormal{Figure 20. The tropical morphism $\varphi_2:\Gamma_2' \to T_2$}
    \end{center}
\end{figure}
\textbf{Case 3.} If the metric graph $\Gamma$ has  $2$ bridges, then $\Gamma$ is one of the metric graphs given in Figure 21 or 23. \par 
\textit{Solution of Case 3.} \par 
\textbf{Case 3.1.} 
 Consider the metric graph $\Gamma$ with essential model $(G,l)$ in Figure 21,  where $a,b,c,d, e$, and $f$ are real positive numbers such that $b>a$. Note that if $b=a$, then $\Gamma$ is a hyperelliptic metric graph. The model $(G_1,l_1)$ that obtained by subdividing $(G,l)$ is shown in Figure 21.1. The tropical modification $\Gamma'$, the metric tree $T$ with model $(G',l')$, $(T',t')$ is given in Figure 21.2, 21.3,  respectively. The construction of the tropical morphism $\varphi: \Gamma' \to T$ of degree $3$ is depicted in Figure 22.
 \begin{figure}[H]
    \centering
    \begin{tikzpicture}[line cap=round,line join=round,>=triangle 45,x=1.0cm,y=1.0cm]
\definecolor{wrwrwr}{rgb}{0.3803921568627451,0.3803921568627451,0.3803921568627451}
\clip(-2.,-1.) rectangle (2.8,1.);
\draw (1.,0.)-- (1.5,0.);
\draw(2.,0.) circle (0.5cm);
\draw (-0.5,0.)-- (0.,0.);
\draw (0.,0.)-- (1.,0.);
\draw [shift={(0.5,0.)}] plot[domain=0.:3.141592653589793,variable=\t]({1.*0.5*cos(\t r)+0.*0.5*sin(\t r)},{0.*0.5*cos(\t r)+1.*0.5*sin(\t r)});
\draw(-1.,0.) circle (0.5cm);
\begin{scriptsize}
\draw [fill=black] (0.,0.) circle (0.5pt);
\draw[color=black] (-0.005769492737813506,-0.18312667562322926) node {$v_1$};
\draw [fill=black] (1.,0.) circle (0.5pt);
\draw[color=black] (0.9993671532817349,-0.18312667562322926) node {$v_2$};
\draw [fill=black] (1.5,0.) circle (0.5pt);
\draw[color=black] (1.7710527073225495,0.07626342657536228) node {$v_4$};
\draw[color=black] (1.2295758689829863,0.17029233862235174) node {$d$};
\draw[color=black] (1.9818071653589062,0.624225017469887) node {$e$};
\draw[color=wrwrwr] (-4.866091532683952,1.4348190868404855) node {$Curve$};
\draw [fill=black] (-0.5,0.) circle (0.5pt);
\draw[color=black] (-0.7450312840038039,0.06329392146543271) node {$v_3$};
\draw[color=black] (-0.2489477135489946,0.15083808095745738) node {$c$};
\draw[color=black] (0.49031407771699587,0.11841431818263343) node {$a$};
\draw[color=black] (0.5032835828269255,0.656648780244711) node {$b$};
\draw[color=black] (-1.0141485150348444,0.6825877904645701) node {$f$};
\end{scriptsize}
\end{tikzpicture}
    \begin{center}
        \textnormal{Figure 21. The essential model $(G,l)$ of $\Gamma$}
    \end{center}
\end{figure}
\begin{figure}[H]
    \centering
    \begin{tikzpicture}[line cap=round,line join=round,>=triangle 45,x=1.0cm,y=1.0cm]
\definecolor{wrwrwr}{rgb}{0.3803921568627451,0.3803921568627451,0.3803921568627451}
\definecolor{ttzzqq}{rgb}{0.2,0.6,0.}
\definecolor{ffqqqq}{rgb}{1.,0.,0.}
\definecolor{qqqqff}{rgb}{0.,0.,1.}
\clip(-3.,0.) rectangle (3.,2.5);
\draw [color=qqqqff] (-1.5,1.)-- (-0.5,1.);
\draw [color=ffqqqq] (-0.5,1.)-- (0.5,1.);
\draw (0.5,1.)-- (1.5,1.);
\draw(-2.,1.) circle (0.5cm);
\draw [color=qqqqff] (2.,1.) circle (0.5cm);
\draw [color=ffqqqq] (-0.5,1.)-- (0.5,2.);
\draw [color=ffqqqq] (-0.5,-1.5)-- (0.5,-1.5);
\draw [color=ttzzqq] (0.5,-1.5)-- (1.,-2.);
\draw (0.5,-1.5)-- (1.5,-1.5);
\draw [color=qqqqff] (1.5,-1.5)-- (2.5,-1.5);
\draw [color=qqqqff] (-0.5,-1.5)-- (-1.5,-1.5);
\draw (-1.5,-1.5)-- (-2.5,-1.5);
\draw [shift={(0.5,1.5)},color=ttzzqq]  plot[domain=-1.5707963267948966:1.5707963267948966,variable=\t]({1.*0.5*cos(\t r)+0.*0.5*sin(\t r)},{0.*0.5*cos(\t r)+1.*0.5*sin(\t r)});
\draw (3.376379223364904,1.16416072678588) node[anchor=north west] {$\Gamma'$};
\draw (3.3219149764113602,-1.2322661391700531) node[anchor=north west] {$T$};
\draw [->] (3.562465400456176,0.34719702248271994) -- (3.562465400456176,-0.7420879165881593);
\begin{scriptsize}
\draw [fill=black] (-1.5,1.) circle (0.5pt);
\draw[color=black] (-1.303007327393744,0.8192204960801017) node {$v_3$};
\draw [fill=black] (-0.5,1.) circle (0.5pt);
\draw[color=black] (-0.4951209975828425,0.8101431215878443) node {$v_1$};
\draw[color=qqqqff] (-0.9807605329186092,1.150544665047494) node {$c$};
\draw [fill=black] (0.5,1.) circle (0.5pt);
\draw[color=black] (0.48523544758094794,0.8101431215878443) node {$v_2$};
\draw[color=ffqqqq] (-0.018558836739333284,0.7783723108649437) node {$a$};
\draw [fill=black] (1.5,1.) circle (0.5pt);
\draw[color=black] (1.7560678764969726,0.9553811134639616) node {$v_4$};
\draw[color=black] (0.9436428594399425,0.7965270598494584) node {$d$};
\draw [color=black] (-2.5,1.)-- ++(-1.0pt,-1.0pt) -- ++(2.0pt,2.0pt) ++(-2.0pt,0) -- ++(2.0pt,-2.0pt);
\draw[color=black] (-2.7463098716626577,1.1006191053400787) node {$v_3'$};
\draw[color=black] (-2.188051340388833,0.9780745496946048) node {$f$};
\draw [color=qqqqff] (2.5,1.)-- ++(-1.0pt,-1.0pt) -- ++(2.0pt,2.0pt) ++(-2.0pt,0) -- ++(2.0pt,-2.0pt);
\draw[color=qqqqff] (2.78181119412205,1.0189227349097627) node {$v_4'$};
\draw[color=qqqqff] (2.296171658786283,0.9417650517255756) node {$e$};
\draw [fill=black] (0.5,2.) circle (0.5pt);
\draw[color=black] (0.5396996945344918,2.271600414841273) node {$v_2''$};
\draw[color=ffqqqq] (-0.009481462247075967,1.7314966325519625) node {$a$};
\draw [fill=black] (-0.5,-1.5) circle (0.5pt);
\draw[color=black] (-0.4951209975828425,-1.7042896127674338) node {$w_1$};
\draw [fill=black] (0.5,-1.5) circle (0.5pt);
\draw[color=black] (0.43984857511966136,-1.7496764852287205) node {$w_2$};
\draw[color=ffqqqq] (-0.04579096021610525,-1.264036949892954) node {$a$};
\draw [fill=black] (1.,-2.) circle (0.5pt);
\draw[color=black] (0.9754136701628432,-2.2761642057796454) node {$w_2'$};
\draw[color=ttzzqq] (1.4610532054986098,-1.9266852878277385) node {$\frac{b-a}{2}$};
\draw [fill=black] (1.5,-1.5) circle (0.5pt);
\draw[color=black] (1.483746641729253,-1.3048851351081117) node {$w_4$};
\draw[color=black] (0.9436428594399425,-1.2458822009084394) node {$d$};
\draw [fill=black] (2.5,-1.5) circle (0.5pt);
\draw[color=black] (2.8725849390446228,-1.4410457524919715) node {$w_4'$};
\draw[color=qqqqff] (1.9149219301114757,-1.264036949892954) node {$\frac{e}{2}$};
\draw [fill=black] (-1.5,-1.5) circle (0.5pt);
\draw[color=black] (-1.5117869407156623,-1.7405991107364631) node {$w_3$};
\draw[color=qqqqff] (-0.9898379074108665,-1.2549595754006966) node {$2c$};
\draw [fill=black] (-2.5,-1.5) circle (0.5pt);
\draw[color=black] (-2.7190777481858857,-1.7133669872596913) node {$w_3'$};
\draw[color=black] (-2.0700454719894874,-1.236804826416182) node {$\frac{f}{2}$};
\draw [fill=black] (1.,-0.5) circle (0.5pt);
\draw[color=black] (1.2477349049305626,-0.5968499247120407) node {$x_2'$};
\draw[color=ttzzqq] (0.5805478797496498,1.5045622702455295) node {$b-a$};
\draw [color=ttzzqq] (1.,1.5)-- ++(-1.0pt,-1.0pt) -- ++(2.0pt,2.0pt) ++(-2.0pt,0) -- ++(2.0pt,-2.0pt);
\draw[color=ttzzqq] (1.2749670284073347,1.4546367105381144) node {$v_2'$};
\draw[color=black] (3.9845633143461447,-0.2110615087911044) node {$\varphi$};
\draw[color=wrwrwr] (-5.265281293264064,2.929710065529929) node {$Curve$};
\end{scriptsize}
\end{tikzpicture}
    \begin{center}
        \textnormal{Figure 21.1. The model $(G_1,l_1)$ of $\Gamma$}
    \end{center}
\end{figure}
\begin{figure}[H]
    \centering
    \begin{tikzpicture}[line cap=round,line join=round,>=triangle 45,x=1.0cm,y=1.0cm]
\definecolor{wrwrwr}{rgb}{0.3803921568627451,0.3803921568627451,0.3803921568627451}
\definecolor{ttzzqq}{rgb}{0.2,0.6,0.}
\definecolor{ffqqqq}{rgb}{1.,0.,0.}
\definecolor{qqqqff}{rgb}{0.,0.,1.}
\clip(-3.,-1.) rectangle (3.,2.8);
\draw [color=qqqqff] (-1.5,1.)-- (-0.5,1.);
\draw [color=ffqqqq] (-0.5,1.)-- (0.5,1.);
\draw (0.5,1.)-- (1.5,1.);
\draw(-2.,1.) circle (0.5cm);
\draw [color=qqqqff] (2.,1.) circle (0.5cm);
\draw [color=ffqqqq] (-0.5,1.)-- (0.5,2.);
\draw (0.5,2.)-- (1.5,2.);
\draw [color=qqqqff] (1.5,2.)-- (2.5,2.);
\draw (1.5,1.)-- (0.5,0.);
\draw [color=ffqqqq] (0.5,0.)-- (-0.5,0.);
\draw [color=qqqqff] (-0.5,0.)-- (-1.5,0.);
\draw (-1.5,0.)-- (-2.5,0.);
\draw [color=ffqqqq] (-0.5,-1.5)-- (0.5,-1.5);
\draw [color=ttzzqq] (0.5,-1.5)-- (1.,-2.);
\draw (0.5,-1.5)-- (1.5,-1.5);
\draw [color=qqqqff] (1.5,-1.5)-- (2.5,-1.5);
\draw [color=qqqqff] (-0.5,-1.5)-- (-1.5,-1.5);
\draw (-1.5,-1.5)-- (-2.5,-1.5);
\draw [color=ttzzqq] (0.5,0.)-- (1.,-0.5);
\draw [shift={(0.5,1.5)},color=ttzzqq]  plot[domain=-1.5707963267948966:1.5707963267948966,variable=\t]({1.*0.5*cos(\t r)+0.*0.5*sin(\t r)},{0.*0.5*cos(\t r)+1.*0.5*sin(\t r)});
\draw (3.376379223364904,1.16416072678588) node[anchor=north west] {$\Gamma'$};
\draw (3.3219149764113602,-1.2322661391700531) node[anchor=north west] {$T$};
\draw [->] (3.562465400456176,0.34719702248271994) -- (3.562465400456176,-0.7420879165881593);
\begin{scriptsize}
\draw [fill=black] (-1.5,1.) circle (0.5pt);
\draw[color=black] (-1.303007327393744,0.8192204960801017) node {$v_3$};
\draw [fill=black] (-0.5,1.) circle (0.5pt);
\draw[color=black] (-0.4951209975828425,0.8101431215878443) node {$v_1$};
\draw[color=qqqqff] (-0.9807605329186092,1.150544665047494) node {$c$};
\draw [fill=black] (0.5,1.) circle (0.5pt);
\draw[color=black] (0.48523544758094794,0.8101431215878443) node {$v_2$};
\draw[color=ffqqqq] (-0.018558836739333284,0.7783723108649437) node {$a$};
\draw [fill=black] (1.5,1.) circle (0.5pt);
\draw[color=black] (1.7560678764969726,0.9553811134639616) node {$v_4$};
\draw[color=black] (0.9436428594399425,0.7965270598494584) node {$d$};
\draw [color=black] (-2.5,1.)-- ++(-1.0pt,-1.0pt) -- ++(2.0pt,2.0pt) ++(-2.0pt,0) -- ++(2.0pt,-2.0pt);
\draw[color=black] (-2.7463098716626577,1.1006191053400787) node {$v_3'$};
\draw[color=black] (-2.188051340388833,0.9780745496946048) node {$f$};
\draw [color=qqqqff] (2.5,1.)-- ++(-1.0pt,-1.0pt) -- ++(2.0pt,2.0pt) ++(-2.0pt,0) -- ++(2.0pt,-2.0pt);
\draw[color=qqqqff] (2.78181119412205,1.0189227349097627) node {$v_4'$};
\draw[color=qqqqff] (2.296171658786283,0.9417650517255756) node {$e$};
\draw [fill=black] (0.5,2.) circle (0.5pt);
\draw[color=black] (0.6032413159802931,2.3714515342561038) node {$v_2''$};
\draw[color=ffqqqq] (-0.009481462247075967,1.7314966325519625) node {$a$};
\draw [fill=black] (1.5,2.) circle (0.5pt);
\draw[color=black] (1.4928240162215105,1.8086543157361499) node {$x_4$};
\draw[color=black] (1.2431962176844342,2.2307522296261153) node {$d$};
\draw [fill=black] (2.5,2.) circle (0.5pt);
\draw[color=black] (2.500412584862073,1.8086543157361499) node {$x_4'$};
\draw[color=qqqqff] (1.8876898066347039,2.3487580980254603) node {$\frac{e}{2}$};
\draw [fill=black] (0.5,0.) circle (0.5pt);
\draw[color=black] (0.3399974557048308,0.3018101500214343) node {$x_2$};
\draw[color=black] (0.9254881104554279,0.10664659843790189) node {$d$};
\draw [fill=black] (-0.5,0.) circle (0.5pt);
\draw[color=black] (-0.4860436230905852,-0.1792906980682038) node {$x_1$};
\draw[color=ffqqqq] (-0.009481462247075967,0.17926559437596046) node {$a$};
\draw [fill=black] (-1.5,0.) circle (0.5pt);
\draw[color=black] (-1.4845548172388903,-0.1792906980682038) node {$x_3$};
\draw[color=qqqqff] (-0.8990641624882934,0.2246524668372471) node {$2c$};
\draw [fill=black] (-2.5,0.) circle (0.5pt);
\draw[color=black] (-2.7100003736936285,-0.17021332357594643) node {$x_3'$};
\draw[color=black] (-2.1154323444507743,-0.24737100676013368) node {$\frac{f}{2}$};
\draw [fill=black] (-0.5,-1.5) circle (0.5pt);
\draw[color=black] (-0.4951209975828425,-1.7042896127674338) node {$w_1$};
\draw [fill=black] (0.5,-1.5) circle (0.5pt);
\draw[color=black] (0.43984857511966136,-1.7496764852287205) node {$w_2$};
\draw[color=ffqqqq] (-0.04579096021610525,-1.264036949892954) node {$a$};
\draw [fill=black] (1.,-2.) circle (0.5pt);
\draw[color=black] (0.9754136701628432,-2.2761642057796454) node {$w_2'$};
\draw[color=ttzzqq] (1.4610532054986098,-1.9266852878277385) node {$\frac{b-a}{2}$};
\draw [fill=black] (1.5,-1.5) circle (0.5pt);
\draw[color=black] (1.483746641729253,-1.3048851351081117) node {$w_4$};
\draw[color=black] (0.9436428594399425,-1.2458822009084394) node {$d$};
\draw [fill=black] (2.5,-1.5) circle (0.5pt);
\draw[color=black] (2.8725849390446228,-1.4410457524919715) node {$w_4'$};
\draw[color=qqqqff] (1.9149219301114757,-1.264036949892954) node {$\frac{e}{2}$};
\draw [fill=black] (-1.5,-1.5) circle (0.5pt);
\draw[color=black] (-1.5117869407156623,-1.7405991107364631) node {$w_3$};
\draw[color=qqqqff] (-0.9898379074108665,-1.2549595754006966) node {$2c$};
\draw [fill=black] (-2.5,-1.5) circle (0.5pt);
\draw[color=black] (-2.7190777481858857,-1.7133669872596913) node {$w_3'$};
\draw[color=black] (-2.0700454719894874,-1.236804826416182) node {$\frac{f}{2}$};
\draw [fill=black] (1.,-0.5) circle (0.5pt);
\draw[color=black] (1.2477349049305626,-0.5968499247120407) node {$x_2'$};
\draw[color=ttzzqq] (0.5805478797496497,-0.5287696160201106) node {$\frac{b-a}{2}$};
\draw[color=ttzzqq] (0.5805478797496498,1.5045622702455295) node {$b-a$};
\draw [color=ttzzqq] (1.,1.5)-- ++(-1.0pt,-1.0pt) -- ++(2.0pt,2.0pt) ++(-2.0pt,0) -- ++(2.0pt,-2.0pt);
\draw[color=ttzzqq] (1.2205027814537908,1.472791459522629) node {$v_2'$};
\draw[color=black] (3.9845633143461447,-0.2110615087911044) node {$\varphi$};
\draw[color=wrwrwr] (-5.265281293264064,2.929710065529929) node {$Curve$};
\end{scriptsize}
\end{tikzpicture}
    \begin{center}
        \textnormal{Figure 21.2. The model $(G',l')$ of $\Gamma'$}
    \end{center}
\end{figure}
\begin{figure}[H]
    \centering
    \begin{tikzpicture}[line cap=round,line join=round,>=triangle 45,x=1.0cm,y=1.0cm]
\definecolor{wrwrwr}{rgb}{0.3803921568627451,0.3803921568627451,0.3803921568627451}
\definecolor{ttzzqq}{rgb}{0.2,0.6,0.}
\definecolor{ffqqqq}{rgb}{1.,0.,0.}
\definecolor{qqqqff}{rgb}{0.,0.,1.}
\clip(-3.,-2.5) rectangle (3.,-0.9);
\draw [color=qqqqff] (-1.5,1.)-- (-0.5,1.);
\draw [color=ffqqqq] (-0.5,1.)-- (0.5,1.);
\draw (0.5,1.)-- (1.5,1.);
\draw(-2.,1.) circle (0.5cm);
\draw [color=qqqqff] (2.,1.) circle (0.5cm);
\draw [color=ffqqqq] (-0.5,1.)-- (0.5,2.);
\draw (0.5,2.)-- (1.5,2.);
\draw [color=qqqqff] (1.5,2.)-- (2.5,2.);
\draw (1.5,1.)-- (0.5,0.);
\draw [color=ffqqqq] (0.5,0.)-- (-0.5,0.);
\draw [color=qqqqff] (-0.5,0.)-- (-1.5,0.);
\draw (-1.5,0.)-- (-2.5,0.);
\draw [color=ffqqqq] (-0.5,-1.5)-- (0.5,-1.5);
\draw [color=ttzzqq] (0.5,-1.5)-- (1.,-2.);
\draw (0.5,-1.5)-- (1.5,-1.5);
\draw [color=qqqqff] (1.5,-1.5)-- (2.5,-1.5);
\draw [color=qqqqff] (-0.5,-1.5)-- (-1.5,-1.5);
\draw (-1.5,-1.5)-- (-2.5,-1.5);
\draw [color=ttzzqq] (0.5,0.)-- (1.,-0.5);
\draw [shift={(0.5,1.5)},color=ttzzqq]  plot[domain=-1.5707963267948966:1.5707963267948966,variable=\t]({1.*0.5*cos(\t r)+0.*0.5*sin(\t r)},{0.*0.5*cos(\t r)+1.*0.5*sin(\t r)});
\draw (3.376379223364904,1.16416072678588) node[anchor=north west] {$\Gamma'$};
\draw (3.3219149764113602,-1.2322661391700531) node[anchor=north west] {$T$};
\draw [->] (3.562465400456176,0.34719702248271994) -- (3.562465400456176,-0.7420879165881593);
\begin{scriptsize}
\draw [fill=black] (-1.5,1.) circle (0.5pt);
\draw[color=black] (-1.303007327393744,0.7556788746343004) node {$v_3$};
\draw [fill=black] (-0.5,1.) circle (0.5pt);
\draw[color=black] (-0.3226508822299535,0.7647562491265578) node {$v_1$};
\draw[color=qqqqff] (-0.9807605329186092,1.150544665047494) node {$c$};
\draw [fill=black] (0.5,1.) circle (0.5pt);
\draw[color=black] (0.6304734394570651,0.7829109981110725) node {$v_2$};
\draw[color=ffqqqq] (-0.009481462247075967,0.6785211914501132) node {$a$};
\draw [fill=black] (1.5,1.) circle (0.5pt);
\draw[color=black] (1.7560678764969726,0.9553811134639616) node {$v_4$};
\draw[color=black] (1.2431962176844342,1.1959315375087805) node {$d$};
\draw [color=black] (-2.5,1.)-- ++(-1.0pt,-1.0pt) -- ++(2.0pt,2.0pt) ++(-2.0pt,0) -- ++(2.0pt,-2.0pt);
\draw[color=black] (-2.7463098716626577,1.1006191053400787) node {$v_3'$};
\draw[color=black] (-2.188051340388833,0.9780745496946048) node {$f$};
\draw [color=qqqqff] (2.5,1.)-- ++(-1.0pt,-1.0pt) -- ++(2.0pt,2.0pt) ++(-2.0pt,0) -- ++(2.0pt,-2.0pt);
\draw[color=qqqqff] (2.78181119412205,1.0189227349097627) node {$v_4'$};
\draw[color=qqqqff] (2.296171658786283,0.9417650517255756) node {$e$};
\draw [fill=black] (0.5,2.) circle (0.5pt);
\draw[color=black] (0.6032413159802931,2.3714515342561038) node {$v_2''$};
\draw[color=ffqqqq] (-0.009481462247075967,1.7314966325519625) node {$a$};
\draw [fill=black] (1.5,2.) circle (0.5pt);
\draw[color=black] (1.76514525098923,1.8903506861664658) node {$x_4$};
\draw[color=black] (1.2431962176844342,2.2307522296261153) node {$d$};
\draw [fill=black] (2.5,2.) circle (0.5pt);
\draw[color=black] (2.736424321660763,1.8903506861664658) node {$x_4'$};
\draw[color=qqqqff] (1.8876898066347039,2.3487580980254603) node {$\frac{e}{2}$};
\draw [fill=black] (0.5,0.) circle (0.5pt);
\draw[color=black] (0.3399974557048308,0.3018101500214343) node {$x_2$};
\draw[color=black] (1.2159640942076622,0.3517357097288496) node {$d$};
\draw [fill=black] (-0.5,0.) circle (0.5pt);
\draw[color=black] (-0.3317282567222108,-0.1792906980682038) node {$x_1$};
\draw[color=ffqqqq] (-0.009481462247075967,0.17926559437596046) node {$a$};
\draw [fill=black] (-1.5,0.) circle (0.5pt);
\draw[color=black] (-1.366548948839545,-0.14298120009917448) node {$x_3$};
\draw[color=qqqqff] (-0.8990641624882934,0.2246524668372471) node {$2c$};
\draw [fill=black] (-2.5,0.) circle (0.5pt);
\draw[color=black] (-2.7100003736936285,-0.17021332357594643) node {$x_3'$};
\draw[color=black] (-1.8612658586675692,-0.23829363226787637) node {$\frac{f}{2}$};
\draw [fill=black] (-0.5,-1.5) circle (0.5pt);
\draw[color=black] (-0.4951209975828425,-1.7042896127674338) node {$w_1$};
\draw [fill=black] (0.5,-1.5) circle (0.5pt);
\draw[color=black] (0.43984857511966136,-1.7496764852287205) node {$w_2$};
\draw[color=ffqqqq] (-0.04579096021610525,-1.264036949892954) node {$a$};
\draw [fill=black] (1.,-2.) circle (0.5pt);
\draw[color=black] (0.9754136701628432,-2.2761642057796454) node {$w_2'$};
\draw[color=ttzzqq] (1.4610532054986098,-1.9266852878277385) node {$\frac{b-a}{2}$};
\draw [fill=black] (1.5,-1.5) circle (0.5pt);
\draw[color=black] (1.483746641729253,-1.3048851351081117) node {$w_4$};
\draw[color=black] (0.9436428594399425,-1.2458822009084394) node {$d$};
\draw [fill=black] (2.5,-1.5) circle (0.5pt);
\draw[color=black] (2.8725849390446228,-1.4410457524919715) node {$w_4'$};
\draw[color=qqqqff] (1.9149219301114757,-1.264036949892954) node {$\frac{e}{2}$};
\draw [fill=black] (-1.5,-1.5) circle (0.5pt);
\draw[color=black] (-1.5117869407156623,-1.7405991107364631) node {$w_3$};
\draw[color=qqqqff] (-0.9898379074108665,-1.2549595754006966) node {$2c$};
\draw [fill=black] (-2.5,-1.5) circle (0.5pt);
\draw[color=black] (-2.7190777481858857,-1.7133669872596913) node {$w_3'$};
\draw[color=black] (-2.0700454719894874,-1.236804826416182) node {$\frac{f}{2}$};
\draw [fill=black] (1.,-0.5) circle (0.5pt);
\draw[color=black] (1.2477349049305626,-0.45161193283592344) node {$x_2'$};
\draw[color=ttzzqq] (0.5805478797496497,-0.5287696160201106) node {$\frac{b-a}{2}$};
\draw[color=ttzzqq] (0.5805478797496498,1.4591753977842428) node {$b-a$};
\draw [color=ttzzqq] (1.,1.5)-- ++(-1.0pt,-1.0pt) -- ++(2.0pt,2.0pt) ++(-2.0pt,0) -- ++(2.0pt,-2.0pt);
\draw[color=ttzzqq] (1.2205027814537908,1.608952076906489) node {$v_2'$};
\draw[color=black] (3.9845633143461447,-0.2110615087911044) node {$\varphi$};
\draw[color=wrwrwr] (-5.265281293264064,2.929710065529929) node {$Curve$};
\end{scriptsize}
\end{tikzpicture}
    \begin{center}
        \textnormal{Figure 21.3. The model $(T',t')$ of $T$}
    \end{center}
\end{figure}
\begin{figure}[H]
    \centering
    \begin{tikzpicture}[line cap=round,line join=round,>=triangle 45,x=1.0cm,y=1.0cm]
\definecolor{cqcqcq}{rgb}{0.7529411764705882,0.7529411764705882,0.7529411764705882}
\definecolor{ttzzqq}{rgb}{0.2,0.6,0.}
\definecolor{ffqqqq}{rgb}{1.,0.,0.}
\definecolor{qqqqff}{rgb}{0.,0.,1.}
\clip(-3.,-2.5) rectangle (4.,2.5);
\draw [color=qqqqff] (-1.5,1.)-- (-0.5,1.);
\draw [color=ffqqqq] (-0.5,1.)-- (0.5,1.);
\draw (0.5,1.)-- (1.5,1.);
\draw(-2.,1.) circle (0.5cm);
\draw [color=qqqqff] (2.,1.) circle (0.5cm);
\draw [color=ffqqqq] (-0.5,1.)-- (0.5,2.);
\draw (0.5,2.)-- (1.5,2.);
\draw [color=qqqqff] (1.5,2.)-- (2.5,2.);
\draw (1.5,1.)-- (0.5,0.);
\draw [color=ffqqqq] (0.5,0.)-- (-0.5,0.);
\draw [color=qqqqff] (-0.5,0.)-- (-1.5,0.);
\draw (-1.5,0.)-- (-2.5,0.);
\draw [color=ffqqqq] (-0.5,-1.5)-- (0.5,-1.5);
\draw [color=ttzzqq] (0.5,-1.5)-- (1.,-2.);
\draw (0.5,-1.5)-- (1.5,-1.5);
\draw [color=qqqqff] (1.5,-1.5)-- (2.5,-1.5);
\draw [color=qqqqff] (-0.5,-1.5)-- (-1.5,-1.5);
\draw (-1.5,-1.5)-- (-2.5,-1.5);
\draw [color=ttzzqq] (0.5,0.)-- (1.,-0.5);
\draw [shift={(0.5,1.5)},color=ttzzqq]  plot[domain=-1.5707963267948966:1.5707963267948966,variable=\t]({1.*0.5*cos(\t r)+0.*0.5*sin(\t r)},{0.*0.5*cos(\t r)+1.*0.5*sin(\t r)});
\draw [line width=0.4pt,dash pattern=on 2pt off 2pt,color=cqcqcq] (-2.5,1.)-- (-2.5,-1.5);
\draw [line width=0.4pt,dash pattern=on 2pt off 2pt,color=cqcqcq] (-1.5,1.)-- (-1.5,-1.5);
\draw [line width=0.4pt,dash pattern=on 2pt off 2pt,color=cqcqcq] (-0.5,1.)-- (-0.5,-1.5);
\draw [line width=0.4pt,dash pattern=on 2pt off 2pt,color=cqcqcq] (0.5,2.)-- (0.5,-1.5);
\draw [line width=0.4pt,dash pattern=on 2pt off 2pt,color=cqcqcq] (1.,1.5)-- (1.,-2.);
\draw [line width=0.4pt,dash pattern=on 2pt off 2pt,color=cqcqcq] (1.5,2.)-- (1.5,-1.5);
\draw [line width=0.4pt,dash pattern=on 2pt off 2pt,color=cqcqcq] (2.5,2.)-- (2.5,-1.5);
\draw (3.376379223364904,1.16416072678588) node[anchor=north west] {$\Gamma'$};
\draw (3.3219149764113602,-1.2322661391700531) node[anchor=north west] {$T$};
\draw [->] (3.562465400456176,0.34719702248271994) -- (3.562465400456176,-0.7420879165881593);
\begin{scriptsize}
\draw [fill=black] (-1.5,1.) circle (0.5pt);
\draw[color=black] (-1.303007327393744,0.7556788746343004) node {$v_3$};
\draw [fill=black] (-0.5,1.) circle (0.5pt);
\draw[color=black] (-0.3226508822299535,0.7647562491265578) node {$v_1$};
\draw[color=qqqqff] (-0.9807605329186092,1.150544665047494) node {$c$};
\draw [fill=black] (0.5,1.) circle (0.5pt);
\draw[color=black] (0.6304734394570651,0.7829109981110725) node {$v_2$};
\draw[color=ffqqqq] (-0.009481462247075967,0.6785211914501132) node {$a$};
\draw [fill=black] (1.5,1.) circle (0.5pt);
\draw[color=black] (1.7560678764969726,0.9553811134639616) node {$v_4$};
\draw[color=black] (1.2431962176844342,1.1959315375087805) node {$d$};
\draw [color=black] (-2.5,1.)-- ++(-1.0pt,-1.0pt) -- ++(2.0pt,2.0pt) ++(-2.0pt,0) -- ++(2.0pt,-2.0pt);
\draw[color=black] (-2.188051340388833,0.9780745496946048) node {$f$};
\draw [color=qqqqff] (2.5,1.)-- ++(-1.0pt,-1.0pt) -- ++(2.0pt,2.0pt) ++(-2.0pt,0) -- ++(2.0pt,-2.0pt);
\draw[color=qqqqff] (2.296171658786283,0.9417650517255756) node {$e$};
\draw [fill=black] (0.5,2.) circle (0.5pt);
\draw[color=ffqqqq] (-0.009481462247075967,1.7314966325519625) node {$a$};
\draw [fill=black] (1.5,2.) circle (0.5pt);
\draw[color=black] (1.2431962176844342,2.2307522296261153) node {$d$};
\draw [fill=black] (2.5,2.) circle (0.5pt);
\draw[color=qqqqff] (1.8876898066347039,2.3487580980254603) node {$\frac{e}{2}$};
\draw [fill=black] (0.5,0.) circle (0.5pt);
\draw[color=black] (1.2159640942076622,0.3517357097288496) node {$d$};
\draw [fill=black] (-0.5,0.) circle (0.5pt);
\draw[color=ffqqqq] (-0.009481462247075967,0.17926559437596046) node {$a$};
\draw [fill=black] (-1.5,0.) circle (0.5pt);
\draw[color=qqqqff] (-0.8990641624882934,0.2246524668372471) node {$2c$};
\draw [fill=black] (-2.5,0.) circle (0.5pt);
\draw[color=black] (-1.8612658586675692,-0.23829363226787637) node {$\frac{f}{2}$};
\draw [fill=black] (-0.5,-1.5) circle (0.5pt);
\draw [fill=black] (0.5,-1.5) circle (0.5pt);
\draw[color=ffqqqq] (-0.04579096021610525,-1.264036949892954) node {$a$};
\draw [fill=black] (1.,-2.) circle (0.5pt);
\draw[color=ttzzqq] (1.4610532054986098,-1.9266852878277385) node {$\frac{b-a}{2}$};
\draw [fill=black] (1.5,-1.5) circle (0.5pt);
\draw[color=black] (1.2250414686999196,-1.2912690733697259) node {$d$};
\draw [fill=black] (2.5,-1.5) circle (0.5pt);
\draw[color=qqqqff] (1.9149219301114757,-1.1460310814936088) node {$\frac{e}{2}$};
\draw [fill=black] (-1.5,-1.5) circle (0.5pt);
\draw[color=qqqqff] (-0.9898379074108665,-1.2549595754006966) node {$2c$};
\draw [fill=black] (-2.5,-1.5) circle (0.5pt);
\draw[color=black] (-1.8431111096830544,-1.7360604234903347) node {$\frac{f}{2}$};
\draw [fill=black] (1.,-0.5) circle (0.5pt);
\draw[color=ttzzqq] (0.5805478797496497,-0.5287696160201106) node {$\frac{b-a}{2}$};
\draw[color=ttzzqq] (0.5805478797496498,1.4591753977842428) node {$b-a$};
\draw [color=ttzzqq] (1.,1.5)-- ++(-1.0pt,-1.0pt) -- ++(2.0pt,2.0pt) ++(-2.0pt,0) -- ++(2.0pt,-2.0pt);
\draw[color=black] (3.9300990673926006,-0.12028776386853117) node {$\varphi$};
\end{scriptsize}
\end{tikzpicture}

    \begin{center}
        \textnormal{Figure 22. The tropical morphism $\varphi:\Gamma' \to T$}
    \end{center}
\end{figure}
\textbf{Case 3.2.} 
 Consider the metric graph $\Gamma_1$ with essential model $(G,l)$ in Figure 23, where $a,b,c,d,e$, and $f$ are real positive numbers. The model $(G_1,l_1)$ that is obtained by subdividing $(G,l)$ is shown in Figure 23.1. The tropical modification $\Gamma_1'$, the metric tree $T_1$ with model $(G_1',l_1')$, $(T_1',t_1')$ is given in Figure 23.2, 23.3, respectively. The construction of the tropical morphism $\varphi_1: \Gamma_1' \to T_1$ of degree $3$ is depicted in Figure 24. 
\begin{figure}[H]
    \centering
    \begin{tikzpicture}[line cap=round,line join=round,>=triangle 45,x=1.0cm,y=1.0cm]
\clip(-1.5,-1.5) rectangle (3.,2.5);
\draw (1.,0.)-- (1.7494319444761144,0.);
\draw (1.7494319444761144,0.)-- (2.,1.);
\draw(2.,1.5) circle (0.5cm);
\draw(2.,-0.5) circle (0.559271267319376cm);
\draw(0.,0.) circle (1.cm);
\begin{scriptsize}
\draw [fill=black] (1.,0.) circle (0.5pt);
\draw[color=black] (0.7282116035578747,0.038220811472343394) node {$v_3$};
\draw [fill=black] (1.7494319444761144,0.) circle (0.5pt);
\draw[color=black] (1.8390126761356087,-0.24843752983803968) node {$v_1$};
\draw [fill=black] (2.,1.) circle (0.5pt);
\draw[color=black] (2.0450483589524464,1.2386026157095724) node {$v_4$};
\draw[color=black] (1.3597557617573122,0.1770709455445602) node {$c$};
\draw[color=black] (2.0405693223694716,0.4458131405230443) node {$d$};
\draw[color=black] (2.013695102871623,2.183679334717242) node {$e$};
\draw[color=black] (2.380976102675551,-0.5485329808973469) node {$b$};
\draw[color=black] (0.02500286003084167,1.2341235791265979) node {$f$};
\end{scriptsize}
\end{tikzpicture}
    \begin{center}
        \textnormal{Figure 23. The essential model $(G,l)$ of $\Gamma_1$ }
    \end{center}
\end{figure}
\begin{figure}[H]
    \centering
    
\begin{tikzpicture}[line cap=round,line join=round,>=triangle 45,x=1.0cm,y=1.0cm]
\definecolor{wrwrwr}{rgb}{0.3803921568627451,0.3803921568627451,0.3803921568627451}
\definecolor{ttzzqq}{rgb}{0.2,0.6,0.}
\definecolor{qqttzz}{rgb}{0.,0.2,0.6}
\definecolor{qqqqff}{rgb}{0.,0.,1.}
\clip(-3.,-0.75) rectangle (3.,3.);
\draw [color=qqqqff] (-1.5,2.)-- (0.5,2.);
\draw (0.5,2.)-- (1.5,0.);
\draw(-2.,2.) circle (0.5cm);
\draw [color=qqqqff] (2.,0.) circle (0.5cm);
\draw (0.5,-1.5)-- (1.5,-1.5);
\draw [color=qqqqff] (1.5,-1.5)-- (2.5,-1.5);
\draw [color=qqqqff] (0.5,-1.5)-- (-1.5,-1.5);
\draw (-1.5,-1.5)-- (-2.5,-1.5);
\draw [color=ttzzqq] (0.8,2.) circle (0.3cm);
\draw [->] (3.905535460266646,-0.04060990164981755) -- (3.9172221450831706,-0.8586778388064921);
\draw (3.671801763936166,0.9644449925712402) node[anchor=north west] {$\Gamma'_1$};
\draw (3.6484283943031186,-1.3728919707335436) node[anchor=north west] {$T_1$};
\draw (4.0107156236153605,-0.07566995609938866) node[anchor=north west] {$\varphi_1$};
\begin{scriptsize}
\draw [fill=black] (-1.5,2.) circle (0.5pt);
\draw[color=black] (-1.3242559951278232,1.852633038627058) node {$v_3$};
\draw [fill=black] (0.5,2.) circle (0.5pt);
\draw[color=qqqqff] (-0.46528466111331274,2.1973902407145136) node {$c$};
\draw [fill=black] (0.5,2.) circle (0.5pt);
\draw[color=black] (0.2651331399194343,1.7941996145444383) node {$v_1$};
\draw [fill=black] (1.5,0.) circle (0.5pt);
\draw[color=black] (1.772715481251024,0.08794363133194621) node {$v_4$};
\draw[color=black] (1.2059112676496124,0.9469149653464543) node {$d$};
\draw [color=black] (-2.5,2.)-- ++(-1.0pt,-1.0pt) -- ++(2.0pt,2.0pt) ++(-2.0pt,0) -- ++(2.0pt,-2.0pt);
\draw[color=black] (-2.796778282009841,2.121426789407108) node {$v_3'$};
\draw[color=black] (-2.1832273291423334,1.9402831747509872) node {$f$};
\draw [color=qqttzz] (2.5,0.)-- ++(-1.0pt,-1.0pt) -- ++(2.0pt,2.0pt) ++(-2.0pt,0) -- ++(2.0pt,-2.0pt);
\draw[color=qqttzz] (2.8128304299216556,0.07625694651542228) node {$v_4'$};
\draw[color=qqqqff] (2.29277295558634,-0.10488666814069846) node {$e$};
\draw [fill=black] (0.5,-1.5) circle (0.5pt);
\draw [fill=black] (0.5,-1.5) circle (0.5pt);
\draw[color=black] (0.37031330326814993,-1.7234925152292613) node {$w_1$};
\draw [fill=black] (1.5,-1.5) circle (0.5pt);
\draw[color=black] (1.503921730470973,-1.6767457759631657) node {$w_4$};
\draw[color=black] (1.0072376257687052,-1.2852418346096144) node {$d$};
\draw [fill=qqttzz] (2.5,-1.5) circle (0.5pt);
\draw[color=qqttzz] (2.777770375472084,-1.4313253948161633) node {$w_4'$};
\draw[color=qqqqff] (2.035665889622813,-1.8228293361697145) node {$\frac{e}{2}$};
\draw [fill=black] (-1.5,-1.5) circle (0.5pt);
\draw[color=black] (-1.5579896914583022,-1.6884324607796894) node {$w_3$};
\draw[color=qqqqff] (-0.45359797629678866,-1.2618684649765666) node {$2c$};
\draw [fill=black] (-2.5,-1.5) circle (0.5pt);
\draw[color=black] (-2.7149714882941733,-1.7118058304127373) node {$w_3'$};
\draw[color=black] (-2.066360480977094,-1.9046361298853822) node {$\frac{f}{2}$};
\draw[color=ttzzqq] (0.8669974079704179,1.9519698595675112) node {$b$};
\draw [color=ttzzqq] (1.0999508514070155,2.0054301694454697)-- ++(-1.0pt,-1.0pt) -- ++(2.0pt,2.0pt) ++(-2.0pt,0) -- ++(2.0pt,-2.0pt);
\draw[color=ttzzqq] (1.3519948278561618,2.133113474223632) node {$v_2'$};
\draw[color=wrwrwr] (-5.081525163640274,4.254246768422724) node {$Curve$};
\end{scriptsize}
\end{tikzpicture}
    \begin{center}
        \textnormal{Figure 23.1. The model $(G_1,l_1)$ of $\Gamma_1$ }
    \end{center}
\end{figure}
\begin{figure}[H]
    \centering
    
\begin{tikzpicture}[line cap=round,line join=round,>=triangle 45,x=1.0cm,y=1.0cm]
\definecolor{wrwrwr}{rgb}{0.3803921568627451,0.3803921568627451,0.3803921568627451}
\definecolor{ttzzqq}{rgb}{0.2,0.6,0.}
\definecolor{qqttzz}{rgb}{0.,0.2,0.6}
\definecolor{qqqqff}{rgb}{0.,0.,1.}
\clip(-3.,-1.) rectangle (3.5,3.);
\draw [color=qqqqff] (-1.5,2.)-- (0.5,2.);
\draw (0.5,2.)-- (1.5,0.);
\draw(-2.,2.) circle (0.5cm);
\draw [color=qqqqff] (2.,0.) circle (0.5cm);
\draw [color=qqqqff] (1.5,1.)-- (2.5,1.);
\draw (1.5,0.)-- (0.5,0.);
\draw [color=qqqqff] (0.5047101786581755,0.)-- (-1.5192111463852906,0.);
\draw (-1.5192111463852906,0.)-- (-2.4850528667591707,0.);
\draw (0.5,-1.5)-- (1.5,-1.5);
\draw [color=qqqqff] (1.5,-1.5)-- (2.5,-1.5);
\draw [color=qqqqff] (0.5,-1.5)-- (-1.5,-1.5);
\draw (-1.5,-1.5)-- (-2.5,-1.5);
\draw (0.5,2.)-- (1.5,1.);
\draw [color=ttzzqq] (0.8,2.) circle (0.3cm);
\draw [color=ttzzqq] (0.5047101786581755,0.)-- (1.0999508514070155,-0.5349682960084712);
\draw [color=ttzzqq] (0.5,-1.5)-- (1.0999508514070155,-2.084143929596917);
\begin{scriptsize}
\draw [fill=black] (-1.5,2.) circle (0.5pt);
\draw[color=black] (-1.3242559951278234,1.8526330386270577) node {$v_3$};
\draw [fill=black] (0.5,2.) circle (0.5pt);
\draw[color=qqqqff] (-0.4652846611133127,2.1973902407145136) node {$c$};
\draw [fill=black] (0.5,2.) circle (0.5pt);
\draw[color=black] (0.2651331399194344,1.794199614544438) node {$v_1$};
\draw [fill=black] (1.5,0.) circle (0.5pt);
\draw[color=black] (1.7727154812510242,0.08794363133194556) node {$v_4$};
\draw[color=black] (1.2059112676496127,0.9469149653464539) node {$d$};
\draw [color=black] (-2.5,2.)-- ++(-1.0pt,-1.0pt) -- ++(2.0pt,2.0pt) ++(-2.0pt,0) -- ++(2.0pt,-2.0pt);
\draw[color=black] (-2.7967782820098415,2.121426789407108) node {$v_3'$};
\draw[color=black] (-2.183227329142334,1.9402831747509872) node {$f$};
\draw [color=qqttzz] (2.5,0.)-- ++(-1.0pt,-1.0pt) -- ++(2.0pt,2.0pt) ++(-2.0pt,0) -- ++(2.0pt,-2.0pt);
\draw[color=qqttzz] (2.812830429921656,0.07625694651542164) node {$v_4'$};
\draw[color=qqqqff] (2.2927729555863405,-0.10488666814069914) node {$e$};
\draw [fill=black] (1.5,1.) circle (0.5pt);
\draw[color=black] (1.5272951001040214,0.8242047747729526) node {$x_4$};
\draw [fill=qqttzz] (2.5,1.) circle (0.5pt);
\draw[color=qqttzz] (2.812830429921656,1.0929985255530028) node {$x_4'$};
\draw[color=qqqqff] (1.9889191503567174,0.6898078993829276) node {$\frac{e}{2}$};
\draw [fill=black] (0.5,0.) circle (0.5pt);
\draw[color=black] (1.0306109954017533,0.2457138763550185) node {$d$};
\draw [fill=black] (0.5047101786581755,0.) circle (0.5pt);
\draw[color=black] (0.39368667290119785,-0.22759685871420032) node {$x_1$};
\draw [fill=black] (-1.5192111463852906,0.) circle (0.5pt);
\draw[color=black] (-1.5346163218252546,-0.21591017389767642) node {$x_3$};
\draw[color=qqqqff] (-0.3717911825811211,0.29246061562111425) node {$2c$};
\draw [fill=black] (-2.4850528667591707,0.) circle (0.5pt);
\draw[color=black] (-2.750031542743746,-0.1925368042646286) node {$x_3'$};
\draw[color=black] (-2.0546737961605706,-0.43211384300336897) node {$\frac{f}{2}$};
\draw [fill=black] (0.5,-1.5) circle (0.5pt);
\draw [fill=black] (0.5,-1.5) circle (0.5pt);
\draw[color=black] (0.37031330326815,-1.7234925152292624) node {$w_1$};
\draw [fill=black] (1.5,-1.5) circle (0.5pt);
\draw[color=black] (1.5039217304709733,-1.6767457759631668) node {$w_4$};
\draw[color=black] (1.0072376257687055,-1.2852418346096155) node {$d$};
\draw [fill=qqttzz] (2.5,-1.5) circle (0.5pt);
\draw[color=qqttzz] (2.7777703754720844,-1.4313253948161644) node {$w_4'$};
\draw[color=qqqqff] (2.0356658896228135,-1.8228293361697159) node {$\frac{e}{2}$};
\draw [fill=black] (-1.5,-1.5) circle (0.5pt);
\draw[color=black] (-1.5579896914583025,-1.6884324607796908) node {$w_3$};
\draw[color=qqqqff] (-0.4535979762967888,-1.2618684649765675) node {$2c$};
\draw [fill=black] (-2.5,-1.5) circle (0.5pt);
\draw[color=black] (-2.714971488294174,-1.7118058304127386) node {$w_3'$};
\draw[color=black] (-2.0663604809770946,-1.9046361298853831) node {$\frac{f}{2}$};
\draw[color=black] (1.2760313765487563,1.53124920617265) node {$d$};
\draw[color=ttzzqq] (0.866997407970418,1.951969859567511) node {$b$};
\draw [color=ttzzqq] (1.0999508514070155,2.0054301694454697)-- ++(-1.0pt,-1.0pt) -- ++(2.0pt,2.0pt) ++(-2.0pt,0) -- ++(2.0pt,-2.0pt);
\draw[color=ttzzqq] (1.3519948278561622,2.133113474223632) node {$v_2'$};
\draw [fill=black] (1.0999508514070155,-0.5349682960084712) circle (0.5pt);
\draw[color=black] (1.3403081430396382,-0.5548240335768702) node {$x_2'$};
\draw[color=ttzzqq] (0.890370777603466,-0.724280963416467) node {$\frac{b}{2}$};
\draw [fill=black] (1.0999508514070155,-2.084143929596917) circle (0.5pt);
\draw[color=black] (1.2117546100578747,-2.307826756055458) node {$w_2'$};
\draw[color=ttzzqq] (0.8319373535208461,-2.196803250298481) node {$\frac{b}{2}$};
\draw[color=wrwrwr] (-4.427070813914933,3.658225842780004) node {$Curve$};
\end{scriptsize}
\end{tikzpicture}
    \begin{center}
        \textnormal{Figure 23.2. The model $(G_1',l_1')$ of $\Gamma_1'$}
    \end{center}
\end{figure}
\begin{figure}[H]
    \centering
    \begin{tikzpicture}[line cap=round,line join=round,>=triangle 45,x=1.0cm,y=1.0cm]
\definecolor{wrwrwr}{rgb}{0.3803921568627451,0.3803921568627451,0.3803921568627451}
\definecolor{ttzzqq}{rgb}{0.2,0.6,0.}
\definecolor{qqttzz}{rgb}{0.,0.2,0.6}
\definecolor{qqqqff}{rgb}{0.,0.,1.}
\clip(-3.,-3.) rectangle (3.5,-1.);
\draw [color=qqqqff] (-1.5,2.)-- (0.5,2.);
\draw (0.5,2.)-- (1.5,0.);
\draw(-2.,2.) circle (0.5cm);
\draw [color=qqqqff] (2.,0.) circle (0.5cm);
\draw [color=qqqqff] (1.5,1.)-- (2.5,1.);
\draw (1.5,0.)-- (0.5,0.);
\draw [color=qqqqff] (0.5047101786581755,0.)-- (-1.5192111463852906,0.);
\draw (-1.5192111463852906,0.)-- (-2.4850528667591707,0.);
\draw (0.5,-1.5)-- (1.5,-1.5);
\draw [color=qqqqff] (1.5,-1.5)-- (2.5,-1.5);
\draw [color=qqqqff] (0.5,-1.5)-- (-1.5,-1.5);
\draw (-1.5,-1.5)-- (-2.5,-1.5);
\draw (0.5,2.)-- (1.5,1.);
\draw [color=ttzzqq] (0.8,2.) circle (0.3cm);
\draw [color=ttzzqq] (0.5047101786581755,0.)-- (1.0999508514070155,-0.5349682960084712);
\draw [color=ttzzqq] (0.5,-1.5)-- (1.0999508514070155,-2.084143929596917);
\begin{scriptsize}
\draw [fill=black] (-1.5,2.) circle (0.5pt);
\draw[color=black] (-1.3593160495773953,1.8993797778931534) node {$v_3$};
\draw [fill=black] (0.5,2.) circle (1.0pt);
\draw[color=qqqqff] (-0.4652846611133127,2.1973902407145136) node {$c$};
\draw [fill=black] (0.5,2.) circle (0.5pt);
\draw[color=black] (0.2885065095524823,1.9110664627096774) node {$v_1$};
\draw [fill=black] (1.5,0.) circle (0.5pt);
\draw[color=black] (1.7142820571684045,0.06457026169889772) node {$v_4$};
\draw[color=black] (1.1942245828330886,0.8534214868142626) node {$d$};
\draw [color=black] (-2.5,2.)-- ++(-1.0pt,-1.0pt) -- ++(2.0pt,2.0pt) ++(-2.0pt,0) -- ++(2.0pt,-2.0pt);
\draw[color=black] (-2.8084649668263655,2.144800159040156) node {$v_3'$};
\draw[color=black] (-2.183227329142334,1.9402831747509872) node {$f$};
\draw [color=qqttzz] (2.5,0.)-- ++(-1.0pt,-1.0pt) -- ++(2.0pt,2.0pt) ++(-2.0pt,0) -- ++(2.0pt,-2.0pt);
\draw[color=qqttzz] (2.871263854004276,0.07625694651542164) node {$v_4'$};
\draw[color=qqqqff] (2.210966161870673,-0.15163340740679485) node {$e$};
\draw [fill=black] (1.5,1.) circle (0.5pt);
\draw[color=black] (1.6792220027188327,0.9293849381216679) node {$x_4$};
\draw [fill=qqttzz] (2.5,1.) circle (0.5pt);
\draw[color=qqttzz] (2.871263854004276,1.0579384711034312) node {$x_4'$};
\draw[color=qqqqff] (2.386266434118532,1.507875836539602) node {$\frac{e}{2}$};
\draw [fill=black] (0.5,0.) circle (1.0pt);
\draw[color=black] (0.9488042016860857,0.2690872459880664) node {$d$};
\draw [fill=black] (0.5047101786581755,0.) circle (0.5pt);
\draw[color=black] (0.32356656400205414,-0.11073001054896112) node {$x_1$};
\draw [fill=black] (-1.5192111463852906,0.) circle (0.5pt);
\draw[color=black] (-1.3359426799443475,-0.13410338018200896) node {$x_3$};
\draw[color=qqqqff] (-0.3717911825811211,0.29246061562111425) node {$2c$};
\draw [fill=black] (-2.4850528667591707,0.) circle (0.5pt);
\draw[color=black] (-2.750031542743746,-0.01723653201676975) node {$x_3'$};
\draw[color=black] (-1.7625066757474719,0.6898078993829276) node {$\frac{f}{2}$};
\draw [fill=black] (0.5,-1.5) circle (0.5pt);
\draw [fill=black] (0.5,-1.5) circle (0.5pt);
\draw[color=black] (0.37031330326815,-1.7234925152292624) node {$w_1$};
\draw [fill=black] (1.5,-1.5) circle (0.5pt);
\draw[color=black] (1.5039217304709733,-1.6767457759631668) node {$w_4$};
\draw[color=black] (1.0072376257687055,-1.2852418346096155) node {$d$};
\draw [fill=qqttzz] (2.5,-1.5) circle (0.5pt);
\draw[color=qqttzz] (2.836203799554704,-1.4313253948161644) node {$w_4'$};
\draw[color=qqqqff] (2.0356658896228135,-1.8228293361697159) node {$\frac{e}{2}$};
\draw [fill=black] (-1.5,-1.5) circle (0.5pt);
\draw[color=black] (-1.5579896914583025,-1.6884324607796908) node {$w_3$};
\draw[color=qqqqff] (-0.4535979762967888,-1.2618684649765675) node {$2c$};
\draw [fill=black] (-2.5,-1.5) circle (0.5pt);
\draw[color=black] (-2.714971488294174,-1.7118058304127386) node {$w_3'$};
\draw[color=black] (-2.0663604809770946,-1.9046361298853831) node {$\frac{f}{2}$};
\draw[color=black] (1.2643446917322323,1.4143823580074106) node {$d$};
\draw[color=ttzzqq] (0.9488042016860857,1.9285964899344632) node {$b$};
\draw [color=ttzzqq] (1.0999508514070155,2.0054301694454697)-- ++(-1.0pt,-1.0pt) -- ++(2.0pt,2.0pt) ++(-2.0pt,0) -- ++(2.0pt,-2.0pt);
\draw[color=ttzzqq] (1.3169347734065904,2.285040376838443) node {$v_2'$};
\draw [fill=black] (1.0999508514070155,-0.5349682960084712) circle (0.5pt);
\draw[color=black] (1.3403081430396382,-0.5548240335768702) node {$x_2'$};
\draw[color=ttzzqq] (0.9254308320530378,-0.6424741697007995) node {$\frac{b}{2}$};
\draw [fill=black] (1.0999508514070155,-2.084143929596917) circle (0.5pt);
\draw[color=black] (1.2117546100578747,-2.307826756055458) node {$w_2'$};
\draw[color=ttzzqq] (0.8319373535208461,-2.196803250298481) node {$\frac{b}{2}$};
\draw[color=wrwrwr] (-4.427070813914933,3.658225842780004) node {$Curve$};
\end{scriptsize}
\end{tikzpicture}
    \begin{center}
        \textnormal{Figure 23.3. The model $(T_1',t_1')$ of $T_1$}
    \end{center}
\end{figure}
\begin{figure}[H]
    \centering
    \begin{tikzpicture}[line cap=round,line join=round,>=triangle 45,x=1.0cm,y=1.0cm]
\definecolor{cqcqcq}{rgb}{0.7529411764705882,0.7529411764705882,0.7529411764705882}
\definecolor{wrwrwr}{rgb}{0.3803921568627451,0.3803921568627451,0.3803921568627451}
\definecolor{ttzzqq}{rgb}{0.2,0.6,0.}
\definecolor{qqttzz}{rgb}{0.,0.2,0.6}
\definecolor{qqqqff}{rgb}{0.,0.,1.}
\clip(-3.5,-3.) rectangle (4.5,3.);
\draw [color=qqqqff] (-1.5,2.)-- (0.5,2.);
\draw (0.5,2.)-- (1.5,0.);
\draw(-2.,2.) circle (0.5cm);
\draw [color=qqqqff] (2.,0.) circle (0.5cm);
\draw [color=qqqqff] (1.5,1.)-- (2.5,1.);
\draw (1.5,0.)-- (0.5,0.);
\draw [color=qqqqff] (0.5047101786581755,0.)-- (-1.5192111463852906,0.);
\draw (-1.5192111463852906,0.)-- (-2.4850528667591707,0.);
\draw (0.5,-1.5)-- (1.5,-1.5);
\draw [color=qqqqff] (1.5,-1.5)-- (2.5,-1.5);
\draw [color=qqqqff] (0.5,-1.5)-- (-1.5,-1.5);
\draw (-1.5,-1.5)-- (-2.5,-1.5);
\draw (0.5,2.)-- (1.5,1.);
\draw [color=ttzzqq] (0.8,2.) circle (0.3cm);
\draw [color=ttzzqq] (0.5047101786581755,0.)-- (1.0999508514070155,-0.5349682960084712);
\draw [color=ttzzqq] (0.5,-1.5)-- (1.0999508514070155,-2.084143929596917);
\draw [line width=0.4pt,dash pattern=on 2pt off 2pt,color=cqcqcq] (-2.5,2.)-- (-2.5,-1.5);
\draw [line width=0.4pt,dash pattern=on 2pt off 2pt,color=cqcqcq] (-1.5,2.)-- (-1.5,-1.5);
\draw [line width=0.4pt,dash pattern=on 2pt off 2pt,color=cqcqcq] (0.5,2.)-- (0.5,-1.5);
\draw [line width=0.4pt,dash pattern=on 2pt off 2pt,color=cqcqcq] (1.0999508514070155,2.0054301694454697)-- (1.0999508514070155,-2.084143929596917);
\draw [line width=0.4pt,dash pattern=on 2pt off 2pt,color=cqcqcq] (1.5,1.)-- (1.5,-1.5);
\draw [line width=0.4pt,dash pattern=on 2pt off 2pt,color=cqcqcq] (2.5,1.)-- (2.5,-1.5);
\draw [->] (3.905535460266646,-0.04060990164981755) -- (3.9172221450831706,-0.8586778388064921);
\draw (3.6718017639361666,0.9644449925712397) node[anchor=north west] {$\Gamma'_1$};
\draw (3.648428394303119,-1.3728919707335447) node[anchor=north west] {$T_1$};
\draw (4.010715623615361,-0.07566995609938935) node[anchor=north west] {$\varphi_1$};
\begin{scriptsize}
\draw [fill=black] (-1.5,2.) circle (0.5pt);
\draw[color=black] (-1.3242559951278234,1.8526330386270577) node {$v_3$};
\draw [fill=black] (0.5,2.) circle (0.5pt);
\draw[color=qqqqff] (-0.4652846611133127,2.1973902407145136) node {$c$};
\draw [fill=black] (0.5,2.) circle (0.5pt);
\draw[color=black] (0.2651331399194344,1.794199614544438) node {$v_1$};
\draw [fill=black] (1.5,0.) circle (0.5pt);
\draw[color=black] (1.7727154812510242,0.08794363133194556) node {$v_4$};
\draw[color=black] (1.2059112676496127,0.9469149653464539) node {$d$};
\draw [color=black] (-2.5,2.)-- ++(-1.0pt,-1.0pt) -- ++(2.0pt,2.0pt) ++(-2.0pt,0) -- ++(2.0pt,-2.0pt);
\draw[color=black] (-2.7967782820098415,2.121426789407108) node {$v_3'$};
\draw[color=black] (-2.183227329142334,1.9402831747509872) node {$f$};
\draw [color=qqttzz] (2.5,0.)-- ++(-1.0pt,-1.0pt) -- ++(2.0pt,2.0pt) ++(-2.0pt,0) -- ++(2.0pt,-2.0pt);
\draw[color=qqttzz] (2.812830429921656,0.07625694651542164) node {$v_4'$};
\draw[color=qqqqff] (2.2927729555863405,-0.10488666814069914) node {$e$};
\draw [fill=black] (1.5,1.) circle (0.5pt);
\draw[color=black] (1.5272951001040214,0.8242047747729526) node {$x_4$};
\draw [fill=qqttzz] (2.5,1.) circle (0.5pt);
\draw[color=qqttzz] (2.812830429921656,1.0929985255530028) node {$x_4'$};
\draw[color=qqqqff] (1.9889191503567174,0.6898078993829276) node {$\frac{e}{2}$};
\draw [fill=black] (0.5,0.) circle (0.5pt);
\draw[color=black] (1.0306109954017533,0.2457138763550185) node {$d$};
\draw [fill=black] (0.5047101786581755,0.) circle (0.5pt);
\draw[color=black] (0.39368667290119785,-0.22759685871420032) node {$x_1$};
\draw [fill=black] (-1.5192111463852906,0.) circle (0.5pt);
\draw[color=black] (-1.5346163218252546,-0.21591017389767642) node {$x_3$};
\draw[color=qqqqff] (-0.3717911825811211,0.29246061562111425) node {$2c$};
\draw [fill=black] (-2.4850528667591707,0.) circle (0.5pt);
\draw[color=black] (-2.750031542743746,-0.1925368042646286) node {$x_3'$};
\draw[color=black] (-2.0546737961605706,-0.43211384300336897) node {$\frac{f}{2}$};
\draw [fill=black] (0.5,-1.5) circle (0.5pt);
\draw [fill=black] (0.5,-1.5) circle (0.5pt);
\draw[color=black] (0.37031330326815,-1.7234925152292624) node {$w_1$};
\draw [fill=black] (1.5,-1.5) circle (0.5pt);
\draw[color=black] (1.5039217304709733,-1.6767457759631668) node {$w_4$};
\draw[color=black] (1.0072376257687055,-1.2852418346096155) node {$d$};
\draw [fill=qqttzz] (2.5,-1.5) circle (0.5pt);
\draw[color=qqttzz] (2.7777703754720844,-1.4313253948161644) node {$w_4'$};
\draw[color=qqqqff] (2.0356658896228135,-1.8228293361697159) node {$\frac{e}{2}$};
\draw [fill=black] (-1.5,-1.5) circle (0.5pt);
\draw[color=black] (-1.5579896914583025,-1.6884324607796908) node {$w_3$};
\draw[color=qqqqff] (-0.4535979762967888,-1.2618684649765675) node {$2c$};
\draw [fill=black] (-2.5,-1.5) circle (0.5pt);
\draw[color=black] (-2.714971488294174,-1.7118058304127386) node {$w_3'$};
\draw[color=black] (-2.0663604809770946,-1.9046361298853831) node {$\frac{f}{2}$};
\draw[color=black] (1.2760313765487563,1.53124920617265) node {$d$};
\draw[color=ttzzqq] (0.866997407970418,1.951969859567511) node {$b$};
\draw [color=ttzzqq] (1.0999508514070155,2.0054301694454697)-- ++(-1.0pt,-1.0pt) -- ++(2.0pt,2.0pt) ++(-2.0pt,0) -- ++(2.0pt,-2.0pt);
\draw[color=ttzzqq] (1.3519948278561622,2.133113474223632) node {$v_2'$};
\draw [fill=black] (1.0999508514070155,-0.5349682960084712) circle (0.5pt);
\draw[color=black] (1.3403081430396382,-0.5548240335768702) node {$x_2'$};
\draw[color=ttzzqq] (0.890370777603466,-0.724280963416467) node {$\frac{b}{2}$};
\draw [fill=black] (1.0999508514070155,-2.084143929596917) circle (0.5pt);
\draw[color=black] (1.2117546100578747,-2.307826756055458) node {$w_2'$};
\draw[color=ttzzqq] (0.8319373535208461,-2.196803250298481) node {$\frac{b}{2}$};
\draw[color=wrwrwr] (-4.427070813914933,3.658225842780004) node {$Curve$};
\end{scriptsize}
\end{tikzpicture}
    \begin{center}
        \textnormal{Figure 24. The tropical morphism $\varphi_1:\Gamma_1' \to T_1$}
    \end{center}
\end{figure}
\textbf{Case 4.} If the metric graph $\Gamma$ has  $3$ bridges, then $\Gamma$ is the metric graph given in Figure 25. \par
\textit{Solution of Case 4.} \par 
Consider the metric graph $\Gamma$ with essential model $(G,l)$ in Figure 25,  where $a,b,c,d,e$, and $f$ are real positive numbers. Note that the metric graph $\Gamma$ is hyperelliptic in the sense of Kawaguchi-Yamaki \cite{KJ} i.e., there is a harmonic morphism from $\Gamma$ to a metric tree, but it is not hyperelliptic in our sense because the harmonic map coming from the unique hyperelliptic involution $\iota$ on $\Gamma$ (see Theorem 3.5, \cite{KJ}) is not a tropical morphism in our sense because it does not satisfy the Riemann-Hurwitz condition. The model $(G_1,l_1)$ that is obtained by subdividing $(G,l)$ is shown in Figure 25.1. The tropical modification $\Gamma'$, the metric tree $T$ with model $(G',l')$, $(T',t')$ is given in Figure 25.2, 25.3, respectively. The construction of the tropical morphism $\varphi: \Gamma' \to T$ of degree $3$ is depicted in Figure 26. This ends our constructive solution of Problem 1.

\begin{figure}[H]
    \centering
    \begin{tikzpicture}[line cap=round,line join=round,>=triangle 45,x=1.0cm,y=1.0cm]
\clip(-2.,-2.) rectangle (3.,1.5);
\draw (1.,0.)-- (-0.5,0.);
\draw (1.,-0.5)-- (1.,0.);
\draw (1.,0.)-- (1.5,0.5);
\draw(2.,0.5) circle (0.5cm);
\draw(1.,-1.) circle (0.5cm);
\draw(-1.,0.) circle (0.5cm);
\begin{scriptsize}
\draw [fill=black] (1.,0.) circle (0.5pt);
\draw[color=black] (0.8787747942569901,0.21080577489341928) node {$v_1$};
\draw [fill=black] (-0.5,0.) circle (0.5pt);
\draw[color=black] (-0.2911421578497567,-0.1568824100544158) node {$v_2$};
\draw[color=black] (0.23114219576932665,0.15648820211703454) node {$a$};
\draw [fill=black] (1.,-0.5) circle (0.5pt);
\draw[color=black] (0.9874099398097594,-0.7334843364498844) node {$v_3$};
\draw[color=black] (1.1252930091651976,-0.2529827311203272) node {$b$};
\draw [fill=black] (1.5,0.5) circle (0.5pt);
\draw[color=black] (1.731142859363334,0.5116415625780116) node {$v_4$};
\draw[color=black] (1.1921454064284402,0.3570453939067628) node {$c$};
\draw[color=black] (2.0110872729031626,1.2344831079868237) node {$f$};
\draw[color=black] (0.364846990295812,-1.0552114982792402) node {$e$};
\draw[color=black] (-0.9889140542848521,0.6745942809071658) node {$d$};
\end{scriptsize}
\end{tikzpicture}
    \begin{center}
        \textnormal{Figure 25. The essential model $(G,l)$ of $\Gamma$ }
    \end{center}
\end{figure}
\begin{figure}[H]
    \centering
    \begin{tikzpicture}[line cap=round,line join=round,>=triangle 45,x=1.0cm,y=1.0cm]
\definecolor{wrwrwr}{rgb}{0.3803921568627451,0.3803921568627451,0.3803921568627451}
\definecolor{qqwuqq}{rgb}{0.,0.39215686274509803,0.}
\definecolor{rvwvcq}{rgb}{0.08235294117647059,0.396078431372549,0.7529411764705882}
\definecolor{dtsfsf}{rgb}{0.8274509803921568,0.1843137254901961,0.1843137254901961}
\clip(-3.,-0.5) rectangle (2.5,3.);
\draw [color=dtsfsf] (-2.,1.5) circle (0.5cm);
\draw (-1.5,1.5)-- (-0.5,1.);
\draw [color=rvwvcq] (-0.5,1.)-- (0.5,2.);
\draw(1.,2.) circle (0.5cm);
\draw (0.4828052989015277,-0.8977665785955626)-- (1.4881393572775417,-0.9143836704695463);
\draw (0.5,-2.)-- (1.5,-2.);
\draw [color=dtsfsf] (-0.5,1.)-- (1.,0.5);
\draw (2.8922836206291644,1.2707639109593092) node[anchor=north west] {$\Gamma'$};
\draw (2.9005921665661565,-1.2550340538862115) node[anchor=north west] {$T$};
\draw [->] (3.056672577778454,0.5664683674959348) -- (3.056672577778454,-0.9335316325040652);
\draw [color=qqwuqq] (1.5,0.5) circle (0.5cm);
\begin{scriptsize}
\draw [color=dtsfsf] (-2.5,1.5)-- ++(-1.0pt,-1.0pt) -- ++(2.0pt,2.0pt) ++(-2.0pt,0) -- ++(2.0pt,-2.0pt);
\draw[color=dtsfsf] (-2.7866075273047657,1.5864886565649992) node {$v_2'$};
\draw[color=dtsfsf] (-2.283940498116759,1.4577061945416254) node {$d$};
\draw [fill=black] (-1.5,1.5) circle (0.5pt);
\draw[color=black] (-1.2993778045832245,1.6695741159349176) node {$v_2$};
\draw [fill=black] (-0.5,1.) circle (0.5pt);
\draw[color=black] (-0.5017573946320066,0.8304109762987415) node {$v_1$};
\draw[color=black] (-1.104126975063916,1.009044713944066) node {$a$};
\draw [fill=black] (-1.5,-1.5) circle (0.5pt);
\draw[color=black] (-1.4821658151970452,-1.662152804798812) node {$w_2$};
\draw [fill=rvwvcq] (-2.5,-1.5) circle (0.5pt);
\draw[color=rvwvcq] (-2.711830613871839,-1.5790673454288935) node {$w_2'$};
\draw [fill=black] (0.5,2.) circle (0.5pt);
\draw[color=black] (0.7445244959167711,2.1016185046584934) node {$v_3$};
\draw[color=rvwvcq] (-0.09048437075091002,1.6238771132814627) node {$b$};
\draw[color=black] (1.3385855304116883,1.922984767013169) node {$e$};
\draw [color=black] (1.5,2.)-- ++(-1.0pt,-1.0pt) -- ++(2.0pt,2.0pt) ++(-2.0pt,0) -- ++(2.0pt,-2.0pt);
\draw[color=black] (1.8246354677257117,2.08500141278451) node {$v_3'$};
\draw [fill=black] (0.4828052989015277,-0.8977665785955626) circle (0.5pt);
\draw[color=black] (0.47865102593303177,-1.122097318894342) node {$x_3$};
\draw [fill=black] (1.4881393572775417,-0.9143836704695463) circle (0.5pt);
\draw[color=black] (1.533836359930997,-1.1304058648313338) node {$x_3'$};
\draw[color=black] (1.0311693307429899,-1.159485775610805) node {$\frac{e}{2}$};
\draw [fill=black] (0.5,-2.) circle (0.5pt);
\draw[color=black] (0.39556556656311326,-2.1440484691443387) node {$w_3$};
\draw [fill=black] (1.5,-2.) circle (0.5pt);
\draw[color=black] (1.525527813994005,-2.1191228313333634) node {$w_3'$};
\draw[color=black] (1.081020606364941,-2.2312882014827533) node {$\frac{e}{2}$};
\draw [fill=black] (1.,0.5) circle (0.5pt);
\draw[color=black] (1.243037252136282,0.5645375063150024) node {$v_4$};
\draw[color=dtsfsf] (0.2750916504767314,0.8927250708261802) node {$c$};
\draw [color=qqwuqq] (2.,0.5)-- ++(-1.0pt,-1.0pt) -- ++(2.0pt,2.0pt) ++(-2.0pt,0) -- ++(2.0pt,-2.0pt);
\draw[color=qqwuqq] (2.2816054942602633,0.5645375063150024) node {$v_4'$};
\draw [fill=dtsfsf] (2.,-1.5) circle (0.5pt);
\draw[color=dtsfsf] (2.173594397079369,-1.5707587994919017) node {$w_4'$};
\draw[color=black] (3.3907963768486757,-0.14584317129780028) node {$\varphi$};
\draw[color=qqwuqq] (1.8204811947572157,0.4606806821026043) node {$f$};
\draw[color=wrwrwr] (-4.253065885183827,3.518225586915603) node {$Curve$};
\end{scriptsize}
\end{tikzpicture}
    \begin{center}
        \textnormal{Figure 25.1. The model $(G_1,l_1)$ of $\Gamma$}
    \end{center}
\end{figure}
\begin{figure}[H]
    \centering
    
\begin{tikzpicture}[line cap=round,line join=round,>=triangle 45,x=1.0cm,y=1.0cm]
\definecolor{wrwrwr}{rgb}{0.3803921568627451,0.3803921568627451,0.3803921568627451}
\definecolor{qqzzqq}{rgb}{0.,0.6,0.}
\definecolor{qqwuqq}{rgb}{0.,0.39215686274509803,0.}
\definecolor{rvwvcq}{rgb}{0.08235294117647059,0.396078431372549,0.7529411764705882}
\definecolor{dtsfsf}{rgb}{0.8274509803921568,0.1843137254901961,0.1843137254901961}
\clip(-3.,-1.4) rectangle (2.5,3.);
\draw [color=dtsfsf] (-2.,1.5) circle (0.5cm);
\draw (-1.5,1.5)-- (-0.5,1.);
\draw (-0.5,1.)-- (-1.5,0.5);
\draw [color=dtsfsf] (-1.5,0.5)-- (-2.5,0.5);
\draw [color=rvwvcq] (-0.5,1.)-- (0.5,2.);
\draw(1.,2.) circle (0.5cm);
\draw [color=rvwvcq] (-0.5,1.)-- (0.4828052989015277,-0.8977665785955626);
\draw (0.4828052989015277,-0.8977665785955626)-- (1.4881393572775417,-0.9143836704695463);
\draw (0.5,-2.)-- (1.5,-2.);
\draw [color=dtsfsf] (-0.5,1.)-- (1.,0.5);
\draw [color=dtsfsf] (-0.5,1.)-- (0.989626601058031,-0.3078598170691431);
\draw [color=qqzzqq] (0.989626601058031,-0.3078598170691431)-- (2.003269205371037,-0.31616836300613493);
\draw (2.8922836206291644,1.2707639109593092) node[anchor=north west] {$\Gamma'$};
\draw (2.9005921665661565,-1.2550340538862115) node[anchor=north west] {$T$};
\draw [->] (3.056672577778454,0.5664683674959348) -- (3.056672577778454,-0.9335316325040652);
\draw [color=qqwuqq] (1.5,0.5) circle (0.5cm);
\begin{scriptsize}
\draw [color=dtsfsf] (-2.5,1.5)-- ++(-1.0pt,-1.0pt) -- ++(2.0pt,2.0pt) ++(-2.0pt,0) -- ++(2.0pt,-2.0pt);
\draw[color=dtsfsf] (-2.7866075273047657,1.5864886565649992) node {$v_2'$};
\draw[color=dtsfsf] (-2.283940498116759,1.4577061945416254) node {$d$};
\draw [fill=black] (-1.5,1.5) circle (0.5pt);
\draw[color=black] (-1.2993778045832245,1.6695741159349176) node {$v_2$};
\draw [fill=black] (-0.5,1.) circle (0.5pt);
\draw[color=black] (-0.5433001243169658,1.2873810028332928) node {$v_1$};
\draw[color=black] (-0.979498786009038,1.4244720107936581) node {$a$};
\draw [fill=black] (-1.5,0.5) circle (0.5pt);
\draw[color=black] (-1.4987829070710288,0.265429852583296) node {$x_2$};
\draw[color=black] (-0.8964133266391197,0.5603832333465064) node {$2a$};
\draw [fill=dtsfsf] (-2.5,0.5) circle (0.5pt);
\draw[color=dtsfsf] (-2.6869049760608634,0.30697258226825525) node {$x_2'$};
\draw[color=dtsfsf] (-2.10946103343993,0.27789267148878377) node {$\frac{d}{2}$};
\draw [fill=black] (-1.5,-1.5) circle (0.5pt);
\draw[color=black] (-1.4821658151970452,-1.662152804798812) node {$w_2$};
\draw [fill=rvwvcq] (-2.5,-1.5) circle (0.5pt);
\draw[color=rvwvcq] (-2.711830613871839,-1.5790673454288935) node {$w_2'$};
\draw [fill=black] (0.5,2.) circle (0.5pt);
\draw[color=black] (0.7445244959167711,2.1016185046584934) node {$v_3$};
\draw[color=rvwvcq] (-0.09048437075091002,1.6238771132814627) node {$b$};
\draw[color=black] (1.3385855304116883,1.922984767013169) node {$e$};
\draw [color=black] (1.5,2.)-- ++(-1.0pt,-1.0pt) -- ++(2.0pt,2.0pt) ++(-2.0pt,0) -- ++(2.0pt,-2.0pt);
\draw[color=black] (1.8246354677257117,2.08500141278451) node {$v_3'$};
\draw [fill=black] (0.4828052989015277,-0.8977665785955626) circle (0.5pt);
\draw[color=black] (0.47865102593303177,-1.122097318894342) node {$x_3$};
\draw[color=rvwvcq] (0.24185746672876401,-0.021214982242922606) node {$2b$};
\draw [fill=black] (1.4881393572775417,-0.9143836704695463) circle (0.5pt);
\draw[color=black] (1.533836359930997,-1.1304058648313338) node {$x_3'$};
\draw[color=black] (1.0311693307429899,-1.159485775610805) node {$\frac{e}{2}$};
\draw [fill=black] (0.5,-2.) circle (0.5pt);
\draw[color=black] (0.39556556656311326,-2.1440484691443387) node {$w_3$};
\draw [fill=black] (1.5,-2.) circle (0.5pt);
\draw[color=black] (1.525527813994005,-2.1191228313333634) node {$w_3'$};
\draw[color=black] (1.081020606364941,-2.2312882014827533) node {$\frac{e}{2}$};
\draw [fill=black] (1.,0.5) circle (0.5pt);
\draw[color=black] (1.243037252136282,0.5645375063150024) node {$v_4$};
\draw[color=dtsfsf] (0.2750916504767314,0.8927250708261802) node {$c$};
\draw [color=qqwuqq] (2.,0.5)-- ++(-1.0pt,-1.0pt) -- ++(2.0pt,2.0pt) ++(-2.0pt,0) -- ++(2.0pt,-2.0pt);
\draw[color=qqwuqq] (2.2816054942602633,0.5645375063150024) node {$v_4'$};
\draw [fill=black] (0.989626601058031,-0.3078598170691431) circle (0.5pt);
\draw[color=black] (1.135026154955388,-0.10014616864434513) node {$x_4$};
\draw[color=dtsfsf] (0.23354892079177214,0.585308871157482) node {$2c$};
\draw [fill=qqwuqq] (2.003269205371037,-0.31616836300613493) circle (0.5pt);
\draw[color=qqwuqq] (2.1403602133314017,-0.10014616864434513) node {$x_4'$};
\draw[color=qqzzqq] (1.7207786435133134,-0.5778875600213762) node {$\frac{f}{2}$};
\draw [fill=dtsfsf] (2.,-1.5) circle (0.5pt);
\draw[color=dtsfsf] (2.173594397079369,-1.5707587994919017) node {$w_4'$};
\draw[color=black] (3.3907963768486757,-0.14584317129780028) node {$\varphi$};
\draw[color=qqwuqq] (1.8204811947572157,0.4606806821026043) node {$f$};
\draw[color=wrwrwr] (-4.253065885183827,3.518225586915603) node {$Curve$};
\end{scriptsize}
\end{tikzpicture}
    \begin{center}
        \textnormal{Figure 25.2. The model $(G',l')$ of $\Gamma'$}
    \end{center}
\end{figure}
\begin{figure}[H]
    \centering
    \begin{tikzpicture}[line cap=round,line join=round,>=triangle 45,x=1.0cm,y=1.0cm]
\definecolor{wrwrwr}{rgb}{0.3803921568627451,0.3803921568627451,0.3803921568627451}
\definecolor{qqwuqq}{rgb}{0.,0.39215686274509803,0.}
\definecolor{rvwvcq}{rgb}{0.08235294117647059,0.396078431372549,0.7529411764705882}
\definecolor{dtsfsf}{rgb}{0.8274509803921568,0.1843137254901961,0.1843137254901961}
\clip(-3.,-2.5) rectangle (2.5,-0.6);
\draw [color=dtsfsf] (-2.,1.5) circle (0.5cm);
\draw (-1.5,1.5)-- (-0.5,1.);
\draw (-0.5,1.)-- (-1.5,0.5);
\draw [color=dtsfsf] (-1.5,0.5)-- (-2.5,0.5);
\draw (-0.5,-1.)-- (-1.5,-1.5);
\draw [color=dtsfsf] (-1.5,-1.5)-- (-2.5,-1.5);
\draw [color=rvwvcq] (-0.5,1.)-- (0.5,2.);
\draw(1.,2.) circle (0.5cm);
\draw [color=rvwvcq] (-0.5,1.)-- (0.49111384483851955,-0.5072649195569464);
\draw (0.49111384483851955,-0.5072649195569464)-- (1.513064995088517,-0.5155734654939388);
\draw [color=rvwvcq] (-0.5,-1.)-- (0.5,-2.);
\draw (0.5,-2.)-- (1.5,-2.);
\draw [color=dtsfsf] (-0.5,1.)-- (1.,0.5);
\draw [color=dtsfsf] (-0.5,1.)-- (0.989626601058031,-0.3078598170691431);
\draw [color=dtsfsf] (0.989626601058031,-0.3078598170691431)-- (2.003269205371037,-0.31616836300613493);
\draw [color=dtsfsf] (-0.5,-1.)-- (1.,-1.5);
\draw [color=qqwuqq] (1.,-1.5)-- (2.,-1.5);
\draw (2.8922836206291644,1.2707639109593083) node[anchor=north west] {$\Gamma'$};
\draw (2.9005921665661565,-1.2550340538862124) node[anchor=north west] {$T$};
\draw [->] (3.056672577778454,0.5664683674959348) -- (3.056672577778454,-0.9335316325040652);
\draw [color=qqwuqq] (1.5,0.5) circle (0.5cm);
\begin{scriptsize}
\draw [color=dtsfsf] (-2.5,1.5)-- ++(-1.0pt,-1.0pt) -- ++(2.0pt,2.0pt) ++(-2.0pt,0) -- ++(2.0pt,-2.0pt);
\draw[color=dtsfsf] (-2.283940498116759,1.4577061945416245) node {$d$};
\draw [fill=black] (-1.5,1.5) circle (0.5pt);
\draw[color=black] (-1.2993778045832245,1.6695741159349167) node {$v_2$};
\draw [fill=black] (-0.5,1.) circle (0.5pt);
\draw[color=black] (-0.5433001243169658,1.2873810028332917) node {$v_1$};
\draw[color=black] (-0.979498786009038,1.4244720107936573) node {$a$};
\draw [fill=black] (-1.5,0.5) circle (0.5pt);
\draw[color=black] (-0.9545731481980624,0.48560631991357894) node {$2a$};
\draw [fill=dtsfsf] (-2.5,0.5) circle (0.5pt);
\draw[color=dtsfsf] (-1.8518961093931825,0.1781901202448807) node {$\frac{d}{2}$};
\draw [fill=black] (-0.5,-1.) circle (0.5pt);
\draw[color=black] (-0.40205484338810443,-0.7565212976667017) node {$w_1$};
\draw [fill=black] (-1.5,-1.5) circle (0.5pt);
\draw[color=black] (-1.4821658151970452,-1.6621528047988128) node {$w_2$};
\draw[color=black] (-0.9961158778830217,-1.0265490406189368) node {$2a$};
\draw [fill=black] (-2.5,-1.5) circle (0.5pt);
\draw[color=black] (-2.711830613871839,-1.5790673454288944) node {$w_2'$};
\draw[color=dtsfsf] (-1.8851302931411498,-1.1844114134217816) node {$\frac{d}{2}$};
\draw [fill=black] (0.5,2.) circle (0.5pt);
\draw[color=black] (0.7445244959167711,2.1016185046584925) node {$v_3$};
\draw[color=rvwvcq] (-0.09048437075091002,1.6238771132814618) node {$b$};
\draw[color=black] (1.3385855304116883,1.9229847670131681) node {$e$};
\draw [color=black] (1.5,2.)-- ++(-1.0pt,-1.0pt) -- ++(2.0pt,2.0pt) ++(-2.0pt,0) -- ++(2.0pt,-2.0pt);
\draw [fill=black] (0.49111384483851955,-0.5072649195569464) circle (0.5pt);
\draw[color=rvwvcq] (0.2501660126657559,0.16157302837089701) node {$2b$};
\draw [fill=black] (1.513064995088517,-0.5155734654939388) circle (0.5pt);
\draw [fill=black] (0.5,-2.) circle (0.5pt);
\draw[color=black] (0.39556556656311326,-2.1440484691443396) node {$w_3$};
\draw[color=rvwvcq] (-0.007398911380991524,-1.7909352668221863) node {$2b$};
\draw [fill=black] (1.5,-2.) circle (0.5pt);
\draw[color=black] (1.525527813994005,-2.1191228313333643) node {$w_3'$};
\draw[color=black] (1.081020606364941,-2.231288201482754) node {$\frac{e}{2}$};
\draw [fill=black] (1.,0.5) circle (0.5pt);
\draw[color=black] (1.243037252136282,0.5645375063150014) node {$v_4$};
\draw[color=dtsfsf] (0.2750916504767314,0.8927250708261792) node {$c$};
\draw [color=qqwuqq] (2.,0.5)-- ++(-1.0pt,-1.0pt) -- ++(2.0pt,2.0pt) ++(-2.0pt,0) -- ++(2.0pt,-2.0pt);
\draw [fill=black] (0.989626601058031,-0.3078598170691431) circle (0.5pt);
\draw[color=dtsfsf] (0.23354892079177214,0.5853088711574811) node {$2c$};
\draw [fill=black] (2.003269205371037,-0.31616836300613493) circle (0.5pt);
\draw[color=dtsfsf] (2.2774512212917677,-0.2538542684786952) node {$\frac{f}{2}$};
\draw [fill=black] (1.,-1.5) circle (0.5pt);
\draw[color=black] (0.9937808740265266,-1.6372271669878373) node {$w_4$};
\draw[color=dtsfsf] (0.3664856557836418,-1.49182761309048) node {$2c$};
\draw [fill=black] (2.,-1.5) circle (0.5pt);
\draw[color=black] (2.173594397079369,-1.5707587994919026) node {$w_4'$};
\draw[color=qqwuqq] (1.7207786435133134,-1.1844114134217816) node {$\frac{f}{2}$};
\draw[color=black] (3.3907963768486757,-0.14584317129780122) node {$\varphi$};
\draw[color=qqwuqq] (1.8204811947572157,0.4606806821026034) node {$f$};
\draw[color=wrwrwr] (-4.327842798616754,3.4766828572306427) node {$Curve$};
\end{scriptsize}
\end{tikzpicture}
    \begin{center}
        \textnormal{Figure 25.3. The model $(T',t')$ of $T$}
    \end{center}
\end{figure}
\begin{figure}[H]
    \centering
    \begin{tikzpicture}[line cap=round,line join=round,>=triangle 45,x=1.0cm,y=1.0cm]

\definecolor{qqwuqq}{rgb}{0.,0.39215686274509803,0.}
\definecolor{eqeqeq}{rgb}{0.8784313725490196,0.8784313725490196,0.8784313725490196}
\definecolor{cqcqcq}{rgb}{0.7529411764705882,0.7529411764705882,0.7529411764705882}
\definecolor{rvwvcq}{rgb}{0.08235294117647059,0.396078431372549,0.7529411764705882}
\definecolor{dtsfsf}{rgb}{0.8274509803921568,0.1843137254901961,0.1843137254901961}
\clip(-3.,-3.) rectangle (4.,3.);
\draw [color=dtsfsf] (-2.,1.5) circle (0.5cm);
\draw (-1.5,1.5)-- (-0.5,1.);
\draw (-0.5,1.)-- (-1.5,0.5);
\draw [color=dtsfsf] (-1.5,0.5)-- (-2.5,0.5);
\draw (-0.5,-1.)-- (-1.5,-1.5);
\draw [color=dtsfsf] (-1.5,-1.5)-- (-2.5,-1.5);
\draw [line width=0.4pt,dash pattern=on 2pt off 2pt,color=cqcqcq] (-2.5,1.5)-- (-2.5,-1.5);
\draw [line width=0.4pt,dash pattern=on 2pt off 2pt,color=eqeqeq] (-2.5,-1.5)-- (-2.5,1.5);
\draw [line width=0.4pt,dash pattern=on 2pt off 2pt,color=cqcqcq] (-1.5,-1.5)-- (-1.5,1.5);
\draw [line width=0.4pt,dash pattern=on 2pt off 2pt,color=cqcqcq] (-0.5,-1.)-- (-0.5,1.);
\draw [color=rvwvcq] (-0.5,1.)-- (0.5,2.);
\draw(1.,2.) circle (0.5cm);
\draw [color=rvwvcq] (-0.5,1.)-- (0.49942239077551187,-0.7482127517297108);
\draw (0.49942239077551187,-0.7482127517297108)-- (1.5047564491515257,-0.7565212976667027);
\draw [color=rvwvcq] (-0.5,-1.)-- (0.5,-2.);
\draw (0.5,-2.)-- (1.5,-2.);
\draw [line width=0.4pt,dash pattern=on 2pt off 2pt,color=cqcqcq] (0.5,2.)-- (0.5,-2.);
\draw [line width=0.4pt,dash pattern=on 2pt off 2pt,color=cqcqcq] (1.5,2.)-- (1.5,-2.);
\draw [color=dtsfsf] (-0.5,1.)-- (1.,0.5);
\draw [color=dtsfsf] (-0.5,1.)-- (0.989626601058031,-0.3078598170691431);
\draw [color=dtsfsf] (0.989626601058031,-0.3078598170691431)-- (2.003269205371037,-0.31616836300613493);
\draw [color=dtsfsf] (-0.5,-1.)-- (1.,-1.5);
\draw [color=qqwuqq] (1.,-1.5)-- (2.,-1.5);
\draw [line width=0.4pt,dash pattern=on 2pt off 2pt,color=cqcqcq] (2.,0.5)-- (2.,-1.5);
\draw [line width=0.4pt,dash pattern=on 2pt off 2pt,color=cqcqcq] (1.,0.5)-- (1.,-1.5);
\draw (2.8922836206291644,1.2707639109593076) node[anchor=north west] {$\Gamma'$};
\draw (2.9005921665661565,-1.2550340538862128) node[anchor=north west] {$T$};
\draw [->] (3.056672577778454,0.5664683674959348) -- (3.056672577778454,-0.9335316325040652);
\draw [color=qqwuqq] (1.5,0.5) circle (0.5cm);
\begin{scriptsize}
\draw [color=dtsfsf] (-2.5,1.5)-- ++(-1.0pt,-1.0pt) -- ++(2.0pt,2.0pt) ++(-2.0pt,0) -- ++(2.0pt,-2.0pt);
\draw[color=dtsfsf] (-2.283940498116759,1.457706194541624) node {$d$};
\draw [fill=black] (-1.5,1.5) circle (0.5pt);
\draw[color=black] (-1.2993778045832245,1.6695741159349162) node {$v_2$};
\draw [fill=black] (-0.5,1.) circle (0.5pt);
\draw[color=black] (-0.5433001243169658,1.2873810028332913) node {$v_1$};
\draw[color=black] (-0.979498786009038,1.4244720107936568) node {$a$};
\draw [fill=black] (-1.5,0.5) circle (0.5pt);
\draw[color=black] (-0.9545731481980624,0.48560631991357844) node {$2a$};
\draw [fill=dtsfsf] (-2.5,0.5) circle (0.5pt);
\draw[color=dtsfsf] (-1.8518961093931825,0.17819012024488023) node {$\frac{d}{2}$};
\draw [fill=black] (-0.5,-1.) circle (0.5pt);
\draw [fill=black] (-1.5,-1.5) circle (0.5pt);
\draw[color=black] (-0.9961158778830217,-1.0265490406189373) node {$2a$};
\draw [fill=rvwvcq] (-2.5,-1.5) circle (0.5pt);
\draw[color=dtsfsf] (-1.8851302931411498,-1.184411413421782) node {$\frac{d}{2}$};
\draw [fill=black] (0.5,2.) circle (0.5pt);
\draw[color=black] (0.7445244959167711,2.101618504658492) node {$v_3$};
\draw[color=rvwvcq] (-0.09048437075091002,1.623877113281461) node {$b$};
\draw[color=black] (1.3385855304116883,1.9229847670131677) node {$e$};
\draw [color=black] (1.5,2.)-- ++(-1.0pt,-1.0pt) -- ++(2.0pt,2.0pt) ++(-2.0pt,0) -- ++(2.0pt,-2.0pt);
\draw [fill=black] (0.49942239077551187,-0.7482127517297108) circle (0.5pt);
\draw[color=rvwvcq] (0.2501660126657559,0.053561931190002576) node {$2b$};
\draw [fill=black] (1.5047564491515257,-0.7565212976667027) circle (0.5pt);
\draw[color=black] (1.1641060657348596,-0.9600806731230024) node {$\frac{e}{2}$};
\draw [fill=black] (0.5,-2.) circle (0.5pt);
\draw[color=rvwvcq] (-0.007398911380991524,-1.7909352668221867) node {$2b$};
\draw [fill=black] (1.5,-2.) circle (0.5pt);
\draw[color=black] (1.081020606364941,-2.2312882014827546) node {$\frac{e}{2}$};
\draw [fill=black] (1.,0.5) circle (0.5pt);
\draw[color=black] (1.243037252136282,0.564537506315001) node {$v_4$};
\draw[color=dtsfsf] (0.2750916504767314,0.8927250708261788) node {$c$};
\draw [color=qqwuqq] (2.,0.5)-- ++(-1.0pt,-1.0pt) -- ++(2.0pt,2.0pt) ++(-2.0pt,0) -- ++(2.0pt,-2.0pt);
\draw [fill=black] (0.989626601058031,-0.3078598170691431) circle (0.5pt);
\draw[color=dtsfsf] (0.23354892079177214,0.5853088711574805) node {$2c$};
\draw [fill=black] (2.003269205371037,-0.31616836300613493) circle (0.5pt);
\draw[color=dtsfsf] (2.2774512212917677,-0.25385426847869563) node {$\frac{f}{2}$};
\draw [fill=black] (1.,-1.5) circle (0.5pt);
\draw[color=dtsfsf] (0.3664856557836418,-1.4918276130904804) node {$2c$};
\draw [fill=dtsfsf] (2.,-1.5) circle (0.5pt);
\draw[color=qqwuqq] (1.7207786435133134,-1.184411413421782) node {$\frac{f}{2}$};
\draw[color=black] (3.3907963768486757,-0.1458431712978017) node {$\varphi$};
\draw[color=qqwuqq] (1.8204811947572157,0.46068068210260293) node {$f$};
\end{scriptsize}
\end{tikzpicture}
    \begin{center}
        \textnormal{Figure 26. The tropical morphism $\varphi:\Gamma' \to T$}
    \end{center}
\end{figure}

\end{document}